\documentstyle[amssymb,10pt]{amsart}
 
\newtheorem{theorem}{Theorem}[section]
\newtheorem{lemma}[theorem]{Lemma}
\newtheorem{proposition}[theorem]{Proposition}

\theoremstyle{remark}  
\newtheorem{remark}[theorem]{Remark}

\newcommand{\nc}{\newcommand}
\nc{\Symm}{{\on{Sym}}}
\nc{\Perm}{{\on{Perm}}}
\nc{\Iff}{\Leftrightarrow}
\newcommand{\on}{\operatorname}   
\newcommand{\eps}{\varepsilon}

 \nc{\cE}{{\cal E}}
\nc{\cL}{{\cal L}}

\nc{\SL}{{\mathfrak sl}}
\nc{\gt}{{\mathfrak gt}}
\nc{\grt}{{\mathfrak grt}}
\nc{\gtm}{{\mathfrak gtm}}
\nc{\grtm}{{\mathfrak grtm}}
\nc{\gtmd}{{\mathfrak gtmd}}
\nc{\grtmd}{{\mathfrak grtmd}}

\renewcommand{\a}{{\mathfrak{a}}}

\nc{\HH}{{\mathfrak H}}
\newcommand{\g}{{\mathfrak{g}}}

\renewcommand{\t}{{\mathfrak{t}}}
\newcommand{\SG}{{\mathfrak{S}}}

\newcommand{\ul}{\underline}

\nc{\wh}{\widehat}\nc{\wt}{\widetilde}

\newcommand{\kk}{{\bf k}}

\newcommand{\ben}{\begin{enumerate}}
\newcommand{\een}{\end{enumerate}}

\newcommand{\Vect}{{\text{Vect}}}

\newcommand{\cO}{{\mathcal O}}

\newcommand{\cS}{{\mathcal S}}

\newcommand{\cX}{{\mathcal X}}
\newcommand{\cY}{{\mathcal Y}}

\newcommand{\QQ}{{\mathbb{Q}}}

\newcommand{\cC}{{\mathcal C}}
\newcommand{\cA}{{\mathcal A}}
\newcommand{\cF}{{\mathcal F}}
\newcommand{\cG}{{\mathcal G}}

\renewcommand{\t}{{\mathfrak t}}
\newcommand{\cB}{{\mathcal B}}
\newcommand{\cD}{{\mathcal D}}

\newcommand{\NN}{{\mathbb N}}

\begin{document}

\title{Quantization of coboundary Lie bialgebras}

\begin{abstract}
We show that any coboundary Lie bialgebra can be quantized. 
For this, we prove that: (a) Etingof-Kazhdan quantization functors 
are compatible with Lie bialgebra twists, and (b) if such a quantization 
functor corresponds to an even associator, then it is also compatible with the 
operation of taking coopposites. We also use the relation between 
the Etingof-Kazhdan construction of quantization functors and the 
alternative approach to this problem, which was established in a
previous work. 
\end{abstract}

\author{Benjamin Enriquez}
\address{IRMA (CNRS), rue Ren\'e Descartes, F-67084 Strasbourg, France}
\email{enriquez@@math.u-strasbg.fr}

\author{Gilles Halbout}
\address{IRMA (CNRS), rue Ren\'e Descartes, F-67084 Strasbourg, France}
\email{halbout@@math.u-strasbg.fr}

\maketitle

Let $\kk$ be a field of characteristic $0$. Unless specified otherwise, 
``algebra'', ``vector space'', etc.,  means ``algebra over $\kk$", etc. 

\section*{Introduction}

In this paper, we solve the problem of quantization of coboundary Lie 
bialgebras. This is one of the quantization problems of Drinfeld's list 
(\cite{Dr:uns}). This result can be viewed as a completion of the result 
of twist quantization of Lie bialgebras (\cite{H}, solving a problem posed 
in \cite{Stolin}.) 

We show that our result, together with a proposition of \cite{Dr:QH}, 
implies that quasi-Poisson manifolds over a pair $(\g,Z)$ ($\g$ a Lie algebra, 
$Z\in \wedge^3(\g)^\g$) can be quantized in the case when the underlying space 
is the group itself (this problem was posed in \cite{EE}). 

To solve the problem of quantization of coboundary Lie bialgebras, 
we show that quantization functors of Lie bialgebras are compatible with Lie 
bialgebra twists. The quantization of all the affine Poisson groups (\cite{DS}) 
of Dazord and Sondaz (i.e., Poisson homogeneous spaces under a 
Poisson-Lie group, which are principal as homogeneous spaces; 
see \cite{Dr:hom}) follows immediately from this result. 
It is also a basic case of the quantization problem of quasi-Lie bialgebras 
(together with their twists) into quasi-Hopf algebras (also a problem 
of Drinfeld's list), which is still open. 

We now describe the problem of quantization of coboundary Lie bialgebras. 
A coboundary Lie bialgebra is a pair $(\a,r_\a)$, where $\a$ is a 
Lie algebra (with Lie bracket denoted by $\mu_\a$) and $r_\a\in \wedge^2(\a)$ 
is such that 
$Z_\a := [r_\a^{12},r_\a^{13}] + [r_\a^{12},r_\a^{23}]
+ [r_\a^{13},r_\a^{23}]\in \wedge^3(\a)^\a$. To $(\a,r_\a)$
is associated a Lie bialgebra with cobracket $\delta_\a : \a\to \wedge^2(\a)$
given by $\delta_\a(x) = [r_\a,x\otimes 1 + 1 \otimes x]$.  

A coboundary QUE algebra is a pair $(U,R_U)$, where 
$(U,m_U,\Delta_U,\eps_U,\eta_U)$ is a QUE (quantized universal 
enveloping) algebra (i.e., a deformation of an enveloping algebra
in the category of topological $\kk[[\hbar]]$-modules), and 
$R_U\in (U^{\otimes 2})^\times$ is such that 
\begin{equation} \label{cond:0}
\Delta_U(x)^{21} = R_U \Delta_U(x) R_U^{-1}, 
\quad R_U R_U^{21} = 1_U^{\otimes 2} 
\end{equation}
\begin{equation} \label{twists:1}
R_U^{12} (\Delta_U \otimes \on{id}_U)(R_U)
= R_U^{23} (\on{id}_U \otimes \Delta_U)(R_U), 
\end{equation}
\begin{equation} \label{twists:2}
R_U = 1_U^{\otimes 2}\on{\ mod\ }\hbar, \quad 
(\eps_U \otimes \on{id}_U)(R_U) = (\on{id}_U \otimes \eps_U)(R_U) = 1_U.  
\end{equation}
$(U,R_U)$ is a quantization of $(\a,r_\a)$ if the classical limit of $U$ is 
$(\a,\mu_\a,\delta_\a)$, and if 
\begin{equation} \label{twists:3}
\big( \hbar^{-1}(R_U^{21} - R_U) \on{\ mod\ }\hbar \big) = 2r_\a
\end{equation} 
($1_U = \eta_U(1)$ is the unit of $U$). 
The problem of quantization of coboundary Lie bialgebras is that of 
constructing a quantization $(U,R_U)$ for each coboundary Lie bialgebra
$(\a,r_\a)$ (\cite{Dr:uns,Dr:QG}). Our solution is 
formulated in the language of props (\cite{McL}). 
Recall that to a prop $P$ and a symmetric tensor category $\cS$, one associates
the category $\on{Rep}_{\cS}(P)$ of $P$-modules in $\cS$. A prop morphism 
$P\to Q$ then gives rise to a functor $\on{Rep}_\cS(Q) \to \on{Rep}_\cS(P)$. 
A quantization problem may often be formulated as the problem of constructing 
a functor $\on{Rep}_{\cS}(P_{\on{class}}) \to \on{Rep}_{\cS}(P_{\on{quant}})$, 
where $\cS = \on{Vect}$ (the category of vector spaces) and 
$P_{\on{class}}$, $P_{\on{quant}}$ are suitable 
``classical" and ``quantum" props. The propic version of the quantization
problem  is then to construct a suitable prop morphism $P_{\on{quant}} 
\to P_{\on{class}}$. 

We construct props COB 
and $\on{Cob}$ of coboundary bialgebras and of coboundary Lie bialgebras. 
Using an even associator defined over $\kk$ (see \cite{Dr:Gal,BN}), we 
construct a prop morphism $\on{COB} \to S({\bf Cob})$ with suitable 
properties (${\bf Cob}$ is a completion of Cob, and $S$ is the 
symmetric algebra Schur functor). This allows to also solve the 
problem of quantization of coboundary Lie bialgebras in symmetric 
tensor categories (when $\cS = \on{Vect}$, this is the original problem). 

Our construction is based on the theory of twists of Lie bialgebras 
(\cite{Dr:QH}). Recall that if $(\a,\mu_\a,\delta_\a)$ is a Lie bialgebra, 
then $f_\a\in \wedge^2(\a)$ is called a twist of $\a$ if 
$(\delta_\a\otimes \on{id}_\a)(f_\a) + [f_\a^{13},f_\a^{23}] +$ 
cyclic permutations $= 0$. If we set $\on{ad}(f_\a)(x) = [f_\a,x^1+x^2]$, 
then $(\a,\mu_\a,\delta_\a + \on{ad}(f_\a))$ is again a Lie bialgebra
(the twisted Lie bialgebra). 

A quantization of $(\a,f_\a)$ is a pair $(U,F_U)$, where
$(U,m_U,\Delta_U,\eps_U,\eta_U)$ is a QUE algebra quantizing 
$(\a,\mu_\a,\delta_\a)$, and $F_U\in (U^{\otimes 2})^\times$
satisfies the above conditions (\ref{twists:1}), (\ref{twists:2}), 
and (\ref{twists:3}), with $(-2r_\a,R_U)$ replaced by $(f_\a,F_U)$. 
Then $(U,m_U,\on{Ad}(F_U)\circ \Delta_U,\eps_U,\eta_U)$ is again a 
QUE algebra (the twisted QUE algebra, denoted ${}^{F_U}U$) and is 
a quantization of 
$(\a,\mu_\a,\delta_\a + \on{ad}(f_\a))$ (here $\on{Ad}(F_U)\in 
\on{Aut}(U^{\otimes 2})$ is $x\mapsto F_U x F_U^{-1}$).  

We notice that if $(\a,r_\a)$ is a coboundary Lie bialgebra, then 
$-2r_\a$ is a twist of $(\a,\mu_\a,\delta_\a = \on{ad}(r_\a))$, and the 
resulting twisted Lie bialgebra is $(\a,\mu_\a,-\delta_\a)$ (which is the 
coopposite of $(\a,\mu_\a,\delta_\a)$). Moreover, a quantization of 
$(\a,r_\a)$ is the same as a quantization $(U,m_U,\Delta_U,\eps_U,\eta_U)$
of $(\a,\mu_\a,\delta_a)$, together with a twist $R_U$ of this QUE algebra, 
satisfying the additional condition (\ref{cond:0}); the second part of 
(\ref{cond:0}) means in particular that the twisted QUE algebra is 
$(U,m_U,\Delta_U^{21},\eps_U,\eta_U)$, i.e., the coopposite of the initial QUE 
algebra. 

On the other hand, Etingof and Kazhdan constructed a quantization functor 
$Q : \on{Bialg} \to \on{S}({\bf LBA})$ for each Drinfeld associator 
defined over $\kk$ (\cite{EK}); here Bialg is the prop of 
bialgebras and ${\bf LBA}$ is a suitable completion of the prop LBA of 
Lie bialgebras. We also denote by $Q : \{$Lie bialgebras over 
Vect$\} \to \{$QUE algebras over Vect$\}$ the functor induced by this prop
morphism. Our construction involves three steps: 

(a) we show that any Etingof-Kazhdan quantization functor $Q$ is compatible 
with twists. This is a propic 
version of the statement that for any $(\a,f_\a)$, 
where $\a$ is a Lie bialgebra and $f_\a$ is a twist of $\a$, there exists 
an element $\on{F}(\a,f_\a)\in Q(\a,\mu_\a,\delta_\a)^{\otimes 2}$
satisfying the twist conditions, such that the twisted QUE algebra 
${}^{\on{F}(\a,f_\a)}Q(\a,\mu_\a,\delta_\a)$ is isomorphic to 
$Q(\a,\mu_\a,\delta_\a + \on{ad}(f_\a))$;  

(b) we show that if $Q$ corresponds to an even associator, then 
$Q$ is compatible with the operation of taking coopposite Lie bialgebras and 
QUE algebras. This is a propic version of the statement that for any 
Lie bialgebra $(\a,\mu_\a,\delta_\a)$, the QUE algebras 
$Q(\a,\mu_\a,-\delta_\a)$ and $Q(\a,\mu_\a,\delta_\a)^{\on{cop}}$
are isomorphic (here $U^{\on{cop}}$ is the coopposite QUE algebra of a 
QUE algebra $U$); 

(c) we are then in the following situation (at the propic level). If 
$(\a,r_\a)$ is a coboundary Lie bialgebra, then $Q(\a,\mu_\a,-\delta_\a)
\simeq Q(\a,\mu_\a,\delta_\a)^{\on{cop}}$ (where cop means the bialgebra with 
the opposite coproduct, and $\simeq$ is an isomorphism of QUE algebras)
and  $Q(\a,\mu_\a,-\delta_\a)
\simeq Q(\a,\mu_\a,\delta_\a)^{\on{F}(\a,-2r_\a)}$, therefore 
$Q(\a,\mu_\a,\delta_\a)^{\on{cop}}
\simeq Q(\a,\mu_\a,\delta_\a)^{\on{F}(\a,-2r_\a)}$. One then proves
(at the propic level) that this implies the existence of a twist  
$R(\a,r_\a)$, such that $R(\a,r_\a)R(\a,r_\a)^{21} = 1_U^{\otimes 2}$
and $Q(\a,\mu_\a,\delta_\a)^{\on{cop}} = 
Q(\a,\mu_\a,\delta_\a)^{R(\a,r_\a)}$. This solves the quantization problem of
coboundary Lie bialgebras.  

Let us now describe the contents of the paper. In Section \ref{sect:props}, 
we recall the formalism of props. We introduce related notions: quasi-props 
and quasi-bi-multiprops. Recall that a prop (e.g., LBA) consists of the 
universal versions $\on{LBA}(F,G)$ of the spaces
of linear maps $F(\a) \to G(\a)$ (where $\a$ is a Lie bialgebra, 
and $F,G$ are Schur functors), constructed from $\mu_\a,\delta_\a$ and avoiding
cycles. The corresponding quasi-biprop consists of universal versions of the  
spaces of maps $F(\a) \otimes F'(\a^*) \to G(\a) \otimes G'(\a^*)$; 
by partial transposition, this identifies with 
$\on{LBA}(F\otimes (G')^*,F'\otimes G^*)$, but due to 
the possible introduction of cycles, the composition is only partially 
defined: it is defined iff the ``trace" of some element is. One 
then constructs a (partially defined) trace map 
$\on{LBA}(F\otimes G,F\otimes G') \to \on{LBA}(G,G')$, which consists in 
closing the graph by connecting $F$ with itself. One can also encounter the
following situation: $F = \otimes_{i=1}^n F_i$, and $x\in 
\on{LBA}(F\otimes G,F\otimes G')$  is such that the element obtained by 
connecting each $F_i$ with itself has no cycle. This defines a trace map 
$\on{LBA}(F\otimes G,F\otimes G') \to \on{LBA}(G,G')$, which depends on the
data of $(F_i)_{i=1,...,n}$ such that $F = \otimes_{i=1}^n F_i$; it 
actually depends on the multi-Schur functor $\boxtimes_{i=1}^n F_i$. 
In the corresponding notion of a prop
(quasi-bi-multiprops), the basic objects are bi-multi-Schur functors
(the ``bi" analogue of a multi-Schur functor). We introduce in the end of
Section \ref{sect:props} the main quasi-bi-multiprops we will be working 
with, $\Pi$ and $\Pi_f$ and their variants. In Section  \ref{sec:2}, we introduce 
the universal algebras ${\bf U}_n$ and ${\bf U}_{n,f}$ (we have morphisms
${\bf U}_n \to U(\a)^{\otimes n}$ if $\a$ is any Lie bialgebra, and 
${\bf U}_{n,f} \to U(\a)^{\otimes n}$ if $\a$ is any Lie bialgebra equipped 
with a Lie bialgebra twist). In Section \ref{sect:inject}, we prove the 
injectivity of a map; this will be crucial for proving the compatibility of
quantization functors with twists (step (a) above). In Section 
\ref{sect:quant}, we present the construction of quantization 
functors of \cite{Enr:coh} (in the framework of 
quasi-bi-multiprops), which can be viewed as an alternative 
to the construction of \cite{EK}. Its basic ingredients are a twist 
$\on{J}$ killing an associator $\Phi$, and a factorization 
result for the corresponding $R$-matrix. 
In Section \ref{sect:twists}, we prove the compatibility of quantization 
with twists (step (a) in the above description). As in \cite{Enr:coh}, the 
proof involves two steps: an ``easy" co-Hochschild cohomology argument, 
and a more involved injectivity result (which was proved in 
Section \ref{sect:inject}). In Section 
\ref{sect:cob}, we perform steps (b) and (c), i.e., we study the 
behavior of quantization functors with the operation of taking coopposites, 
and ``correct" the twist $\on{F}(\a,-2r_\a)$ into a quantization of 
coboundary Lie bialgebras. 
Finally, in Subsection \ref{sec:qpoisson}, we show how quantization 
of coboundary Lie bialgebra implies that of certain quasi-Poisson 
homogeneous spaces. 

\subsection*{Notation}

If $A = \oplus_{n\geq 0} A_n$ is a graded vector space, we 
denote by $\wh A = \wh\oplus_{n\geq 0} A_n$ its completion w.r.t. 
the grading. If $A$ is an algebra, we denote by $A^\times$
the group of its invertible elements. If the algebra $A$ is equipped with a
character $\chi$, then we denote by $A^\times_1$ the kernel of 
$\chi : A^\times \to \kk^\times$. If $A$ is a graded and connected algebra, 
then a graded character $\chi$ is unique; we will use it for defining 
$A^\times_1$ and $\wh A^\times_1$. 

\section{Props and (quasi)(multi)(bi)props} \label{sect:props}

In this section, we define various ``Schur categories'', which are all symmetric 
monoidal categories. We then define monoidal quasi-categories and show how 
they can be constructed using partial traces on monoidal categories. 
We then define (quasi)(multi)(bi)props, and show that variants of the 
prop of Lie bialgebras yield examples of these structures. 

\subsection{Schur categories}

If $\cO$ is a category, we denote by $\on{Ob}(\cO)$ its set of objects and 
by $\on{Irr}(\cO)$ the set of 
isomorphism classes of irreducible objects of $\cO$. We denote by 
$\on{Vect}$ the category of finite dimensional $\kk$-vector spaces. 

For $n\geq 0$, let $\wh\SG_{n}$ denote the set of isomorphism classes 
of irreductible representations of $\SG_{n}$ (by convention, 
$\SG_{0}=\{1\}$). We view  $\sqcup_{n\geq 0}\wh\SG_{n}$ 
as the set of pairs $(n,\pi)$, where $n\geq 0$ and $\pi\in\wh\SG_{n}$. 
For $\rho = (n,\pi)$, we set $|\rho|:= n$ and $\pi_{\rho}:= \pi$, so 
$\rho = (|\rho|,\pi_{\rho})$. If $\sigma,\tau$ are finite dimensional 
representations of $\SG_{n}$, $\SG_{m}$, then $\sigma*\tau$ is defined as 
$\on{Ind}_{\SG_{n}\times\SG_{m}}^{\SG_{n+m}}(\sigma\otimes\tau)$; 
we have then an identification $(\sigma*\sigma')*\sigma'' 
\simeq \sigma*(\sigma'*\sigma'')$. The dual representation of $\rho$ is 
denoted $\rho^{*}$. 

\subsubsection{The category $\on{Sch}$} \label{sect:Delta}

Define the Schur category $\on{Sch}$ as follows. $\on{Ob}(\on{Sch})
:= \on{Ob(Vect)}^{(\sqcup_{n\geq 0}\wh\SG_{n})} = \{$finitely supported 
families $F = (F_{\rho})$ of finite dimensional vector spaces,  
indexed by $\rho\in \sqcup_{n\geq 0}\wh\SG_{n}\}$. 
For $F=(F_{\rho})$
and $G= (G_{\rho})$ in $\on{Ob(Sch)}$, we set $\on{Sch}(F,G):= 
\oplus_{\rho}\on{Vect}(F_{\rho},G_{\rho})$; $F\oplus G = 
(F_{\rho}\oplus G_{\rho})$; $F^{*} := (F^{*}_{\rho^{*}})$; 
and $(F\otimes G)_{\rho}:= \oplus_{\rho',\rho''\in\sqcup_{n\geq 0}
\wh\SG_{n}} F_{\rho'}\otimes G_{\rho''}\otimes 
\mu_{\rho'\rho''}^{\rho}$, where
for $\rho,\rho',\rho''\in \sqcup_{n\geq 0}\wh\SG_{n}$, we set 
$\mu_{\rho'\rho''}^{\rho}= \on{Hom}_{\SG_{|\rho|}}
(\pi_{\rho'}*\pi_{\rho''},\pi_{\rho})$ if $|\rho| = |\rho'|+|\rho''|$
and $0$ otherwise. 
The direct sum, and the involution $\rho\mapsto \rho^{*}$ followed by the 
transposition induce canonical maps
$\on{Sch}(F,G)\oplus \on{Sch}(F',G')\to \on{Sch}(F\oplus F',G\oplus G')$ 
and $\on{Sch}(F,G)\to \on{Sch}(G^{*},F^{*})$, and 
for $f=(f_{\rho})\in \on{Sch}(F,F')$ and $g=
(g_{\rho})\in \on{Sch}(G,G')$, we define $(f\otimes g)_{\rho}:=
\oplus_{\rho',\rho''\in\sqcup_{n\geq 0}\wh\SG_{n}}
 f_{\rho'}\otimes g_{\rho''}\otimes \on{id}_{\mu_{\rho'\rho''}^{\rho}}$. 

Then $\on{Sch}$ is a symmetric additive strict monoidal category with an 
anti-automorphism\footnote{An anti-automorphism of a category $\cC$
is the data of a permutation $X\mapsto X^{*}$ of $\on{Ob}(\cC)$, and of maps 
$\cC(X,Y)\to \cC(Y^{*},X^{*})$, $x\mapsto x^{*}$, such that $(y\circ x)^{*}=
x^{*} \circ y^{*}$; if $\cC$ is additive, we require compatibility with direct
sums and the linear structure of the $\cC(X,Y)$; if $\cC$ is monoidal, we require
$(X\otimes Y)^{*}=X^{*}\otimes Y^{*}$, ${\bf 1}^{*}={\bf 1}$ 
and $(x\otimes y)^{*} = x^{*}\otimes y^{*}$.}; it is also Karoubian
(i.e., every projector has a kernel and a cokernel). 
We have a canonical 
bijection $\on{Irr(Sch)}\simeq \sqcup_{n\geq 0}\wh\SG_{n}$, with inverse 
given by $\rho\mapsto Z_{\rho}$, where $(Z_{\rho})_{\rho'} = \kk$ if
$\rho'=\rho$ and $0$ otherwise. We denote by 
${\bf 1}$, ${\bf id}$, $S^{n}$, $\wedge^{n}$ the elements of $\on{Irr(Sch)}$
corresponding to the elements of $\wh\SG_{0}$, $\wh\SG_{1}$, the trivial and the 
signature characters of $\SG_{n}$; ${\bf 1}$ is the unit object of $\on{Sch}$. 
$\on{Sch}$ has the following universal property: if 
${\mathcal C}$ is a Karoubian additive symmetric strict monoidal category with a 
distinguished object $M$, 
then there exists a unique tensor functor $F_{(\cC,M)}:\on{Sch}\to \cC$ such that 
$F({\bf id}) = M$. In particular, for $G\in\on{Ob(Sch)}$, we get an endofunctor
$F_{(\on{Sch},G)}:\on{Sch}\to \on{Sch}$, which we denote by $F\mapsto F\circ G$
(or $F(G)$)
at the level of objects and $f\mapsto f\circ G$ (or $f(G)$) at the level of morphisms. 
We say that $F=(F_{\rho})\in\on{Ob(Sch)}$ is homogeneous of degree $n$ iff 
$F_{\rho}=0$ for $|\rho|\neq n$. If $F$ is a homogeneous 
Schur functor, we denote by $|F|$ its degree. 

Let $\on{End(Vect)}$ be the symmetric additive strict monoidal 
category where objects are endofunctors $F:\on{Vect}
\to \on{Vect}$, and morphisms $F\to G$ are natural transformations, i.e., assignments 
$\on{Vect}\ni V\mapsto f_{V}\in \Vect(F(V),G(V))$, such that $f_{W}\circ 
F(\phi) = G(\phi)\circ f_{V}$ for $\phi\in\on{Vect}(V,W)$. We define
direct sums in $\on{End(Vect)}$ by $(F\oplus F')(V) := F(V)\oplus F'(V)$
and  $(V\mapsto f_{V})\oplus (V\mapsto f'_{V}):= 
(V\mapsto f_{V}\oplus f'_{V})$. We define a tensor 
product in $\on{End(Vect)}$ by $(F\otimes F')(V):= F(V)\otimes F'(V)$ 
and  $(V\mapsto f_{V})\otimes (V\mapsto f'_{V}):= 
(V\mapsto f_{V}\otimes f'_{V})$. We define an 
anti-automorphism
of $\on{End(Vect)}$ by $F^{*}(V):= F(V^{*})^{*}$ 
and $(V\mapsto f_{V})^{*}:= (V\mapsto f_{V^{*}}^{*})$, where
$(-)^{*}$ is the transposed endomorphism. Each $G\in \on{End(Vect)}$
gives rise to an endomorphism of $\on{End(Sch)}$, $F\mapsto F\circ G$, 
where $F\circ G(V):= F(G(V))$. 

We then have a tensor functor $\on{Sch}\to \on{End(Vect)}$, compatible with 
the anti-automorphisms and with the endomorphisms
$F\mapsto F\circ G$, defined at the level of objects by $F = (F_{\rho})
\mapsto (V\mapsto \oplus_{\rho\in\sqcup_{n\geq 0}\wh\SG_{n}} F_{\rho}
\otimes Z_{\rho}(V)$, where $Z_{\rho}(V) := 
\on{Hom}_{\SG_{|\rho|}}(\pi_{\rho},V^{\otimes|\rho|})$,  
and at the level 
of morphisms by $f = (f_{\rho}) \mapsto (V\mapsto f_{V})$, where 
$f_{V} = \oplus_{\rho\in\sqcup_{n\geq 0}\wh\SG_{n}} f_{\rho}\otimes 
\on{id}_{Z_{\rho}(V)}$. 

For later use, we define\footnote{For $I$ a finite set, 
we define $\on{Ob}(\on{Sch}_{I})$ similarly, where $(\rho_{1},..,
\rho_{K})$ is replaced by a map $I\to\sqcup_{n\geq 0}\wh\SG_{n}$.}
 the set $\on{Ob}(\on{Sch}_{k}) := 
\on{Ob(Vect)}^{((\sqcup_{n\geq 0}\wh\SG_{n})^{k})} = 
\{$finitely supported families $F=(F_{\rho_{1},...,\rho_{k}})$ of finite 
dimensional vector spaces, indexed by $(\rho_{1},...,\rho_{k})
\in (\sqcup_{n\geq 0}\wh\SG_{n})^{k}\}$. 
The direct sums and duality are defined component-wise as before.
Note that $\on{Ob}(\on{Sch}_{0}) = \on{Ob}(\on{Vect})$. 
If\footnote{We set $[n]:= \{1,...,n\}$.} $\phi:[m]\to [n]$ is a 
partially defined map (often identified with the collection of 
preimages $\phi^{-1}(1),...,\phi^{-1}(n)$), we define an additive map 
$\Delta^{\phi}:\on{Ob}(\on{Sch}_{n}) \to \on{Ob}(\on{Sch}_{m})$, 
taking $(F_{\rho_{1},...,\rho_{n}})$ to 
$(\Delta^{\phi}(F)_{\pi_{1},...,\pi_{m}})$, where 
$\Delta^{\phi}(F)_{\pi_{1},...,\pi_{m}} = \oplus_{\rho_{1},...,\rho_{m}}
(\mu^{\rho_{1}}_{\pi_{i},i\in\phi^{-1}(1)})^{*}\otimes ... \otimes 
(\mu^{\rho_{n}}_{\pi_{i},i\in\phi^{-1}(n)})^{*})\otimes 
F_{\rho_{1}...\rho_{n}}$. 
We set $\Delta:= \Delta^{\{1,2\}} : \on{Ob(Sch)} \to \on{Ob}(\on{Sch}_{2})$. 
Let $\on{Fun}(\on{Vect}^{n},\on{Vect})$ be the set of 
functors $\on{Vect}^{n}\to \on{Vect}$; the direct sum is defined by 
$(F\oplus G)(V_{1},...,V_{n}):= F(V_{1},...,V_{n})
\oplus G(V_{1},...,V_{n})$ and the duality by 
$F^{*}(V_{1},...,V_{n}):= F(V_{1}^{*},...,V_{n}^{*})^{*}$; 
we also define $\Delta^{\phi}:\on{Fun}(\on{Vect}^{n},\on{Vect})
\to \on{Fun}(\on{Vect}^{m},\on{Vect})$ by 
$(\Delta^{\phi}F)(V_{1},...,V_{m}):= 
F(\oplus_{i\in \phi^{-1}(1)}V_{i},...,\oplus_{i\in \phi^{-1}(n)}V_{i})$. 
Then the map $\on{Sch}_{n}\to \on{Fun}(\on{Vect}^{n},\on{Vect})$
taking $(F_{\rho_{1}...\rho_{n}})$ to $F:(V_{1},...,V_{n})\mapsto 
\oplus_{\rho_{1},...,\rho_{n}} F_{\rho_{1},...,\rho_{n}} \otimes 
Z_{\rho_{1}}(V_{1}) \otimes ... \otimes Z_{\rho_{n}}(V_{n})$ is compatible 
with the direct sums, the duality and the maps $\Delta^{\phi}$. 

\subsubsection{The category $\on{Sch}_{1+1}$}

We define the symmetric additive strict monoidal category of Schur bifunctors 
$\on{Sch}_{1+1}$ as follows. $\on{Ob}(\on{Sch}_{1+1}):= 
\on{Ob}(\on{Sch}_{2})$.  For $F,G\in\on{Ob}(\on{Sch}_{1+1})$, 
$\on{Sch}_{1+1}(F,G):= \oplus_{\rho_{1},\rho_{2}}
\on{Vect}(F_{\rho_{1},\rho_{2}},G_{\rho_{1},\rho_{2}})$. 
We also define $(F\otimes G)_{\rho_{1},\rho_{2}}:= 
\oplus_{\rho'_{i},\rho''_{i}} F_{\rho'_{1},\rho'_{2}} \otimes 
G_{\rho''_{1},\rho''_{2}}\otimes \mu_{\rho'_{1}\rho''_{1}}^{\rho_{1}}
\otimes \mu_{\rho'_{1}\rho''_{2}}^{\rho_{2}}$. 
The direct sums and tensor products of morphisms are then defined
component-wise. An anti-automorphism of $\on{Sch}_{1+1}$ is defined 
by $(F_{\rho,\sigma})^{*} = (F^{*}_{\sigma^{*},\rho^{*}})$ 
and by $((\rho,\sigma)\mapsto f_{\rho,\sigma})^{*}=
((\rho,\sigma)\mapsto f_{\sigma^{*},\rho^{*}}^{*})$. 
We define a tensor morphism 
$\boxtimes:\on{Sch}^{2}\to \on{Sch}_{1+1}$
at the level of objects by $(F\boxtimes G)_{\rho,\sigma}:= 
F_{\rho}\otimes G_{\sigma}$, and component-wise at the level of morphisms. 
We then have for $F,...,G'\in\on{Ob(Sch)}$, 
\begin{equation} \label{eq:boxtimes}
\on{Sch}_{1+1}(F\boxtimes G^{*},F'\boxtimes G^{\prime *})
\simeq \on{Sch}(F,F') \otimes \on{Sch}(G',G) 
\end{equation}
and $(F\boxtimes G)^{*} = G^{*}\boxtimes F^{*}$. 
As $\on{Sch}_{1+1}$ is Karoubian, any $G\in\on{Sch}_{1+1}$
gives rise to a unique tensor functor $\on{Sch}\to \on{Sch}_{1+1}$
taking ${\bf id}$ to $G$, which we denote by $F\mapsto F\circ G$. 

Let $\on{Fun}(\on{Vect}^{2},\on{Vect})$ be the symmetric 
additive strict monoidal category 
where objects are functors $\on{Vect}^{2}\to \on{Vect}$ and morphisms
are natural transformations; the direct sum is defined by 
$(F\oplus G)(V,W):= F(V,W)\oplus G(V,W)$, the tensor products by 
$(F\otimes G)(V,W):= F(V,W)\otimes G(V,W)$ and the duality by 
$F^{*}(V,W):= F(W^{*},V^{*})^{*}$. Then we have a tensor 
functor $\on{Sch}_{1+1}\to \on{Fun}(\on{Vect}^{2},\on{Vect})$
taking $F$ to $((V,W)\mapsto \oplus_{\rho_{1},\rho_{2}}
F_{\rho_{1},\rho_{2}}\otimes Z_{\rho_{1}}(V) \otimes 
Z_{\rho_{2}}(W))$, compatible with the dualities.  
It is also compatible with the tensor functor $\on{End(Sch)}\to 
\on{Fun}(\on{Vect}^{2},\on{Vect})$, $G\mapsto F\circ G$, 
where $F\circ G(V,W) = F(G(V,W))$. 
We define a tensor functor $\boxtimes : \on{End}(\on{Vect})^{2}
\to \on{Fun}(\on{Vect}^{2},\on{Vect})$ at the level of objects by 
$(F\boxtimes G)(V,W):= F(V)\otimes G(W)$. Then the morphisms 
$\on{Sch}\to \on{End}(\on{Vect})$, $\on{Sch}_{1+1}\to 
\on{Fun}(\on{Vect}^{2},\on{Vect})$
intertwine the functors $\boxtimes : \on{Sch}^{2}\to \on{Sch}_{1+1}$ and 
$\boxtimes : \on{End}(\on{Vect})^{2}\to \on{Fun}(\on{Vect}^{2},\on{Vect})$. 

\subsubsection{The category $\on{Sch}_{(1)}$}

We now define the additive symmetric strict monoidal category $\on{Sch}_{(1)}$ as follows. 
We set $\on{Ob}(\on{Sch}_{(1)}):= \on{Ob(Vect)}^{(\sqcup_{k\geq 0}
(\sqcup_{n\geq 0}\SG_{n})^{k})} = \prod'_{k\geq 0} \on{Ob}(\on{Sch}_{k})
= \{$finitely supported collections $(F_{k})_{k\geq 0}$, where $F_{k}
\in\on{Ob}(\on{Sch}_{k})$ is a family 
$F_{k}=(F_{\rho_{1},...,\rho_{k}})\}$. 
The direct sum of objects is defined by component-wise addition. The
tensor product of objects is defined by $(F_{k})\boxtimes (G_{k}) := 
((F\boxtimes G)_{k})$, where $(F\boxtimes G)_{k} := \oplus_{k',k''|
k'+k'' = k} F_{k'}\boxtimes G_{k''}$, and if $F_{k'} = (F_{\rho_{1},...,\rho_{k'}})
\in \on{Ob}(\on{Sch}_{k'})$, $G_{k''} = (G_{\rho_{1},...,\rho_{k''}})
\in \on{Ob}(\on{Sch}_{k''})$, then $F_{k'}\boxtimes G_{k''} = 
((F\boxtimes G)_{\rho_{1},...,\rho_{k'+k''}})\in \on{Ob}(\on{Sch}_{k'+k''})$, 
where $(F\boxtimes G)_{\rho_{1},...,\rho_{k'+k''}} := F_{\rho_{1},...,\rho_{k'}}
\otimes G_{\rho_{k'+1},...,\rho_{k'+k''}}$. 

In order to define the morphisms, we first define a ``contraction'' map 
$c: \on{Ob}(\on{Sch}_{k})\to \on{Ob(Sch)}$, $F_{k} = (F_{\rho_{1},...,\rho_{k}})
\mapsto c(F_{k})$
by $c(F_{k})_{\rho}:= \oplus_{\rho_{1},...,\rho_{k}}
F_{\rho_{1},...,\rho_{k}} \otimes \mu_{\rho_{1},...,\rho_{k}}^{\rho}$, 
where 
$\mu_{\rho_{1},...,\rho_{k}}^{\rho}:= \on{Hom}_{\SG_{|\rho|}}
(\rho_{1}* ... * \rho_{k},\rho)$ if $\sum_{i}|\rho_{i}| = |\rho|$ 
and $0$ otherwise. For $F = (F_{k})$, we then set $c(F):= \oplus_{k} c(F_{k})$
and for $F,G\in\on{Ob}(\on{Sch}_{(1)})$, we set $\on{Sch}_{(1)}(F,G)
:= \on{Sch}(c(F),c(G))$.
We define the direct sum and the tensor product of morphisms using the identifications 
$c(F\oplus G)\simeq c(F)\oplus c(G)$ and 
$c(F\boxtimes G) \simeq c( F)\otimes c(G)$.
The symmetry constraint in $\on{Sch}_{(1)}(F\boxtimes G,G
\boxtimes F) = 
\on{Sch}(c(F\boxtimes G),c(G\boxtimes F))$ is then given by 
the identifications $c(X\boxtimes Y) \simeq c(X)\otimes c(Y)$
and the symmetry constraint for $\on{Sch}$. The unit object of 
$\on{Ob}(\on{Sch}_{(1)})$ is ${\mathfrak 1}$, whose only nonzero component
is ${\mathfrak 1}_{0} = \kk\in \on{Ob}(\on{Sch}_{0}) = \on{Ob(Vect)}$. 

We define an additive symmetric strict monoidal category 
$\prod'_{k\geq 0}\on{Fun}(\on{Vect}^{k},\on{Vect})$
as follows. The objects are finitely supported families $(F_{k})_{k\geq 0}$, 
where $F_{k}\in \on{Fun}(\on{Vect}^{k},\on{Vect})$. 
The direct sum of objects is defined component-wise, where for $F_{k},G_{k}\in 
\on{Fun}(\on{Vect}^{k},\on{Vect})$, $F_{k}\oplus 
G_{k}\in \on{Fun}(\on{Vect}^{k},\on{Vect})$ is given by 
$(F_{k}\oplus G_{k})(V_{1},...,V_{k}):= 
F_{k}(V_{1},...,V_{k})\oplus G_{k}(V_{1},...,V_{k})$. The tensor 
product of objects is $(F_{k})\boxtimes (G_{k}):= ((F\boxtimes G)_{k})$, 
where $(F\boxtimes G)_{k}:= \oplus_{k'+k''=k} F_{k'}\boxtimes_{k',k''} G_{k''}$, 
and $\boxtimes_{k',k''} : \on{Fun}(\on{Vect}^{k'},\on{Vect})
\times \on{Fun}(\on{Vect}^{k''},\on{Vect}) \to 
\on{Fun}(\on{Vect}^{k'+k''},\on{Vect})$ is given by 
$(F_{k'}\boxtimes_{k',k''}G_{k''})(V_{1},...,V_{k'+k''}):= 
F_{k'}(V_{1},...,V_{k'})\otimes G_{k''}(V_{k'+1},...,V_{k'+k''})$. 
The contraction $c:\prod'_{k\geq 0}\on{Fun}(\on{Vect}^{k},\on{Vect})
\to \on{End}(\on{Fun})$ is defined by $c((F_{k})):= \oplus_{k\geq 0}c(F_{k})$, 
where $c(F_{k})(V):= F_{k}(V,...,V)$. The space of morphisms $F\to G$
is then defined as $\on{Fun}(\on{Vect})(c(F),c(G))$. 

There is a unique tensor morphism $\on{Sch}_{(1)}\to 
\prod'_{k\geq 0}\on{Fun}(\on{Vect}^{k},\on{Vect})$, taking 
$F=(F_{k})_{k\geq 0}$, where $F_{k}=(F_{\rho_{1},...,\rho_{k}})$,
to the collection $(\tilde F_{k})_{k\geq 0}$, where 
$\tilde F_{k}(V_{1},...,V_{k}):= \oplus_{\rho_{1},...,\rho_{k}}
F_{\rho_{1},...,\rho_{k}}\otimes Z_{\rho_{1}}(V_{1})\otimes ... 
\otimes Z_{\rho_{k}}(V_{k})$. 

\subsubsection{The category $\on{Sch}_{(1+1)}$}

We now define an additive symmetric strict monoidal category $\on{Sch}_{(1+1)}$
as follows. $\on{Ob}(\on{Sch}_{(1+1)}):= \on{Ob}(\on{Vect})^{(\sqcup_{k,l\geq 0}
(\sqcup_{n\geq 0}\wh\SG_{n})^{k+l})} = \{$finitely supported 
collections $(F_{k,l})$, where 
$F_{k,l} = (F_{\rho_{1},...,\rho_{k}; \sigma_{1},...,\sigma_{l}})\in 
\on{Ob}(\on{Sch}_{k+l})\}$. The direct sum of objects is defined component-wise, 
and the tensor product is given by $(F\boxtimes F')_{k,l}:= 
\oplus_{(k_{1},l_{2})+(k_{2},l_{2})=(k,l)} 
F_{k_{1},l_{1}}\boxtimes_{k_{1},l_{1},k_{2},l_{2}} F'_{k_{2},l_{2}}$
if $F= (F_{k,l})$ and $F' = (F'_{k,l})$, where 
$\boxtimes_{k,l,k',l'} : \on{Ob}(\on{Sch}_{k+l})
\times \on{Ob}(\on{Sch}_{k'+l'}) \to \on{Ob}(\on{Sch}_{k+k'+l+l'})$ is given
by \\
\noindent $(F\boxtimes_{k,l,k',l'} F')_{\rho_{1},...,\rho_{k+k'};
\sigma_{1},...,\sigma_{l+l'}}:= F_{\rho_{1},...,\rho_{k};\sigma_{1},...,\sigma_{l}}
\otimes F'_{\rho_{k+1},...,\rho_{k+k'};\sigma_{l+1},...,\sigma_{l+l'}}$. 
An involution is defined by 
$F^{*} = ((k,l,\rho_{1},...,\rho_{k},\sigma_{1},...,\sigma_{l})\mapsto
F(l,k,\sigma_{1}^{*},...,\sigma_{l}^{*},\rho_{1}^{*},...,\rho_{l}^{*})^{*})$
for $F = ((k,l,\rho_{1},...,\rho_{k},\sigma_{1},...,\sigma_{l})\mapsto 
F(k,l,\rho_{1},...,\rho_{k},\sigma_{1},...,\sigma_{l}))$. 

In order to define the morphisms, we define a map $c :\on{Ob}(\on{Sch}_{(1+1)})
\to \on{Ob}(\on{Sch}_{1+1})$, $F\mapsto c(F)$, by 
$c(F) := \oplus_{k,l} c(F_{k,l})$ for $F=(F_{k,l})$, and if
$F_{k,l} = (F_{\rho_{1},...,\rho_{k};\sigma_{1},...,\sigma_{l}})$, then 
$c(F_{k,l})_{\rho,\sigma} = \oplus_{\rho_{1},...,\rho_{k};
\sigma_{1},...,\sigma_{k}} F_{\rho_{1},...,\rho_{k};\sigma_{1},...,\sigma_{l}}
\otimes \mu_{\rho_{1}...\rho_{k}}^{\rho} \otimes 
\mu_{\sigma_{1}...\sigma_{l}}^{\sigma}$. 
We then set $\on{Sch}_{(1+1)}(F,G):= \on{Sch}_{1+1}(c(F),c(G))$. 
The direct sum, tensor product and duality of morphisms are then induced by those
of $\on{Sch}_{1+1}$ and the identifications $c(F\oplus G)\simeq c(F)\oplus c(G)$, 
$c(F\boxtimes G)\simeq c(F)\boxtimes c(G)$, $c(F^{*}) = c(F)^{*}$. 

We define a tensor morphism $\ul\boxtimes:(\on{Sch}_{(1)})^{2}\to 
\on{Sch}_{(1+1)}$; at the level of objects, it is defined by $(F_{k})\ul\boxtimes
(G_{l}):= (F_{k}\boxtimes_{k,l} G_{l})$; at the level of morphisms, it is induced 
by the tensor morphism $\on{Sch}^{2}\to \on{Sch}_{1+1}$. The unit object of 
$\on{Sch}_{(1+1)}$ is ${\mathfrak 1}\ul\boxtimes{\mathfrak 1}$. We then 
have for $F,...,G'\in\on{Ob}(\on{Sch}_{(1)})$
$$\on{Sch}_{(1+1)}(F\ul\boxtimes G^{*},F'\ul\boxtimes G^{\prime*})
= \on{Sch}_{(1)}(F,F')\otimes \on{Sch}_{(1)}(G',G) 
$$
and $(F\ul\boxtimes G)^{*} = G^{*}\ul\boxtimes F^{*}$. 

As before, we define an additive symmetric strict monoidal 
category $\prod'_{k,l}\on{Fun}(\on{Vect}^{k+l},\on{Vect})$; 
objects are fiinitely supported families $(F_{k,l})_{k,l\geq 0}$; 
$(F_{k,l}\oplus G_{k,l})(V_{1},...,V_{k};W_{1},...,W_{l}):= 
F_{k,l}(V_{1},...,W_{l})\oplus G_{k,l}(V_{1},...,W_{l})$; 
$\boxtimes_{k',l',k'',l''}:\on{Fun}(\on{Vect}^{k'+l'},\on{Vect})
\times \on{Fun}(\on{Vect}^{k''+l''},\on{Vect})\to 
\on{Fun}(\on{Vect}^{k'+l'+k''+l''},
\on{Vect})$ is $(F\boxtimes_{k',l',k'',l''}G)(V_{1},...,W_{l'+l''}):= 
F(V_{1},...,W_{l'})\otimes G(V_{k'+1},...,W_{l'+l''})$. We define 
\\ \noindent $c:\prod'_{k,l}\on{Fun}(\on{Vect}^{k+l},\on{Vect})\to 
\on{Fun}(\on{Vect}^{2},\on{Vect})$ by $c(F) = \oplus_{k,l}c(F_{k,l})$
and $c(F_{k,l})(V,W):= F_{k,l}(V,...,V;W,...,W)$ and the space of 
morphisms $F\to G$ as $\on{Fun}(\on{Vect}^{2},\on{Vect})(c(F),c(G))$. 
We also define a tensor morphism $\underline\boxtimes:
(\prod'_{k\geq 0}\on{Fun}(\on{Vect}^{k},
\on{Vect}))^{2}\to \prod'_{k,l\geq 0} \on{Fun}(\on{Vect}^{k+l},\on{Vect})$
at the level of objects by $(F_{k})\ul\boxtimes(G_{k}):= 
(F_{k}\boxtimes_{k,l} G_{l})$. 

Then we have a tensor morphism $\on{Sch}_{(1+1)}\to \prod'_{k,l\geq 0}
\on{Fun}(\on{Vect}^{k+l},\on{Vect})$, taking $(F_{k,l})$ to $(\tilde F_{k,l})$, 
where for $F_{k,l}=(F_{\rho_{1},...,\sigma_{l}})$, $\tilde F_{k,l}(V_{1},...,W_{l})
:= \oplus_{\rho_{1},...,\sigma_{l}} F_{\rho_{1},...,\sigma_{l}}\otimes
\rho_{1}(V_{1})\otimes ... \otimes \sigma_{l}(W_{l})$. This morphism 
is compatible with the morphism $\on{Sch}_{(1)}\to \prod'_{k\geq 0}
\on{Fun}(\on{Vect}^{k},\on{Vect})$ and the morphisms 
$(\on{Sch}_{(1)})^{2}\to \on{Sch}_{(1+1)}$, 
$(\prod'_{k\geq 0}\on{Fun}(\on{Vect}^{k},\on{Vect}))^{2} \to 
\prod'_{k,l\geq 0}\on{Fun}(\on{Vect}^{k+l},\on{Vect})$. 

\subsubsection{Completions}

If in the definition of $\on{Sch}$, we forget the condition that $(F_{\rho})$ is 
finitely supported, we get a symmetric additive strict monoidal category with duality
${\bf Sch}$. Infinite sums of objects of increasing degrees are
defined in ${\bf Sch}$. 
For each $G = \oplus_{i\geq 1} G_{i}\in {\bf Sch}$ ($|G_{i}|=i$), we have an 
endofunctor of ${\bf Sch}$, $F\mapsto F\circ G$, $f\mapsto f\circ G$
(also written $F(G),f(G)$). 
We also 
define $\on{Ob}({\bf Sch}_{k})$ by dropping the finite support condition. 
The maps $\Delta^{\phi}$ extend to these sets. 
We define ${\bf Sch}_{1+1}$, ${\bf Sch}_{(1)}$ and ${\bf Sch}_{(1+1)}$ similarly to 
$\on{Sch}_{1+1}$, $\on{Sch}_{(1)}$, $\on{Sch}_{(1+1)}$, namely 
$\on{Ob}({\bf Sch}_{1+1}) = \on{Ob}({\bf Sch}_{2})$, 
$\on{Ob}({\bf Sch}_{(1)}) = \{$finitely supported families 
$(F_{k})_{k\geq 0}$, where $F_{k}\in\on{Ob}({\bf Sch}_{k})$, 
and $\on{Ob}({\bf Sch}_{(1+1)}) = \{$finitely supported 
families $(F_{k,l})_{k,l\geq 0}$, where $F_{k,l}\in
\on{Ob}({\bf Sch}_{k+l})\}$. Then $\on{Sch}_{1+1}$, 
$\on{Sch}_{(1)}$ and $\on{Sch}_{(1+1)}$ have structures of additive
symmetric strict monoidal categories. The map $c$ and the bifunctor $\ul\boxtimes$ 
extend to these categories; the duality extends to ${\bf Sch}_{1+1}$ and 
${\bf Sch}_{(1+1)}$.

\medskip 

\noindent{\bf Examples.} Let $T_{n}\in\on{Ob(Sch)}$ be such that 
$(T_{n})_{\rho'}=\pi_{\rho'}$ if $|\rho'|=n$ and $0$ otherwise. 
The corresponding endofunctor of $\on{Vect}$ is 
$V\mapsto T_{n}(V) = V^{\otimes n} = \oplus_{\rho\in\wh\SG_{n}}\pi_{\rho}
\otimes Z_{\rho}(V)$. Using the obvious module category structure of 
$\on{Sch}$ over $\on{Vect}$, we write 
\begin{equation} \label{schur:recip}
T_{n} = \oplus_{Z\in\on{Irr(Sch)},|Z|=n} \pi_{Z}\otimes Z, 
\end{equation}
where $Z\mapsto (|Z|,\pi_{Z})$ is the inverse to $\sqcup_{n\geq 0}\wh\SG_{n}
\to \on{Irr(Sch)}$, $\rho\mapsto Z_{\rho}$. 

The endofunctors of $\on{Vect}$ corresponding to $S^{n}$ and $\wedge^{n}$
are the $n$th symmetric and exterior powers functors. 
The symmetric and exterior algebra functors 
$S := \oplus_{n\geq 0} S^n$ and $\wedge := \oplus_{n\geq 0} \wedge^n$
are objects in ${\bf Sch}$. We then have $\Delta(S) = S \boxtimes S$, 
$\Delta(\wedge) = \wedge \boxtimes \wedge$. Note that while the map 
$\on{Ob(Sch)}\to \on{Ob(End(Vect))}$ is injective, it is not 
surjective, e.g. the exterior algebra functor is not in the image of this map. 
 

\begin{remark} For any $F,G\in \on{Ob(Sch)}$, we have 
\begin{equation} \label{eq:sch}
\on{Sch}(F,G) = \oplus_{Z\in \on{Irr}(\on{Sch})}
\on{Sch}(F,Z) \otimes \on{Sch}(Z,G);  
\end{equation}
for any $B,B'\in\on{Ob}(\on{Sch}_{1+1})$, we have 
$$ 
\on{Sch}_{1+1}(B,B') = \oplus_{Z,Z'\in\on{Irr}(\on{Sch})}
\on{Sch}_{1+1}(B,Z\boxtimes Z') \otimes \on{Sch}_{1+1}(Z\boxtimes Z',B'). 
$$ 
\end{remark}

\begin{remark} We have 
$\on{Irr}(\on{Sch}_{(1)}) = \{Z_1\boxtimes ... \boxtimes 
Z_k | k\geq 0, Z_1,...,Z_k\in\on{Irr(Sch)}\}$; and 
$\on{Irr}(\on{Sch}_{(2)}) = \{(Z_1\boxtimes ... \boxtimes 
Z_k) \underline\boxtimes (W_1\boxtimes...\boxtimes W_\ell)| k,\ell\geq 0, 
Z_1,...,W_\ell\in\on{Irr(Sch)}\}$.  Then $c(Z_{1}\boxtimes...\boxtimes Z_{k})
= Z_{1}\otimes...\otimes Z_{k}$ and $c((Z_{1}\boxtimes ... \boxtimes 
Z_{k})\ul\boxtimes(W_{1}\boxtimes ... \boxtimes W_{l})) = 
(Z_{1}\otimes ... \otimes Z_{k})\boxtimes (W_{1}\otimes ... \otimes W_{l})$. 
\end{remark}

\subsection{Quasi-categories}

We define a quasi-category $\cC$ to be the data of: (a) a set of objects
$\on{Ob}(\cC)$; (b) for any $X,Y\in \on{Ob}(\cC)$, a set of morphisms
$\cC(X,Y)$, and for any $X\in \on{Ob}(\cC)$, an element 
$\on{id}_{X}\in \cC(X,X)$; (c) for $X_{i}\in \on{Ob}(\cC)$
($i=1,2,3$), a subset $\cC(X_{1},X_{2},X_{3})\subset \cC(X_{1},X_{2})
\times \cC(X_{2},X_{3})$ and a map $\cC(X_{1},X_{2},X_{3})
\stackrel{\circ}{\to}\cC(X_{1},X_{3})$, $(x_{1},x_{2})\mapsto 
x_{2}\circ x_{1}$, such that:  

(identity axiom) if $X,Y\in\on{Ob}(\cC)$, and $x\in \cC(X,Y)$, then
$\on{id}_{Y}\circ x\in \cC(X,Y,Y)$, $x\circ \on{id}_{X}\in 
\cC(X,X,Y)$, and $\on{id}_{Y}\circ x = x \circ \on{id}_{X}=x$; 

(associativity axiom) if $X_{i}\in \on{Ob}(\cC)$ ($i=1,...,4$) and 
$x_{i}\in \cC(X_{i},X_{i+1})$ ($i=1,2,3$), then if: $(x_{1},x_{2})
\in \cC(X_{1},X_{2},X_{3})$, $(x_{2}\circ x_{1},x_{3})\in 
\cC(X_{1},X_{3},X_{4})$, $(x_{2},x_{3})\in \cC(X_{2},X_{3},X_{4})$
and $(x_{1},x_{3}\circ x_{2})\in \cC(X_{1},X_{2},X_{4})$, then 
$x_{3}\circ (x_{2}\circ x_{1}) = (x_{3}\circ x_{2})\circ x_{1}$. 

We then define inductively a diagram $\cC(X_{1},X_{2})\times...\times
\cC(X_{n-1},X_{n}) \supset \cC(X_{1},...,X_{n})\stackrel{\circ}{\to}
\cC(X_{1},X_{n})$, as follows: $(x_{1},...,x_{n-1})\in \cC(X_{1},...,X_{n})$
iff for any $k=2,...,n-1$, $(x_{1},...,x_{k-1})\in \cC(X_{1},...,X_{k})$, 
$(x_{k},...,x_{n-1})\in \cC(X_{k},...,X_{n})$, and 
$(x_{k-1}\circ...\circ x_{1},x_{n-1}\circ...\circ x_{k})\in 
\cC(X_{1},X_{k},X_{n})$; if $(x_{1},...,x_{n-1})$ satisfies these conditions, 
then the $(x_{n-1}\circ...\circ x_{k})\circ (x_{k-1}\circ...\circ x_{1})$
all coincide; this defines the map 
$\cC(X_{1},...,X_{n})\to \cC(X_{1},X_{n})$. 

If $1<n_{1}<...<n_{k}<n$ and $x = (x_{1},...,x_{n-1})
\in \cC(X_{1},X_{2})\times...
\times \cC(X_{n-1},X_{n})$, then $x\in \cC(X_{1},...,X_{n})$ iff: (a)
$(x_{1},...,x_{n_{1}-1})\in \cC(X_{1},...,X_{n_{1}})$, 
$(x_{n_{1}},...,x_{n_{2}-1})\in \cC(X_{n_{1}},...,X_{n_{2}})$, ..., and 
$(x_{n_{k-1}},...,x_{n-1})\in \cC(X_{n_{k-1}},...,X_{n})$;  
(b) moreover, $(x_{n_{1}-1}\circ ... \circ x_{1},...,x_{n-1}\circ ... 
\circ x_{n_{k-1}})\in \cC(X_{1},X_{n_{1}},X_{n_{2}},...,X_{n})$.
If these conditions are satisfied, then 
$x_{n-1}\circ ... \circ x_{1} = (x_{n-1}\circ ... \circ x_{n_{k-1}+1})
\circ ... \circ (x_{n_{1}-1}\circ ... \circ x_{1})$. 

The quasi-category $\cC$ is called strict monoidal if it is equipped
with: (a) a map $\otimes : \on{Ob}(\cC)^{2}\to \cC$, 
$(X,Y)\mapsto X\otimes Y$ and an object ${\bf 1}\in 
\on{Ob}(\cC)$, such that $(X\otimes Y)\otimes Z = 
X\otimes (Y\otimes Z)$, ${\bf 1}\otimes X = X \otimes {\bf 1} = X$; 
(b) a map $\otimes : \cC(X,Y)\times \cC(X',Y')\to \cC(X\otimes X',
Y\otimes Y')$ such that $(f\otimes f')\otimes f'' = f\otimes
(f'\otimes f'')$, $f\otimes \on{id}_{{\bf 1}} 
= \on{id}_{{\bf 1}}\otimes f = f$; (c) a map 
$\otimes : \cC(X_{1},X_{2},X_{3}) \times \cC(X'_{1},X'_{2},X'_{3})
\to \cC(X_{1}\otimes X'_{1},X_{2}\otimes X'_{3},X_{3}\otimes X'_{3})$, 
such that 
$$\begin{matrix}
\cC(X_{1},X_{2},X_{3}) \times \cC(X'_{1},X'_{3},X'_{3})
&\to &\cC(X_{1}\otimes X'_{1},X_{2}\otimes X'_{2},X_{3}\otimes X'_{3}) \\
 \downarrow & & \downarrow \\
 \cC(X_{1},X_{2})\times\cC(X_{2},X_{3})\times
 \cC(X'_{1},X'_{2})\times\cC(X'_{2},X'_{3})
& \to &  
\cC(X_{1}\otimes X'_{1},X_{2}\otimes X'_{2})
\times\cC(X_{2}\otimes X'_{2},X_{3}\otimes X'_{3})
\end{matrix}
$$
commutes. Then we have maps $\cC(X_{1},...,X_{n})\otimes 
\cC(X'_{1},...,X'_{n})\to \cC(X_{1}\otimes X'_{1},...,X_{n}\otimes X'_{n})$, 
such that the analogous diagram (with $3$ replaced by $n$) commutes. 

\medskip  
\noindent{\bf Example.} $\cG$ is the category where objects are pairs
$(I,J)$ of finite sets and $\cG((I,J),(I',J'))$ is the set of oriented acyclic 
graphs with vertices $i_{in},j_{in},i'_{out},j'_{out}$, $i\in I$, $j\in J$, 
$i'\in I'$, $j'\in J'$, where each edge has its origin in $\{i_{in},j'_{out}| 
i\in I, j'\in J'\}$ and its end in $\{i'_{out},j_{in}|i'\in I', j\in J\}$, 
and there is at most one edge through two given vertices. 
Equivalently, a graph is a subset of $(I\sqcup J')\times (I'\sqcup J)$. 
If $X_{\alpha}=(I_{\alpha},J_{\alpha})$ and $x_{\alpha}\in \cG(X_{\alpha},
X_{\alpha+1})$
($\alpha=1,...,k-1$), we obtain a composed graph with edges 
$x_{in}$, $y_{out}$, $x\in I_{1}\sqcup J_{1}$, $y\in I_{n}\sqcup J_{n}$, 
by declaring that 
two edges are connected if there exists an oriented path in the juxtaposition of 
$x_{1},...,x_{k-1}$ relating them. Then $\cG(X_{1},...,X_{k})\subset 
\cG(X_{1},X_{2})\times ... \times \cG(X_{k-1},X_{k})$ is the set of 
tuples of graphs whose composed graph is acyclic, which is then their composition. 
The tensor product is 
given by $(I,J)\otimes (I',J'):= (I\sqcup J,I'\sqcup J')$ at the level of 
objects, and by the disjoint union of graphs at the level of morphisms.
Note that $\cG$ contains subcategories $\cG^{\on{left}}$ and $\cG^{\on{right}}$, 
where $\cG^{\on{left}}((I,J),(I',J')) = \{S\in\cG((I,J),(I',J')) | S\cap (J'\times I')
=\emptyset\}$, and $\cG^{\on{right}}((I,J),(I',J')) = \{S|S\cap (I\times J) = 
\emptyset\}$.  \hfill \qed \medskip 

A $\kk$-additive quasi-category $\cC$ is the data of: (a) a set of objects 
$\on{Ob}(\cC)$, (b) for any $X,Y\in\on{Ob}(\cC)$, a vector space 
$\cC(X,Y)$, and for any $X_{1},...,X_{n}\in \on{Ob}(\cC)$, 
a vector subspace $\cC(X_{1},...,X_{n})\subset \cC(X_{1},X_{2})
\otimes ... \otimes \cC(X_{n-1},X_{n})$, and a linear map 
$\cC(X_{1},...,X_{n})\to \cC(X_{1},X_{n})$, satisfying the 
axioms of a quasi-category (with products replaced by tensor products); 
(c) an associative direct sum map 
$\oplus : \on{Ob}(\cC)^{2}\to \on{Ob}(\cC)$, 
$(X,Y)\mapsto X\oplus Y$, an object ${\bf 0}\in\on{Ob}(\cC)$, and
isomorphisms  $\cC(Z,X\oplus Y)\simeq \cC(Z,X)\oplus \cC(Z,Y)$ and 
$\cC(X\oplus Y,Z)\simeq \cC(X,Z)\oplus \cC(Y,Z)$, 
such that: 
 
$\cC(X_{1}\oplus X'_{1},X_{2},X_{3})\simeq \cC(X_{1},X_{2},X_{3})
\oplus \cC(X'_{1},X_{2},X_{3})$, $\cC(X_{1},X_{2},X_{3}\oplus X'_{3})
\simeq \cC(X_{1},X_{2},X_{3})\oplus\cC(X_{1},X_{2},X'_{3})$, 
$\cC(X_{1},X_{2}\oplus X'_{2},X'_{3}) \simeq 
\cC(X_{1},X_{2},X_{3})\oplus \cC(X_{1},X'_{2},X_{3})
\oplus \cC(X_{1},X_{2})\otimes \cC(X_{2},X_{3})\oplus 
\cC(X_{1},X'_{2})\otimes \cC(X_{2},X_{3})$, and the composition 
map on left sides coincides with the sum of compositions on the right sides, 
and of the zero maps on the two last summands in the last case (this statement 
then generalizes to $\cC(X_{1},...,X_{i}\oplus X'_{i},...,X_{n})$);  

$X\oplus {\bf 0} = X = {\bf 0}\oplus X$ and 
$\cC(X,{\bf 0}) = \cC({\bf 0},X) = 0$ for any $X$, 
and the composed isomorphisms $\cC(X,Y) = \cC(X\oplus{\bf 0},Y)\simeq 
\cC(X,Y)$, $\cC(X,Y) = \cC({\bf 0}\oplus X,Y)\simeq 
\cC(X,Y)$, $\cC(X,Y) = \cC(X,Y\oplus{\bf 0})\simeq 
\cC(X,Y)$ and $\cC(X,Y) = \cC(X,{\bf 0}\oplus Y)\simeq 
\cC(X,Y)$ are the identity. 

Such a $\cC$ is called strict monoidal if it satisfies the above
axioms of a strict monoidal quasi-category, where $\otimes$ is bilinear and 
biadditive. 

A functor $F:\cC\to \cD$ between quasi-categories
is defined as the data of a 
map $F:\on{Ob}(\cC)\to \on{Ob}(\cD)$, and a collection of 
maps $F(X,Y):\cC(X,Y)\to \cD(F(X),F(Y))$, such that 
$\times_{i=1}^{n-1}F(X_{i},X_{i+1})$ restricts to a map 
$\cC(X_{1},...,X_{n})\to \cD(F(X_{1}),...,F(X_{n}))$ and the 
natural diagrams commute; natural additional axioms are imposed if
the categories are strict monoidal and/or additive. 

\medskip 
\noindent{\bf Example.} $\kk\cG$ is the category with $\on{Ob}(\kk\cG) = 
\on{Ob}(\cG)$ and $(\kk\cG)((I,J),(I',J')) = \kk \cG((I,J),(I',J'))$; then 
$\kk\cG$ is an additive strict monoidal quasi-category. \hfill \qed \medskip 

\subsection{Partial traces and quasi-categories} \label{sect:part:traces}

If $\cC_{0}$ is a symmetric strict monoidal category with symmetry 
constraint $\beta_{X,Y}\in \cC_{0}(X\otimes Y,Y\otimes X)$, 
a partial trace on $\cC_{0}$ is the data of
diagrams $\cC_{0}(X\otimes Z,Y\otimes Z)\supset \cC_{0}(X,Y|Z)
\stackrel{\on{tr_{Z}}}{\to} \cC_{0}(X,Y)$ for $X,Y,Z\in\on{Ob}(\cC)$, 
such that:  $\cC_{0}(X,Y|Z\otimes Z') \subset \cC_{0}(X\otimes Z,Y\otimes Z|Z')
\cap \on{tr}_{Z'}^{-1}(\cC_{0}(X,Y|Z))$, and $\on{tr}_{Z\otimes Z'} = \on{tr}_{Z}
\circ \on{tr}_{Z'}$; the map
$x\mapsto x':= (\on{id}_{Y}\otimes \beta_{Z',Z}) \circ 
x \circ (\on{id}_{X}\otimes \beta_{Z,Z'})$ induces an isomorphism 
$\cC_{0}(X,Y|Z\otimes Z') \to \cC_{0}(X,Y|Z'\otimes Z)$, and 
$\on{tr}_{Z'\otimes Z}(x') = \on{tr}_{Z\otimes Z'}(x)$;  
the composition takes $\cC_{0}(X,Y|T)\times 
\cC_{0}(Y,Z)$ to $\cC_{0}(X,Z|T)$, and $\on{tr}_{T}(
(y\otimes \on{id}_{T}) \circ x) = y \circ 
\on{tr}_{T}(x)$, and similarly it takes $\cC_{0}(X,Y)\times 
\cC_{0}(Y,Z|T)$ to $\cC_{0}(X,Z|T)$, and $\on{tr}_{T}(y \circ
(x\otimes \on{id}_{T})) = \on{tr}_{T}(y)\circ x$; 
the map $\times_{i=1}^{2}\cC_{0}(X_{i}\otimes T_{i},Y_{i}\otimes T_{i})
\to \cC_{0}(X_{1}\otimes X_{2}\otimes T_{1}\otimes T_{2},
Y_{1}\otimes Y_{2}\otimes T_{1}\otimes T_{2})$, 
$(x_{1}, x_{2})\mapsto 
x:= (\on{id}_{Y_{1}}
\otimes \beta_{T_{1},Y_{2}}\otimes \on{id}_{T_{2}})
\circ (x_{1}\otimes x_{2})\circ (\on{id}_{X_{1}}
\otimes \beta_{X_{2},T_{1}}\otimes \on{id}_{T_{2}})$
takes $\cC_{0}(X_{1},Y_{1}|T_{1}) \times \cC_{0}(X_{2},Y_{2}|T_{2})$
to $\cC_{0}(X_{1}\otimes X_{2},Y_{1}\otimes Y_{2}|T_{1}\otimes T_{2})$, 
and $\on{tr}_{T_{1}\otimes T_{2}}(x) 
= \on{tr}_{T_{1}}(x_{1})\otimes \on{tr}_{T_{2}}(x_{2})$; 
$\cC_{0}(X,Y|{\bf 1}) = \cC_{0}(X,Y)$ and $\on{tr}_{{\bf 1}}(x)=x$. 

Set
$\on{Ob}(\cC) := \on{Ob}(\cC_{0})^{2}$, $(X,Y)\otimes (X',Y'):= 
(X\otimes X',Y\otimes Y')$, 
$\cC((X,Y),(X',Y')):= 
\cC_{0}(X\otimes Y', X'\otimes Y)$, 
$\cC((X_{1},Y_{1}),(X_{2},Y_{2}),(X_{3},Y_{3})):=
\{(x_{1},x_{2}) | x_{2}*x_{1}\in \cC_{0}(X_{1}\otimes Y_{3},
X_{3}\otimes Y_{1}|Y_{2})\}$, where $x_{2}*x_{1} = 
(\on{id}_{X_{3}}\otimes \beta_{Y_{2}Y_{1}}) \circ (x_{2}
\otimes \on{id}_{Y_{1}}) \circ (\on{id}_{X_{2}}\otimes 
\beta_{Y_{1}Y_{3}}) \circ (x_{1}\otimes\on{id}_{Y_{3}})
\circ (\on{id}_{X_{1}}\otimes \beta_{Y_{3}Y_{2}})$; 
then $x_{2}\circ x_{1}:= \on{tr}_{Y_{2}}(x_{2}*x_{1})$. 
The tensor product of morphisms
is defined as $\times_{i=1}^{2}\cC_{0}(X_{i}\otimes Y'_{i},X'_{i}\otimes Y_{i})
\ni (x_{1},x_{2})\mapsto (\on{id}_{X'_{1}}
\otimes \beta_{Y_{1},X'_{2}} \otimes \on{id}_{Y_{2}})
 \circ (x_{1}\otimes x_{2}) \circ (\on{id}_{X_{1}}
\otimes \beta_{X_{2},Y'_{1}} \otimes \on{id}_{Y'_{2}})
\in \cC_{0}(X_{1}\otimes X_{2}\otimes Y'_{1}\otimes Y'_{2},
X'_{1}\otimes X'_{2}\otimes Y_{1}\otimes Y_{2})$; the unit of $\cC$
is $({\bf 1},{\bf 1})$. 

\begin{proposition}
$\cC$ is a strict monoidal quasi-category. 
\end{proposition}

{\em Proof.} Let $U_{i}=(X_{i},Y_{i})$; let $x_{i}\in \cC(U_{i},U_{i+1})$ 
($i=1,2,3$); assume that $(x_{1},x_{2})\in \cC(U_{1},U_{2},U_{3})$ and 
$(x_{2}\circ x_{1},x_{3})\in \cC(U_{1},U_{3},U_{4})$; define
$x_{3}*x_{2}*x_{1}$ by formula (\ref{formula*}) below. 
Let us show that $x_{3}*x_{2}*x_{1}\in \cC(X_{1}\otimes Y_{4},X_{4}\otimes 
Y_{1}|Y_{2}\otimes Y_{3})$ and that 
$x_{3}\circ (x_{2}\circ x_{1}) = \on{tr}_{Y_{2}\otimes Y_{3}}
(x_{3}*x_{2}*x_{1})$. 

$x_{3}\circ (x_{2}\circ x_{1}) = \on{tr}_{Y_{3}}(x_{3} * 
\on{tr}_{Y_{2}}(x_{2}*x_{1}))$. Using the fact that $x_{2}*x_{1}$ may 
as well be expressed as 
$x_{2}* x_{1} = (\on{id}_{X_{3}}\otimes \beta_{Y_{2}Y_{1}})
\circ (x_{2}\otimes \on{id}_{Y_{1}})\circ (\beta_{Y_{3}X_{2}}
\otimes \on{id}_{Y_{1}})\circ (\on{id}_{Y_{3}}\otimes x_{1})
\circ (\beta_{X_{1}Y_{3}}\otimes \on{id}_{Y_{2}})$, 
we write $x_{3}*\on{tr}_{Y_{2}}(x_{2}*x_{1})
=(\on{id}_{X_{4}}\otimes \beta_{Y_{3},Y_{1}})
\circ  (x_{3}\otimes \on{id}_{Y_{1}})
\circ (\beta_{Y_{4},X_{3}}\otimes\on{id}_{Y_{1}})
\circ (\on{id}_{Y_{4}}\otimes \on{tr}_{Y_{2}}(x_{2}*x_{1}))
\circ (\beta_{X_{1},Y_{4}}\otimes \on{id}_{Y_{3}})$. Now 
$\on{id}_{Y_{4}}\otimes (x_{2}*x_{1})\in \cC_{0}(Y_{4}\otimes
X_{1}\otimes Y_{3},Y_{4}\otimes X_{3}\otimes Y_{1}|Y_{2})$, and 
$\on{id}_{Y_{4}}\otimes \on{tr}_{Y_{2}}(x_{2}*x_{1}) = 
\on{tr}_{Y_{2}}(\on{id}_{Y_{4}}\otimes (x_{2}*x_{1}))$. 
We have then 
$ [\{(\on{id}_{X_{4}}\otimes \beta_{Y_{3},Y_{1}})
\circ (x_{3}\otimes \on{id}_{Y_{1}})
\circ (\beta_{Y_{4},X_{3}}\otimes \on{id}_{Y_{1}})\}\otimes \on{id}_{Y_{2}}] 
\circ [\on{id}_{Y_{4}}\otimes (x_{2}*x_{1})]
\circ [\beta_{X_{1},Y_{4}}\otimes 
\on{id}_{Y_{3}}\otimes \on{id}_{Y_{2}}]
\in \cC_{0}(X_{1}\otimes Y_{4}\otimes Y_{3},X_{4}\otimes Y_{1}\otimes
Y_{3}|Y_{2})$, and  
$x_{3}*\on{tr}_{Y_{2}}(x_{2}*x_{1})
= \on{tr}_{Y_{2}} \{ [\{(\on{id}_{X_{4}}\otimes \beta_{Y_{3},Y_{1}})
\circ (x_{3}\otimes \on{id}_{Y_{1}})
\circ (\beta_{Y_{4},X_{3}}\otimes \on{id}_{Y_{1}})\}\otimes \on{id}_{Y_{2}}] 
\circ [\on{id}_{Y_{4}}\otimes (x_{2}*x_{1})]
\circ [\beta_{X_{1},Y_{4}}\otimes 
\on{id}_{Y_{3}}\otimes \on{id}_{Y_{2}}]\}$. As the right side is in 
the domain of $\on{tr}_{Y_{3}}$, the argument of $\on{tr}_{Y_{2}}$ in the
right side is in the domain of $\on{tr}_{Y_{3}\otimes Y_{2}}$, and 
$x_{3}\circ (x_{2}\circ x_{1}) = \on{tr}_{Y_{3}\otimes Y_{2}}
\{ [\{(\on{id}_{X_{4}}\otimes \beta_{Y_{3},Y_{1}})
\circ (x_{3}\otimes \on{id}_{Y_{1}})
\circ (\beta_{Y_{4},X_{3}}\otimes \on{id}_{Y_{1}})\}\otimes \on{id}_{Y_{2}}] 
\circ [\on{id}_{Y_{4}}\otimes (x_{2}*x_{1})]
\circ [\beta_{X_{1},Y_{4}}\otimes 
\on{id}_{Y_{3}}\otimes \on{id}_{Y_{2}}]\}$. On the other hand, this
argument also expressed as $(\on{id}_{X_{1}}\otimes \on{id}_{Y_{4}}\otimes 
\beta_{Y_{2},Y_{3}}) \circ (x_{3}*x_{2}*x_{1}) \circ (\on{id}_{X_{1}}
\otimes \on{id}_{Y_{4}}\otimes \beta_{Y_{3},Y_{2}})$, therefore 
$x_{3}*x_{2}*x_{1}$ is in the domain of $\on{tr}_{Y_{2}\otimes Y_{3}}$
and $x_{3}\circ (x_{2}\circ x_{1}) = \on{tr}_{Y_{2}\otimes Y_{3}}(x_{3}
*x_{2}*x_{1})$. One proves in the same way that $(x_{3}\circ x_{2})\circ x_{1}
= \on{tr}_{Y_{2}\otimes Y_{3}}(x_{3}*x_{2}*x_{1})$, which proves the 
asociativity identity. 
\hfill \qed \medskip 

More generally, one shows that for any 
$(x_{1},...,x_{n-1})\in \cC((X_{1},Y_{1}),...,(X_{n},Y_{n}))$, 
we have  $x_{n-1}*...*x_{1}\in 
\cC_{0}(X_{1}\otimes Y_{n},X_{n}\otimes Y_{1}|Y_{2}\otimes ... \otimes Y_{n-1})$, 
where $x_{n-1}*...*x_{1}\in \cC_{0}(X_{1}\otimes Y_{n}\otimes 
Y_{2}\otimes ... \otimes Y_{n-1},X_{n}\otimes Y_{1}\otimes 
Y_{2}\otimes ... \otimes Y_{n-1})$ is defined inductively by 
\begin{align} \label{formula*}
& \nonumber x_{n}*... * x_{1}:= 
(\on{id}_{X_{n+1}}\otimes\beta_{Y_{n},Y_{1}\otimes...\otimes Y_{n-1}})
\circ (x_{n}\otimes\on{id}_{Y_{1}\otimes ...\otimes Y_{n-1}}) \circ 
(\on{id}_{X_{n}}\otimes \beta_{Y_{1}\otimes...\otimes Y_{n-1},Y_{n+1}})
\\ 
 & \circ [(x_{n-1}*...*x_{1})\otimes\on{id}_{Y_{n+1}}] \circ 
(\on{id}_{X_{1}} \circ\beta_{Y_{n+1},Y_{2}\otimes...\otimes Y_{n-1},Y_{n}}),  
\end{align}
where $\beta_{X,Y,Z}\in \cC_{0}(X\otimes Y\otimes Z,Z\otimes Y\otimes X)$ 
is $\beta_{X\otimes Y,Z} \circ 
(\beta_{X,Y}\otimes\on{id}_{Z})$, and that 
$x_{n-1}\circ ... \circ x_{1}:= \on{tr}_{Y_{2}\otimes ... 
\otimes Y_{n-1}}(x_{n-1}*...*x_{1})$. 

If $X\mapsto X^*$ is an involution of $\cC_{0}$, another symmetric 
strict monoidal quasi-category $\cC'$ may be defined by 
$\on{Ob}(\cC') = \on{Ob}(\cC_{0})^{2}$, 
$\cC'((X,Y),(X',Y')):= \cC_{0}(X\otimes Y^{\prime *},
X'\otimes Y^{*})$. 

A functor between categories with partial traces is a tensor functor 
$F:\cC_{0}\to \cD_{0}$, such that $F(\cC_{0}(X,Y|Z)) \subset 
\cD_{0}(F(X),F(Y)|F(Z))$, and such that $\on{tr}_{F(Z)} \circ F = 
F\circ \on{tr}_{Z}$ (equality of maps $\cC_{0}(X,Y|Z)\to 
\cD_{0}(F(X),F(Y))$). Such a functor induces a functor $\cC\to\cD$
between the corresponding quasi-categories. 

If now $\cC_{0}$ is additive and is a free module category over $\on{Vect}$, 
this construction can be extended as follows. In the definition of a trace, 
the maps are now linear and products are replaced by tensor products. 
Let $\on{Ob}'(\cC_{0})\subset \on{Ob}(\cC_{0})$ be a set of 
generators, i.e., each $F\in\on{Ob}(\cC_{0})$ has the form 
$\oplus_{X\in\on{Ob}'(\cC_{0})} F_{X}\otimes X$, where 
$X\mapsto F_{X}$ is a finitely supported map $\on{Ob}'(\cC_{0})
\to \on{Ob(Vect)}$. Then $\on{Ob}(\cC):= \{$finitely supported
maps $\on{Ob}'(\cC_{0})^{2}\to \on{Vect}$, 
$F = [(X,Y)\mapsto F_{X,Y}]\}$. 
We then set $\cC(F,G):= \oplus_{(X,Y),(X',Y')}
\on{Vect}(F_{X,Y},G_{X',Y'})\otimes 
\cC_{0}(X\otimes Y',X'\otimes Y)$ and extend the above composition and 
tensor product operations by linearity. 
In the case of $\cC'$, we replace 
$\cC_{0}(X\otimes Y',X'\otimes Y)$ by 
$\cC_{0}(X\otimes Y^{\prime*},X'\otimes Y^{*})$. 

\medskip\noindent
{\bf Example.} Let $\cG_{0}$ be the category where objects are finite sets, and 
$\cG_{0}(I,J) = \{$subsets of $I\times J\}$, and composition given by 
$S'\circ S:=$ the image in 
$I\times I''$ of $S\times_{I'}S'$, for $S\subset I\times I'$
and $S'\subset I'\times I''$;  to $S\subset I\times J$, we associate
the oriented graph with vertices $I\sqcup J$ and egdes $i\to j$ if 
$(i,j)\in S$, composition then corresponds to the composition of graphs. 
The tensor product is $I\otimes I':= 
I\sqcup I'$ and for $S\in \cG_{0}(I,I')$, $T\in \cG_{0}(J,J')$, 
$S\otimes T:= S\sqcup T\subset (I\times I')\sqcup 
(J\times J')\subset (I\sqcup I')\times (J\sqcup J')$. 
Then $\cG_{0}$ is a strict monoidal category. It has a partial 
trace defined as follows. For $I,J,K$ finite sets, 
let $\cG_{0}(I,J|K)\subset \cG_{0}(I\sqcup K,J\sqcup K)$ 
be the set of graphs, such that the introduction of 
the edges $k_{out}\to k_{in}$ ($k\in K$) does not introduce cycles
(alternatively, the set of $S\subset (I\sqcup K)\times (J\sqcup K)$, 
such that the relation in $K$ defined by $u\prec v$ if
$(u,v)\in S$, has no cycle), and if $x$ is such a graph, then 
$\on{tr}_{K}(x)\in \cG_{0}(I,J)$ corresponds to 
$\{(i,j)\in I\times J |$ there exists $s\geq 0$ and 
a sequence $(k_{1},...,k_{s})$ of elements of $K$, such that 
$i\prec k_{1}\prec ... \prec k_{s}\prec j\}$, where the relation 
$\prec$ is extended to $I\sqcup K\sqcup J$ by $u\prec v$ iff 
$(u,v)\in S$. 
Then the strict monoidal 
quasi-category constructed from $\cG_{0}$, equipped with its partial trace, coincides
with $\cG$. 

Here is another description of $\on{tr}_{K}(x)$. As the relation $\prec$
on $K$ is acyclic, we may extend it to a total order relation $<$ on $K$. 
Extend it to $I\sqcup K\sqcup J$ by $i<k<j$ for any 
$i,j,k\in I,J,K$. The relation $<$ induces
a numbering $K=\{k_{1},...,k_{|K|}\}$ for $K$, where $k_{1}<...<k_{|K|}$. 
For $\alpha\in [|K|]$, let $K_{\alpha}:= \{k_{\alpha}\}\sqcup 
\{(u,v)\in (I\sqcup K\sqcup J)^{2}|u\prec v, 
u<k_{\alpha}<v\}\in \on{Ob}(\cG_{0})$. Then 
$\on{tr}_{K}(x) = x_{K_{|K|}J}\circ...\circ x_{K_{1}K_{2}}\circ x_{IK_{1}}$, 
where: 

$\bullet$ $x_{K_{\alpha}K_{\alpha+1}}\in \cG_{0}(K_{\alpha},K_{\alpha+1})$
is defined as follows: we have identifications
$K_{\alpha} \simeq \{k_{\alpha}\}\sqcup K'_{\alpha,\alpha+1}\sqcup 
K_{\alpha,\alpha+1}$ and $K_{\alpha+1} \simeq \{k_{\alpha+1}\}\sqcup
K''_{\alpha,\alpha+1}\sqcup K_{\alpha,\alpha+1}$, where
$K'_{\alpha,\alpha+1} := \{s\in I\sqcup K|s<k_\alpha, s\prec k_{\alpha+1}\}$, 
$K''_{\alpha,\alpha+1} := \{t\in K\sqcup J|t>k_{\alpha+1}, t\succ k_\alpha\}$, 
$K_{\alpha,\alpha+1} := \{(s,t)\in (I\sqcup K)\times (K\sqcup J)|s
<k_\alpha,k_{\alpha+1}<t, s\prec t\}$; let $\diamond$ be a one-element set if 
$k_\alpha\prec k_{\alpha+1}$ 
and $\emptyset$ otherwise; then we define $\kappa_{\alpha,\alpha+1}
:= \{k_{\alpha}\}\times (\diamond
\sqcup K''_{\alpha+1})\in \cG_{0}(\{k_{\alpha}\},\diamond
\sqcup K''_{\alpha,\alpha+1})$ and 
$\lambda_{\alpha,\alpha+1}:= (\diamond\sqcup K'_{\alpha,\alpha+1})\times
\{k_{\alpha+1}\}\in \cG_{0}(\diamond\sqcup K'_{\alpha,\alpha+1}, 
\{k_{\alpha+1}\})$; then $x_{K_{\alpha}K_{\alpha+1}}:=
[(\lambda_{\alpha,\alpha+1}\otimes \on{id}_{K''_{\alpha,\alpha+1}})
\circ (\on{id}_{\diamond}\otimes \beta_{K''_{\alpha,\alpha+1},
K'_{\alpha,\alpha+1}}) \circ (\kappa_{\alpha,\alpha+1}
\otimes\on{id}_{K'_{\alpha,\alpha+1}})]\otimes 
\on{id}_{K_{\alpha,\alpha+1}}$;  

$\bullet$ $x_{IK_{1}}\in \cG_{0}(I,K_{1})$ is defined as follows: $K_{1} \simeq 
\{k_{1}\}\sqcup (\sqcup_{i\in I}K''_{i})$, where $K''_{i}:= \{t\in K\sqcup J|
t\neq k_{1}, t\succ i\}$; set $\diamond:= \{i\in I|i\prec k_{1}\}$ and $\diamond_{i}:= 
\diamond\cap \{i\}$ so $\diamond = \sqcup_{i\in I}\diamond_{i}$; let $\kappa_{i}:= 
\{i\}\times (\diamond_{i}\sqcup K''_{i})\in \cG_{0}(\{i\},\diamond_{i}
\sqcup K''_{i})$, 
$br\in \cG_{0}(\sqcup_{i}(*_{i}\sqcup K''_{i}),*\sqcup (\sqcup_{i\in I}
K''_{i}))$ be the canonical braiding morphism; and let 
$\lambda_{01}:= \diamond\times \{k_{1}\}
\in \cG_{0}(\diamond,\{k_{1}\})$; then $x_{IK_{1}}:= (\lambda_{01}\otimes 
(\otimes_{i\in I}\on{id}_{K''_{i}})) \circ br \circ 
(\otimes_{i\in I}\kappa_{i})$; 

$\bullet$ $x_{K_{|K|}J}\in \cG_{0}(K_{|K|},J)$ is defined as follows: 
$K_{|K|} \simeq \{k_{|K|}\}\sqcup (\sqcup_{j\in J}K'_{j})$, where $K'_{j}:=
\{s\in I\sqcup K| s\neq k_{|K|}, s\prec j\}$, set $\diamond:= 
\{j\in J|j\succ k_{|K|}\}$; 
let $\diamond_{J}:= \diamond\cap\{j\}$, then $\diamond=\sqcup_{j\in J}
\diamond_{j}$; 
let $\kappa_{|K|,|K|+1}:= \{k_{|K|}\}\times \diamond\in\cG_{0}(\{k_{|K|}\},
\diamond)$; 
let $\lambda_{j}:= K'_{j}\times \{j\}\in \cG_{0}(K'_{j},\{j\})$; 
let $br\in \cG_{0}(\diamond\sqcup
(\sqcup_{j\in J}K'_{j}),\sqcup_{j\in J}(\diamond_{j}\sqcup K'_{j}))$
be the canonical braiding map; then $x_{K_{|K|}J}:= (\otimes_{j\in J}\lambda_{j})
\circ br \circ 
(\kappa_{|K|,|K|+1}\otimes (\otimes_{j\in J}\on{id}_{K'_{j}}))$.

\subsection{Props and (quasi)(bi)(multi)props}

A prop $P$ is a symmetric additive strict monoidal category, equipped with a tensor
functor $i_P : \on{Sch} \to P$, inducing a bijection on the sets of  
objects; so $\on{Ob}(P) = \on{Ob}(\on{Sch})$ (see, e.g., \cite{Tam}). 
It is easy to check that this definition is equivalent to the original one 
(\cite{McL}). 
For $\phi\in\on{Sch}(F,G) \to P(F,G)$, we sometimes write 
$\phi$ instead of $i_P(\phi)$. A prop morphism $f:P\to Q$ is a 
tensor functor, inducing a bijection on the sets of objects, and such that $f\circ
i_{P}=i_{Q}$. 

A biprop (resp., multiprop, bi-multiprop) is a symmetric additive monoidal category 
$\pi$ (resp., $\Pi^{0},\Pi$), equipped with a 
tensor functor $\on{Sch}_{1+1} \to \pi$ (resp., 
$\on{Sch}_{(1)}\to \Pi^{0}$, $\on{Sch}_{(1+1)}\to\Pi$), 
which induces a bijection on the sets of objects. Morphisms betweens these structures
are defined as above. 

A quasi-prop is a symmetric additive strict monoidal quasi-category $P$, 
equipped with a morphism $i_{P}:\on{Sch}\to P$, inducing a bijection on 
the sets of
objects. Quasi(bi)(multi)props are defined in the same way, as well
as morphisms between these structures. 

A topological (quasi)(bi)(multi)prop is defined in the same way as its
non-topological analogue, replacing $\on{Sch}_{*}$ by ${\bf Sch}_{*}$. 
E.g., a topological prop ${\bf P}$ is a symmetric tensor category, 
equipped with a morphism ${\bf Sch} \to {\bf P}$, which is the identity 
on objects.

\subsection{Operations on props}

If $H\in \on{Ob(Sch)}$ (resp., $\on{Ob}(\on{Sch}_{1+1})$) 
and $P$ is a (bi)prop, then we define a prop $H(P)$ by
$H(P)(F,G):= P(F\circ H,G\circ H)$. 
A (bi)prop morphism $P\to Q$ 
gives rise to a prop morphism $H(P)\to H(Q)$. 
Similarly, if ${\bf P}$ is a topological (bi)prop, 
then for any $H\in \on{Ob}({\bf Sch})$ 
(resp., $\on{Ob}({\bf Sch}_{1+1})$), 
we get a prop $H({\bf P})$, such that 
$H({\bf P})(F,G) := {\bf P}(F\circ H,G\circ H)$. 
A morphism of topological (bi)props ${\bf P} \to {\bf Q}$ then gives
rise to a prop morphism $H({\bf P}) \to H({\bf Q})$. 
 
To each (quasi)(bi)multiprop $\Pi$, one associates a (quasi)(bi)prop 
$\pi$ by $\pi(F,F') := \Pi(F,F')$, i.e., using the injections 
$\on{Ob}(\on{Sch}) = \on{Ob}(\on{Sch}_{1})
\subset \on{Ob}(\on{Sch}_{(1)})$ in the ``non-bi'' case, and
$\on{Ob}(\on{Sch}_{1+1}) 
\subset \on{Ob}(\on{Sch}_{(1+1)})$, $F\mapsto (F_{k,l})$, where
$F_{k,l}=0$ if $(k,l)\neq (1,1)$ and $F_{1,1}=F$, in the  ``bi'' case. 

 If $P$ is a prop, then one defines a 
multiprop $\Pi^{0}_P$ by $\Pi^{0}_P(F,G) := P(c(F),c(G))$. 
The tensor product is induced by the tensor product of $P$
and the identity $c(F\boxtimes F') = c(F)\otimes c(F')$.

\subsection{Presentation of a prop} 

If $P$ is a prop, then a prop ideal $I_P$ of $P$
is a collection of vector subspaces $I_P(F,G) \subset P(F,G)$, such that 
$(F,G) \mapsto P(F,G)/I_P(F,G)$ is a prop, which we denote by $P/I_P$. 
Then $P\to P/I_P$ is a prop morphism. 

If $P$ is a prop, $(F_i,G_i)_{i\in I}$ is a collection of pairs 
of Schur functors 
and $V_i \subset P(F_i,G_i)$ are vector subspaces, then $(V_i,i\in I)$
is the smallest of all prop ideals $I_P$ of $P$, such that 
$V_i \subset I_P(F_i,G_i)\subset P(F_i,G_i)$ for any $i\in I$.  

Let $(F_i,G_i)_{i\in I}$ be a collection of Schur functors, 
and let $(V_i)_{i\in I}$ be a collection of vector spaces.  Then there exists 
a unique (up to isomorphism) prop $\cF = \on{Free}(V_i,F_i,G_i,i\in I)$, 
which is initial in the category of
all props $P$ equipped with linear maps $V_i \to P(F_i,G_i)$. We call it the 
free prop generated by $(V_i,F_i,G_i)$. 
 
If $(F'_\alpha,G'_\alpha)$ is a collection of Schur functors
and $R_\alpha \subset \cF(F'_\alpha,G'_\alpha)$ is a collection of vector
spaces, then the prop with generators $(V_i,F_i,G_i)$ with relations 
$(R_\alpha,F'_\alpha,G'_\alpha)$ is the quotient of $\cF$ with the prop ideal
generated by $R_\alpha$. 

\subsection{Topological props} \label{1:4}

Let $P$ be a prop equipped with a filtration 
$P(F,G) = P^0(F,G) \supset P^1(F,G) \supset ...$ for any 
$F,G\in\on{Ob(Sch)}$, compatible with direct sums and such that: 

(a) $\circ$ and $\otimes$ induce maps 
$\circ : P^i(F,G) \otimes P^j(G,H) \to P^{i+j}(F,H)$, and 
$\otimes : P^i(F,G) \otimes P^{i'}(F',G') \to P^{i+i'}(F\otimes F',
G\otimes G')$

(b) if $F,G\in \on{Ob(Sch)}$ are homogeneous, then 
$P(F,G) = P^{||F| - |G||}(F,G)$. 

For $F,G\in \on{Ob(Sch)}$, we then define $\wh P(F,G) = \on{lim}_{\leftarrow}
P(F,G)/P^n(F,G)$ as 
the completed separated of $P(F,G)$ w.r.t. the filtration 
$P^n(F,G)$.  Then $\wh P$ is a prop.   

If $F,G\in \on{Ob}({\bf Sch})$, define ${\bf P}(F,G)$ as follows: 
for $F = \widehat\oplus_{i\geq 0} F_i$, $G = \widehat\oplus_{i\geq 0} G_i$
the decompositions of $F,G$ into sums of homogeneous components, 
we set ${\bf P}(F,G) = \wh\oplus_{i,j\geq 0} \wh P(F_i,G_j)$
(where $\wh\oplus$ is the direct product). 

\begin{proposition}
${\bf P}$ is a symmetric additive strict monoidal category, equipped with a morphism 
${\bf Sch} \to {\bf P}$, which is the identity on objects. 
\end{proposition} 

Recall that ${\bf P}$ is called a topological prop. 

\medskip 
{\em Proof.} Let $F = \wh\oplus_i F_i$, $G = \wh\oplus_i G_i$, $H = \oplus_i H_i$
be in ${\bf Sch}$. We define a map $\circ : {\bf P}(F,G) \otimes {\bf P}(G,H)
\to {\bf P}(F,H)$ as follows. We first define a map 
${\bf P}(F_i,G) \otimes {\bf P}(G,H_k) \to \wh P(F_i,H_k)$. 
The left vector space injects in $\wh\oplus_j \wh P(F_i,G_j)
\otimes \wh P(G_j,H_k)$. The composition takes the $j$th summand 
to $P^{|j-i| + |j-k|}(F_i,H_k)$. As $|j-i| + |j-k| \to\infty$ as 
$j\to\infty$, we have a well-defined map  
${\bf P}(F_i,G) \otimes {\bf P}(G,H_k) \to \wh P(F_i,H_k)$. 
The direct product of these maps then induces a map 
$\circ : {\bf P}(F,G) \otimes {\bf P}(G,H)
\to {\bf P}(F,H)$. 

Let now $F' = \oplus_i F'_i$, $G' = \oplus_i G'_i$ be in ${\bf Sch}$. 
We define a map $\otimes : {\bf P}(F,G) \otimes {\bf P}(F',G')
\to {\bf P}(F\otimes F',G\otimes G')$ as the direct product of the 
maps $\wh P(F_i,G_j) \otimes \wh P(F'_{i'},G'_{j'}) \to \wh P(F_i \otimes 
F'_{i'}, G_j \otimes G'_{j'})$. This is well-defined since 
$(F\otimes F')_i = \oplus_{j=0}^n F_j \otimes F'_{i-j}$ is the sum of 
a finite number of tensor products, and the same holds for $(G\otimes G')_j$. 
\hfill \qed \medskip 

A grading of $P$ by an abelian semigroup $\Gamma$ is a decomposition 
$P(F,G) = \oplus_{\gamma\in \Gamma} P_\gamma(F,G)$, such that the prop 
operations are compatible with the semigroup structure of $\Gamma$. 
Then if $P$ is graded by $\NN$ and if we set 
$P^n(F,G) = \oplus_{i\geq n} P_i(F,G)$, the descending filtration 
$P = P^0 \supset...$ satisfies condition (a) above. 

If $P \to Q$ is a surjective prop morphism (i.e., the maps 
$P(F,G) \to Q(F,G)$ are all surjective), and if $P$ 
is equipped with a filtration as above, then so is $Q$
(we define $Q^n(F,G)$ as the image of $P^n(F,G)$). Then we get a 
morphism ${\bf P} \to {\bf Q}$ of topological props, i.e., a morphism of
tensor categories such that the 
morphisms ${\bf Sch} \to {\bf P} \to {\bf Q}$ and ${\bf Sch} \to 
{\bf Q}$ coincide. 

If $P$, $R$ are props equipped with a filtration as above, and 
$P\to R$ is a prop morphism compatible with the filtration (i.e., 
$P^n(F,G)$ maps to $R^n(F,G)$), then we get a morphism   
of topological props ${\bf P} \to {\bf R}$. 

\subsection{Modules over props}

If $\cS$ is an additive symmetric strict monoidal category, 
and $V\in \on{Ob}(\cS)$, then we have a prop $\on{Prop}(V)$, s.t. 
$\on{Prop}(V)(F,G)
= \on{Hom}_\cS(F(V),G(V))$. Then a $P$-module (in the category $\cS$) 
is a pair $(V,\rho)$, where $V\in \on{Ob}(\cS)$ and $\rho : P \to \on{Prop}(V)$ 
is a tensor functor. Then $P$-modules in the category $\cS$
form a category. The tautological $P$-module is $\cS=P$, 
$V={\bf id}$. 
  
\subsection{Examples of props} \label{1:2:4}
 
We will define several props by generators and relations.  

\subsubsection{The prop $\on{Bialg}$} This is the prop with generators 
$m\in \on{Bialg}(T_2,{\bf id})$, 
$\Delta\in \on{Bialg}({\bf id},T_2)$, $\eta\in \on{Bialg}({\bf 1},{\bf id})$, 
$\eps\in \on{Bialg}({\bf id},{\bf 1})$, and relations  
$$
m \circ (m\otimes \on{id}_{{\bf id}})
= m \circ (\on{id}_{{\bf id}} \otimes m), 
\quad 
(\Delta\otimes \on{id}_{{\bf id}})
\circ \Delta  = 
(\on{id}_{{\bf id}} \otimes \Delta)
\circ \Delta, 
$$
$$
\Delta \circ m = (m\otimes m) \circ (1324) \circ (\Delta \otimes\Delta), 
$$
$$
m \circ (\eta \otimes \on{id}_{{\bf id}})
= m \circ (\on{id}_{{\bf id}} \otimes \eta) 
= \on{id}_{{\bf id}},  
\quad 
(\eps \otimes \on{id}_{{\bf id}}) \circ \Delta  
= (\on{id}_{{\bf id}} \otimes \eps) \circ \Delta
= \on{id}_{{\bf id}}.    
$$
When $\cS = \on{Vect}$, the category of $\on{Bialg}$-modules is
that of bialgebras. 

\subsubsection{ The prop $\on{COB}$} This is the prop with generators 
$m\in \on{COB}(T_2,{\bf id})$, $\Delta \in \on{COB}({\bf id},T_2)$, 
$\eta \in \on{COB}({\bf 1},{\bf id})$, $\eps \in \on{COB}({\bf id},{\bf 1})$ 
and $R\in \on{COB}({\bf 1},T_2)$, and relations: 
$m,\Delta,\eta,\eps$ satisfy the relations of Bialg, 
$$
(m \otimes m) \circ (1324) \circ (R \otimes ((21)\circ R))   
= (m \otimes m) \circ (1324) \circ ((21)\circ R) \otimes R)   
= \eta\otimes \eta, 
$$
$$
(21) \circ (m\otimes m) \circ (1324) \circ (\Delta \otimes R) 
= (m \otimes m) \circ (1324) \circ (R \otimes \Delta), 
$$
$$
m^{\otimes 3} \circ (142536) \circ \big( (R \otimes \eta) \otimes 
((\Delta \otimes \on{id}_{{\bf id}}) \circ R) \big)  
= 
m^{\otimes 3} \circ (142536) \circ \big( (\eta \otimes R) \boxtimes 
(\on{id}_{{\bf id}} \otimes \Delta) \circ R) \big)  
$$
The category of COB-modules over $\cS = \on{Vect}$ is that of 
coboundary bialgebras, i.e., pairs $(A,R_A)$, where $A$ is a bialgebra, 
and $R_A\in A^{\otimes 2}$ satisfies 
$R_A R_A^{21} = R_A^{21} R_A = 1_A^{\otimes 2}$, 
$\Delta_A^{21}(x)R_A^{21} = R_A\Delta_A(x)$, and 
$$
(R_A\otimes 1_A) \big( (\Delta_A \otimes \on{id}_A)(R_A) \big) = 
(1_A \otimes R_A) \big( (\on{id}_A \otimes \Delta_A)(R_A) \big). 
$$

\subsubsection{The prop $\on{LA}$} 
This is the prop with generator the bracket 
$\mu\in\on{LA}(\wedge^2,{\bf id})$ and relation the Jacobi identity
\begin{equation} \label{Jacobi}
\mu \circ (\mu\otimes \on{id}_{{\bf id}}) \circ 
((123) + (231) + (312)) = 0. 
\end{equation}
When $\cS = \on{Vect}$, the category of $\on{LA}$-modules is
that of Lie algebras.

\subsubsection{The prop $\on{LCA}$} This is the prop with 
generator the cobracket $\delta\in \on{LCA}({\bf id},\wedge^2)$
and relation the co-Jacobi identity 
$$
((123) + (231) + (312)) \circ 
(\delta\otimes\on{id}_{{\bf id}}) \circ \delta = 0. 
$$
When $\cS = \on{Vect}$, the category of $\on{LCA}$-modules is
that of Lie coalgebras.

\subsubsection{The prop $\on{LBA}$} This is the prop with generators 
$\mu\in \on{LBA}(\wedge^2,{\bf id})$, 
$\delta\in \on{LBA}({\bf id},\wedge^2)$; relations are the Jacobi and the 
co-Jacobi identities, and the cocycle relation  
$$
\delta \circ \mu = ((12) - (21)) \circ  (\mu \otimes 
\on{id}_{{\bf id}}) \circ 
(\on{id}_{{\bf id}} \otimes \delta) \circ ((12)-(21)).   
$$
When $\cS = \on{Vect}$, the category of $\on{LBA}$-modules is
that of Lie bialgebras. 

\subsubsection{The prop $\on{LBA}_f$}
This is the prop with generators $\mu\in \on{LBA}_f(\wedge^2,{\bf id})$,
$\delta\in \on{LBA}_f({\bf id},\wedge^2)$, $f\in \on{LBA}_f({\bf 1},\wedge^2)$
and relations: $\mu,\delta$ satisfy the relations of LBA, and
$$
\big( (123) + (231) +(312) \big)  \circ 
\Big( (\delta\otimes \on{id}_{{\bf id}}) \circ f
+ (\mu \otimes \on{id}_{{\bf id}^{\otimes 2}})\circ (1324) \circ 
(f\otimes f) \Big) = 0. 
$$
The category of $\on{LBA}_f$-modules 
is the category of pairs $(\a,f_\a)$ 
where $\a$ is a Lie bialgebra and $f_\a$ is a  twist of $\a$.

\subsubsection{The prop $\on{Cob}$}
This is the prop with generators 
$\mu\in \on{Cob}(\wedge^2,{\bf id})$
and $\rho\in \on{Cob}({\bf 1},\wedge^2)$ and relations: 
$\mu$ satisfies the Jacobi identity (\ref{Jacobi}), 
and the element $Z\in \on{Cob}({\bf 1},\wedge^3)$ defined 
by 
$$
Z:= ((123) + (231) + (312)) \circ (\on{id}_{{\bf id}} \otimes 
\mu \otimes \on{id}_{{\bf id}}) \circ (\rho \otimes \rho)
$$
is invariant, i.e., it satisfies 
$$
\Big( (\mu\otimes \on{id}_{{\bf id}^{\otimes 2}}) \circ (1423)
+ (\on{id}_{{\bf id}} \otimes \mu \otimes \on{id}_{{\bf id}}) \circ (1243)
+ (\on{id}_{{\bf id}^{\otimes 2}} \otimes \mu) \Big) 
\circ (Z\otimes \on{id}_{{\bf id}}) = 0.
$$

The category of $\on{Cob}$-modules over $\cS = \on{Vect}$ is that of 
coboundary Lie bialgebras, i.e., pairs $(\a,\rho_\a)$, where $\a$  is a 
Lie algebra and $\rho_\a\in \wedge^2(\a)$ is such that $Z_\a := 
[\rho_\a^{12},\rho_\a^{13}] 
+ [\rho_\a^{12},\rho_\a^{23}] + [\rho_\a^{13},\rho_\a^{23}]$ is $\a$-invariant.  

\subsection{Some prop morphisms}

We have unique prop morphisms $\on{Cob} \to \on{Sch}$,
$\on{LBA} \to \on{Sch}$ and $\on{LBA}_f \to \on{Sch}$, respectively defined 
by $(\mu,\rho)\mapsto (0,0)$, $(\mu,\delta)\mapsto (0,0)$ and 
$(\mu,r)\mapsto (0,0)$. 

If $\on{LA} \to P$ is a prop morphism, and $\alpha\in P({\bf 1},T_2)$, 
define 
$$
\on{ad}(\alpha) := 
\Big( ((\mu\circ \on{Alt})\otimes 
\on{id}_{{\bf id}}) \circ (132)
+ \on{id}_{{\bf id}} \otimes (\mu\circ \on{Alt})
\Big) \circ (\alpha\otimes 
\on{id}_{{\bf id}}); 
$$
(here $\on{Alt} : T_2 \to \wedge^2$ is the alternation map); 
this is a propic version of the map $x\mapsto [\alpha_\a,x^1 + x^2]$, 
where $\a\in\on{Rep}(P)$. If $\alpha\in P({\bf 1},\wedge^2)$, then 
$\on{ad}(\alpha)\in P({\bf 1},\wedge^2)$. 
 
Using the presentations of $\on{LBA}$ and $\on{LBA}_f$, we get:  
\begin{proposition}
We have unique prop morphisms 
$$
\kappa_1,\kappa_2 : \on{LBA} \to \on{LBA}_f, \quad such\ that \quad  
\kappa_1(\mu) = \kappa_2(\mu) = \mu, \; \kappa_1(\delta) = \delta, \; 
\kappa_2(\delta) = \delta + \on{ad}(f)  
$$
and a unique prop morphism 
$$
\kappa_0 : \on{LBA}_f \to \on{LBA}, \quad such\  that \quad 
\kappa_0(\mu) = \mu, \; \kappa_0(\delta) = \delta, \; \kappa_0(f) = 0.  
$$
\end{proposition}

We also have a prop morphism 
$$
\kappa : \on{LBA}_f \to \on{Cob},
$$ 
such that 
$\mu\mapsto \mu$, $\delta \mapsto \on{ad}(\rho)$, $f\mapsto -2\rho$
and $\tau_{\on{LBA}} : \on{LBA} \to \on{LBA}$, defined by $(\mu,\delta)\mapsto 
(\mu,-\delta)$. 

\subsubsection{The prop $\on{Sch}$} $\on{Sch}$ it itself a prop (with no
generator and relation). The corresponding category of modules over $\cS$
is $\cS$ itself. 

\subsection{Examples of topological props}

\subsubsection{The prop ${\bf Sch}$}

We set $\on{Sch}^0(F,G) = \on{Sch}(F,G)$, 
$\on{Sch}^1(F,G)= ... = 0$. This filtration satisfies conditions (a) and 
(b) of Section \ref{1:4}, since for $F,G$ homogeneous, $\on{Sch}(F,G) = 0$ 
unless $F$ and $G$ have
the same degree. The corresponding completion of $\on{Sch}$ coincides with 
${\bf Sch}$. 

\subsubsection{The props ${\bf LA}$ and ${\bf LCA}$}

Since the relation in $\on{LA}$ is homogeneous in $\mu$, 
the prop $\on{LA}$ has a grading $\on{deg}_\mu$. If $F,G\in \on{Ob(Sch)}$
and $x\in \on{LA}(F,G)$ are homogeneous, then $|G| - |F| = 
-\on{deg}_\mu(x)$, which implies that the filtration induced by $\on{deg}_\mu$
satisfies conditions (a) and (b) above. We denote by ${\bf LA}$ the
corresponding topological prop. 

In the same way, $\on{LCA}$ has a grading $\on{deg}_\delta$, 
$|G| - |F| = \on{deg}_\delta(x)$, so the filtration 
induced by $\on{deg}_\delta$ satisfies conditions (a) and (b) above. 
We denote by ${\bf LCA}$ the corresponding topological prop. 

\subsubsection{The prop ${\bf LBA}$}

Since the relations in $\on{LBA}$ are homogeneous in both $\mu$
and $\delta$, the prop $\on{LBA}$ is equipped with a grading 
$(\on{deg}_\mu,\on{deg}_\delta)$ by $\NN^2$. Moreover, if 
$F,G\in\on{Ob(Sch)}$ and $x\in \on{LBA}(F,G)$ are homogeneous,  
then 
\begin{equation} \label{homogeneity}
|G| - |F| = \on{deg}_\delta(x) - \on{deg}_\mu(x).
\end{equation} 

Then $\on{deg}_\mu + \on{deg}_\delta$ is a grading of $\on{LBA}$ by $\NN$. 
The corresponding filtration therefore satisfies condition (a) above. 
(\ref{homogeneity}) also implies that it satisfies condition (b),
since $\on{deg}_\mu$ and $\on{deg}_\delta$ are $\geq 0$. We denote by 
${\bf LBA}$ the resulting topological prop. 

Let $\cL\cB\cA$ be the category of Lie bialgebras over $\kk$.
Let $\cS_1$ be the category of topological $\kk[[\hbar]]$-modules
(i.e., quotients of modules of the form $V[[\hbar]]$, where $V\in\on{Vect}$
and the topology is 
given by the images of $\hbar^n V[[\hbar]]$) and let $\cS_2$ be the category of 
modules of the same form, where $V$ is a complete separated $\kk$-vector space. 

Then we have a functor $\cL\cB\cA \to \{S({\bf LBA})$-modules over 
$\cS_1\}$, taking $\a$ to $S(\a)[[\hbar]]$; the representation of 
$S({\bf LBA})$ is given by $\mu\mapsto \mu_\a$, 
$\delta \mapsto \hbar\delta_\a$. 

We also have a functor $\cL\cB\cA \to \{S({\bf LBA})$-modules over 
$\cS_2\}$, taking $\a$ to $\wh S(\a)[[\hbar]]$; the representation of 
$S({\bf LBA})$ is given by $\mu\mapsto \hbar \mu_\a$, 
$\delta \mapsto \delta_\a$. 

\subsubsection{The prop ${\bf LBA}_f$} 

Define $\wt{\on{LBA}}_f$  as the prop with generators 
$\mu,f,\delta$ and only relations: $\mu$ and $\delta$ satisfy the relations 
of LBA. Then $\wt{\on{LBA}}_f$ has a grading $(\on{deg}_\mu,\on{deg}_\delta,
\on{deg}_f)$ by $\NN^3$. For $F,G\in\on{Ob(Sch)}$ and 
$x\in \wt{\on{LBA}}_f(F,G)$ homogeneous, we have 
\begin{equation} \label{homog:2}
|G| - |F| = \on{deg}_\delta(x) - \on{deg}_\mu(x) 
+ 2 \on{deg}_f(x). 
\end{equation}
Then $\on{deg}_\mu + \on{deg}_\delta + 2 \on{deg}_f$ is a grading of 
$\on{LBA}_f$ by $\NN$. The corresponding filtration therefore satisfies 
condition (a). Since $\on{deg}_\mu$, $\on{deg}_\delta$ and 
$\on{deg}_g$ are $\geq 0$, (\ref{homog:2}) implies that it also satisfies
conditions (b). Since the morphism $\wt{\on{LBA}}_f$ to $\on{LBA}_f$ 
is surjective, the filtration of $\wt{\on{LBA}}_f$ induces a 
filtration of $\on{LBA}_f$
satisfying (a) and (b). We denote by ${\bf LBA}_f$ 
the corresponding completion of $\on{LBA}_f$. 

As before, if $\cL\cB\cA_f$ is the category of pairs $(\a,f_\a)$ of 
Lie bialgebras with twists, we have 
(a) a functor $\cL\cB\cA_f \to \{S({\bf LBA}_f)$-modules over 
$\cS_1\}$, taking $(\a,f_\a)$ to $S(\a)[[\hbar]]$; the representation of 
$S({\bf LBA}_f)$ is given by $\mu\mapsto \mu_\a$, 
$\delta \mapsto \hbar\delta_\a$, $f\mapsto \hbar f_a$; and (b)  
a functor $\cL\cB\cA_f \to \{S({\bf LBA}_f)$-modules over 
$\cS_2\}$, taking $(\a,f_\a)$ to $\wh S(\a)[[\hbar]]$; the representation of 
$S({\bf LBA}_f)$ is given by $\mu\mapsto \hbar \mu_\a$, 
$\delta \mapsto \delta_\a$, $f\mapsto f_a$. 

\subsubsection{The prop ${\bf Cob}$}

$\on{Cob}$ has a grading $(\on{deg}_\mu,\on{deg}_r)$ by $\NN^2$. 
If $F,G\in \on{Ob(Sch)}$ and $x\in \on{Cob}(F,G)$ are homogeneous, then 
$|G| - |F| = 2 \on{deg}_r(x) - \on{deg}_\mu(x)$. 
Then the $\NN$-grading of $\on{Cob}$ by $\on{deg}_\mu + 2 \on{deg}_r$
induces a filtration satisfying (a) and (b). We denote the resulting   
topological prop by ${\bf Cob}$. 

If $\cC ob$ is the category of coboundary Lie bialgebras $(\a,r_\a)$,
we have 
(a) a functor $\cC ob \to \{S($\boldmath$\on{Cob}$\unboldmath$)$-modules over 
$\cS_1\}$, taking $(\a,r_\a)$ to $S(\a)[[\hbar]]$; the representation of 
$S($\boldmath$\on{Cob}$\unboldmath$)$ is given by $\mu\mapsto \mu_\a$, 
$r\mapsto \hbar r_a$; and (b)  
a functor $\cC ob \to \{S($\boldmath$\on{Cob}$\unboldmath$)$-modules over 
$\cS_2\}$, taking $(\a,r_\a)$ to $\wh S(\a)[[\hbar]]$; the representation of 
$S($\boldmath$\on{Cob}$\unboldmath$)$ is given by $\mu\mapsto \hbar\mu_\a$, 
$r\mapsto r_a$. 

\subsubsection{Morphisms between completed props}

The above morphisms $\on{Cob} \to \on{Sch}$, $\on{LBA} \to \on{Sch}$
and $\on{LBA}_f\to \on{Sch}$ are compatible with the filtrations, so they 
induce topological prop morphisms ${\bf Cob} \to {\bf Sch}$, 
${\bf LBA} \to {\bf Sch}$ and ${\bf LBA}_f\to {\bf Sch}$.  

Since $\kappa_1$ preserves the $\NN$-grading, it extends to a 
morphism ${\bf LBA}\to {\bf LBA}_f$ of completed props. 

$\kappa_2$ takes a monomial in $(\mu,\delta)$ of bidegree 
$(a,b)$ to a sum of monomials in $(\mu,\delta,f)$ of degrees
$(a+b'',b',b'')$, where $b'+b'' = b$. The $\NN$-degree of 
such a monomial is $a+b'+3b''\geq a+b$. So $\kappa_2$ preserves the
descending filtrations of both props and extends to a morphism 
${\bf LBA}\to {\bf LBA}_f$.  
 
$\kappa_0$ takes a monomial in $(\mu,\delta,f)$ either to $0$
if the $f$-degree is $>0$, or to the same monomial (which has the 
same $\NN$-degree) otherwise. So $\kappa_0$ preserves the descending filtration 
and extends to a morphism ${\bf LBA}_f \to {\bf LBA}$.  
 
Finally, $\kappa$ takes a monomial in $(\mu,\delta,f)$ of degree 
$(a,b,c)$ (of $\NN$-degree $a+b+2c$) to a monomial in 
$(\mu,\rho)$ of degree $(a+b,b+c)$, of $\NN$-degree $a+3b+2c$.
Since the $\NN$-degree increases, $\kappa$ preserves the descending filtration
and extends to a morphism ${\bf LBA}_f\to {\bf Cob}$. 

\subsection{The props $P_\alpha$} \label{P:alpha}

Let $C$ be a coalgebra in ${\bf Sch}$. This means that 
$C = \wh\oplus_i C_i \in {\bf Sch}$ ($|C_{i}|=i$), and we have 
prop morphisms $C \to C^{\otimes 2}$, $C \to {\bf 1}$ in ${\bf Sch}$, 
such that the two morphisms $C \to C^{\otimes 3}$ coincide, 
and the composed morphisms $C \to C^{\otimes 2} \to C \otimes {\bf 1}
\simeq C$ and $C \to C^{\otimes 2} \to {\bf 1} \otimes C
\simeq C$ are the identity. 
  
Let $P$ be a prop. For $F = \wh\oplus_i F_i$, $G = \wh\oplus_i G_i$ in ${\bf Sch}$, 
we set $P(F,G) = \oplus_{i,j} P(F_i,G_j)$. Then the operations $\circ : 
P(F,G) \otimes P(G,H) \to P(F,H)$ and $\otimes : P(F,G) \otimes P(F',G')
\to P(F\otimes F',G\otimes G')$ are well-defined. Moreover, operations 
$\circ : P(F,G) \otimes {\bf Sch}(G,H) \to P(F,H)$ and 
$\circ : {\bf Sch}(F,G) \otimes P(G,H) \to P(F,H)$ are also well-defined. 

We define a prop $P_C$ by $P_C(F,G) := P(C \otimes F,G)$ for $F,G\in\on{Sch}$. 
The composition of $P_C$ is then defined as the map 
\begin{align*}
& P_C(F,G) \otimes P_C(G,H) \simeq 
P(C \otimes F,G) \otimes P(C \otimes G,H)
\stackrel{(P(\on{id}_{C}) \otimes -)\otimes \on{id}}{\to} 
\\ & P(C^{\otimes 2} \otimes F, C \otimes G) \otimes P(C \otimes G,H)
\to P(C^{\otimes 2} \otimes F, H)
\to P(C \otimes F, H)
\simeq P_C(F,H)
\end{align*}
and the tensor product is defined by  
$$
P_C(F,G) \otimes P_C(F',G') \simeq 
P(C \otimes F,G) \otimes P(C \otimes F',G')
\stackrel{\otimes}{\to}
P(C^{\otimes 2} \otimes F \otimes F',G \otimes G')
\to P(C \otimes F \otimes F',G \otimes G'). 
$$

We then have an isomorphism $P \simeq P_{\bf 1}$ and 
a prop morphism $P \to P_C$ induced by $C\to {\bf 1}$. 

Let us define a $P$-coideal $D$ of $C$ to be the data of 
$D = \wh\oplus_i D_i\in {\bf Sch}$, and morphisms $\alpha \in 
\wh\oplus_i P(C_i,D)$, $\beta \in \wh\oplus_i P(D_i,C\otimes D)$ 
and $\gamma \in \wh\oplus_i P(D_i,D \otimes C)$, such that the 
diagrams
$$
\begin{matrix}
C_i & \stackrel{\Delta_{|C_{i}}}{\to}  & C^{\otimes 2}\\ 
\scriptstyle{\alpha_{|C_{i}}}\downarrow & & 
\downarrow\scriptstyle{\alpha\otimes\on{id}_{C}} \\ 
D & \stackrel{\gamma}{\to} & D \otimes C
\end{matrix}
\quad \on{and} \quad 
\begin{matrix}
C_i & \stackrel{\Delta_{|C_{i}}}{\to}  & C^{\otimes 2}\\ 
\scriptstyle{\alpha_{|C_{i}}}\downarrow & & \downarrow  
\scriptstyle{\on{id}_{C}\otimes\alpha} \\ 
D & \stackrel{\beta}{\to} & C \otimes D
\end{matrix}
$$
commute for each $i$. 
A $P$-coideal $D$ of $C$ may be constructed as follows. 
Let $D'\in \on{Sch}$, let $\alpha'\in P(C,D')$. 
Set $D := D' \otimes C$ and define $\alpha \in \wh\oplus_i 
P(C_i,D)$ 
as the composed morphism $C\stackrel{\Delta}{\to} C^{\otimes 2} 
\stackrel{\alpha' \otimes {\on{id}}}{\to} D' \otimes C = D$. 
We also define the morphism $\gamma\in \wh\oplus_i P(D_i,D \otimes C)$ 
as the composition  
$D = D' \otimes C \stackrel{\on{id} \otimes \Delta}{\to} 
D' \otimes C^{\otimes 2} = D \otimes C$ and $\beta\in \wh\oplus_i
\on{LBA}(D,C\otimes D)$
as the composed morphism $D\to D\otimes C \to C \otimes D$. 

If $D$ is a $P$-coideal of $C$, set $P_D(F,G) := P(D\otimes F,G)$. 
Then for each $(F,G)$, 
we have a morphism $P_D(F,G) \to P_C(F,G)$, such that the collection of 
all $\on{Im}(P_D(F,G) \to P_C(F,G))$ is an ideal of $P_C$. 

We denote by $P_\alpha$ as the corresponding quotient prop. 
We then have $P_\alpha(F,G) = \on{Coker}(P_D(F,G) \to P_C(F,G))$
for any $(F,G)$. 

\subsection{Automorphisms of props} 

For $\xi\in P({\bf id},{\bf id})$, 
$\xi^{\otimes n}\in P(T_{n},T_{n}) = 
\oplus_{\rho,\rho'\in \wh\SG_{n}} \on{Hom}(\pi_{\rho},\pi_{\rho'})
\otimes P(Z_{\rho},Z_{\rho'})$; as $\xi^{\otimes n}$
is $\SG_{n}^{diag}$-invariant ($\SG_{n}^{diag}$ being the diagonal subgroup 
of $\SG_{n}\times\SG_{n}$), we have $\xi^{\otimes n} = 
\oplus_{\rho\in\wh\SG_{n}}\on{id}_{\pi_{\rho}} \otimes \xi_{\rho}$, 
for some $\xi_{\rho}\in P(Z_{\rho},Z_{\rho})$. For $F = (F_{\rho})_{\rho\in
\sqcup_{n\geq 0}\wh\SG_{n}}$, we set $\xi_{F}:= \oplus_{\rho\in
\sqcup_{n\geq 0}\wh\SG_{n}}\on{id}_{\pi_{\rho}} \otimes \xi_{\rho}
\in P(F,F)$. One can prove that $\xi_{F\oplus G}
= \xi_{F}\otimes \xi_{G}$, $\xi_{F\otimes G} = \xi_{F}\otimes \xi_{G}$, 
$(\xi\circ \eta)_{F} = \xi_{F}\circ\eta_{F}$. So if $\xi$
is invertible, so are the $\xi_{F}$, and there is a unique prop automorphism 
$\theta(\xi)$ of $P$, taking $x\in P(F,G)$ to $\xi_{G}\circ x \circ 
\xi_{F}^{-1}$. The map $P({\bf id},{\bf id})^{\times}\to \on{Aut}(P)$
is a group morphism with normal image $\on{Inn}(P)$. We call the elements 
of this image the inner automorphisms of $P$. 

\subsection{Structure of the prop $\on{LBA}$}

\begin{lemma} \label{prop:str}
If $F,G\in\on{Ob(Sch)}$, then we have an isomorphism 
$$
\on{LBA}(F,G) \simeq \oplus_{N\geq 0}
(\on{LCA}(F,T_N) \otimes \on{LA}(T_N,G))_{\SG_N},
$$
with inverse given by $f \otimes g \mapsto g \circ f$ 
(the prop morphisms $\on{LCA} \to \on{LBA}$, 
$\on{LA} \to \on{LBA}$ are understood). 
\end{lemma}

{\em Proof.} This has been proved in the case $F = T_n$, $G = T_m$
in \cite{Enr:univ,Pos}. We then pass to the case of $F = Z_{\rho}$, 
$G = Z_{\sigma}$, $\rho\in\wh\SG_{n}$, $\sigma\in\wh\SG_{m}$
by identifying the isotypic components of this identity under the 
action of $\SG_{n}\times\SG_{m}$. The general case follows by 
linearity. \hfill \qed \medskip 

According to (\ref{schur:recip}), this result may be expressed 
as the isomorphism   
\begin{equation} \label{eq:lba} 
\on{LBA}(F,G) \simeq 
\oplus_{Z\in\on{Irr}(\on{Sch})} \on{LCA}(F,Z) \otimes \on{LA}(Z,G). 
\end{equation} 
We also have: 

\begin{lemma} \label{str:prop:LA}
If $A,B$ are finite sets, then 
$\on{LA}(T_A,T_B) = \oplus_{f : A\to B \on{\ surjective}}
\otimes_{b\in B} \on{LA}(T_{f^{-1}(b)},T_{\{b\}})$, where the inverse map 
is given by the tensor product.  
\end{lemma}

We now prove: 

\begin{lemma} \label{decomp:LBA}
Let $F_1,...,F_n,G_1,...,G_p\in\on{Ob(Sch)}$, then we have a decomposition 
$$
\on{LBA}(\otimes_{i=1}^n F_i, \otimes_{j=1}^p G_j)  =  
\bigoplus_{(Z_{ij})_{i,j}\in \on{Irr}(\on{Sch})^{[n]\times [p]}}
\on{LBA}((F_{i})_{i},(G_{j})_{j})_{(Z_{ij})_{i,j}}, 
$$
where  
$$
\on{LBA}((F_{i})_{i},(G_{j})_{j})_{(Z_{ij})_{i,j}} =  
\bigotimes_{i=1}^n \on{LCA}(F_i, 
\otimes_{\beta = 1}^p Z_{i\beta}) \big) \otimes 
\big( \bigotimes_{j=1}^p \on{LA}(\otimes_{\alpha=1}^n 
Z_{\alpha j},G_j) , 
$$
with inverse given by $(\otimes_i f_i) \otimes (\otimes_j g_j)
\mapsto (\otimes_j g_j) \circ \sigma_{n,p} \circ (\otimes_i f_i)$; 
here $\sigma_{n,p}$ is the braiding isomorphism $\otimes_i(\otimes_j
Z_{ij}) \to \otimes_j(\otimes_i Z_{ij})$. 
\end{lemma}

{\em Proof.} The l.h.s. is equal to 
$$
\oplus_{N\geq 0} \oplus_{(n_{ij}) | \sum_{i,j} n_{ij} = N}
\big( (\otimes_i \on{LCA}(F_i,T_{\sum_j n_{ij}})) 
\otimes (\otimes_j \on{LA}(T_{\sum_i n_{ij}},G_j)) 
\big)_{\prod_{i,j} \SG_{n_{ij}}}. 
$$
(\ref{schur:recip}) then implies the result.
\hfill \qed \medskip 

\subsection{Structure of the prop $\on{LBA}_f$} 

In the construction of Subsection \ref{P:alpha}, we set 
$C := S \circ \wedge^2$, 
$D' = \wedge^3$, and $\alpha' \in \on{LBA}(C,D') = \oplus_{k\geq 0}
\on{LBA}(S^k \circ \wedge^2,\wedge^3)$ has only nonzero components for 
$k = 1,2$; for $k = 1$, this component specializes to $\wedge^2(\a)
\to \wedge^3(\a)$, 
$$
f \mapsto (\delta\otimes\on{id})(f) 
+ (\delta\otimes\on{id})(f)^{231} + (\delta\otimes\on{id})(f)^{312}
$$
and for $k = 2$ it specializes to $S^2(\wedge^2(\a)) \to \wedge^3(\a)$, 
$$
f^{\otimes 2} \mapsto [f^{12},f^{13}] + [f^{12},f^{23}] + [f^{13},f^{23}]. 
$$

Then $\alpha : C\to D$ is a $\on{LBA}$-coideal. 
We denote by $\on{LBA}_\alpha$ the corresponding quotient 
prop $P_\alpha$. 

\begin{proposition} There exists a prop isomorphism 
$\on{LBA}_f \stackrel{\sim}{\to} \on{LBA}_\alpha$. 
\end{proposition} 

{\em Proof.} 
Using the presentation of $\on{LBA}_f$, one checks that there is 
a unique prop morphism $\on{LBA}_f \to \on{LBA}_\alpha$, taking $\mu$
to the class of $\mu\in \on{LBA}(\wedge^2,{\bf id}) \subset
\oplus_{k\geq 0} \on{LBA}((S^k\circ\wedge^2)\otimes \wedge^2,{\bf id})
= \on{LBA}((S\circ\wedge^2)\otimes \wedge^2,{\bf id})$, taking
$\delta$ to the class of $\delta\in \on{LBA}({\bf id},\wedge^2) \subset
\oplus_{k\geq 0} \on{LBA}((S^k\circ\wedge^2)\otimes {\bf id},\wedge^2)
= \on{LBA}((S\circ\wedge^2)\otimes {\bf id},\wedge^2)$,
and taking $f$ to the class of 
$\on{id}_{\wedge^2} \in \on{LBA}(\wedge^2,\wedge^2) \subset
\oplus_{k\geq 0} \on{LBA}(S^k\circ\wedge^2,\wedge^2)
= \on{LBA}((S\circ\wedge^2)\otimes {\bf 1},\wedge^2)$.  

We now construct a prop morphism $\on{LBA}_\alpha \to \on{LBA}_f$. 

We construct a linear map 
$\on{LBA}((S^k\circ \wedge^2) \otimes F,G) \to \on{LBA}_f(F,G)$
as follows: using the prop morphism $\on{LBA} \to \on{LBA}_f$ given by 
$\mu,\delta\mapsto \mu,\delta$,  
we get a linear map $\on{LBA}((S^k\circ \wedge^2) \otimes F,G) \to 
\on{LBA}_f( (S^k \circ \wedge^2) \otimes F,G)$. We have an element 
$f^{\otimes k} \in \on{LBA}_f(S^k \circ {\bf 1},S^k \circ \wedge^2)$
so the operation $x\mapsto x\circ (f^{\otimes k} \otimes \on{id}_F)$ 
is a linear map 
$$
\on{LBA}_f( (S^k \circ \wedge^2) \otimes F,G) \to 
\on{LBA}_f( (S^k \circ {\bf 1}) \otimes F,G) \simeq \on{LBA}_f(F,G).  
$$ 
The composition of these maps is a linear map 
$\on{LBA}((S^k\circ \wedge^2) \otimes F,G) \to \on{LBA}_f(F,G)$. 
Summing up these maps, we get a linear map 
$\on{LBA}((S\circ \wedge^2) \otimes F,G) \to \on{LBA}_f(F,G)$, and 
one checks that 
it factors through a linear map  $\on{LBA}_\alpha(F,G) \to \on{LBA}_f(F,G)$. 
One also checks that this map is compatible with the prop operations, so it 
is a prop morphism.  

We now show that the composed morphisms $\on{LBA}_f \to \on{LBA}_\alpha 
\to \on{LBA}_f$ and $\on{LBA}_\alpha \to \on{LBA}_f \to \on{LBA}_\alpha$ 
are both the identity. 

In the case of $\on{LBA}_f \to \on{LBA}_\alpha \to \on{LBA}_f$, one 
shows that the composed map takes each generator of $\on{LBA}_f$ to 
itself, hence is the identity. 

Let us show that $\on{LBA}_\alpha \to \on{LBA}_f \to \on{LBA}_\alpha$ 
is the identity. We already defined the prop $\wt{\on{LBA}}_f$. 
Then we have a canonical prop morphism $\wt{\on{LBA}}_f \to \on{LBA}_f$. 
We also have prop morphisms 
$\on{LBA}_{S\circ\wedge^2} \to \wt{\on{LBA}}_f$ and $\wt{\on{LBA}}_f\to
\on{LBA}_{S\circ\wedge^2}$, defined similarly to 
$\on{LBA}_\alpha \to \on{LBA}_f$ and $\on{LBA}_f\to
\on{LBA}_\alpha$. We then have commuting squares 
$$
\begin{matrix}
\on{LBA}_{S\circ\wedge^2} & \to & \wt{\on{LBA}}_f \\
\downarrow  & & \downarrow\\
\on{LBA}_\alpha & \to & \on{LBA}_f 
\end{matrix} \quad \on{and} \quad 
\begin{matrix}
\wt{\on{LBA}}_f & \to & \on{LBA}_{S\circ\wedge^2} \\
\downarrow  & & \downarrow \\
\on{LBA}_f & \to & \on{LBA}_\alpha
\end{matrix}
$$
One checks that the composed morphism $\on{LBA}_{S\circ\wedge^2} 
\to \wt{\on{LBA}}_f \to \on{LBA}_{S\circ\wedge^2}$ 
is the identity, which implies that 
$\on{LBA}_\alpha \to \on{LBA}_f \to \on{LBA}_\alpha$ is the identity. 
\hfill \qed\medskip 

In what follows, we will use the above isomorphism to identify 
$\on{LBA}_\alpha$ with $\on{LBA}_f$. The main output of this 
identification is the construction of a grading of 
$\on{LBA}_{f}(\otimes_{i}F_{i},\otimes_{j}G_{j})$
by families $(Z_{ij})$ in $\on{Irr(Sch)}$, since as we now 
show, props of the form $\on{LBA}_{\alpha}$ all give rise to
such a grading. 

Let $C,D\in \on{Ob}({\bf Sch})$ and let 
$\alpha \in \wh\oplus_i \on{LBA}(C_i,D)$. 

\begin{proposition} \label{str:LBAf}
Let $F_1,...,F_n,G_1,...,G_p\in \on{Ob(Sch)}$. 
Set $F:= \otimes_i F_i$, $G := \otimes_j G_j$.  
For $Z = (Z_{ij})_{i\in [n],j\in [p]}$ a map $[n]\times [p]\to\on{Irr(Sch)}$, 
set 
$$
\on{LBA}_{C}((F_{i})_{i},(G_{j})_{j})_{Z} := 
\oplus_{\tilde Z\in \on{Irr(Sch)}^{(\{0\}\cup [n])\times[p]}|
\tilde Z_{|[n]\times[p]}=Z}
\on{LBA}([C,(F_{i})_{i}],(G_{j})_{j})_{\tilde Z}, 
$$
where $[C,(F_{i})_{i}]\in \on{Ob(Sch)}^{\{0\}\cup [n]}$
is the extension of $(F_{i})_{i}$ defined by $0\mapsto C$. 
 
Then  $\on{LBA}_C(F,G) = 
\oplus_{Z\in\on{Irr(Sch)}^{[n]\times [p]}}
\on{LBA}_{C}((F_{i})_{i},(G_{j})_{j})_{Z}$. 

Moreover, the map $\on{LBA}_D(F,G) \to \on{LBA}_C(F,G)$ preserves the 
grading by $\on{Irr(Sch)}^{[n]\times[p]}$. 
The cokernel of this map therefore inherits a grading 
$$
\on{LBA}_\alpha(F,G) = \oplus_{Z\in\on{Irr(Sch)}^{[n]\times[p]}}
\on{LBA}_{\alpha}((F_{i})_{i},(G_{j})_{j})_{Z}. 
$$
\end{proposition}

{\em Proof.} The first statement follows from Lemma \ref{decomp:LBA} 
(with $F_{0}=C$). Let us prove the second statement. 
Consider the sequence of maps  
\begin{align*}
& \on{LCA}(D,\otimes_{j\in [p]}Z_{0j}) \stackrel{\alpha\otimes -}{\to}
\on{LBA}(C,D) \otimes \on{LBA}(D,\otimes_{j\in [p]} Z_{0j})
\stackrel{\circ}{\to} \on{LBA}(D,\otimes_{j\in [p]} Z_{0j}) 
\\ & \simeq 
\oplus_{(Z'_{0j})_{j}\in\on{Irr(Sch)}^{[p]}} 
\on{LCA}(C,\otimes_{j\in[p]} Z'_{0j})
\otimes (\otimes_{j\in[p]}\on{LA}(Z'_{0j},Z_{0j})), 
\end{align*}
$$
\kappa\mapsto \oplus_{(Z'_{0j})_{j}\in \on{Irr(Sch)}^{[p]}} 
\sum_{\alpha}\kappa'_{\alpha}\otimes (\otimes_{j\in [p]}\lambda'_{j,\alpha}),  
$$
where the first map uses the prop morphism $\on{LCA}\to\on{LBA}$. 

For $\tilde Z \in \on{Irr(Sch)}^{(\{0\}\cup[n])\times [p]}$,  
$\tilde Z_{|[n]\times [p]}$ is its restriction to $[n]\times [p]$. 
The map $\on{LBA}(D\otimes F,G)\simeq 
\on{LBA}_D(F,G) \to \on{LBA}_C(F,G)$ restricts to 
$$
\on{LBA}([D,(F_{i})_{i}],(G_{j})_{j})_{\tilde Z} 
\to \on{LBA}_C((F_{i})_{i},(G_{j})_{j})_{\tilde Z_{|[n]\times [p]}} 
$$
$$
(\otimes_{j\in [p]}\lambda_j) \circ \sigma_{n+1,p} \circ 
(\kappa \otimes (\otimes_{i\in [n]}\kappa_i))
\mapsto
\oplus_{(Z'_{0j})_{j}\in \on{Irr(Sch)}^{[p]}} 
\sum_{\alpha}\Big( \otimes_{j\in[p]} \big( \lambda_j \circ 
(\lambda'_{j,\alpha} \otimes \on{id}_{\otimes_{i\in [n]}Z_{ij}})\big) 
\Big) \circ \sigma_{n+1,p} \circ \kappa'_{\alpha}.  
$$
Summing up over $(Z_{0j})_{j}\in\on{Irr(Sch)}^{[p]}$, we get the result. 
\hfill \qed \medskip 

\subsection{Partial traces on $\Pi^{0}_{\on{LBA}}$, $\Pi^{0}_{\on{LBA}_{f}}$}

Recall that for $F,G\in\on{Ob}(\on{Sch}_{(1)})$, $\Pi^{0}_{\on{LBA}}(F,G)=
\on{LBA}(c(F),c(G))$. For $F,G\in\on{Ob}(\on{Sch}_{p,q})$, we introduce
a grading of $\Pi^{0}_{\on{LBA}}(F,G)$ by $\cG_{0}([p],[q])$ as follows. 
Assume first that $F,G$ are simple, so $F=\boxtimes_{i=1}^{n}Z_{\rho_{i}}$, 
$G= \boxtimes_{j=1}^{p}Z_{\sigma_{p}}$. If $Z = (Z_{ij})_{i,j}\in 
\on{Irr(Sch)}^{[n]\times[p]}$, we define the support of $Z$ as
$\on{supp}(Z):= \{(i,j)|Z_{ij}\neq {\bf 1}\}$. Then 
$$
\Pi^{0}_{\on{LBA}}(\boxtimes_{i=1}^{n}Z_{\rho_{i}},
\boxtimes_{j=1}^{p}Z_{\sigma_{j}})_{S}:= 
\oplus_{Z\in \on{Irr(Sch)}^{[n]\times [p]}|
\on{supp}(Z)=S}
\on{LBA}((Z_{\rho_{i}}),(Z_{\sigma_{j}}))_{Z}. 
$$
If $F = (F_{\rho_{1},...,\rho_{n}})$, $G = (G_{\sigma_{1},...,\sigma_{p}})$, 
$$
\Pi^{0}_{\on{LBA}}(F,G)_{S}:= \oplus_{(\rho_{i}),(\sigma_{j})}
\on{Vect}(F_{(\rho_{i})},G_{(\sigma_{j})})\otimes \Pi^{0}_{\on{LBA}}
(\boxtimes_{i=1}^{n} Z_{\rho_{i}},\boxtimes_{j=1}^{p} Z_{\sigma_{j}})_{S}.
$$ 

\begin{proposition} \label{compat:G0}
This grading is compatible with the monoidal category structure of $\cG_{0}$, namely: 
for $F,G,H\in\on{Ob}(\on{Sch}_{p,q,r})$, 
$$
\Pi^{0}_{\on{LBA}}(G,H)_{S'} \circ \Pi^{0}_{\on{LBA}}(F,G)_{S}\subset
\Pi^{0}_{\on{LBA}}(F,H)_{S'\circ S}, 
$$
for $F_{i},G_{i}\in \on{Ob}(\on{Sch}_{p_{i},q_{i}})$, 
$$
\Pi^{0}_{\on{LBA}}(F_{1},G_{1})_{S_{1}} \boxtimes 
\Pi^{0}_{\on{LBA}}(F_{2},G_{2})_{S_{2}} \subset 
\Pi^{0}_{\on{LBA}}(F_{1}\boxtimes F_{2},
G_{1}\boxtimes G_{2})_{S_{1}\otimes S_{2}}
$$
(here $\boxtimes$ denotes the tensor product operation of $\Pi^{0}_{\on{LBA}}$), 
and $\beta_{F,G}\in \Pi^{0}_{\on{LBA}}(F\boxtimes G,
G\boxtimes F)_{\beta_{[n],[p]}}$. 
\end{proposition}

{\em Proof.} The only nontrivial statement is the first one. Let $Z,Z'$ be such that 
$\on{supp}(Z)=S$, $\on{supp}(Z')=S'$. For $F,G,H$ simple, the composition factors as 
\begin{align*}
& \Pi^{0}_{\on{LBA}}(F,G)_{Z}\otimes \Pi^{0}_{\on{LBA}}(G,H)_{Z'}
\\
 & \to (\otimes_{i}\on{LCA}(Z_{\rho_{i}},\otimes_{j}Z_{ij}))
\otimes (\otimes_{j}\on{LA}(\otimes_{i}Z_{ij},Z_{\sigma_{j}}))
\otimes (\otimes_{j}\on{LCA}(Z_{\sigma_{j}},\otimes_{k}Z'_{jk}))
\otimes (\otimes_{k}\on{LA}(\otimes_{j}Z'_{jk},Z_{\tau_{k}}))
\\
 & 
 \to \oplus_{(Z''_{ijk})}
 (\otimes_{i}\on{LCA}(Z_{\rho_{i}},\otimes_{j}Z_{ij}))
\otimes (\otimes_{j,i}\on{LCA}(Z_{ij},\otimes_{k}Z''_{ijk}))
\\
 & \otimes (\otimes_{j,k}\on{LA}(\otimes_{i}Z''_{ijk},Z'_{jk}))
\otimes (\otimes_{k}\on{LA}(\otimes_{j}Z'_{jk},Z_{\tau_{k}}))
\\
& 
\to \oplus_{(Z''_{ijk})} (\otimes_{i}\on{LA}(Z_{\rho_{i}},\otimes_{j,k}Z''_{ijk}))
\otimes (\otimes_{j}\on{LCA}(\otimes_{i,j}Z''_{ijk},Z_{\sigma_{j}}))
\\
& \simeq
\oplus_{(Z''_{ijk})} (\otimes_{i}\on{LA}(Z_{\rho_{i}},\otimes_{k}Z'''_{ik}))
\otimes (\otimes_{j}\on{LCA}(\otimes_{i}Z'''_{ik},Z_{\sigma_{j}}))
\to \Pi^{0}_{\on{LBA}}(F,H),  \end{align*}
where the first map is the decomposition map, the second map is the tensor product over 
$j$ of the $j$th exchange map (composition followed by decomposition)
$\on{LA}(\otimes_{i}Z_{ij},Z_{\sigma_{j}}) \otimes \on{LCA}(Z_{\sigma_{j}},
\otimes_{k}Z'_{jk}) \to \on{LBA}(\otimes_{i}Z_{ij},\otimes_{k}Z'_{jk})
\to \oplus_{(Z''_{ijk})} (\otimes_{i}\on{LCA}(Z_{ij},\otimes_{k}Z''_{ijk}))
\otimes (\otimes_{k}\on{LA}(\otimes_{i}Z''_{ijk},Z'_{jk}))$, the third
map is composition in LA and LCA, the fourth map is obtained by 
$Z'''_{ik}:= \otimes_{j}Z''_{ijk}$, the fifth map is composition.

If $(i,k)\in S'\circ S$, then for some $j\in J$, $(i,j)\in S$ and $(j,k)\in S'$; 
so $Z_{ij},Z'_{jk}$ are $\neq {\bf 1}$. So the component of the target of the 
$j$th exchange map corresponding to $(i,k)\mapsto Z'''_{ijk}={\bf 1}$, is zero; 
hence the component of the composition of the 3 first maps of the diagram, 
coresponding to  $(i,k)\mapsto Z'''_{ijk}={\bf 1}$, is zero. It follows
that $Z''_{ik}$ in the forelast vector space are sums of objects of 
$\on{Irr(Sch)}$ with degree $>0$. 

On the other hand, if $(i,k)\notin S'\circ S$, then for any $j\in J$, 
$Z_{ij}$ or $Z'_{jk}={\bf 1}$. Then the component of the 
target of the $j$th exchange map corresponding to any 
$(i,k)\mapsto Z'''_{ijk}$ different from $(i,k)\mapsto {\bf 1}$, is zero; 
si $Z''_{ik}$ is the forelast vector space are sums of copies of ${\bf 1}$. 

It follows that the image of the overall map is contained in 
$\Pi_{\on{LBA}}(F,G)_{S'\circ S}$. 
\hfill \qed \medskip 

Let $F,G,H\in\on{Ob}(\on{Sch}_{p,q,r})$, and let us
define the diagram 
$$\Pi^{0}_{\on{LBA}}(F\boxtimes H,G\boxtimes H)\supset 
\Pi^{0}_{\on{LBA}}(F,G|H) \stackrel{\on{tr}_{H}}{\to}
\Pi^{0}_{\on{LBA}}(F,G)
$$
(the general case is then derived by linearity). Let us first assume that $F,G,H$ are 
simple, so $F = Z_{(\rho_{i})} := \boxtimes_{i\in I}Z_{\rho_{i}}$, 
$G=Z_{(\sigma_{j})} = \boxtimes_{j\in J}Z_{\sigma_{j}}$, 
$H=Z_{(\tau_{k})} = \boxtimes_{k\in K}Z_{\tau_{k}}$ ($I,J,K$ are ordered
sets of cardinality $p,q,r$, and $(\rho_{i}),(\sigma_{j}),(\tau_{k})$ are 
maps $I,J,K\to\sqcup_{n\geq 0}\wh\SG_{n}$). 
Recall that $\Pi^{0}_{\on{LBA}}(F\boxtimes H,G\boxtimes H)= 
\oplus_{S\in\cG_{0}\Pi^{0}_{\on{LBA}}
(I\otimes K,J\otimes K)}(F\boxtimes H,G\boxtimes H)_{S}$. 
We then set 
$$
\Pi^{0}_{\on{LBA}}(F,G|H):= \oplus_{S\in
\cG_{0}(I,J|K)}\Pi^{0}_{\on{LBA}}(F\boxtimes H,G\boxtimes H)_{S}.
$$ 
We then define the linear map 
$$
\on{tr}_{H} : 
\Pi^{0}_{\on{LBA}}(F,G|H)_{S} =\oplus_{Z|\on{supp}(Z)=S}
\Pi^{0}_{\on{LBA}}(F\boxtimes H,G\boxtimes H)_{Z}
\to \Pi^{0}_{\on{LBA}}(F,G)_{\on{tr}_{K}(S)}
$$
as follows. Recall from section \ref{sect:part:traces}
the order relation $\prec$ on $K$, the total order relation 
$<$ on $K$, its extension to a relation on $I\sqcup K\sqcup J$, 
the numbering $K = \{k_{1},...,k_{|K|}\}$, the sets $K_{\alpha}$. 

Let $Z=(Z_{uv})_{(u,v)\in (I\sqcup K)\times (K\sqcup J)}$ 
be such that $\on{supp}(Z)=S$. For $k\in K$, set 
$H_{k}:= \boxtimes_{x\in K_{k}}Z(x)$, where 
$Z(k):= Z_{\tau_{k}}$, $Z(u,v):= Z_{\rho_{u}\sigma_{v}}$ for 
$(u,v)\in (I\sqcup K)\times (K\sqcup J)$
(we extend $\rho$ to $I\sqcup K$ by $\rho_{k}:= \tau_{k}$ and 
$\sigma$ to $K\sqcup J$ by $\sigma_{k}:= \tau_{k}$). Set 
also $H_{k,k+1}:= Z_{k,k+1}\boxtimes (\boxtimes_{(u,v)\in 
K'_{\alpha,\alpha+1}\sqcup K''_{\alpha,\alpha+1}\sqcup K_{\alpha,\alpha+1}}
Z(u,v))$. 

Set also $H_{0}:= F$, $H_{01}:= Z(k_{1}) \boxtimes 
(\boxtimes_{(u,v)\in I \times (K\sqcup J), u\prec v}Z(u,v))$, and $H_{|K|+1}=G$, 
$H_{|K|,|K|+1}=Z(k_{|K|})\boxtimes (\boxtimes_{(u,v)\in (I\sqcup K)\times J}
Z(u,v))$. 

Then $\on{tr}_{H}$ is the sum over $Z$ with $\on{supp}(Z)=S$ of the compositions
\begin{align} \label{def:tr}
& \nonumber \Pi^{0}_{\on{LBA}}(F\boxtimes H,G\boxtimes H)_{Z} = 
(\otimes_{u\in I\sqcup K}\on{LCA}(Z_{\rho_{u}},\otimes_{v\in K\sqcup J}
Z_{uv})) \otimes (\otimes_{v\in K\sqcup J}\on{LA}(\otimes_{u\in I\sqcup K}
Z_{uv},Z_{\sigma_{v}})) 
\\ & \to 
 \bigotimes_{k=0}^{|K|} \Pi^{0}_{\on{LCA}}(H_{k},H_{k,k+1}) \otimes 
\Pi^{0}_{\on{LA}}(H_{k,k+1},
 H_{k+1}) \stackrel{\circ}{\to} \Pi^{0}_{\on{LBA}}(F,G)
\end{align}
(one checks that this map is independent on the ordering of $K$). 

The sum of these maps takes $\Pi^{0}_{\on{LBA}}(F\boxtimes H,G\boxtimes H)_{S}$
to $\otimes_{k=0}^{|K|} \Pi^{0}_{\on{LBA}}(H_{k},H_{k+1})_{S_{K_{k},K_{k+1}}}$, 
where $K_{0}:= I$, $K_{k+1}:= J$ (with the notation of Section \ref{sect:part:traces}), 
therefore the image of the above map is contained in 
$\Pi^{0}_{\on{LBA}}(F,G)_{\on{tr}_{K}(S)}$. 

If now $F,G,H$ are arbitrary elements of $\on{Ob}(\on{Sch}_{p,q,r})$, namely 
$F = (F_{\rho_{1},...,\rho_{p}})$, $G = (G_{\sigma_{1},...,\sigma_{q}})$, 
$H = (H_{\tau_{1},...,\tau_{r}})$, then 
$$
\Pi^{0}_{\on{LBA}}(F\boxtimes H,G\boxtimes H) = 
\oplus_{(\rho_{i}),(\sigma_{j}),(\tau_{k}),(\tau'_{k})}
\on{Vect}(F_{(\rho_{i})}\otimes H_{(\tau_{k})},G_{(\sigma_{j})}\otimes
H_{(\tau'_{k})}) \otimes \Pi^{0}_{\on{LBA}}(Z_{(\rho_{i})} \boxtimes 
Z_{(\tau_{k})},Z_{(\sigma_{j})}\boxtimes Z_{(\tau'_{k})}); 
$$
here $Z_{(\rho_{i})} = \boxtimes_{i\in I}Z_{\rho_{i}}$, etc. 
Then $\Pi^{0}_{\on{LBA}}(F,G|H)$ is homogeneous w.r.t. this decomposition; 
its components for $(\tau_{k})\neq (\tau'_{k})$ coincide with those of 
$\Pi^{0}_{\on{LBA}}(F\boxtimes H,G\boxtimes H)$ and the component for 
$(\tau'_{k})=(\tau_{k})$ is 
$\on{Vect}(F_{(\rho_{i})}\otimes H_{(\tau_{k})},G_{(\sigma_{j})}\otimes
H_{(\tau_{k})}) \otimes \Pi^{0}_{\on{LBA}}(Z_{(\rho_{i})} , Z_{(\sigma_{j})} | 
Z_{(\tau'_{k})})$; the restriction of $\on{tr}_{H}$ to a components of the 
first kind is $0$ and its restriction to the last component is the tensor product
the partial trace with $\on{tr}_{(Z_{\tau_{k}})}$. 
Then one checks that the diagrams $\Pi^{0}_{\on{LBA}}(F\boxtimes H,G\boxtimes H)
\supset \Pi^{0}_{\on{LBA}}(F,G|H)\to \Pi^{0}_{\on{LBA}}(F,G)$ define a partial trace on 
$\Pi^{0}_{\on{LBA}}$. 

Let us define now a partial trace on $\Pi^{0}_{\on{LBA}_{f}}$ (and more generally on 
the $\Pi^{0}_{\on{LBA}_{\alpha}}$). Let $C$ be a coalgebra in ${\bf Sch}$, 
we have for $F,G\in\on{Ob}(\on{Sch}_{p,q})$, 
$\Pi^{0}_{\on{LBA}_C}(F,G) = \Pi^{0}_{\on{LBA}}(C\boxtimes F,G)$. A partial trace is 
then defined on $\Pi^{0}_{\on{LBA}_{C}}$ as follows: 
$\Pi^{0}_{\on{LBA}_{C}}(F,G|H) \subset \Pi^{0}_{\on{LBA}_{C}}(F\boxtimes H,
G\boxtimes H)$ is $\Pi^{0}_{\on{LBA}_{C}}(F,G|H) := \Pi^{0}_{\on{LBA}}
(C\boxtimes F,G|H)$
and $\Pi^{0}_{\on{LBA}_{C}}(F,G|H) \stackrel{\on{tr}_{H}}{\to}
\Pi^{0}_{\on{LBA}_{C}}(F,G)$ coincides with  
$\Pi^{0}_{\on{LBA}}(C\boxtimes F,G|H) \stackrel{\on{tr}_{H}}{\to}
\Pi^{0}_{\on{LBA}}(C\boxtimes F,G)$. One checks that this defines 
a partial trace on $\Pi^{0}_{\on{LBA}_{C}}$.

If $Z\in \on{Irr(Sch)}^{[n]\times [p]}$, we set $\Pi^{0}_{\on{LBA}_{C}}(F,G)_{Z}:= 
\oplus_{\tilde Z}\Pi^{0}_{\on{LBA}}(C\boxtimes F,G)_{\tilde Z}$, where the sum is
over the $\tilde Z:(\{0\}\sqcup [n])\times [p]\to \on{Irr(Sch)}$, such that 
$\tilde Z_{|[n]\times [p]}=Z$. We also set, for $S\subset [n]\times [p]$,
$\Pi^{0}_{\on{LBA}_{C}}(F,G)_{S}:= \oplus_{Z|\on{supp}(Z)=S}
\Pi^{0}_{\on{LBA}_{C}}(F,G)_{Z}$. 
Then: 

\begin{lemma} 
The properties of $\Pi^{0}_{\on{LBA}}$ extend to $\Pi^{0}_{\on{LBA}_{C}}$, 
namely $\Pi^{0}_{\on{LBA}_{C}}(G,H)_{S'}\circ \Pi^{0}_{\on{LBA}_{C}}(F,G)_{S}
\subset \Pi^{0}_{\on{LBA}_{C}}(F,H)_{S'\circ S}$, $\Pi^{0}_{\on{LBA}_{C}}(F,G)_{S}
\boxtimes \Pi^{0}_{\on{LBA}_{C}}(F',G')_{S'}\subset 
\Pi^{0}_{\on{LBA}_{C}}(F\boxtimes
F',G\boxtimes G')_{S\otimes S'}$, 
$\beta_{F,G}\in\Pi^{0}_{\on{LBA}_{C}}(F\boxtimes G,G
\boxtimes F)_{\beta_{[n],[p]}}$, and for $S\in \cG_{0}(I,J|K)$, 
$\on{tr}_{H}(\Pi^{0}_{\on{LBA}_{C}}(F,G|H)_{S}) \subset 
\Pi^{0}_{\on{LBA}_{C}}(F,G)_{\on{tr}_{K}(S)}$. 
\end{lemma}

{\em Proof.} $\Pi_{\on{LBA}_{C}}^{0}(F,G)_{S} = \oplus_{\tilde S\in 
(\{0'\}\sqcup [n])\times [p] | \tilde S \cap ([n]\times [p]) = S} 
\Pi^{0}_{\on{LBA}}(C\boxtimes F,G)_{\tilde S}$. 
In the same way, 
$\Pi_{\on{LBA}_{C}}^{0}(G,H)_{S} = \oplus_{\tilde S'\in 
(\{0''\}\sqcup [p])\times [q] | \tilde S' \cap ([p]\times [q]) = S'} 
\Pi^{0}_{\on{LBA}}(C\boxtimes G,H)_{S'}$. Then 
for $\tilde S\subset (\{0'\}\sqcup [n])\times [p]$, $\tilde S'
\subset (\{0''\}\sqcup [p])\times [q]$, let 
$\tilde S' * \tilde S := 
\{ (i,k)\in (\{0',0''\}\sqcup [n])\times p|
(i,k)\in (\{0'\}\sqcup [n])\times [q]$ and there exists $k\in [p]$ 
with $(i,j)\in S$, $(j,k)\in S'$ or 
$i=0''$ and $(i,k)\in S'\}$. Then the composition  
$\Pi_{\on{LBA}}^{0}(C\boxtimes F,G)\otimes 
\Pi_{\on{LBA}}^{0}(C\boxtimes G,H) \to \Pi_{\on{LBA}}^{0}(C\boxtimes 
F,H)$ maps $\Pi_{\on{LBA}}^{0}(C\boxtimes F,G)_{\tilde S} \otimes 
\Pi_{\on{LBA}}^{0}(C\boxtimes G,H)_{\tilde S'}$ to 
$\Pi_{\on{LBA}}^{0}(C\boxtimes F,H)_{(\tilde S'*\tilde S) \circ \emptyset}
\oplus 
\Pi_{\on{LBA}}^{0}(C\boxtimes F,H)_{(\tilde S'*\tilde S) \circ \Delta_{0,0'0''}}$, 
where $\emptyset,\Delta_{0,0'0''} 
\in \cG_{0}(\{0\}\sqcup [n],\{0',0''\}\sqcup [n])$ are 
$\emptyset$ and $\Delta_{00'} = \{(0,0'),(0,0'')\}$.
Now both $(\tilde S'*\tilde S)\circ \emptyset$ and
$(\tilde S'\circ \tilde S)\circ \Delta_{0,0'0''}$ are elements of 
$\cG_{0}(\{0\}\sqcup [n],[q])$, with their intersection with $[n]\times [p]$
equal to $S'\circ S$. This proves the first statement. The other 
statements are proved in the same way. \hfill \qed \medskip

If $D\in\on{Ob}({\bf Sch})$, we similarly set $\Pi^{0}_{D}(F,G):= 
\Pi^{0}_{\on{LBA}}(D\boxtimes F,G)$. The diagram $\Pi^{0}_{D}(F\boxtimes H,
G\boxtimes H)\supset \Pi^{0}_{D}(F,G|H)\stackrel{\on{tr}_{H}}{\to}
\Pi^{0}_{D}(F,G)$ is defined as above. We define $\Pi^{0}_{D}(F,G)_{Z}$ and 
$\Pi^{0}_{D}(F,G)_{S}$
as above. The following properties generalize to this more general setup: 
$\beta_{F,G}\in\Pi^{0}_{D}(F\boxtimes G,G
\boxtimes F)_{\beta_{[n],[p]}}$, and for $S\in \cG_{0}(I,J|K)$, 
$\on{tr}_{H}(\Pi^{0}_{D}(F,G|H)_{S}) \subset 
\Pi^{0}_{D}(F,G)_{\on{tr}_{K}(S)}$. 

\begin{lemma}
$\alpha\in \wh\oplus_{i}\on{LBA}(C_{i},D)$ induces a linear map 
$\Pi^{0}_{D}(F,G)\to \Pi^{0}_{\on{LBA}_{C}}(F,G)$, which is 
compatible with the gradings by $\on{Irr(Sch)}^{[n]\times [p]}$ (and therefore
also by $\cG_{0}([n],[p])$). Then $\Pi^{0}_{\on{LBA}_\alpha}(F,G)
= \on{Coker}(\Pi^{0}_{D}(F,G)
\to \Pi^{0}_{\on{LBA}_{C}}(F,G))$, and $\Pi^{0}_{\on{LBA}_\alpha}(F,G) 
= \oplus_{S\in
\cG_{0}(I,J)} \Pi_{\on{LBA}_\alpha}(F,G)_{S}$. 

For each $S\in \cG_{0}(I,J|K)$, the diagram 
$$
\begin{matrix}
\Pi^{0}_{D}(F\boxtimes H,G\boxtimes H)_{S} 
& \stackrel{\on{tr}_{H}}{\to}& \Pi^{0}_{D}(F,G)_{\on{tr}_{K}(S)} \\
\downarrow & & \downarrow \\
\Pi^{0}_{\on{LBA}_{C}}(F\boxtimes H,G\boxtimes H)_{S} 
& \stackrel{\on{tr}_{H}}{\to}& 
\Pi^{0}_{\on{LBA}_{C}}(F,G)_{\on{tr}_{K}(S)} 
\end{matrix}
$$
commutes; its vertical cokernel is a linear map $\on{tr}_{H}: 
\Pi^{0}_{\on{LBA}_\alpha}(F\boxtimes H,G\boxtimes H)_{S}\to 
\Pi^{0}_{\on{LBA}_\alpha}(F,G)_{\on{tr}_{K}(S)}$. We set 
$\Pi^{0}_{\on{LBA}_\alpha}(F,G|H):= \oplus_{S\in \cG_{0}(I,J|K)}
\Pi^{0}_{\on{LBA}_{\alpha}}(F\boxtimes H,G\boxtimes H)_{S}$, 
then we have a diagram $\Pi_{\on{LBA}_{\alpha}}^{0}(F\boxtimes H,
G\boxtimes H)\supset \Pi_{\on{LBA}_{\alpha}}^{0}(F,G|H)
\stackrel{\on{tr}_{H}}{\to}\Pi_{\on{LBA}_{\alpha}}^{0}(F,G)$. 

If $\alpha:C\to D$ is a $\on{LBA}$-coideal, then the multiprop 
$\Pi^{0}_{\on{LBA}_{\alpha}}$
is graded by $\cG_{0}$, and $(\on{tr}_{H})$ is a partial trace on 
$\Pi^{0}_{\on{LBA}_{\alpha}}$, compatible with this grading. 
\end{lemma}

{\em Proof.} The first statement is a consequence of Proposition \ref{str:LBAf}. 
The commutativity of the diagram follows from the fact that for 
$x\in\Pi^{0}_{\on{LBA}_{D}}(F,G|H)$, $x\circ (\alpha\boxtimes
\on{id}_{F\boxtimes H})\in \Pi^{0}_{\on{LBA}_{C}}(F,G|H)$
and $\on{tr}_{H}(x\circ (\alpha\boxtimes\on{id}_{F\boxtimes H})) 
= \on{tr}_{H}(x) \circ (\alpha\boxtimes\on{id}_{F})$. The remaining
properties follow from those of $\Pi^{0}_{\on{LBA}_{C}}$. 
\hfill \qed \medskip 
 
\begin{remark}
One also checks that the partial trace on $\cG_{0}$, as well as 
its counterparts on $\Pi^{0}_{\on{LBA}_{\alpha}}$, 
have the following cyclicity properties. If $S\in \cG_{0}(U\otimes I,V'\otimes J)$
and $S'\in \cG_{0}(V\otimes J,U'\otimes I)$, then   
$S'\circ (\beta_{V,V'}\otimes \on{id}_{J})\circ S \in 
\cG_{0}(V\otimes U,V'\otimes U'|I)$ iff 
$S\circ (\beta_{U,U'}\otimes\on{id}_{I})\circ S'
\in \cG_{0}(U\otimes V,U'\otimes V'|J)$, and we then have 
$\on{tr}_{I}(S'\circ (\beta_{V,V'}\otimes \on{id}_{J})\circ S)=
\beta_{V',U'}\circ \on{tr}_{J}(S\circ (\beta_{U,U'}\otimes \on{id}_{I})\circ S')
\circ \beta_{U,V}$. 
If now $S,S'$ are as above, $F_{U} = \boxtimes_{u\in I}Z_{\rho_{u}}$, etc., 
and $x\in \Pi_{\on{LBA}}^{0}(F_{U}\boxtimes F_{I},F_{V'}\boxtimes 
F_{J})_{S}$, $x'\in \Pi_{\on{LBA}}^{0}(F_{V}\boxtimes F_{J},F_{U'}
\boxtimes F_{I})_{S'}$, then $x'\circ (\beta_{F_V,F_{V'}}\boxtimes 
\on{id}_{F_{J}})\circ x\in \Pi_{\on{LBA}}^{0}(F_{V}\boxtimes F_{U},
F_{V'}\boxtimes F_{U'}|F_{I})$, $x\circ (\beta_{F_{U},F_{U'}}\boxtimes
\on{id}_{F_{I}})\circ x'\in\Pi_{\on{LBA}}(F_{U}\boxtimes F_{V},F_{U'}
\boxtimes F_{V'}|F_{J})$, and 
$$
\on{tr}_{F_I}(x'\circ (\beta_{F_V,F_{V'}}\otimes \on{id}_{F_J})\circ x)=
\beta_{F_{V'},F_{U'}}\circ \on{tr}_{F_J}(x\circ (\beta_{F_U,F_{U'}}
\otimes \on{id}_{F_I})\circ x')
\circ \beta_{F_U,F_V}. 
$$
\hfill \qed \medskip 
\end{remark}

\subsection{Morphisms of multiprops with partial traces}
 
The prop morphisms $\kappa_{1,2}:\on{LBA}\to \on{LBA}_{f}$, $\kappa_{0}:
\on{LBA}_{f}\to \on{LBA}$ and $\tau:\on{LBA}\to \on{LBA}$ induce 
morphisms between the corresponding multiprops $\Pi^{0}_{\on{LBA}_{f}}$
and $\Pi^{0}_{\on{LBA}}$ (still denoted $\kappa_{i}$, etc.). 

We will prove: 

\begin{proposition} \label{prop:intertw}
These morphisms intertwine the traces. 
\end{proposition}

{\em Proof.} Let $\kappa : \on{LBA}_{\alpha}\to \on{LBA}_{\beta}$ be any 
for these morphisms. We will prove that for any simple $F,G,H\in
\on{Ob}(\on{Sch}_{I,J,K})$, 
$\kappa(\Pi^{0}_{\on{LBA}_{\alpha}}(F,G|H))
\subset \Pi^{0}_{\on{LBA}_{\beta}}(F,G|H)$, and then that 
the diagram 
$$\begin{matrix} \Pi^{0}_{\on{LBA}_{\alpha}}(F,G|H)
& \stackrel{\kappa}{\to}& \Pi^{0}_{\on{LBA}_{\beta}}(F,G|H) \\
\scriptstyle{\on{tr}_{H}}
\downarrow &  & \downarrow\scriptstyle{\on{tr}_{H}}\\
\Pi^{0}_{\on{LBA}_{\alpha}}(F,G)
& \stackrel{\kappa}{\to}& \Pi^{0}_{\on{LBA}_{\beta}}(F,G) 
\end{matrix}
$$
commutes. 

The case of $\tau$ is clear. In the case of $\kappa_{0}$, we argue as follows: 
let $C: S\circ\wedge^{2}$, $D:= \wedge^{3}\otimes (S\circ\wedge^{2})$, 
then $C$ is a coalgebra in ${\bf Sch}$ and $\alpha:C\to D$ is a $\on{LBA}$-coideal
in $C$. The coalgebra morphism ${\bf 1}\to C$ induces a multiprop morphism
$\Pi_{\on{LBA}_{C}}^{0}\to \Pi_{\on{LBA}_{\bf 1}}^{0} 
= \Pi_{\on{LBA}}^{0}$; the composed morphism $\Pi_{\on{LBA}_{D}}^{0}
\to \Pi_{\on{LBA}_{C}}^{0}\to \Pi_{\on{LBA}}^{0}$ so we get a 
multiprop morphism $\Pi_{\on{LBA}_{f}}^{0}\to \Pi_{\on{LBA}}^{0}$, 
compatible with the traces. The maps $\Pi_{\on{LBA}_{f}}^{0}(F,G)\to 
\Pi_{\on{LBA}}^{0}(F,G)$ are the maps induced by $\kappa_{0}$, which proves 
the statement in the case of $\kappa_{0}$. 

The coalgebra morphism $C\to{\bf 1}$ induces a morphism of multiprops 
$\Pi_{\on{LBA}}^{0}\simeq \Pi_{\on{LBA}_{\bf 1}}^{0}\to 
\Pi_{\on{LBA}_{C}}$, compatible with the traces, 
which we compose with the projection $\Pi_{\on{LBA}_{C}}\to 
\Pi_{\on{LBA}_{f}}$. The maps $\Pi_{\on{LBA}}^{0}(F,G)\to 
\Pi_{\on{LBA}_{f}}^{0}(F,G)$ are the maps induced by $\kappa_{1}$, 
which proves the statement in the case of $\kappa_{1}$. 

We now treat the case of $\kappa_{2}$. 
For $F=\boxtimes_{i\in I}Z_{\rho_i}$, $G = \boxtimes_{j\in J}
Z_{\sigma_{j}}$, where $(\rho_{i})_{i}$, $(\sigma_{j})_{j}$ are maps 
$I,J\to \sqcup_{n\geq 0}\wh\SG_{n}$. 
For $Z\in\on{Irr(Sch)}^{I\times J}$, set 
$$
\ul\Pi_{\on{LBA}}^{0}(F,G)_{Z} := 
(\otimes_{i\in I} \on{LBA}(Z_{\rho_{i}},\otimes_{j\in J} Z_{ij}))
\otimes (\otimes_{j\in J}\on{LBA}(\otimes_{i\in I}Z_{ij},Z_{\sigma_{j}}))
$$
and for $S\in\cG_{0}(I,J)$, 
$$
\ul\Pi_{\on{LBA}}^{0}(F,G)_{S}:= 
\oplus_{Z|\on{supp}(Z)=S} \ul\Pi_{\on{LBA}}^{0}(F,G)_{Z}. 
$$
The operations of $\on{LBA}$ (tensor products, composition, braidings)
give rise to a natural map $\ul\Pi_{\on{LBA}}^{0}(F,G)_{Z}\to 
\Pi_{\on{LBA}}^{0}(F,G)$, which add up to a map 
$\ul\Pi_{\on{LBA}}^{0}(F,G)_{S}\to 
\Pi_{\on{LBA}}^{0}(F,G)$. 

\begin{lemma}
The image of this map is equal to $\Pi_{\on{LBA}_{\alpha}}^{0}(F,G)_{S}$.  
\end{lemma}

{\em Proof.} This image contains $\Pi_{\on{LBA}}^{0}(F,G)_{S}$, as 
$\Pi^0_{\on{LBA}}(F,G)_{S}$ is the subspace of 
$\Pi^{0}_{\on{LBA}}(F,G)_{S}$, where the successive $\on{LBA}$ are 
replaced by $\on{LCA},\on{LA}$. Let us prove the opposite inclusion. 

For each $Z$, the map $\ul\Pi_{\on{LBA}}^{0}(F,G)_{Z}\to 
\Pi_{\on{LBA}}^{0}(F,G)$ factors as 
\begin{align*}
& (\otimes_{i=1}^{n}\on{LBA}(Z_{\rho_{i}},\otimes_{j}Z_{ij}))
\otimes (\otimes_{j=1}^{p} \on{LBA}(\otimes_{i}Z_{ij},Z_{\sigma_{j}}))
\\
 & \to \oplus_{Z',Z''\in \on{Irr(Sch)}^{[n]\times [p]}}
(\otimes_{i}\on{LCA}(Z_{\rho_{i}},\otimes_{j}Z'_{ij})) 
\otimes (\otimes_{i,j}\on{LA}(Z'_{ij},Z_{ij})) \otimes 
(\otimes_{i,j} \on{LCA}(Z_{ij},Z''_{ij})) \\
 & \otimes 
(\otimes_{j}\on{LA}(\otimes_{i}Z''_{ij},Z_{\sigma_{j}})) 
\\
& \to 
\oplus_{Z',Z'',Z'''\in \on{Irr(Sch)}^{[n]\times [p]}}
(\otimes_{i}\on{LCA}(Z_{\rho_{i}},\otimes_{j}Z'_{ij})) 
\otimes (\otimes_{i,j}\on{LCA}(Z'_{ij},Z'''_{ij})) \otimes 
(\otimes_{i,j} \on{LA}(Z'''_{ij},Z''_{ij})) \\
 & \otimes 
(\otimes_{j}\on{LA}(\otimes_{i}Z''_{ij},Z_{\sigma_{j}})) 
\\
& \to 
\oplus_{Z'''\in \on{Irr(Sch)}^{[n]\times [p]}}
(\otimes_{i}\on{LCA}(Z_{\rho_{i}},\otimes_{j}Z'''_{ij})) 
(\otimes_{j}\on{LA}(\otimes_{i}Z'''_{ij},Z_{\sigma_{j}})) 
 \to \Pi_{\on{LBA}}^{0}
 (\boxtimes_{i}Z_{\rho_{j}},\boxtimes_{j}Z_{\sigma_{j}}), 
\end{align*}
where the first map is a tensor product of  decompositions of 
$\on{LBA}(\otimes_{i}F_{i}, \otimes_{j}G_{j})$, the second map is a tensor product of 
exchange maps (composition followed by decomposition) 
$\on{LA}(Z'_{ij},Z_{ij})\otimes \on{LCA}(Z_{ij},Z''_{ij}) \to 
\on{LBA}(Z'_{ij},Z''_{ij})\to \oplus_{Z'''_{ij}}\on{LCA}(Z'_{ij},Z'''_{ij})
\otimes \on{LA}(Z''_{ij},Z''_{ij})$, the third map is a tensor product of compositions
in $\on{LA}$ and $\on{LCA}$. Now if $Z'_{ij}={\bf 1}$ (resp., $\neq {\bf 1}$), 
the components 
of the exchange map corresponding to any $Z'''_{ij}\neq {\bf 1}$ (resp., 
$Z'''_{ij}={\bf 1}$), are zero. Therefore the components of the composition
of the three first maps, where 
$\on{supp}(Z''')\neq \on{supp}(Z)$, are zero. It follows that if $\on{supp}(Z)=S$, 
the image of the 
overall map is contained in $\Pi_{\on{LBA}}^{0}(\boxtimes_{i}Z_{\rho_{i}},
\boxtimes_{j}Z_{\sigma_{j}})_{S}$, as wanted. 
\hfill \qed \medskip 

We also define $\ul\Pi_{\on{LBA}_{C}}^{0}(F,G)_{Z}:= 
\oplus_{Z'|Z'_{|I\times J}} \Pi^{0}_{\on{LBA}}(C\boxtimes F,G)$, and 
$\ul\Pi_{\on{LBA}_{f}}^{0}(F,G)_{Z}$ as the image of 
$\ul\Pi_{\on{LBA}_{C}}^{0}(F,G)_{Z}$ in $\ul\Pi_{\on{LBA}_{f}}^{0}(F,G)$. 
We define similarly $\ul\Pi_{\on{LBA}_{f}}^{0}(F,G)_{S}$ for $S\in\cG_{0}(I,J)$. 
Arguing as in the above Lemma, one show that the image of the natural map 
$\ul\Pi^{0}_{\on{LBA}_{f}}(F,G)_{S}\to\Pi^{0}_{\on{LBA}_{f}}(F,G)$
is contained in $\Pi^{0}_{\on{LBA}_{f}}(F,G)_{S}$. 

\begin{lemma}
The map $\kappa_{2}:\Pi_{\on{LBA}}^{0}(F,G)_{S}\to 
\Pi_{\on{LBA}_f}^{0}(F,G)$ factors through 
$\Pi_{\on{LBA}}^{0}(F,G)_{S}\to 
\ul\Pi_{\on{LBA}_f}^{0}(F,G)_{S} \to 
\Pi_{\on{LBA}_f}^{0}(F,G)_{S}$, so it is compatible with
the gradings by $\cG_{0}(I,J)$. 
\end{lemma}

{\em Proof of Lemma.} For each $Z\in \on{Irr(Sch)}^{I\times J}$, 
the restriction of this map factors as 
\begin{align*} & \Pi^{0}_{\on{LBA}}(F,G)_{Z} = 
(\otimes_{i\in I}\on{LCA}(Z_{\rho_{i}},\otimes_{j\in J}Z_{ij}))
\otimes (\otimes_{j\in J}\on{LA}(\otimes_{i\in I}Z_{ij},Z_{\sigma_{j}}))
\\
 & \stackrel{\kappa_{2}}{\to}
(\otimes_{i\in I}\on{LBA}_{f}(Z_{\rho_{i}},\otimes_{j\in J}Z_{ij}))
\otimes (\otimes_{j\in J}\on{LBA}_{f}(\otimes_{i\in I}Z_{ij},Z_{\sigma_{j}}))
\to\Pi_{\on{LBA}_{f}}^{0}(F,G). \end{align*}
The statement follows from the fact that the image of the last map is contained in 
$\Pi^{0}_{\on{LBA}_{f}}(F,G)_{Z}$. 
\hfill \qed\medskip 

If $H = \boxtimes_{k\in K} Z_{\tau_{k}}$, and $S\in\cG_{0}(I,J|K)$, 
then we define 
$$\ul{\on{tr}}_{H} : \ul\Pi^{0}_{\on{LBA}_{f}}(F,G|H)_{S} \to 
\Pi^{0}_{\on{LBA}_{f}}(F,G)_{\on{tr}_{K}(S)}
$$
similarly to $\on{tr}_{H}$ (see (\ref{def:tr}), where $\on{LCA},\on{LA}$
are replaced by $\on{LBA}$ and $\Pi^{0}_{\on{LCA}}$, $\Pi^{0}_{\on{LA}}$
are replaced by $\Pi^{0}_{\on{LBA}}$). 

\begin{lemma}
The diagram 
$$
\begin{matrix}
\ul\Pi^{0}_{\on{LBA}_{f}}(F,G|H)_{S} & 
& 
\\
\downarrow  & \searrow{\scriptstyle{\ul{\on{tr}}}_{H}} & 
\\
\Pi^{0}_{\on{LBA}_{f}}(F,G|H)_{S} & \stackrel{{\on{tr}}_{H}}{\to} 
& \Pi^{0}_{\on{LBA}_{f}}(F,G)_{\on{tr}_{K}(S)}
\end{matrix}
$$
commutes. 
\end{lemma}

{\em Proof.} We first prove that the similar diagram commutes,
where $\on{LBA}_{f}$ is replaced by $\on{LBA}$. For $Z=(Z_{uv})_{(u,v)
\in (I\sqcup K)\times (K\sqcup J)}$ such that 
$\on{supp}(Z)=S$, the vertical map restricts to 
\begin{align*}
& \ul\Pi^{0}_{\on{LBA}}(F,G|H)_{Z} = 
(\otimes_{i\in I\sqcup K}\on{LBA}(Z_{\rho_{i}},
\otimes_{j\in J\sqcup K}Z_{ij})) \otimes 
(\otimes_{j\in K\sqcup J}\on{LBA}(\otimes_{i\in I\sqcup K}Z_{ij},
Z_{\sigma_{j}}))
\\
 & \simeq 
\oplus_{Z_{iij},Z_{ijj}\in \on{Irr(Sch)}}
[\otimes_{i\in I\sqcup K} \{\on{LCA}(Z_{\rho_{i}},
\otimes_{j\in K\sqcup J}Z_{iij}) \otimes
\otimes_{j\in J\sqcup K}\on{LA}(Z_{iij},Z_{ij})\}] \\
 & \otimes 
[\otimes_{j\in K\sqcup J}\{\otimes_{i\in I\sqcup K}
\on{LCA}(Z_{ij},Z_{ijj})\otimes\on{LA}
(\otimes_{i\in I\sqcup K}Z_{ijj},Z_{\sigma_{j}})\}]
\\
 & \to \oplus_{Z_{iij},Z_{ijj},Z'_{ij}\in \on{Irr(Sch)}}
[\otimes_{i\in I\sqcup K}\{\on{LCA}(Z_{\rho_{i}},
\otimes_{j\in K\sqcup J}Z_{iij}) \otimes
\otimes_{j\in J\sqcup K}\on{LCA}(Z_{iij},Z'_{ij})\}] \\
 & \otimes 
[\otimes_{j\in K\sqcup J}\{\otimes_{i\in I\sqcup K}
\on{LA}(Z'_{ij},Z_{ijj})\otimes\on{LA}
(\otimes_{i\in I\sqcup K}Z_{ijj},Z_{\sigma_{j}})\}]
\\
& \to \oplus_{Z'_{ij}\in\on{Irr(Sch)}}
[\otimes_{i\in I\sqcup K} \on{LCA}(Z_{\rho_{i}},
\otimes_{j\in K\sqcup J}Z'_{ij})] 
\otimes [\otimes_{j\in K\sqcup J}
\on{LA}(\otimes_{i\in I\sqcup K}Z'_{ij},Z_{\sigma_{j}})]
\\
& \simeq \Pi^{0}_{\on{LBA}}(F,G|H)_{S}
\end{align*}
Define $(H_{k})_{k}$ and $(H_{k,k+1})_{k}$ as above; define also 
$(H'_{k})_{k}$ and $(H'_{k,k+1})_{k}$ similarly, replacing $(Z_{ij})$
by $(Z'_{ij})$, and $(H''_{k,k+1})_{k}, (H'''_{k,k+1})_{k}$ by 
replacing $(Z_{ij})$ by $(Z_{iij})$, $(Z_{ijj})$.

Then each square of the following diagram commutes 
$$\begin{matrix}
\begin{matrix} 
{\scriptstyle (\otimes_{i\in I\sqcup K}\on{LBA}(Z_{\rho_{i}},
\otimes_{j\in J\sqcup K}Z_{ij}))} \\
{\scriptstyle \otimes 
(\otimes_{j\in K\sqcup J}\on{LBA}(\otimes_{i\in I\sqcup K}Z_{ij},
Z_{\sigma_{j}}))} \end{matrix}
 & \to & {\scriptstyle \otimes_{k=0}^{|K|}\Pi^{0}_{\on{LBA}}(H_{k},H_{k+1})
 \otimes \Pi^{0}_{\on{LBA}}(H_{k,k+1},H_{k+1})}\\
\downarrow  & & \downarrow \\
\begin{matrix} {\scriptstyle \oplus_{Z_{iij},Z_{ijj}}
[\otimes_{i\in I\sqcup K} \{\on{LCA}(Z_{\rho_{i}},
\otimes_{j\in K\sqcup J}Z_{iij}) }\\
{\scriptstyle \otimes
\otimes_{j\in J\sqcup K}\on{LA}(Z_{iij},Z_{ij})\}] }
\\
{\scriptstyle
\otimes[\otimes_{j\in K\sqcup J}\{\otimes_{i\in I\sqcup K}
\on{LCA}(Z_{ij},Z_{ijj})} \\
{\scriptstyle 
\otimes\on{LA}
(\otimes_{i\in I\sqcup K}Z_{ijj},Z_{\sigma_{j}})\}]}
\end{matrix}
   & \to & \begin{matrix} {\scriptstyle\oplus_{Z_{iij},Z_{ijj}}
   \otimes_{k=0}^{|K|} \Pi^{0}_{\on{LCA}}(H_{k},H''_{k,k+1}) \otimes 
   \Pi^{0}_{\on{LA}}(H''_{k,k+1},H_{k,k+1}) }\\
   \scriptstyle{ \otimes   \Pi^{0}_{\on{LCA}}(H_{k,k+1},H'''_{k,k+1})\otimes
   \Pi^{0}_{\on{LA}}(H'''_{k,k+1},H_{k+1})} \end{matrix}\\
\downarrow   & & \downarrow \\
\begin{matrix} {\scriptstyle \oplus_{Z_{iij},Z_{ijj},Z'_{ij}\in \on{Irr(Sch)}}}
\\ {\scriptstyle
[\otimes_{i\in I\sqcup K}\{\on{LCA}(Z_{\rho_{i}},
\otimes_{j\in K\sqcup J}Z_{iij}) \otimes}\\
{\scriptstyle
\otimes_{j\in J\sqcup K}\on{LCA}(Z_{iij},Z'_{ij})\}] }\\
{\scriptstyle\otimes 
[\otimes_{j\in K\sqcup J}\{\otimes_{i\in I\sqcup K}
\on{LA}(Z'_{ij},Z_{ijj})}\\
{\scriptstyle \otimes\on{LA}
(\otimes_{i\in I\sqcup K}Z_{ijj},Z_{\sigma_{j}})\}]}
\end{matrix}
  & \to & \begin{matrix} {\scriptstyle \oplus_{Z_{iij},Z_{ijj},Z'_{ij}}
   \otimes_{k=0}^{|K|} \Pi^{0}_{\on{LCA}}(H_{k},H''_{k,k+1}) \otimes 
   \Pi^{0}_{\on{LCA}}(H''_{k,k+1},H'_{k,k+1})} \\
   \scriptstyle{\otimes 
   \Pi^{0}_{\on{LA}}(H'_{k,k+1},H'''_{k,k+1})\otimes
   \Pi^{0}_{\on{LA}}(H'''_{k,k+1},H_{k+1})} \end{matrix} \\
\downarrow & & \downarrow \\
\begin{matrix}{\scriptstyle \oplus_{Z'_{ij}\in\on{Irr(Sch)}}
[\otimes_{i\in I\sqcup K} \on{LCA}(Z_{\rho_{i}},
\otimes_{j\in K\sqcup J}Z'_{ij})]}\\
{\scriptstyle \otimes [\otimes_{j\in K\sqcup J}
\on{LA}(\otimes_{i\in I\sqcup K}Z'_{ij},Z_{\sigma_{j}})]}
\end{matrix}
 & \to & \scriptstyle{\oplus_{Z'_{ij}}
   \otimes_{k=0}^{|K|} \Pi^{0}_{\on{LCA}}(H_{k},H'_{k,k+1}) \otimes 
      \Pi^{0}_{\on{LA}}(H'_{k,k+1},H_{k+1})} \\
\downarrow & & \downarrow \\
{\scriptstyle \Pi^{0}_{\on{LBA}}(F,G|H)_{S}}
& \to & \scriptstyle{\Pi^{0}_{\on{LBA}}(F,G)}
\end{matrix}$$ 
which implies the commutativity in the case of $\Pi_{\on{LBA}}^{0}$. 
The proof is the same in the case of $\Pi_{\on{LBA}_{f}}^{0}$. 
\hfill \qed \medskip 

{\em End of proof of Proposition \ref{prop:intertw}.} The proposition now follows
from the commutativity of the above diagram, together with that of 
$$\begin{matrix}
\Pi^{0}_{\on{LBA}}(F,G|H)_{S} & \stackrel{\on{tr}_{H}}{\to} & 
\Pi^{0}_{\on{LBA}}(F,G)_{\on{tr}_{K}(S)} 
\\
 \downarrow  &  & \downarrow \\
\ul\Pi_{\on{LBA}_{f}}^{0}(F,G|H)_{S}
  & \stackrel{\ul{\on{tr}}_{H}}{\to} & 
\ul\Pi^{0}_{\on{LBA}_{f}}(F,G)_{\on{tr}_{K}(S)} 
 \end{matrix}$$
 \hfill \qed \medskip

\subsection{The quasi-bi-multiprops $\Pi$ and $\Pi_{f}$}

Let $\Pi,\Pi_{f}$ be the quasi-bi-multiprops associated to the 
multiprops with traces $\Pi_{\on{LBA}}^{0}$, 
$\Pi_{\on{LBA}_{f}}^{0}$, and the involution $F\mapsto F^{*}$
of $\on{Sch}_{(1)}$. Explicitly, we have 
$\Pi(F\ul\boxtimes G,F'\ul\boxtimes G'):= 
\on{LBA}(c(F)\otimes c(G')^{*},c(F')\otimes c(G)^{*})$
and $\Pi_{f}(F\ul\boxtimes G,F'\ul\boxtimes G'):= 
\on{LBA}_{f}(c(F)\otimes c(G')^{*},c(F')\otimes c(G)^{*})$. 

Since $\kappa_{1,2}:
\Pi_{\on{LBA}}^{0}\to \Pi_{\on{LBA}_{f}}^{0}$, 
$\kappa_{0}:\Pi^{0}_{\on{LBA}}\to \Pi^{0}_{\on{LBA}}$
and $\tau:\Pi^{0}_{\on{LBA}}\to \Pi^{0}_{\on{LBA}}$ 
are morphisms of multiprops with traces, they induce morphisms
$\kappa_{1,2}^{\Pi}$, etc., between the corresponding 
quasi-bi-multiprops. 

We now define a degree on $\Pi$ as follows. For $F\in
\on{Ob}(\on{Sch}_{(1)})$ of the form $F = \boxtimes_{i\in[n]}Z_{\rho_{i}}$, 
we set $|F| = \sum_{i\in [n]}|Z_{\rho_{i}}|$. For $F,...,G'\in\on{Ob}(\on{Sch}_{(1)})$
and $x\in \Pi(F\ul\boxtimes G,F'\ul\boxtimes G') = \on{LBA}(c(F)\otimes 
c(G')^{*},c(F')\otimes c(G)^{*})$ homogeneous, we set 
$\on{deg}_{\Pi}(x):= \on{deg}_{\delta}(x)+|G'|-|G|$. 
If $x\in \Pi(F\ul\boxtimes G,F'\ul\boxtimes G')_{Z} = 
\Pi^{0}_{\on{LBA}}(F\boxtimes (G')^{*},F'\boxtimes G^{*})_{Z}$, 
then $\on{deg}_{\Pi}(x)  = \sum_{(s,t)\in ([n]\sqcup[m'])
\times ([n']\sqcup[m])}|Z_{st}| 
- |F| - |G|$, so $\on{deg}_{\Pi}(x)\geq 0$. 
One checks that $\on{deg}_{\Pi}$ is a degree on $\Pi$, i.e., it is 
additive under composition and tensor product. 

We now define a degree on $\Pi_{f}$. We first define a degree on 
$\on{LBA}_{f}$ as follows: $f$ and $\delta$ have degree $1$
and $\mu$ has degree $0$. If now $x\in\Pi_{f}(F\ul\boxtimes G,
F'\ul\boxtimes G')$, we set $\on{deg}_{\Pi_{f}}(x):= 
\on{deg}_{\on{LBA}_{f}}(x) + |G'|-|G|$. Then $\on{deg}_{\Pi_{f}}(x)\geq 0$, 
and $\deg_{\Pi_{f}}$ defines a degree on $\Pi_{f}$

We define completions of $\Pi$ and $\Pi_{f}$ as follows: 
for $B,B'\in \on{Ob}({\bf Sch}_{(1+1)})$, 
$\hbox{\boldmath$\Pi$\unboldmath}(B,B')$ 
(resp., $\hbox{\boldmath$\Pi$\unboldmath}_{f}(B,B')$)
is the degree completion of $\Pi(B,B')$ (resp., $\Pi_{f}(B,B')$). 
The morphisms $\kappa_{1,2}^{\Pi}$ of quasi-bi-multiprops 
are of degree $0$, and induce therefore morphisms between 
their completions. 

It follows from the cyclicity of the trace on $\cG_{0}$
that we have an involution of $\cG$, defined as follows: 
it acts on objects by $(I,J)\mapsto (J,I)$, and on morphisms by
$\cG((I,J),(I',J')) = \cG_{0}(I\sqcup J',I'\sqcup J)
\ni x \mapsto \beta_{I,J'} \circ x \circ \beta_{J,I'} \in 
\cG_{0}(J'\sqcup I,J\sqcup I')\in \cG((J',I'),(J,I))$. 

The cyclicity of the trace on $\Pi^{0}_{\on{LBA}}$ implies that 
the bi-multiprop $\Pi$ is equipped with a compatible involution, 
described as follows: its acts on objects as $F\ul\boxtimes G\mapsto
G^{*}\ul\boxtimes F^{*}$, and on morphisms by 
$\Pi(F\ul\boxtimes G,F'\ul\boxtimes G') = 
\on{LBA}(c(F)\otimes c(G')^{*},c(F')\otimes c(G)^{*})
\ni x \mapsto \beta_{c(F'),c(G^{*})} \circ x \circ 
\beta_{c(G^{\prime *}),c(F)} \in 
\on{LBA}(c(G^{*})\otimes c(F'),c(G^{\prime *})\otimes c(F)) 
= \Pi(G^{*}\ul\boxtimes F^{*},
(G')^{*}\ul\boxtimes (F')^{*})$. 
This involution has degree zero, hence extends to 
$\hbox{\boldmath$\Pi$\unboldmath}$.  

If $B\in \on{Ob}(\on{Sch}_{(1+1)})$, we define $\on{can}_{B}\in 
 \Pi({\mathfrak 1}\ul\boxtimes {\mathfrak 1},B\boxtimes B^{*})$
as follows. If $B = Z_{\rho_{1},...,\rho_{n}} \ul\boxtimes 
Z_{\sigma_{1},...,\sigma_{p}}$, where $Z_{\rho_{1},...,\rho_{n}} = 
(\boxtimes_{i=1}^{n} Z_{\rho_{i}})$, then 
$B\boxtimes B^{*} = (Z_{\rho_{1},...,\rho_{n}}\ul\boxtimes 
Z_{\sigma_{1},...,\sigma_{p}})\boxtimes (Z_{\sigma_{1}^{*},...,\sigma_{p}^{*}}
\ul\boxtimes Z_{\rho_{1}^{*},...,\rho_{n}^{*}}) = 
Z_{\rho_{1},...,\sigma_{p}^{*}}\ul\boxtimes Z_{\sigma_{1},...,\rho_{n}^{*}}$, 
so $\Pi({\mathfrak 1}\ul\boxtimes {\mathfrak 1},
Z_{\rho_{1},...,\sigma_{p}^{*}}\ul\boxtimes
Z_{\sigma_{1},...,\rho_{n}^{*}}) = \Pi_{\on{LBA}}^{0}
(Z^{*}_{\sigma_{1},...,\rho_{n}^{*}},Z_{\rho_{1},...,\sigma_{p}^{*}})
= \Pi_{\on{LBA}}^{0}(Z_{\sigma_{1}^{*},...,\sigma_{n}^{*}}\boxtimes 
Z_{\rho_{1},...,\rho_{n}},Z_{\rho_{1},...,\rho_{n}}\boxtimes 
Z_{\sigma_{1}^{*},...,\sigma_{n}^{*}})$, and $\on{can}_{B}$ corresponds to 
$\beta_{Z_{\rho_{1},...,\rho_{n}},Z_{\sigma_{1}^{*},...,\sigma_{n}^{*}}}$. 
If now $B = (B_{\rho_{1},...,\rho_{n};\sigma_{1},...,\sigma_{p}})$, 
we set $\on{can}_{B} = \oplus \on{id}_{B_{\rho_{1},...,\sigma_{p}}}
\otimes \on{can}_{Z_{\rho_{1},...,\rho_{n}}\ul\boxtimes 
Z_{\sigma_{1},...,\sigma_{p}}} \in \oplus\on{End}(B_{\rho_{1},...,\sigma_{p}})
\otimes \Pi({\mathfrak 1}\ul\boxtimes{\mathfrak 1},
(Z_{\rho_{1},...,\rho_{n}}\ul\boxtimes Z_{\sigma_{1},...,\sigma_{p}})
\boxtimes (Z_{\rho_{1},...,\rho_{n}}\ul\boxtimes Z_{\sigma_{1},...,\sigma_{p}}
)^{*})\subset \Pi({\mathfrak 1}\ul\boxtimes {\mathfrak 1},B\boxtimes B^{*})$. 

We define similarly $\on{can}_{B}^{*}\in \Pi(B\boxtimes B^{*},{\mathfrak 1}
\ul\boxtimes {\mathfrak 1})$ as follows. If $B = Z_{\rho_{1},...,\rho_{n}}
\ul\boxtimes Z_{\sigma_{1},...,\sigma_{p}}$, then 
$\Pi(B\boxtimes B^{*},{\mathfrak 1}\ul\boxtimes{\mathfrak 1}) = 
\Pi_{\on{LBA}}^{0}(Z_{\rho_{1},...,\sigma_{p}^{*}},Z^{*}_{\sigma_{1},...,
\rho_{n}^{*}})$ an $\on{can}_{B}^{*}$ corresponds to 
$\beta_{Z_{\rho_{1},...,\rho_{n}},Z_{\sigma_{1}^{*},...,\sigma_{p}^{*}}}$; 
we then extend this definition by linearity as above. 

The involution of $\Pi$ can then be described as follows: for 
$x\in \Pi(B,C)$, $x^{*}\in \Pi(C^{*},B^{*})$ can be expressed as 
$x^{*} = (\on{can}_{C}^{*}\boxtimes\on{id}_{B^{*}})\circ 
(x\boxtimes\beta_{B^{*},C^{*}})\circ (\on{can}_{B}\boxtimes 
\on{id}_{C^{*}})$.

The quasi-bi-multiprops $\Pi,\Pi_{f}$ give rise to quasi-biprops
$\pi,\pi_{f}$ by the inclusion $\on{Ob}(\on{Sch}_{1+1})\subset 
\on{Ob}(\on{Sch}_{(1+1)})$. Their topological versions 
$\hbox{\boldmath$\Pi,\Pi$\unboldmath}_{f}$ give rise to 
topological quasi-biprops  $\hbox{\boldmath$\pi,\pi$\unboldmath}_{f}$. 

We define sub-bimultiprops $\Pi^{\on{left},\on{right}}$ and 
$\Pi_{f}^{\on{left},\on{right}}$ as follows. For $F = \boxtimes_{i\in 
I}Z_{\rho_{i}}$, $G = \boxtimes_{j\in J} Z_{\sigma_{j}}$, etc., 
we set $\Pi^{\on{left},\on{right}}(F\ul\boxtimes G,F'\ul\boxtimes G'):= 
\oplus_{S\in \cG^{\on{left},\on{right}}\Pi((I,J),(I',J'))}(F\ul\boxtimes G,
F'\ul\boxtimes G')_{S}$, and a similar definition in the case of $\Pi_{f}$. 
These bimultiprops have topological versions
$\hbox{\boldmath$\Pi,\Pi$\unboldmath}_{f}^{\on{left},\on{right}}$. 

The sub-bimultiprops give rise to sub-biprops $\pi^{\on{left},\on{right}}$
of $\pi$ and $\pi_{f}^{\on{left},\on{right}}$ of $\pi_{f}$, as well as 
to topological sub-biprops 
$\hbox{\boldmath$\pi$\unboldmath}^{\on{left},\on{right}}$
and $\hbox{\boldmath$\pi$\unboldmath}_{f}^{\on{left},\on{right}}$ 
of $\hbox{\boldmath$\pi$\unboldmath}$
and $\hbox{\boldmath$\pi$\unboldmath}_{f}$. 


\subsection{Cokernels in $\on{LA}$}

Let $\on{i}_{\on{LA}}\in {\bf LA}(T\otimes T_2 \otimes T,T)$ 
be the prop morphism $m_T^{(2)} \circ \big( \on{id}_T \otimes  \big( 
(12) - (21) - \mu \big) \otimes \on{id}_T \big)$, where 
$m_T\in {\bf Sch}(T^{\otimes 2},T)$ is the propic version of the product in the 
tensor algebra and $m_T^{(2)} \in{\bf Sch}(T^{\otimes 3},T)$ is its
2-fold iterate. 

Let $p_{\on{LA}} \in {\bf LA}(T,S)$ be the direct sum, for $n\geq 0$, 
of $p_{\on{LA},n} \in {\bf LA}({\bf id}^{\otimes n},S)$ given by 
$m_{\on{PBW}}^{(n)} \circ \on{inj}_1^{\boxtimes n}$, 
where $m_{\on{PBW}}\in {\bf LA}(S^{\otimes 2},S)$ is the propic version of
the PBW star-product,\footnote{The PBW star-product is 
the map $S(\a)^{\otimes 2} \to S(\a)$ obtained from the product map 
$U(\a)^{\otimes 2} \to U(\a)$ by the symmetrization map, where 
$\a$ is a Lie algebra.} $m_{\on{PBW}}^{(n)} = m_{\on{PBW}}
 \circ ... \circ (m_{\on{PBW}}
\boxtimes \on{id}_{S^{\otimes n-2}})\in {\bf LA}(S^{\otimes n},S)$
and $\on{inj}_1 : {\bf id} \to S$ is the canonical morphism. 

Then the composed morphism $p_{\on{LA}} \circ \on{i}_{\on{LA}}$ is zero. 
Moreover, one proves using the image of $\on{sym} \in {\bf Sch}(S,T)
\to {\bf LA}(S,T)$, where sym is the symmetrization map, that 
$T\stackrel{p_{\on{LA}}}{\to} S$ is a cokernel for 
$T\otimes T_2\otimes T \stackrel{\on{i}_{\on{LA}}}{\to} T$. 

The following diagram then commutes in ${\bf LA}$: 
$$
\begin{matrix}
T^{\otimes 2} & \stackrel{m_T}{\to}& T \\
\downarrow & & \downarrow\\ 
S^{\otimes 2} & \stackrel{m_{\on{PBW}}}{\to}& S 
\end{matrix}
$$

Let $\Delta_0^S\in {\bf Sch}(S,S^{\otimes 2}) \subset 
{\bf LA}(S,S^{\otimes 2})$ be the propic version 
of the coproduct of the symmetric algebra $S(V)$, where $V$ is a vector space. 
Then 
\begin{equation} \label{LA:m:Delta}
\Delta_0^S \circ m_{\on{PBW}} = (m_{\on{PBW}} \otimes 
m_{\on{PBW}}) \circ (1324) \circ 
\Big( \Delta_0^S \otimes \Delta_0^S \Big);  
\end{equation}

In the same way, $T^{\otimes n} \stackrel{p_{\on{LA}}^{\otimes n}}
{\to} S^{\otimes n}$ is a cokernel for 
$\oplus_{i=0}^{n-1} T^{\otimes i} \otimes \Big( T \otimes T_2 \otimes T\Big) 
\otimes T^{\otimes n-1-i} \to T^{\otimes n} $. 

Observe that $\on{LA}$ is not an abelian category, as some morphisms
(e.g., $\mu\in\on{LA}(\wedge^2,{\bf id})$) do not admit cokernels.

\subsection{The element $\delta_S \in {\bf LBA}(S,S^{\otimes 2})$}

There is a unique element $\delta_T\in {\bf LBA}(T,T^{\otimes 2})$, 
such that $\delta_T \circ \on{inj}_1 = \on{can} \circ \delta$, where 
$\on{inj}_1 : {\bf id} \to T$ and $\on{can} : \wedge^2 \hookrightarrow  
{\bf id}^{\otimes 2} \hookrightarrow T^{\otimes 2}$ are the canonical 
injections, and such that 
$$
\delta_T \circ m_T = m_T^{\otimes 2} \circ (1324) \circ 
(\delta_T \otimes \Delta_0^T + \Delta_0^T \otimes \delta_T), 
$$
where $\Delta_0^T \in {\bf Sch}(T,T^{\otimes 2})$ is the propic 
version of the coproduct of $T(V)$, where $V$ is primitive. There exists
$\delta_{T\otimes T_2 \otimes T} \in {\bf LBA}(T\otimes T_2 \otimes T,
T \otimes (T\otimes T_2 \otimes T) \oplus (T\otimes T_2 \otimes T) \otimes T)$,
such that the diagram 
$$
\begin{matrix}
T\otimes T_2 \otimes T & \stackrel{\delta_{T\otimes T_2 \otimes T}}{\to} & 
T\otimes (T\otimes T_2 \otimes T) \oplus (T\otimes T_2 \otimes T) 
\otimes T \\
\scriptstyle{\on{i}_{\on{LA}}}\downarrow  & & \downarrow
\scriptstyle{\on{id}_T \otimes \on{i}_{\on{LA}} + \on{i}_{\on{LA}} \otimes  
\on{id}_T}\\ T & \stackrel{\delta_T}{\to} & T\otimes T
\end{matrix}
$$
commutes. Taking cokernels, we get a morphism $\delta_S \in{\bf
LBA}(S,S^{\otimes 2})$, such that 
$$
\delta_S \circ m_{\on{PBW}} = m_{\on{PBW}}^{\otimes 2} \circ (1324) \circ 
\big( \delta_S \otimes \Delta_0^S + \Delta_0^S \otimes \delta_S\big),  
$$
and $\delta_S \circ \on{inj}_1 = \on{can} \circ \delta$, where $\on{can} : 
\wedge^2 \subset {\bf id}^{\otimes 2} \subset S^{\otimes 2}$ is the canonical 
inclusion. We also have $((12) + (21)) \circ \delta_S = 0$. 

Thus $\delta_S$ is the propic version of the image under the symmetrization map
of the co-Poisson map $\delta_{U(\a)} : U(\a) \to U(\a)^{\otimes 2}$, where 
$\a $ is a Lie bialgebra.    

\subsection{The morphisms 
$m_\Pi\in\hbox{\boldmath$\Pi$\unboldmath}
((S\underline\boxtimes S)^{\otimes 2},S\underline\boxtimes S)$
and $\Delta_0\in\hbox{\boldmath$\Pi$\unboldmath}
(S\underline\boxtimes S,(S\underline\boxtimes S)^{\otimes 2})$}

We introduce the propic version $m_\Pi$ of the product map $U(\g)^{\otimes 2}
\to U(\g)$, where $\g$ is the double of a Lie bialgebra $\a$, 
transported via the isomorphism  $U(\g) \simeq S(\a) \otimes S(\a^*)$
induced by the symmetrizations and the product map $U(\a) \otimes
U(\a^*) \to U(\g)$. 

We construct a prop morphism $\on{LA} \to ({\bf id} \boxtimes 
{\bf 1} \oplus {\bf 1} \boxtimes {\bf id}) 
(\pi^{\on{left}})$, taking $\mu\in \on{LA}(\wedge^2,{\bf id})$
to the sum of $\mu\in \on{LBA}(\wedge^2,{\bf id}) \simeq 
\pi(\wedge^2 \boxtimes {\bf 1},{\bf id} \boxtimes {\bf 1})$, 
$\delta\in \on{LBA}({\bf id},\wedge^2) \subset 
\on{LBA}({\bf id},T_2) \simeq \pi({\bf id} \boxtimes 
{\bf id},{\bf id} \boxtimes {\bf 1})$, 
$\mu\in \on{LBA}(\wedge^2,{\bf id}) \subset 
\on{LBA}(T_2,{\bf id}) \simeq \pi({\bf id} \boxtimes 
{\bf id},{\bf 1} \boxtimes {\bf id})$, 
and $\delta \in \on{LBA}({\bf id},\wedge^2) \simeq 
\pi({\bf 1} \underline\boxtimes \wedge^2,{\bf 1}
\boxtimes {\bf id})$. This image is a morphism  
$$
\hbox{\boldmath$\mu$\boldmath} : \wedge^2( {\bf id} \boxtimes 
{\bf 1} \oplus {\bf 1} \boxtimes {\bf id})  
\simeq (\wedge^2 \boxtimes {\bf 1}) 
\oplus ({\bf id} \boxtimes {\bf id}) 
\oplus ({\bf 1} \boxtimes \wedge^2) 
\to ({\bf id} \boxtimes {\bf 1}) \oplus 
({\bf 1}\boxtimes {\bf id}) 
$$
of $\pi^{\on{left}}$. It has $\Pi$-degree $0$. 

Let us denote by 
$$
m_\pi \in ({\bf id} \boxtimes {\bf 1} \oplus {\bf 1} \boxtimes {\bf id}) 
(\hbox{\boldmath$\pi$\unboldmath}^{\on{left}})(S^{\otimes 2},S)
\simeq \hbox{\boldmath$\pi$\unboldmath}^{\on{left}}(S^{\otimes 2}  
({\bf id} \boxtimes {\bf 1} \oplus {\bf 1} \boxtimes {\bf id}), 
S\circ ({\bf id} \boxtimes 
{\bf 1} \oplus {\bf 1} \boxtimes {\bf id}))
\simeq 
\hbox{\boldmath$\pi$\unboldmath}^{\on{left}}(
(S\boxtimes S)^{\otimes 2}, S\boxtimes S) 
$$
the image of $m_{\on{PBW}}\in {\bf LA}(S^{\otimes 2},S)$.  

We denote by $m_\Pi \in \hbox{\boldmath$\Pi$\unboldmath}^{\on{left}}
((S\underline\boxtimes S)^{\boxtimes 2},S\underline\boxtimes S)
\subset \hbox{\boldmath$\Pi$\unboldmath}
((S\underline\boxtimes S)^{\boxtimes 2},S\underline\boxtimes S)$
the image of $m_\pi$. Then $m_\Pi$ has $\Pi$-degree $0$.

Then $m_\pi$ is associative, therefore 
\begin{equation} \label{m:assoc}
m_\Pi \circ (m_\Pi \boxtimes \on{id}_{S\underline\boxtimes S}) 
= m_\Pi \circ (\on{id}_{S\underline\boxtimes S} \boxtimes m_\Pi).  
\end{equation}
We denote by $m_\Pi^{(2)}\in \hbox{\boldmath$\Pi$\unboldmath} 
((S\underline\boxtimes S)^{\otimes 3},
S\underline\boxtimes S) = \hbox{\boldmath$\Pi$\unboldmath}
(S^{\boxtimes 3}\underline\boxtimes S^{\boxtimes 3}, 
S\underline\boxtimes S)$ the common value of both sides, and more
generally by $m_{\Pi}^{(n)}$ the $n$fold iterate of $m_{\Pi}$. 

Let us define $m_{ba}\in 
\hbox{\boldmath$\pi$\unboldmath}^{\on{left}}
(({\bf 1}\boxtimes S)\otimes (S\boxtimes {\bf 1}),S\boxtimes S)$
as $m_\pi \circ ((\on{inj}_0 \boxtimes \on{id}_S) \otimes (\on{id}_S 
\boxtimes \on{inj}_0))$ (where $\on{inj}_0 : {\bf 1} \to S$ is 
the canonical morphism). 

Since $m_{\pi}\circ \on{can}_{S\boxtimes S} = \on{id}_{S\boxtimes S}$,  
where $\on{can}_{S\boxtimes S}
\in {\bf Sch}(S\boxtimes S,(S\boxtimes {\bf 1})\otimes 
({\bf 1}\boxtimes S))$ is the canonical map, and since 
we have commutative diagrams 
\begin{equation} \label{diag:subalg}
\begin{matrix}
(S\boxtimes {\bf 1})^{\otimes 2} & \to & 
(S\boxtimes S)^{\otimes 2}\\
\scriptstyle{m_{\on{PBW}} \boxtimes \on{id}_{{\bf 1}}} 
\downarrow & & \downarrow 
\scriptstyle{m_\pi}\\
S\boxtimes {\bf 1}&\to & S\boxtimes S 
\end{matrix}
\quad \on{and} \quad 
\begin{matrix}
({\bf 1}\boxtimes S)^{\otimes 2} & \to & 
(S\boxtimes S)^{\otimes 2}\\
\scriptstyle{\on{id}_{{\bf 1}}\boxtimes m_{\on{PBW}}^*} \downarrow & & 
\downarrow \scriptstyle{m_\pi}\\
{\bf 1}\boxtimes S &\to & S\boxtimes S 
\end{matrix}
\end{equation}
we have 
$$
m_{\Pi} = (\tilde m_{\on{PBW}}\boxtimes \tilde m^{*}_{\on{PBW}})
\circ (\on{id}_{S\ul\boxtimes {\mathfrak 1}}
\boxtimes \tilde m_{ba}\boxtimes \on{id}_{{\mathfrak 1}\ul\boxtimes S})
\circ \on{can}_{S\ul\boxtimes S}^{\boxtimes 2}, 
$$
where
$\on{can}_{S\ul\boxtimes S}\in \hbox{\boldmath$\Pi$\unboldmath}
(S\ul\boxtimes S,
(S\ul\boxtimes {\mathfrak 1})\boxtimes ({\mathfrak 1}
\ul\boxtimes S))$ is the canonical morphism, 
$\tilde m_{ba}\in \hbox{\boldmath$\Pi$\unboldmath}
(({\mathfrak 1}\ul\boxtimes S)
\boxtimes (S\ul\boxtimes{\mathfrak 1}),(S\ul\boxtimes {\mathfrak 1})
\boxtimes ({\mathfrak 1}\ul\boxtimes S))$
is the morphism derived from $m_{ba}$, 
$\tilde m_{\on{PBW}}\in \hbox{\boldmath$\Pi$\unboldmath}
((S\ul\boxtimes{\mathfrak 1})^{\boxtimes 2},
S\ul\boxtimes{\mathfrak 1})$ and 
$\tilde m_{\on{PBW}}^{*}\in \hbox{\boldmath$\Pi$\unboldmath}
(({\mathfrak 1}\ul\boxtimes S)^{\boxtimes 2},
{\mathfrak 1}\ul\boxtimes S)$ are the morphisms derived from $m_{\on{PBW}}$
and $m^{*}_{\on{PBW}}$. 
It follows that a graph\footnote{For $x\in \Pi(F\ul\boxtimes G,
F'\ul\boxtimes G')$, $F=\boxtimes_{i\in I}F_{i}$,..., 
$G' = \boxtimes_{j'\in J'} G'_{j'}$, and $S\subset (I\sqcup J')\times 
(I'\sqcup J)$, we say that $x$ admits the graph $S$ if $x\in \oplus_{S'\subset S}
\Pi(F\ul\boxtimes G,F'\ul\boxtimes G')_{S'}$.} 
for $m_\Pi$ is as follows: set $F_1 = ... = G_2 = S$, 
so that it belongs to $\hbox{\boldmath$\Pi$\unboldmath}
((F_1\boxtimes F_2)\underline\boxtimes (G_1\boxtimes G_2),
F'\underline\boxtimes G')$, then the edges of the graph are 
$F_i\to F'$, $G'\to G_i$, $F_{2}\to G_1$. 

It follows that if we view $m_\Pi^{(n-1)}$ as an element of 
$\hbox{\boldmath$\Pi$\unboldmath}((\boxtimes_{i=1}^n F_i) 
\underline\boxtimes (\boxtimes_{i=1}^n G_i),F'\underline\boxtimes G')$, 
where $F_1 = ... = G' = S$, a graph for $m_{\Pi}^{(n-1)}$ is 
$F_{i}\to F'$, $G'\to G_{i}$ ($i=1,...,n$), $F_{j}\to G_{i}$ ($1\leq i<j\leq n$). 

We observe for later use that the  morphism 
$$
T({\bf id} \boxtimes {\bf 1} \oplus 
{\bf 1} \boxtimes {\bf id}) \to S \boxtimes S
$$
in $\hbox{\boldmath$\pi$\unboldmath}^{\on{left}}$, 
given by the direct sum over $n\geq 0$ of all compositions 
$({\bf id}\boxtimes {\bf 1} \oplus 
{\bf 1} \boxtimes {\bf id})^{\otimes n} \to 
(S\boxtimes S)^{\otimes n}
\stackrel{m_\pi^{(n-1)}}{\to} S \boxtimes S$, is the cokernel of the 
morphism 
$$
T({\bf id} \boxtimes {\bf 1} \oplus 
{\bf 1} \boxtimes {\bf id}) \otimes 
T_2({\bf id} \boxtimes {\bf 1} \oplus 
{\bf 1} \boxtimes {\bf id})
\otimes 
T({\bf id} \boxtimes {\bf 1} \oplus 
{\bf 1} \boxtimes {\bf id})
\to T({\bf id} \boxtimes {\bf 1} \oplus 
{\bf 1} \boxtimes {\bf id})
$$
given by  
$m_T^{(2)} \circ \Big( 
\on{id}_T\otimes 
\big( (12) - (21) 
- \hbox{\boldmath$\bar\mu$\unboldmath} \big) \otimes 
\on{id}_T \Big)$, where 
$\hbox{\boldmath$\bar\mu$\unboldmath}$ is the composed morphism  
\begin{align*}
& \hbox{\boldmath$\bar\mu$\unboldmath} : 
T_2( {\bf id} \boxtimes {\bf 1} \oplus 
{\bf 1} \boxtimes {\bf id}) \to 
\wedge^2( {\bf id} \boxtimes {\bf 1} \oplus 
{\bf 1} \boxtimes {\bf id}) \simeq 
\wedge^2({\bf id}) \boxtimes {\bf 1} \oplus 
{\bf id} \boxtimes {\bf id}
\oplus {\bf 1} \boxtimes \wedge^2({\bf id})
\\ & 
\stackrel{\hbox{\boldmath$\mu$\unboldmath}}{\to} 
{\bf id} \boxtimes {\bf 1} \oplus 
{\bf 1} \boxtimes {\bf id} = 
T_1({\bf id} \boxtimes {\bf 1} \oplus 
{\bf 1} \boxtimes {\bf id}).  
\end{align*}

Let us define $\Delta_0\in \hbox{\boldmath$\Pi$\unboldmath}
(S\underline\boxtimes S,(S\underline\boxtimes S)^{\boxtimes 2})$. 
Recall that $\Delta_0^S\in {\bf Sch}(S,S^{\otimes 2}) \subset 
{\bf LBA}(S,S^{\otimes 2}) \subset \hbox{\boldmath$\Pi$\unboldmath}
(S\underline\boxtimes {\mathfrak 1},(S\boxtimes S)\underline\boxtimes 
{\mathfrak 1})$ and let $m_0^S\in {\bf Sch}(S^{\otimes 2},S) \subset 
{\bf LBA}(S^{\otimes 2},S) =  
\hbox{\boldmath$\Pi$\unboldmath}({\mathfrak 1}\underline\boxtimes S, 
{\mathfrak 1} \underline\boxtimes (S\boxtimes S))$ be the propic version 
of the product of the symmetric algebra $S(V)$. Set 
$$
\Delta_0 := \Delta_0^S \boxtimes m_0^S \in 
\hbox{\boldmath$\Pi$\unboldmath}(S\underline\boxtimes S, 
(S\boxtimes S) \underline\boxtimes (S\boxtimes S) )
\simeq 
\hbox{\boldmath$\Pi$\unboldmath}(S\underline\boxtimes S, 
(S\underline\boxtimes S)^{\boxtimes 2}). 
$$
A graph for this element is as follows. Set $F = G = ... = G'_2 = S$, 
then $\Delta_0 \in \hbox{\boldmath$\Pi$\unboldmath}(F\underline\boxtimes G, 
(F'_1\boxtimes F'_2) \underline\boxtimes (G'_1\boxtimes G'_2))$, and a graph 
is $F\to F'_1$, $F\to F'_2$, $G'_1 \to G$, $G'_2 \to G$. 
 
Then both sides of the following equality are defined, and the equality holds: 
\begin{equation} \label{Delta_0:m}
\Delta_0 \circ m_\Pi = (m_\Pi \boxtimes m_\Pi)\circ (1324) \circ  
(\Delta_0 \boxtimes \Delta_0).   
\end{equation}

This follows from (\ref{LA:m:Delta}). 

We also have commutative diagrams 
\begin{equation} \label{diag:subcoalg}
\begin{matrix}
S\underline\boxtimes {\mathfrak 1} & \to & S \underline\boxtimes S
\\
\scriptstyle{\Delta_0^S\underline\boxtimes 1}\downarrow  & & 
\downarrow\scriptstyle{\Delta_0} \\ 
(S\underline\boxtimes {\mathfrak 1})^{\otimes 2} & \to & 
(S \underline\boxtimes S)^{\otimes 2}
\end{matrix}
\quad \on{and} \quad 
\begin{matrix}
{\mathfrak 1}\underline\boxtimes S & \to & S \underline\boxtimes S
\\
\scriptstyle{1\underline\boxtimes (\Delta_0^S)^t}\downarrow  & & 
\downarrow\scriptstyle{\Delta_0} \\ 
({\mathfrak 1}\underline\boxtimes S)^{\otimes 2} & \to & 
(S \underline\boxtimes S)^{\otimes 2}
\end{matrix}
\end{equation}

\subsection{The element $\delta_{S\underline\boxtimes S}\in 
\hbox{\boldmath$\Pi$\unboldmath}(S\underline\boxtimes S,
(S\underline\boxtimes S)^{\otimes 2})$} \label{def:r}

We have $\hbox{\boldmath$\Pi$\unboldmath}(S\underline\boxtimes S,
(S\underline\boxtimes S)^{\otimes 2})
\simeq \hbox{\boldmath$\pi$\unboldmath}
(S\boxtimes S,S^{\otimes 2} \boxtimes S^{\otimes 2})$. 
Define 
\begin{equation} \label{def:delta:ss}
\delta_{S\boxtimes S} := \delta_S \boxtimes m_0^S + \Delta_0^S \boxtimes
\delta_S^t
\end{equation}
and $\delta_{S\underline\boxtimes S}$ as the image of this element in 
$\hbox{\boldmath$\Pi$\unboldmath} (S\underline\boxtimes S,
(S\underline\boxtimes S)^{\otimes 2})$. 
Then $((12) + (21)) \circ \delta_{S\underline\boxtimes S} = 0$. 

Let $r \in \Pi({\mathfrak 1}\underline\boxtimes {\mathfrak 1},
(S\underline\boxtimes S)^{\boxtimes 2})$ be the image of the 
composition $S\to {\bf id}\to S$ in ${\bf Sch}(S,S) \subset 
{\bf LBA}(S,S) \simeq  \hbox{\boldmath$\pi$\unboldmath}
({\bf 1}\boxtimes {\bf 1},(S\otimes{\bf 1})\boxtimes
({\bf 1}\otimes S))\subset \hbox{\boldmath$\pi$\unboldmath}
({\bf 1}\boxtimes {\bf 1},
S^{\otimes 2}\boxtimes S^{\otimes 2}) \simeq 
\hbox{\boldmath$\Pi$\unboldmath}({\mathfrak 1}
\ul\boxtimes{\mathfrak 1},S^{\boxtimes 2}\ul\boxtimes S^{\boxtimes 2})$. 
Then one checks that 
$$
\delta_{S\underline\boxtimes S} = m_\Pi^{\boxtimes 2} 
\circ (r\boxtimes \Delta_0- \Delta_0 \boxtimes r), 
$$
which is a propic version of the statement that the co-Poisson structure 
on $U(D(\a))$ is quasitriangular, with $r$-matrix $r_\a$. 
(\ref{def:delta:ss}) also implies that the diagrams 
$$
\begin{matrix}
S\underline\boxtimes {\mathfrak 1} & \to & S\underline\boxtimes S\\ 
\scriptstyle{\delta_S \underline\boxtimes 1}\downarrow  & & 
\downarrow\scriptstyle{\delta_{S\underline\boxtimes S}} \\ 
(S\underline\boxtimes {\mathfrak 1})^{\otimes 2} & \to & 
(S\underline\boxtimes S)^{\otimes 2}
\end{matrix}
\quad \on{and} \quad 
\begin{matrix}
{\mathfrak 1}\underline\boxtimes S & \to & S\underline\boxtimes S\\ 
\scriptstyle{1\underline\boxtimes \delta_S^t}\downarrow  & & 
\downarrow\scriptstyle{\delta_{S\underline\boxtimes S}} \\ 
({\mathfrak 1}\underline\boxtimes S)^{\otimes 2} & \to & 
(S\underline\boxtimes S)^{\otimes 2}
\end{matrix}
$$
commute. 

\subsection{The morphism $\Xi_f\in \hbox{\boldmath$\Pi$\unboldmath}_f
(S\underline\boxtimes S,S\underline\boxtimes S)^\times$}

If $\a$ is a Lie bialgebra and $f\in\wedge^2(\a)$ is a 
twist (we denote by $\a_f$ the Lie bialgebra $(\a,\delta + \on{ad}(f))$), 
then the doubles $D(\a)\simeq\a\oplus\a^*$ and $D(\a_f)\simeq\a\oplus\a^*$
are Lie algebra isomorphic, the isomorphism $D(\a) \to D(\a_f)$
being given by the automorphism of $\a\oplus \a^*$, $(a,0)\mapsto (a,0)$
and $(0,\alpha) \mapsto ((\on{id}_\a\otimes \alpha)(f),\alpha)$. 
The composed isomorphism 
$S(\a) \otimes S(\a^*) \to U(D(\a)) \simeq U(D(\a_f)) \simeq 
S(\a) \otimes S(\a^*)$ has a propic version $\Xi_f\in 
\hbox{\boldmath$\Pi$\unboldmath}_f
(S\underline\boxtimes S,S\underline\boxtimes S)^\times$, which we now  
construct. 

We define 
$$
\iota\in 
\pi_f^{\on{left}}({\bf id} \boxtimes {\bf 1}, 
{\bf id}\boxtimes {\bf 1}) 
\oplus
\pi_f^{\on{left}}({\bf 1} \boxtimes {\bf id}, 
{\bf id} \boxtimes {\bf 1}) 
\oplus 
\pi_f^{\on{left}}({\bf 1} \boxtimes {\bf id},
{\bf 1} \boxtimes {\bf id}) 
\subset 
\pi_f^{\on{left}}({\bf id} \boxtimes {\bf 1}
\oplus {\bf 1} \boxtimes {\bf id}, 
{\bf id} \boxtimes {\bf 1} \oplus 
{\bf 1} \boxtimes {\bf id}) 
$$
as the sum of $\pi_f(\on{id}_{{\bf id} \boxtimes {\bf 1}})$, 
of the element corresponding to $f$, and of $\pi_f(\on{id}_{
{\bf 1} \boxtimes {\bf id}})$. Then $\iota$ is homogeneous 
of degree $0$ (in the case of the middle element, the degree of $f$ is
compensated by the fact that the source and the target have different 
degrees). 

Then the diagram 
$$
\begin{matrix}
$$
(T\otimes T_2 \otimes T)({\bf id}\boxtimes{\bf 1}\oplus 
{\bf 1} \boxtimes {\bf id}) 
& \stackrel{T(\iota)\otimes T_2(\iota) \otimes
T(\iota)}{\to} & 
(T \otimes T_2 \otimes T)({\bf id}\boxtimes{\bf 1}\oplus 
{\bf 1} \boxtimes {\bf id}) 
\\
\begin{matrix}
\scriptstyle{ \kappa_1^\pi\big( m_T^{(2)} \circ 
(\on{id}_T \otimes}\\ 
\scriptstyle{ 
 ((12) - (21) - \bar\mu)
\otimes \on{id}_T)\big) } 
\end{matrix} 
\downarrow & & \downarrow 
\begin{matrix}
\scriptstyle{ \kappa_2^\pi\big( m_T^{(2)} \circ 
(\on{id}_T \otimes}\\ 
\scriptstyle{ 
 ((12) - (21) - \bar\mu)
\otimes \on{id}_T)\big) } 
\end{matrix} \\
T({\bf id}\boxtimes{\bf 1}\oplus 
{\bf 1} \boxtimes {\bf id})
 & \stackrel{T(\iota)}{\to} & 
T({\bf id}\boxtimes{\bf 1}\oplus 
{\bf 1} \boxtimes {\bf id})  
$$
\end{matrix}
$$
commutes in $\hbox{\boldmath$\pi$\unboldmath}_f^{\on{left}}$. Taking cokernels, 
we get a morphism $\xi_f\in \hbox{\boldmath$\pi$\unboldmath}_f^{\on{left}}(
S\boxtimes S,S\boxtimes S)$. We denote by $\Xi_f$ its image in 
$\hbox{\boldmath$\Pi$\unboldmath}_f
(S\underline\boxtimes S,S\underline\boxtimes S)$.

A lift of $\Xi_f$ to $\hbox{\boldmath$\Pi$\unboldmath} 
\Big( ((S\circ \wedge^2) \boxtimes S) \underline\boxtimes S,
S\underline\boxtimes S \Big)$ (which we write
$\hbox{\boldmath$\Pi$\unboldmath}\Big( 
((S\circ \wedge^2) \boxtimes F) \underline\boxtimes G,
F'\underline\boxtimes G'\Big)$) admits the graph 
$F\to F'$, $G'\to G$, $F'\to G'$, 
$S\circ \wedge^2 \to G$, $S\circ\wedge^2 \to F'$. 

Since $\iota$ is invertible, so is $\Xi_f$. 

\subsection{Relations between $\kappa_i^\Pi$, $\Xi_f$ and $m_\Pi$, $\Delta_0$}

Let us now study the relations of $\Xi_f$ with $m_\Pi$. 
In the case of a Lie bialgebra with twist $(\a,f)$, 
since $D(\a) \to D(\a_f)$ is a Lie algebra isomorphism, 
the diagram
$$
\begin{matrix}
(S(\a) \otimes S(\a^*))^{\otimes 2} & \to & S(\a) \otimes S(\a^*)\\ 
\downarrow & & \downarrow\\
(S(\a) \otimes S(\a^*))^{\otimes 2} & \to & S(\a) \otimes S(\a^*) 
\end{matrix}
$$
commutes, where the upper (resp., lower) arrow is induced by the product in 
$U(D(\a))$ (resp., $U(D(\a_f))$) and the vertical arrows are given by the 
above automorphism of $S(\a) \otimes S(\a^*)$. The propic version of this 
statement is that both terms of the following identity are defined, and 
the identity holds: 
\begin{equation} \label{kappa1:kappa2:1} \label{24}
\kappa_2^\Pi(m_\Pi) = \Xi_f \circ 
\kappa^\Pi_1(m_\Pi) \circ (\Xi_f^{-1})^{\boxtimes 2}. 
\end{equation}
The proof of this statement relies on the properties of a cokernel
and on the commutativity of the diagram 
$$
\begin{matrix}
T^{\otimes 2}({\bf id} \boxtimes {\bf 1}
\oplus {\bf 1} \boxtimes {\bf id}) 
& \stackrel{T^{\otimes 2}(\iota)}{\to} & 
T^{\otimes 2}({\bf id} \boxtimes {\bf 1}
\oplus {\bf 1} \boxtimes {\bf id}) \\
\scriptstyle{\kappa_1^\pi(m_T^{(2)})} \downarrow & & 
\downarrow \scriptstyle{\kappa_2^\pi(m_T^{(2)})}\\ 
T({\bf id} \boxtimes {\bf 1}
\oplus {\bf 1} \boxtimes {\bf id}) 
& \stackrel{T(\iota)}{\to} & 
T({\bf id} \boxtimes {\bf 1}
\oplus {\bf 1} \boxtimes {\bf id})  
\end{matrix}
$$

We now study the relation of $\Xi_f$ with 
$\Delta_0$. In the case of a Lie bialgebra with twist $(\a,f)$, 
the isomorphism $U(D(\a)) \to U(D(\a_f))$ is also compatible with the 
(cocommutative) bialgebra structures, as it is induced by a Lie 
algebra isomorphism. The propic version of this statement is that
both sides of the following identity are defined, as the identity holds  
\begin{equation} \label{kappa1:kappa2:2} \label{25}
\kappa^\Pi_2(\Delta_0) = \Xi_f^{\boxtimes 2} \circ\kappa_1^\Pi(\Delta_0)
\circ\Xi_f^{-1}. 
\end{equation}
This identity follows from the commutativity of 
$$
\begin{matrix}
T({\bf id}\boxtimes {\bf 1} \oplus 
{\bf 1} \boxtimes {\bf id}) & \stackrel{T(\iota)}{\to} & 
T({\bf id}\boxtimes {\bf 1} \oplus 
{\bf 1} \boxtimes {\bf id}) 
\\
\scriptstyle{
\pi(\Delta_T({\bf id} \boxtimes {\bf 1} \oplus {\bf 1}
\boxtimes {\bf id}))}
\downarrow & & \downarrow
\scriptstyle{
\pi(\Delta_T({\bf id} \boxtimes {\bf 1} \oplus {\bf 1}
\boxtimes {\bf id}))} \\
T^{\otimes 2}({\bf id}\boxtimes {\bf 1} \oplus 
{\bf 1}\boxtimes {\bf id}) 
& \stackrel{T^{\otimes 2}(\iota)}{\to} & 
T^{\otimes 2}({\bf id}\boxtimes {\bf 1} \oplus 
{\bf 1} \boxtimes {\bf id}) 
\end{matrix}
$$ 

\subsection{Relations between $\tau_\Pi$, $m_\Pi$ and $\Delta_0$}

Define $\omega_S\in{\bf Sch}(S,S)^\times$ as $\wh\oplus_{n\geq 0}
(-1)^n \on{id}_{S^n}$. 

\begin{lemma} \label{tau:pi:m:pi}
We have 
$$
\tau_\Pi(m_\Pi) = (\on{id}_S \underline\boxtimes \omega_S) \circ m_\Pi 
\circ ((\on{id}_S \underline\boxtimes \omega_S)^{\otimes 2})^{-1}, 
\quad 
\tau_\Pi(\Delta_0) = (\on{id}_S \underline\boxtimes \omega_S)^{\otimes 2} 
\circ \Delta_0 
\circ (\on{id}_S \underline\boxtimes \omega_S)^{-1}.  
$$
\end{lemma}

{\em Proof.} The first statement follows the the commutativity of the diagram 
$$
\begin{matrix}
\wedge^2({\bf id}\boxtimes {\bf 1} \oplus {\bf 1} \boxtimes {\bf id}) 
& \stackrel{\hbox{\boldmath$\mu$\unboldmath}}{\to} 
& {\bf id} \boxtimes {\bf 1} \oplus {\bf 1} \boxtimes 
{\bf id}\\
\scriptstyle{\wedge^2(\on{id}_{{\bf id}\boxtimes{\bf 1}} 
\oplus (-\on{id}_{{\bf 1}\boxtimes {\bf id}}))}
\downarrow &  & \downarrow 
\scriptstyle{\on{id}_{{\bf id}\boxtimes{\bf 1}} \oplus (-\on{id}_{{\bf 1}
\boxtimes {\bf id}})}
\\
\wedge^2({\bf id}\boxtimes {\bf 1} \oplus {\bf 1} \boxtimes {\bf id}) 
& \stackrel{\tau_{\pi}(\hbox{\boldmath$\mu$\unboldmath})}{\to} 
& {\bf id} \boxtimes {\bf 1} \oplus {\bf 1} \boxtimes {\bf id}
\end{matrix}
$$
The second statement follows from $m_0^S = \omega_S \circ m_0^S \circ 
(\omega_S^{\boxtimes 2})^{-1}$. \hfill \qed \medskip 

\begin{remark} Lemma \ref{tau:pi:m:pi} 
is the propic version of the following statement. 
Let $\a$ be a Lie bialgebra and let $\a' := \a^{\on{cop}}$ be $\a$ 
with opposite coproduct. 
Its double $\g'$ is Lie algebra isomorphic to $\g$, using the 
automorphism $\on{id}_\a \oplus (-\on{id}_{\a^*})$ of $\a\oplus \a^*$. 
The bialgebras $U(\g)$ and $U(\g')$ are therefore isomorphic, 
the isomorphism being given by $U(\g) \simeq S(\a\oplus \a^*)
\stackrel{S(\on{id}_\a\oplus (-\on{id}_{\a^*}))}{\to} 
S(\a\oplus \a^*) \simeq U(\g')$. \end{remark}

\section{The $\cX$-algebras ${\bf U}_n$ and ${\bf U}_{n,f}$} \label{sec:2}

\subsection{The category $\cX$}

Let $\cX$ be the category where objects are finite sets and 
morphisms are partially defined functions. A $\cX$-vector space 
(resp., algebra) is a contravariant functor $\cX \to \on{Vect}$ 
(resp., $\cX \to \on{Alg}$, where $\on{Alg}$ is the category 
of algebras). A $\cX$-vector space (resp., algebra) is the same 
as a collection $(V_s)_{s\geq 0}$ of vector spaces (resp., algebras), 
together with a collection of morphisms (called insertion-coproduct morphisms)
$V_s \to V_t$, $x\mapsto x^\phi$, 
for each function $\phi : [t] \to [s]$, satisfying the chain rule. 
Instead of $x^\phi$, we
often write $x^{\phi^{-1}(1),...,\phi^{-1}(s)}$. An example of a  
$\cX$-vector space (resp., $\cX$-algebra) is 
$V_s = H^{\otimes s}$, where $H$ is a cocommutative coalgebra 
(resp., bialgebra). Then $x^\phi$ is obtained from $x$ by applying 
the $(|\phi^{-1}(\alpha)|-1)$th iterated coproduct to the
component $\alpha$ of $x$ and plugging the result 
in the factors $\phi^{-1}(\alpha)$, for $\alpha = 1,...,s$.

\subsection{The $\cX$-algebra ${\bf U}_n$}

Let us set ${\bf U}_n := \Pi({\mathfrak 1} \underline\boxtimes
{\mathfrak 1},(S\underline\boxtimes S)^{\boxtimes n})$. Let us equip 
it with the $\Pi$-degree: the $\Pi$-degree of a homogeneous element of 
$\Pi(F\underline\boxtimes G,F'\underline\boxtimes G')$ is $\geq -|F|-|G|$, 
therefore ${\bf U}_n$ is $\NN$-graded. 

For $x,y\in {\bf U}_n$, the composition 
$m_\Pi^{\boxtimes n} \circ (1,n+1,2,n+2,...) \circ (x\boxtimes y)$ is
well-defined. We set 
$$
xy := m_\Pi^{\boxtimes n} \circ (1,n+1,2,n+2,...) \circ (x\boxtimes y)
\in {\bf U}_n. 
$$
It follows from (\ref{m:assoc}) and from $\on{deg}_\Pi(m_\Pi) = 0$
that the map $x\otimes y\mapsto xy$ 
defines an associative product of degree $0$ 
on ${\bf U}_n$.  

Let $\on{Coalg}_{\on{coco}}$ be the prop of cocommutative bialgebras, 
let $\phi : [m] \to [n]$ and let $\Delta^\phi \in \on{Coalg}_{\on{coco}}
(T_n,T_m)$ be the corresponding element. We have a prop morphism 
$\on{Coalg}_{\on{coco}}\to (S\underline\boxtimes S)(\Pi^{\on{right}})$, 
induced by $(\on{coproduct}) \mapsto \Delta_0$. We denote by $\Delta_0^\phi$
the image of $\Delta^\phi$ by this morphism, and we set 
$$
x^\phi := \Delta_0^\phi \circ x. 
$$ 
It then follows from (\ref{Delta_0:m}) that the family of all 
$({\bf U}_n)_{n\geq 0}$, equipped with these insertion-coproduct morphisms, 
is a $\cC$-algebra. 

\begin{remark}
The element $r$ defined in Subsection \ref{def:r}  belongs to 
${\bf U}_2$; one checks 
that it satisfies the classical Yang-Baxter equation 
$[r^{1,2},r^{1,3}] + [r^{1,2},r^{2,3}] + [r^{1,3},r^{2,3}] = 0$
in ${\bf U}_3$. \hfill \qed \medskip 
\end{remark}

\begin{lemma}
The map ${\bf U}_n \to {\bf U}_n$, $x\mapsto (\on{id}_S \underline\boxtimes 
\omega_S)^{\boxtimes n} \circ \tau_\Pi(x)$ is a automorphism of $\NN$-graded 
$\cX$-algebra. 
\end{lemma}

{\em Proof.} This follows from Lemma \ref{tau:pi:m:pi} and from the fact
that $\tau_\Pi$ has degree $0$. \hfill \qed \medskip 

\begin{remark} One checks that 
$(\on{id}_S \underline\boxtimes \omega_S)^{\boxtimes 2} \circ \tau_\Pi(r) = -r$, 
so the above automorphism 
will be denoted $x\mapsto x(-r)$. It is involutive, i.e., $(x(-r))(-r) = x$.  
\end{remark}

\subsection{The $\cX$-algebra ${\bf U}_{n,f}$}

We set ${\bf U}_{n,f} := \Pi_f({\mathfrak 1}\underline\boxtimes {\mathfrak 1},
(S\underline\boxtimes S)^{\boxtimes n})$. The $\Pi_f$-degree induces a 
grading on ${\bf U}_{n,f}$. 

\begin{lemma}
${\bf U}_{n,f}$ is a $\NN$-graded vector space. 
\end{lemma}

{\em Proof.} Let $x$ be homogeneous in the image of 
$\on{LBA}((S^k\circ \wedge^2) \otimes F\otimes G', 
F'\otimes G) \to \on{LBA}_f(F\otimes G',F'\otimes G)$. 
Let us show that $\on{deg}_{\on{LBA}_f}(x) 
\geq k - |F| - |G'|$. Indeed, if  $x$ belongs to 
$$
\on{LCA}(S^k\circ \wedge^2,W_1\otimes W_1) \otimes \on{LCA}(F,Z_1\otimes Z_2)
\otimes \on{LCA}(G',Z_3\otimes Z_4) \otimes \on{LA}(W_1\otimes 
Z_1\otimes Z_3,F) \otimes \on{LA}(W_2 \otimes Z_2\otimes Z_4,G),  
$$
its degree is $\on{deg}_{\on{LBA}_f}(x) = 
k + \on{deg}_\delta(x) = k + |W_1| + |W_2| + |Z_1| + ... + |Z_4|
- |F| - |G'| - 2k = (|W_1| + |W_2| - 2k) + |Z_1| + ... + |Z_4| - |F| - |G'|+k$. 
Now $|W_1| + |W_2| \geq 2k$, so $\on{deg}_{\on{LBA}_f}(x) \geq 
k - |F| - |G'|$. Therefore $\on{deg}_{\Pi_f}(x) \geq k - |F| - |G|$. 
In our case, $F = G = {\mathfrak 1}$, so $\on{deg}_{\on{LBA}_f}(x) 
\geq 0$.  
\hfill \qed \medskip 

For $x,y\in {\bf U}_{n,f}$, we define $xy$ and $x^\phi$ 
as above, replacing $m_\Pi$ and $\Delta_0^\phi$ by 
$\kappa_1^\Pi(m_\Pi)$ and $\kappa_1^\Pi(\Delta_0^\phi)$. 
This makes ${\bf U}_{n,f}$ into a $\NN$-graded $\cX$-algebra. 

\begin{lemma} 
The maps ${\bf U}_n \to {\bf U}_{n,f}$, $x\mapsto 
\kappa_1^\Pi(x)$ and $x\mapsto (\Xi_f^{-1})^{\boxtimes n} \circ 
\kappa_2^\Pi(x)$ are morphisms of $\NN$-graded $\cC$-algebras. 
\end{lemma}
 
{\em Proof.} This follows from (\ref{kappa1:kappa2:1}), 
(\ref{kappa1:kappa2:2}) and the fact that $\kappa_i^\Pi(m_\Pi)$, 
$\kappa_i^\Pi(\Delta_0)$ have degree $0$. 
\hfill \qed \medskip 

\begin{remark} Let $f$ be the image of $(\on{inj}_1^{\otimes 2}
\circ \on{can}) \otimes \on{pr}_0^{\otimes 2} \in \on{LBA}(\wedge^2 
\otimes S^{\otimes 2},S^{\otimes 2}) \subset \on{LBA}((S\circ \wedge^2)
\otimes S^{\otimes 2},S^{\otimes 2}) \to \on{LBA}_f(S^{\otimes 2},
S^{\otimes 2}) \simeq \Pi_f({\mathfrak 1} \underline\boxtimes {\mathfrak 1},
(S\underline\boxtimes S)^{\otimes 2}) = {\bf U}_{n,f}$, where $\on{can} : 
\wedge^2 \to {\bf id}^{\otimes 2}$ is the canonical morphism. We then have 
$(\Xi_f^{-1})^{\boxtimes 2} \circ \kappa_2^\Pi(r) = \kappa_1^\Pi(r) + f$. 
The morphisms ${\bf U}_n \to {\bf U}_{n,f}$, $x\mapsto 
\kappa_1^\Pi(x)$ and $x\mapsto (\Xi_f^{-1})^{\boxtimes 2} \circ 
\kappa_2^\Pi(x)$ will be denoted $x\mapsto x(r)$ and $x\mapsto x(r+f)$. 
\end{remark}

\subsection{The algebras ${\bf U}_{n,f}^{c_1...c_n}$}

For $c_1,...,c_n\in \{a,b\}$, we set ${\bf U}_{n,f}^{c_1...c_n}
:= \Pi_f({\mathfrak 1} \underline\boxtimes {\mathfrak 1},
F_{c_1} \boxtimes ... \boxtimes F_{c_n})$, where $F_a = S\underline\boxtimes 
{\bf 1}$ and $F_b = {\bf 1} \underline\boxtimes S$. Then 
${\bf U}_{n,f}^{c_1...c_n} \subset {\bf U}_{n,f}$ is a graded 
subspace. The diagrams (\ref{diag:subalg}) imply that it is also a 
subalgebra. 

Diagrams (\ref{diag:subcoalg}) also imply that for $\phi : [m]\to [n]$
partially defined, $\Delta^\phi$ takes
${\bf U}_{n,f}^{c_1...c_n}$ to ${\bf U}_{m,f}^{c'_1...c'_m}$, where 
$c'_1,...,c'_m$ are such that $c'_k = c_{\phi(k)}$ for any $k$ in the 
domain of $\phi$.  

In particular, $({\bf U}_{n,f}^{a...a})_{n\geq 0}$ is a $\NN$-graded 
$\cX$-algebra. 

\subsection{Hochschild cohomology of ${\bf U}_{n,f}$ and 
${\bf U}_{n,f}^{a...a}$} \label{cohoch}

The co-Hochschild complex of a $\cX$-vector space $(V_n)_{n\geq 0}$
is given by the differentials $d^1 : V_1 \to V_2$, $x\mapsto x^{12} - 
x^1 - x^2$, $d^2 : V_2 \to V_3$, $x\mapsto x^{12,3} - x^{1,23} - x^{2,3}
+ x^{1,2}$, etc. We denote the corresponding cohomology groups by $H^n(V_*)$. 

If $B = (B_{\rho,\sigma}) \in \on{Ob}({\bf Sch}_2)$ and $C\in \on{Sch}$, we set 
$B_C := \oplus_{\rho,\sigma} B_{\rho,\sigma} \otimes 
\on{LBA}(C\otimes Z_{\sigma},Z_{\rho})$. If $\alpha : C \to D$ is a morphism in 
LBA, then we set $B_\alpha := \on{Coker}(B_{D} \to B_{C}) = 
\oplus_{\rho,\sigma} B_{\rho,\sigma} \otimes \on{LBA}_{\alpha}(Z_{\sigma},
Z_{\rho})$. In the case of the above morphism 
$\alpha : \wedge^3 \otimes (S\circ \wedge^2) \to S\circ \wedge^2$, we define 
in this way spaces $B_f$.    

In particular, for $F,G\in \on{Ob}({\bf Sch}_{(1)})$, we have isomorphisms 
$\Pi({\mathfrak 1}\underline\boxtimes {\mathfrak 1},F\underline\boxtimes 
G) \simeq (c(F) \boxtimes c(G))_{\bf 1}$ and   
$\Pi_f({\mathfrak 1}\underline\boxtimes {\mathfrak 1},F\underline\boxtimes 
G) \simeq (c(F) \boxtimes c(G))_f$. 

\begin{lemma}
1) $H^n({\bf U}^{a...a}_{*,f}) \simeq (\wedge^n \boxtimes 
{\bf 1})_f$. If we set $C^n_a := \on{Ker}(d : {\bf U}^{a...a}_{n,f} \to 
{\bf U}^{a...a}_{n+1,f})$, 
then $\on{Alt} = (n!)^{-1}\sum_{\sigma\in \SG_n} \eps(\sigma) \sigma: 
{\bf U}^{a...a}_{n,f} \to {\bf U}^{a...a}_{n,f}$ restricts to a map $C_a^n
\to (\wedge^n \boxtimes {\bf 1})_f$, which factors through the above 
isomorphism. 

2) $H^n({\bf U}_{*,f}) \simeq (\Delta(\wedge^n))_f 
= \oplus_{p,q | p+q = n} (\wedge^p\boxtimes \wedge^q)_f$. 
If we set $C^n_a := \on{Ker}(d : {\bf U}^{a...a}_{n,f} \to 
{\bf U}^{a...a}_{n+1,f})$, 
then $\on{Alt} : 
{\bf U}_{n,f} \to {\bf U}_{n,f}$ restricts to a map $C^n
\to (\Delta(\wedge^n))_f$, which factors through the above 
isomorphism. 
\end{lemma}

{\em Proof.} 1) We have a co-Hochschild complex $S\to S^{\otimes 2} \to
S^{\otimes 3} \to ...$ in ${\bf Sch}$. It is defined as above, 
where $(x\mapsto x^{12,3})$ is replaced by the element of 
$\Delta_0^S \boxtimes \on{id}_S \in 
{\bf Sch}(S^{\otimes 2},S^{\otimes 3})$, etc. We express it as the 
sum of an acyclic complex $\Sigma_1 \to \Sigma_2 \to ...$ and a complex
with zero differential $\wedge^1 \to \wedge^2 \to ...$

The inclusion $\wedge^n\to S^{\otimes n}$ is given by the 
composition $\wedge^n \subset {\bf id}^{\otimes n} \subset 
S^{\otimes n}$. 

The inclusion $\Sigma_n \subset S^{\otimes n}$ is defined as follows: 
$\Sigma_n := \oplus_{k_1,...,k_n\neq 1} (\otimes_{i=1}^n S^{k_i}) \oplus 
\rho_n$, where $\rho_n \subset (S^1)^{\otimes n} = {\bf id}^{\otimes n}$
is the sum of all the images of the pairwise symmetrization maps 
$\on{id} + (ji) : {\bf id}^{\otimes n} \to {\bf id}^{\otimes n}$, 
where $i<j\in [n]$. Then we have a direct sum decomposition 
$S^{\otimes n} = \Sigma_n \oplus \wedge^n$. One checks that this 
is a decomposition of complexes, where $\wedge^n$ has zero differential. 

In particular, when $V$ is a vector space, the co-Hochschild 
complex $V \to S^2(V)\to...$ decomposes as the sum of the complexes
$\wedge^n(V)$ and $\Sigma_n(V)$. Since the cohomology is reduced to 
$\wedge^n(V)$, the complex $\Sigma_n(V)$ is acyclic. It has therefore a 
homotopy $\Sigma_n(V) \stackrel{K_n(V)}{\to} \Sigma_{n-1}(V)$, which 
has a propic version $\Sigma_n \stackrel{K_n}{\to} \Sigma_{n-1}$.  

Recall that  
${\bf U}_{n,f}^{a...a} \simeq (S^{\otimes n} \boxtimes {\bf 1})_f$. 
The co-Hochschild complex for the latter space decomposes as the sum of 
$(\wedge^n\boxtimes {\bf 1})_f$ with zero differential and 
$(\Sigma_n \boxtimes {\bf 1})_f$, which admits a homotopy and is therefore
acyclic. It follows that $H^n({\bf U}^{a...a}_{*,f}) = 
(\wedge^n\boxtimes {\bf 1})_f$. Then $C^n_a = (\wedge^n\boxtimes {\bf 1})_f
\oplus d((\Sigma_{n-1} \boxtimes {\bf 1})_f)$. The restriction of 
$\on{Alt}$ to $C^n_a$ is then the projection on the first summand of this 
decomposition, which implies the second result. This proves 1). 

Let us prove 2). We have ${\bf U}_{n,f} \simeq (S^{\otimes n} \boxtimes 
S^{\otimes n})_f \simeq (\Delta(S^{\otimes n}))_f$, where $\Delta : 
\on{Ob(Sch)} \to \on{Ob}(\on{Sch}_2)$ has been defined in Section 
\ref{sect:Delta}. 

We then have a decomposition ${\bf U}_{n,f} \simeq (\Delta(\wedge^n))_f
\oplus (\Delta(\Sigma_n))_f$, where the first complex has zero differential and 
the second complex admits a homotopy and is therefore acyclic. Therefore 
$H^n({\bf U}_{*,f}) = (\Delta(\wedge^n))_f$. 
As before, $C^n = (\Delta(\wedge^n))_f \oplus d(\Delta(\Sigma_{n-1})_f)$, 
and the restriction of $\on{Alt}$ to $C^n$ is the projection on the first
summand. This proves 2). \hfill \qed \medskip 

\begin{remark}
One can prove that for $B\in \on{Ob}({\bf Sch}_2)$, we have 
$$
B_{{\bf 1}} = \oplus_{N\geq 0} \Big( 
B(\on{Lie}(a_1,...,a_N) \oplus \on{Lie}(b_1,...,b_N)
)_{\sum_{i=1}^N (\alpha_i + \beta_i)}\Big)_{\SG_N}, 
$$
where $\on{Lie}(x_1,...,x_N)$ is the free Lie algebra with generators 
$x_1,...,x_N$, the generators $a_i,b_i$ have degrees $\alpha_i,\beta_i
\in \oplus_{i=1}^N (\NN\alpha_i \oplus \NN \beta_i)$. Here the index
$\sum_{i=1}^N (\alpha_i + \beta_i)$ means the part of degree $\sum_{i=1}^N
(\alpha_i + \beta_i)$, and the index $\SG_N$ means the space of coinvariants
w.r.t. the diagonal action of $\SG_N$ on generators $a_i,b_i$, $i =1,...,N$. 

Using the symmetrization map, we then get 
$$
{\bf U}_n \simeq \oplus_{N\geq 0}
\Big( \big( 
(\kk\langle a_1,...,a_N\rangle  
\kk\langle b_1,...,b_N\rangle)^{\otimes n} 
\big)_{\sum_{i=1}^n (\alpha_i+\beta_i)} \Big)_{\SG_N}, 
$$ 
where $\kk\langle x_1,...,x_N \rangle$ is the free algebra 
with generators $x_1,...,x_N$, and $\kk\langle a_1,...,a_N\rangle  
\kk\langle b_1,...,b_N\rangle$ is the image of the product map 
$\kk\langle a_1,...,a_N \rangle 
\otimes \kk\langle b_1,...,b_N\rangle
\to \kk\langle a_1,...,b_N \rangle$. So 
${\bf U}_n$ identifies with 
$(U(\g)^{\otimes n})_{\on{univ}}$ (see \cite{Enr:shuffles}). 

Then ${\bf U}_n \subset \oplus_{N\geq 0}
(\kk\langle a_1,...,b_N\rangle^{\otimes n})_{\SG_N}$. 
This inclusion is compatible with the $\cX$-structure on the 
right induced by the coalgebra structure of 
$\kk\langle a_1,...,b_N\rangle$. The space 
${\bf U}_n^{c_1,...,c_n}$ identifies with 
$\oplus_{N\geq 0} ((F^N_{c_1} \otimes ... \otimes F^N_{c_n}
)_{\sum_{i=1}^N (\alpha_i + \beta_i)})_{\SG_N}$, where $F^N_a = 
\kk\langle a_1,...,a_N \rangle$ and $F^N_b = 
\kk\langle b_1,...,b_N \rangle$, and if $(c_i,c'_i) \neq (b,a)$
for any $i$, then the product ${\bf U}_n^{c_1...c_n} 
\otimes {\bf U}_n^{c'_1...c'_n} \to {\bf U}_n$
is induced by the maps $F^N_c \otimes F^M_{c'}\to \kk\langle 
a_1,...,b_{N+M}\rangle$, $x(c_1,...,c_N) \otimes x'(c'_1,...,c'_M)
\mapsto x(c_1,...,c_N)x'(c'_{N+1},...,c'_{N+M})$. 
\end{remark}

\section{Injectivity of a map} \label{sect:inject}

Let $n,m$ be integers $\geq 0$. 
Define $\mu_m\in\on{LA}(T_m\otimes {\bf id},T_m)$ as the propic version 
of the map $x_1\otimes ... \otimes x_m \otimes x \mapsto \sum_{i=1}^m 
x_1\otimes ... \otimes [x_i,x] \otimes ... \otimes x_m$. Define a linear map  
$$
i_{n,m} : \on{LA}(T_n\otimes T_m,{\bf id}) \to \on{LA}(T_n \otimes 
T_m \otimes {\bf id},{\bf id}), 
$$
$$
\lambda\mapsto \lambda \circ (\on{id}_{T_n} \otimes \mu_m). 
$$
Define also a linear map 
$$
c_{n',n,m} : \on{LA}(T_{n'},T_n) \otimes \on{LA}(T_n \otimes T_m,{\bf id}) 
\to \on{LA}(T_{n'} \otimes T_m,{\bf id}),   
$$
$$
\lambda_0 \otimes \lambda \mapsto \lambda \circ (\lambda_0 \boxtimes
\on{id}_{T_m}). 
$$

Then $i_{n,m}$ is $\SG_n\times\SG_m$-equivariant, and the diagram 
$$
\begin{matrix}
\on{LA}(T_{n'},T_n) \otimes \on{LA}(T_n \otimes T_m,{\bf id}) & 
\stackrel{\on{id} \otimes i_{n,m}}{\to} & \on{LA}(T_{n'},T_n) 
\otimes \on{LA}(T_n \otimes T_m \otimes {\bf id},{\bf id})\\
\scriptstyle{c_{n',n,m}}\downarrow & & \downarrow\scriptstyle{c_{n',n,m+1}}\\
\on{LA}(T_{n'} \otimes T_m,{\bf id}) & \stackrel{i_{n',m}}{\to} & 
\on{LA}(T_{n'} \otimes T_m \otimes {\bf id},{\bf id})
\end{matrix}
$$
commutes.

\begin{lemma} \label{lemma:T}
There exists a map $p_{n,m} : \on{LA}(T_n\otimes T_m\otimes {\bf id},{\bf id})
\to \on{LA}(T_n \otimes T_m,{\bf id})$, such that $p_{n,m} \circ i_{n,m} = 
\on{id}$, which is $\SG_n \times
\SG_m$-equivariant, and such that the diagram 
$$
\begin{matrix}
\on{LA}(T_{n'},T_n) \otimes \on{LA}(T_n \otimes T_m \otimes {\bf id},{\bf id}) 
& \stackrel{\on{id} \otimes p_{n,m}}{\to} & \on{LA}(T_{n'},T_n) 
\otimes \on{LA}(T_n \otimes T_m,{\bf id})\\
\scriptstyle{c_{n',n,m+1}}\downarrow & & \downarrow\scriptstyle{c_{n',n,m}}\\
\on{LA}(T_{n'} \otimes T_m \otimes {\bf id},{\bf id}) & 
\stackrel{p_{n',m}}{\to} & 
\on{LA}(T_{n'} \otimes T_m,{\bf id})
\end{matrix}
$$
commutes. 
\end{lemma}

{\em Proof.} Let us first recall some results on free Lie algebras. 
Let $\on{L}(x_1,...,x_s)$ (resp., $\on{A}(x_1,...,x_s)$) 
the multilinear part of the free Lie (resp., associative) algebra
generated by $x_1,...,x_s$. Then $\on{L}(x_1,...,x_s) \subset 
\on{A}(x_1,...,x_s)$. For any $i = 1,...,s$, we have an isomorphism 
$\on{L}(x_1,...,x_s) \stackrel{\sim}{\to} \on{A}(x_1,...,x_{i-1},x_{i+1},...,
x_s)$, given by $\on{L}(x_1,...,x_s) \ni P(x_1,...,x_s) \mapsto P_{x_i}
(x_1,...,x_{i-1},x_{i+1},...,x_s)$, where $P_{x_i}$ is the element such that  
$P$ decomposes as $P_{x_i}x_i$ + sum of terms not ending with $x_i$. 
The inverse isomorphism is given by $\on{A}(x_1,...,x_{i-1},x_{i+1},...,x_s)
\ni Q\mapsto \on{ad}(Q)(x_s)$, where $\on{ad}(x_{i_1}...x_{i_{s-1}})(x_i) = 
[x_{i_1},[x_{i_2},...,[x_{i_{s-1}},x_i]]]$. 

Let us now prove the lemma. We have an isomorphism $\on{LA}(T_n,{\bf id})
\simeq \on{L}(x_1,...,x_n)$. The map $i_{n,m}$ is given by 
$$
i_{n,m} : \on{L}(x_1,...,x_{n+m}) \to \on{L}(x_1,...,x_{n+m+1}), 
$$
$$ 
P(x_1,...,x_{n+m}) \mapsto \sum_{i=1}^m P(x_1,...,[x_{n+i},x_{n+m+1}],
...,x_{n+i}).  
$$

Define 
$$
p_{n,m} : \on{L}(x_1,...,x_{n+m+1}) \to \on{L}(x_1,...,x_{n+m}), 
$$
$$
Q(x_1,...,x_{n+m+1}) \mapsto {1\over m} \sum_{i=1}^m 
\on{ad}(Q_{x_{n+i}x_{n+m+1}})(x_{n+i}),  
$$
where $Q_{x_{n+i}x_{n+m+1}}$ is the element of 
$\on{A}(x_1,...,x_{n+i-1},x_{n+i+1},...,x_{n+m})$ such that 
$$
Q = Q_{x_{n+i}x_{n+m+1}}x_{n+i}x_{n+m+1} 
\on{\ +\ terms\ not\ ending\ by\ } x_{n+i}x_{n+m+1}.
$$ 

Let us show that $p_{n,m} \circ i_{n,m} = \on{id}$. 
If $Q_i := P(x_1,...,[x_{n+i},x_{n+m+1}],...,x_{n+m})$, 
then $(Q_i)_{x_{n+j}x_{n+m+1}} = 0$ if $i\neq j$, and 
equals $P_{x_{n+j}}$ if $i=j$. So $(i_{n,m}(P))_{x_{n+j}x_{n+m+1}}
= P_{x_{n+j}}$. So 
$\on{ad}\Big( \big( i_{n,m}(P)\big)_{x_{n+j}x_{x_{n+m+1}}}
\Big) (x_{n+j}) = \on{ad}(P_{x_{n+j}})(x_{n+j}) = P$. 
Averaging over $i \in [m]$, we get $p_{n,m}(i_{n,m}(P)) = P$. 

Let us show that $p_{n,m}$ is $\SG_n\times\SG_m$-equivariant. 
The $\SG_n$-equivariance is clear.  Let us show the $\SG_m$-equivariance. 
Let $\tau\in\SG_m$. We have $Q^\tau(x_1,...,x_{n+m+1}) = 
Q(x_1,...,x_n,x_{n+\tau(1)},...,x_{n+\tau(m)},x_{n+m+1})$, so 
$(Q^\tau)_{x_{n+\tau(i)}x_{n+m+1}} = (Q_{x_{n+i}x_{n+m+1}})^\tau$. 
Then 
\begin{align*}
p_{n,m}(Q^\tau) & = {1\over m} \sum_{i=1}^m \on{ad}\big( 
(Q^\tau)_{x_{n+i}x_{n+m+1}}\big) (x_{n+i})
= {1\over m} \sum_{i=1}^m \on{ad}\big( 
(Q^\tau)_{x_{n+\tau(i)}x_{n+m+1}}\big) (x_{n+\tau(i)})
\\ & = {1\over m} \sum_{i=1}^m \on{ad}\big( 
(Q_{x_{n+i}x_{n+m+1}})^\tau\big) (x_{n+\tau(i)})
= (p_{n,m}(Q))^\tau, 
\end{align*}
which proves the $\SG_m$-equivariance. 

Let us prove that the announced diagram commutes. 
Let $\lambda_0\in \on{LA}(T_{n'},T_n)$ and $P\in 
\on{LA}(T_n \otimes T_m \otimes {\bf id},{\bf id})$.  
We must show that the images of $\lambda_0 \otimes P$
in $\on{LA}(T_{n'} \otimes T_m,{\bf id})$ by two maps coincide. 
By linearity, we may assume that $\lambda_0$ has the form 
$x_1\otimes .... \otimes x_{n'} \mapsto P_1(x_i,i\in f^{-1}(1)) 
\otimes ... \otimes P_n(x_i,i\in f^{-1}(n))$, where $f : [n'] \to [n]$
is a map and $P_i\in \on{F}(x_{i'},i'\in f^{-1}(i))$. The commutativity of 
the diagram then follows from the equality 
\begin{align*}
& P_{x_{n+i}x_{n+m+1}}
\big( P_1(x_i,i\in f^{-1}(1)),...,P_n(x_i,i\in f^{-1}(n)),
x_{n'+1},...,x_{n'+m+1}  \big)
\\ &  = \Big( P(P_1(x_i,i\in f^{-1}(1)),...,P_n(x_i,i\in f^{-1}(n)),
x_{n'+1},...,x_{n'+m+1}) \Big)_{x_{n'+i}x_{n'+m+1}}. 
\end{align*}
\hfill \qed \medskip 

If $Z\in\on{Irr(Sch)}$, we now define 
$$
\mu_Z \in \on{LA}(Z\otimes {\bf id},Z)
$$ as follows. Let $n$ be an integer $\geq 0$. 
The decomposition $T_n = 
\oplus_{Z\in \on{Irr(Sch)}, |Z| = n} Z\otimes \pi_Z$
gives rise to an isomorphism $\on{LA}(T_n\otimes {\bf id},T_n) \simeq 
\oplus_{Z,W\in\on{Irr(Sch)},|Z| = |W| = n}
\on{LA}(Z\otimes {\bf id},W) \otimes \on{Vect}(\pi_Z,\pi_W)$. On the 
other hand, $\mu_n$ has the $\SG_n$-equivariance property 
$\mu_n \circ (\sigma \otimes \on{id}_{\bf id}) = 
\sigma \circ \mu_n$ for any $\sigma\in\SG_n$. It follows that 
$\mu_n$ decomposes as $\oplus_{Z\in\on{Irr(Sch)},|Z| = n} \mu_Z \otimes 
\on{id}_{\pi_Z}$. This defines $\mu_Z$ for any $Z\in\on{Irr(Sch)}$ with 
$|Z| = n$. 

For $W,Z\in\on{Irr(Sch)}$, define 
$$
i_{W,Z} : \on{LA}(W\otimes Z,{\bf id}) \to 
\on{LA}(W\otimes Z \otimes {\bf id},{\bf id}), \quad 
\lambda \mapsto \lambda \circ (\on{id}_W \otimes \mu_Z). 
$$

For $Z,W,W'\in\on{Irr(Sch)}$, define 
$$
c_{W',W,Z} : \on{LA}(W',W) \otimes \on{LA}(W\otimes Z,{\bf id})
\to \on{LA}(W'\otimes Z,{\bf id}), \quad 
\lambda_0 \otimes \lambda \mapsto \lambda \circ (\lambda_0 \otimes 
\on{id}_Z).   
$$

Then the diagram 
$$
\begin{matrix}
\on{LA}(W',W) \otimes \on{LA}(W\otimes Z,{\bf id}) & 
\stackrel{\on{id} \otimes i_{W,Z}}{\to} & \on{LA}(W',W) \otimes \on{LA}(W\otimes Z \otimes {\bf id},
{\bf id}) \\
\scriptstyle{c_{W',W,Z}}\downarrow & & \downarrow 
\scriptstyle{c_{W',W,Z\otimes {\bf id}}}\\
\on{LA}(W'\otimes Z,{\bf id}) & \stackrel{i_{W',Z}}{\to}& \on{LA}(W'\otimes Z
\otimes {\bf id},{\bf id}) 
\end{matrix}
$$
commutes. 

For $W,Z\in\on{Irr(Sch)}$, define a linear map 
$$
p_{W,Z} : \on{LA}(W\otimes Z\otimes {\bf id},{\bf id}) \to 
\on{LA}(W\otimes Z,{\bf id})
$$
as follows. For $n,m$ integers $\geq 0$, the decompositions 
$T_n = \oplus_{W\in\on{Irr(Sch)},|W| = n} W \otimes \pi_W$, 
$T_m = \oplus_{Z\in\on{Irr(Sch)},|Z| = m} Z \otimes \pi_Z$
give rise to a decomposition 
\begin{align*}
& \on{Vect}(\on{LA}(T_n \otimes T_m,{\bf id}),\on{LA}(T_n \otimes T_m \otimes 
{\bf id},{\bf id})) 
\\ & 
\simeq \oplus_{\begin{matrix} \scriptstyle{W,W',Z,Z'||W| = |W'| = n,} \\ 
\scriptstyle{|Z| = |Z'| = m }\end{matrix}}
\on{Vect}(\on{LA}(W\otimes Z,{\bf id}),\on{LA}(W' \otimes Z' \otimes 
{\bf id},{\bf id})) \otimes \on{Vect}(\pi_{W}\otimes \pi_{Z},\pi_{W'}\otimes 
\pi_{Z'})
\end{align*}
which is $\SG_n \times \SG_m$-equivariant. Then $p_{n,m}\in 
\on{Vect}(\on{LA}(T_n \otimes T_m,{\bf id}),\on{LA}(T_n \otimes T_m \otimes 
{\bf id},{\bf id}))$ is $\SG_n \times \SG_m$-invariant, which implies that 
it decomposes as $\sum_{W,Z\in\on{Irr(Sch)} | |W|=n,|Z|=m} p_{W,Z} \otimes
\on{id}_{\pi_W \otimes \pi_Z}$. This defines $p_{W,Z}$ for $W,Z\in 
\on{Irr(Sch)}$. 

\begin{proposition}
We have $p_{W,Z} \circ i_{W,Z} = \on{id}$, and the diagram 
\begin{equation} \label{comm:p:c}
\begin{matrix}
\on{LA}(W',W) \otimes \on{LA}(W\otimes Z \otimes {\bf id},{\bf id}) & 
\stackrel{\on{id} \otimes p_{W,Z}}{\to} & 
\on{LA}(W',W) \otimes \on{LA}(W\otimes Z,{\bf id})\\
\scriptstyle{c_{W',W,Z\otimes {\bf id}}}\downarrow  & & 
\downarrow\scriptstyle{c_{W',W,Z}}\\ 
\on{LA}(W'\otimes Z \otimes {\bf id},{\bf id}) & \stackrel{p_{W',Z}}{\to} & 
\on{LA}(W'\otimes Z,{\bf id}) 
\end{matrix}
\end{equation}
commutes. 
\end{proposition}   

{\em Proof.} This is obtained by taking the isotypic components of 
the statements of Lemma \ref{lemma:T}, and using that $i_{n,m} = 
\oplus_{W,Z\in\on{Irr(Sch)} | |W| = n, |Z| = m} i_{W,Z} \otimes 
\on{id}_{\pi_W \otimes \pi_Z}$. 
\hfill \qed \medskip 

We now prove: 

\begin{proposition} \label{prop:inject}
The map $\on{LBA}_f({\bf id},{\bf id}) \to \on{LBA}_f(\wedge^2,{\bf id})$, 
$x\mapsto x\circ \mu$, is injective. 
\end{proposition}

{\em Proof.} Let $\alpha : C \to D$ be a morphism in $\on{LBA}$. We will 
prove that $i_\alpha : \on{LBA}_\alpha({\bf id},{\bf id}) \to 
\on{LBA}_\alpha(\wedge^2, {\bf id}) \subset 
\on{LBA}_\alpha({\bf id}^{\otimes 2},{\bf id})$, $x\mapsto x\circ \mu$ is 
injective. For this, we will construct a map $p_\alpha : 
\on{LBA}_\alpha({\bf id}^{\otimes 2},{\bf id}) \to 
\on{LBA}_\alpha({\bf id},{\bf id})$, such that 
$p_\alpha \circ i_\alpha = \on{id}$. 

The first map is the vertical cokernel of the commutative diagram 
$$
\begin{matrix}
\on{LBA}(D\otimes {\bf id},{\bf id}) & \stackrel{i_D}{\to}& 
\on{LBA}(D\otimes {\bf id}^{\otimes 2},{\bf id})\\
\scriptstyle{- \circ (\alpha \otimes \on{id}_{{\bf id}})}
\downarrow 
& & \downarrow
\scriptstyle{- \circ (\alpha \otimes \on{id}_{{\bf id}})}
 \\
\on{LBA}(C\otimes {\bf id},{\bf id}) & 
\stackrel{i_C}{\to} & 
\on{LBA}(C\otimes {\bf id}^{\otimes 2},{\bf id})  
\end{matrix}
$$
where $i_X(x) = x \circ (\on{id}_X \otimes \mu)$ for $X = C,D$. 

We will construct a commutative diagram 
\begin{equation} \label{comm:CD}
\begin{matrix}
\on{LBA}(D\otimes {\bf id}^{\otimes 2},{\bf id}) & \stackrel{p_D}{\to}& 
\on{LBA}(D\otimes {\bf id},{\bf id}) \\
\scriptstyle{- \circ (\alpha \otimes \on{id}_{{\bf id}^{\otimes 2}})}
\downarrow 
 & & \downarrow
\scriptstyle{- \circ (\alpha \otimes \on{id}_{{\bf id}^{\otimes 2}})}
\\
\on{LBA}(C\otimes {\bf id}^{\otimes 2},{\bf id}) & \stackrel{p_C}{\to}& 
\on{LBA}(C\otimes {\bf id},{\bf id})  
\end{matrix}
\end{equation}
such that $p_C \circ i_C = \on{id}$ and $p_D \circ i_D = \on{id}$; then 
we will define $p_\alpha$ as the vertical cokernel of this diagram. 

Set $A_C := \on{LBA}(C\otimes {\bf id},{\bf id})$, $A'_C := \on{LBA}
(C\otimes {\bf id}^{\otimes 2},{\bf id})$. 
Let us study the map 
$$
i_C : A_C \to A'_C. 
$$
We have 
$A_C = \oplus_{W,Z\in\on{Irr(Sch)}} A_C(W,Z)$
and $A'_C = \oplus_{W,Z',Z''\in\on{Irr(Sch)}} A'_C(W,Z',Z'')$, 
where 
$$
A_C(W,Z) := \on{LCA}(C,W) \otimes \on{LCA}({\bf id},Z) \otimes 
\on{LA}(W\otimes Z,{\bf id}), 
$$
$$
A'_C(W,Z',Z'') := \on{LCA}(C,W) \otimes \on{LCA}({\bf id},Z') \otimes
\on{LCA}({\bf id},Z'') \otimes 
\on{LA}(W\otimes Z'\otimes Z'',{\bf id}). 
$$  
Set $A''_C := \oplus_{W,Z\in\on{Irr(Sch)}} A'_C(W,Z,{\bf id})$. 
We have a natural projection map  $A'_C \to A''_C$. 

Then the composition 
$A_C \stackrel{i_C}{\to} A'_C \to A''_C$ is the direct sum over 
$W,Z$ of the maps $A_C(W,Z) \to A'_C(W,Z,{\bf id})$, given by 
\begin{align*}
& \on{LCA}(C,W) \otimes \on{LCA}({\bf id},Z) \otimes 
\on{LA}(W\otimes Z, {\bf id}) \to 
\on{LCA}(C,W) \otimes \on{LCA}({\bf id},Z) \otimes
\on{LCA}({\bf id},{\bf id}) 
\\ & \otimes 
\on{LA}(W\otimes Z\otimes {\bf id},{\bf id}),  
\quad 
\kappa_C \otimes \kappa_{{\bf id}} \otimes \lambda
\mapsto \kappa_C \otimes \kappa_{{\bf id}} \otimes 1 
\otimes i_{W,Z}(\lambda). 
\end{align*}

Define the map 
$$
p_C : A'_C \to A_C
$$
as the composition $A'_C \to A''_C \to A_C$, where the first map 
is the natural projection and the second 
map is the direct sum over $W,Z$ of the maps $A_C(W,Z)\to 
A'_C(W,Z,{\bf id})$, given by 
\begin{align*}
& \on{LCA}(C,W) \otimes \on{LCA}({\bf id},Z) \otimes \on{LCA}({\bf id},{\bf id}) 
\otimes \on{LA}(W\otimes Z\otimes {\bf id},{\bf id})
\to 
\on{LCA}(C,W) \otimes \on{LCA}({\bf id},Z) \\ & \otimes 
\on{LA}(W\otimes Z,{\bf id}), 
\quad 
\kappa_C \otimes \kappa_{{\bf id}} \otimes 1 \otimes \lambda'
\mapsto \kappa_C \otimes \kappa_{{\bf id}} \otimes p_{W,Z}(\lambda').  
\end{align*}
Then $p_{W,Z}\circ i_{W,Z} = \on{id}$ implies that $p_C \circ i_C = \on{id}$. 

Let us prove that (\ref{comm:CD}) commutes. For this, we will prove that 
$$
\begin{matrix}
\on{LBA}(C,D) \otimes \on{LBA}(D\otimes {\bf id}^{\otimes 2},{\bf id}) 
& \stackrel{\on{id} \otimes p_D}{\to}& 
\on{LBA}(C,D) \otimes \on{LBA}(D\otimes {\bf id},{\bf id}) \\
\downarrow  & & \downarrow \\
\on{LBA}(C\otimes {\bf id}^{\otimes 2},{\bf id}) 
& \stackrel{p_D}{\to}& 
\on{LBA}(C\otimes {\bf id},{\bf id}) 
\end{matrix}
$$
commutes, where the vertical maps are $\alpha\otimes x\mapsto x \circ (\alpha
\otimes \on{id}_{{\bf id}^{\otimes 2}})$ (right map) and 
$\alpha\otimes x\mapsto x \circ (\alpha \otimes \on{id}_{{\bf id}})$ 
(left map). 

This diagram is the same as 
$$
\begin{matrix}
\oplus_{U,W,Z',Z''\in\on{Irr(Sch)}} \on{LCA}(C,U) \otimes \on{LA}(U,D)
& \stackrel{(1)}{\to}
& \oplus_{U,W,Z\in\on{Irr(Sch)}} \on{LCA}(C,U) \\
\otimes \on{LCA}(D,W) \otimes \on{LCA}({\bf id},Z') \otimes \on{LCA}({\bf
id},Z'') & & \otimes \on{LA}(U,D) \otimes \on{LCA}(D,W) \\
\otimes \on{LA}(W\otimes Z'\otimes Z'',{\bf id})& 
& \otimes \on{LCA}({\bf id},Z) \otimes \on{LA}(W\otimes
Z,{\bf id})\\
\scriptstyle{(3)}\downarrow & & \downarrow\scriptstyle{(4)}\\
\oplus_{V,Z',Z''\in\on{Irr(Sch)}} \on{LCA}(C,V) \otimes \on{LCA}({\bf id},Z')
& \stackrel{(2)}{\to}& 
\oplus_{V,Z\in\on{Irr(Sch)}} \on{LCA}(C,V)  \\ 
\otimes \on{LCA}({\bf id},Z'') \otimes \on{LA}(V\otimes Z'\otimes Z'',{\bf id})
& &  \otimes \on{LCA}({\bf id},Z) \otimes \on{LA}(V\otimes Z,{\bf id}) 
\end{matrix}
$$
where (1) is zero on the components with $Z''\neq {\bf id}$; it takes the 
component $(U,W,Z,{\bf id})$ to the component $(U,W,Z)$ by the map 
$\on{id} \otimes \on{id} \otimes \on{id} \otimes 1 \otimes p_{W,Z}$; 

(2) is zero on the components with $Z''\neq {\bf 1}$; its takes the 
component $(U,Z,{\bf id})$ to the component $(U,Z)$ by the map 
$\on{id} \otimes 1 \otimes p_{U,Z}$; 

(3) is the composition of the natural map 
$$
\on{LA}(U,D) \otimes 
\on{LCA}(D,W) \to \on{LBA}(U,W)  \simeq 
\oplus_{V\in\on{Irr(Sch)}} \on{LCA}(U,V) \otimes 
\on{LA}(V,W),
$$ 
of the composition $\on{LCA}(C,U) \otimes \on{LCA}(U,V) 
\to \on{LCA}(C,V)$ and of the map 
$\on{LA}(V,W) \otimes \on{LA}(W\otimes Z'\otimes Z'',{\bf id}) 
\to \on{LA}(V,\otimes Z'\otimes Z'',{\bf id})$, $\alpha \otimes \beta 
\mapsto \beta \circ (\alpha\otimes \on{id}_{Z'\otimes Z''})$; 

(4) is the composition of same maps, where in the last step $Z'\otimes Z''$
is replaced by $Z$. 

The commutativity of the diagram formed by these maps 
then follows from that of (\ref{comm:p:c}). 
\hfill \qed \medskip

\section{Quantization functors} \label{sect:quant}

\subsection{Definition}
A quantization functor is a prop morphism $Q : \on{Bialg} \to 
S({\bf LBA})$, such that: 

(a) the composed morphism 
$\on{Bialg} \stackrel{Q}{\to} S({\bf LBA}) \to 
S({\bf Sch})$ (where the second morphism is given by the specialization $\mu =
\delta = 0$) is the propic version of the bialgebra structure of the symmetric
algebras $S(V)$, where the elements of $V$ are primitive, and 

(b) (classical limits) 
$\on{pr}_1 \circ Q(m) \circ (\on{inj}_1^{\otimes 2}
\circ \on{can})
\in {\bf LBA}(\wedge^2,{\bf id})$ has the form $\mu$ + terms of positive 
$\delta$-degree, and 
$(\on{Alt} \circ \on{pr}_1^{\otimes 2}) \circ Q(\Delta)
\circ \on{inj}_1
\in {\bf LBA}({\bf id},\wedge^2) = \delta$ + terms of 
positive $\mu$-degree.  

Here $\on{inj}_1 : {\bf id} \to S$ and $\on{pr}_1 : 
S \to {\bf id}$ are the canonical 
injection and projection maps, and $\on{inc} : \wedge^2 \to T_2$, 
$\on{Alt} : T_2 \to \wedge^2$ are the 
inclusion and alternation maps. 

Note that (a) implies that $Q(\eta) = \on{inj}_0\in {\bf LBA}({\bf 1},S)$, 
and $Q(\eps) = \on{pr}_0\in {\bf LBA}(S,{\bf 1})$, where $\on{inj}_0 : {\bf 1}
\to S$ and $\on{pr}_0 : S \to {\bf 1}$ are the natural injection and 
projection.  

The quantization functors $Q,Q'$ are called equivalent iff there exists an 
inner automorphism $\theta(\xi_0)$ of $S({\bf LBA})$, such that $Q' = 
\theta(\xi_0) \circ Q$. 

\subsection{Construction of quantization functors}

In \cite{EK}, Etingof and Kazhdan constructed a quantization functor
corresponding to each associator $\Phi$. This construction can be 
described as follows (\cite{Enr:coh}). 

Let $\t_n$ be the Lie algebra with generators $t_{ij}$, $1\leq i\neq j \leq n$
and relations $t_{ij} = t_{ji}$, $[t_{ij},t_{ik} + t_{jk}] = 0$, 
$[t_{ij},t_{kl}] = 0$ ($i,j,k,l$ distinct). It is graded by 
$\on{deg}(t_{ij}) = 1$. We have a graded  algebra morphism 
$U(\t_n) \to {\bf U}_n$, taking $t_{ij}$ to $t^{i,j}$, 
where $t\in {\bf U}_2$ is $r + r^{2,1}$.  

The family $(\t_n)_{n\geq 0}$ is a $\cC$-Lie algebra, and 
$U(\t_n) \to {\bf U}_n$ is a morphism of $\cC$-algebras. 

An associator $\Phi$ is an element of $\wh{U(\t_3)}^\times_1$, 
satisfying certain relations (see \cite{Dr:Gal}, where it is proved that 
associators exist over $\kk$). We fix an associator $\Phi$; we also 
denote by $\Phi$ its image in $(\wh{\bf U}_3)^\times_1$.  
 
One constructs ${\on J}\in (\wh{\bf U}_2)^\times_1$, 
such that ${\on J} = 1 - r/2 + ...$, and 
\begin{equation} \label{eq:J}
\on{J}^{1,2} \on{J}^{12,3} = \on{J}^{2,3} 
\on{J}^{1,23} \Phi. 
\end{equation}
Then one sets 
$$
\on{R} := \on{J}^{2,1} e^{t/2} \on{J}^{-1} 
\in (\wh{\bold U}_2)^\times_1. 
$$

Using $\on{J}$ and $\on{R}$, we will define elements of 
the quasi-bi-multiprop $\hbox{\boldmath$\Pi$\unboldmath}$. 

We define 
$$
\Delta_\Pi\in 
\hbox{\boldmath$\Pi$\unboldmath}
(S\underline\boxtimes S,(S\underline\boxtimes S)^{\boxtimes 2}), \quad 
\on{Ad}(\on{J}) \in 
\hbox{\boldmath$\Pi$\unboldmath}
((S\underline\boxtimes S)^{\boxtimes 2},
(S\underline\boxtimes S)^{\boxtimes 2}).  
$$

One checks that the elements $m_\Pi^{(2)} \boxtimes m_\Pi^{(2)} \in 
\hbox{\boldmath$\Pi$\unboldmath}
((S\underline\boxtimes S)^{\boxtimes 6},(S\underline\boxtimes S)^{\boxtimes 2})$, 
$(142536)\in \hbox{\boldmath$\Pi$\unboldmath}
((S\underline\boxtimes S)^{\boxtimes 6},(S\underline\boxtimes S)^{\boxtimes 6})$
and $\on{J} \boxtimes \on{id}_{(S\underline\boxtimes S)^{\boxtimes 2}}
\boxtimes \on{J}^{-1} \in \hbox{\boldmath$\Pi$\unboldmath}
((S\underline\boxtimes S)^{\boxtimes 2},(S\underline\boxtimes S)^{\boxtimes 6})$
are composable, and we set 
$$
\on{Ad}(\on{J}) := (m_\Pi^{(2)} \boxtimes m_\Pi^{(2)})
\circ (142536)\circ \big( \on{J} \boxtimes 
\on{id}_{(S\underline\boxtimes S)^{\boxtimes 2}} \boxtimes 
\on{J}^{-1}\big) 
\in \hbox{\boldmath$\Pi$\unboldmath}
((S\underline\boxtimes S)^{\boxtimes 2}, 
(S\underline\boxtimes S)^{\boxtimes 2}). 
$$
A graph for this element is as follows. Set $F_1 = ... = G'_2 = S$, 
then this is an element of $\hbox{\boldmath$\Pi$\unboldmath} 
((F_1\boxtimes F_2) \underline\boxtimes 
(G_1\boxtimes G_2), (F'_1\boxtimes F'_2) \underline\boxtimes 
(G'_1\boxtimes G'_2))$, and the edges are 
$F_i \to F'_j$, $G'_i \to G_j$, $G'_j \to F'_i$ 
($i,j = 1,2$).  

Now $\on{Ad}(\on{J})$ and $\Delta_0$ can be composed, and we set 
$$ 
\Delta_\Pi := \on{Ad}(\on{J}) \circ \Delta_0 \in
\hbox{\boldmath$\Pi$\unboldmath}
(S\underline\boxtimes S, 
(S\underline\boxtimes S)^{\boxtimes 2}). 
$$
A graph for this element is as follows. If we set 
$F = ... = G'_2 = S$, then this is an element of 
$\hbox{\boldmath$\Pi$\unboldmath}(F\underline\boxtimes G,
(F'_1\boxtimes F'_2)\underline\boxtimes (G'_1\boxtimes G'_2))$. The 
vertices are then $F\to F'_i$, $G'_i \to G$, $G'_i \to F'_j$
($i,j = 1,2$). 

The elements $m_\Pi,\Delta_\Pi$ then satisfy (\ref{m:assoc}); moreover, 
the following elements make sense, and the identities hold:   
\begin{equation} \label{other:ids}
\Delta_\Pi\circ m_\Pi = (m_\Pi\boxtimes m_\Pi) \circ (1324) \circ 
(\Delta_\Pi\boxtimes\Delta_\Pi), 
\quad   
(\Delta_\Pi \boxtimes 
\on{id}_{S\underline\boxtimes S})\circ \Delta_\Pi
= (\on{id}_{S\underline\boxtimes S}
\boxtimes \Delta_\Pi)\circ \Delta_\Pi.  
\end{equation}

In particular, $\overline{m}_\Pi := \Delta_\Pi^* \circ (21)\in 
\hbox{\boldmath$\Pi$\unboldmath}
((S\underline\boxtimes S)^{\boxtimes 2},
S\underline\boxtimes S)$ and $\overline{\Delta}_\Pi := m_\Pi^*\in 
\hbox{\boldmath$\Pi$\unboldmath}(S\underline\boxtimes S,
(S\underline\boxtimes S)^{\boxtimes 2})$ satisfy relations (\ref{m:assoc}), 
(\ref{other:ids}). 

Moreover, $\on{R}\in \hbox{\boldmath$\Pi$\unboldmath}
({\mathfrak 1}\underline\boxtimes{\mathfrak 1},
(S\underline\boxtimes S)^{\boxtimes 2})$ satisfies the 
quasitriangular identities 
\begin{equation} \label{QT:1}
\big( \Delta_\Pi \boxtimes
\on{id}_{S\underline\boxtimes S}
\big) \circ \on{R} = \big(
\on{id}_{(S\underline\boxtimes S)^{\boxtimes 2}}
\boxtimes m_\Pi \big) \circ (1324) \circ (\on{R}
\boxtimes \on{R})
\end{equation}
and
\begin{equation} \label{QT:2}
\big( \on{id}_{S\underline\boxtimes S}
\boxtimes \Delta_\Pi \big) \circ \on{R} = (132) \circ \big(
m_\Pi \boxtimes \on{id}_{(S\underline\boxtimes S)^{\boxtimes 2}}\big) 
\circ (1324) \circ (\on{R} \boxtimes \on{R}). 
\end{equation}

Define 
$$
\ell := 
\big( \on{id}_{S\underline\boxtimes S} \boxtimes
\on{can}^*_{S\underline\boxtimes S}\big) \circ 
\big( \on{R}\boxtimes  
\on{id}_{S\underline\boxtimes S}\big)
\in \hbox{\boldmath$\Pi$\unboldmath}
(S\underline\boxtimes S,S\underline\boxtimes S). 
$$

The following proposition is a consequence of the quasitriangular identities 
(\ref{QT:1}), (\ref{QT:2}): 

\begin{proposition} 
The following elements are defined, and the equations hold: 
$$
m_\Pi \circ \ell^{\boxtimes 2} = \ell \circ \overline{m}_\Pi, \quad 
\Delta_\Pi \circ \ell = \ell^{\boxtimes 2} \circ \overline{\Delta}_\Pi. 
$$
\end{proposition}

The flatness statement of \cite{Enr:coh} can be restated as follows: 

\begin{proposition} There exist elements $\on{R}_+\in 
\hbox{\boldmath$\Pi$\unboldmath}
(S\underline\boxtimes {\mathfrak 1},S\underline\boxtimes S)$
and $\on{R}_-\in \hbox{\boldmath$\Pi$\unboldmath}
({\mathfrak 1}\underline\boxtimes S,S\underline\boxtimes S)$, such that 
\begin{equation} \label{fact:R}
\on{R} = (\on{R}_+ \boxtimes \on{R}_-) \circ \on{can}_{S\underline\boxtimes
{\mathfrak 1}}. 
\end{equation}
Moreover, $\on{R}_\pm$ are right-invertible, i.e., there exist 
$\on{R}_+^{(-1)}\in \hbox{\boldmath$\Pi$\unboldmath}  
(S\underline\boxtimes S,S\underline\boxtimes
{\mathfrak 1})$ and $\on{R}_-^{(-1)}\in
\hbox{\boldmath$\Pi$\unboldmath}
(S\underline\boxtimes S,{\mathfrak 1}\underline\boxtimes S)$, 
with graphs $F\to F'$ and $G'\to G$
[where $\on{R}_+^{(-1)}$ is viewed as an 
element of $\hbox{\boldmath$\Pi$\unboldmath}
(F\underline\boxtimes G,F'\underline\boxtimes {\mathfrak 1})$
and $\on{R}_-^{(-1)}$ as an 
element of $\hbox{\boldmath$\Pi$\unboldmath}
(F\underline\boxtimes G,{\mathfrak 1}
\underline\boxtimes G')$], such that $\on{R}_+^{(-1)} \circ \on{R}_+ 
= \on{id}_{S\underline\boxtimes {\mathfrak 1}}$ and 
$\on{R}_-^{(-1)} \circ \on{R}_-
= \on{id}_{{\mathfrak 1}\underline\boxtimes S}$
(where the compositions are well-defined). 
\end{proposition}

Notice that $(\on{R}_+,\on{R}_-)$ is uniquely defined only up to 
a transformation 
$$
(\on{R}_+,\on{R}_-) \to (\on{R}_+ \circ \on{R}'',
\on{R}_- \circ ((\on{R}'')^*)^{-1}), 
$$
where $\on{R}''\in \hbox{\boldmath$\Pi$\unboldmath}(S\underline\boxtimes
{\mathfrak 1}, S\underline\boxtimes {\mathfrak 1})^\times$.  
This transformation will not change the equivalence class of $Q$. 

(\ref{fact:R}) implies that 
$$
\ell = \on{R}_+ \circ \on{R}_-^*. 
$$

\begin{proposition}
The following elements are defined, and the equations hold: 
\begin{equation} \label{pre:ma}
\on{R}_+^{(-1)} \circ m_\Pi \circ \on{R}_+^{\boxtimes 2} 
= \on{R}_-^* \circ \overline{m}_\Pi \circ (\on{R}_-^{(-1)*})^{\boxtimes 2}, 
\quad 
(\on{R}_+^{(-1)})^{\boxtimes 2} \circ \Delta_\Pi \circ \on{R}_+ 
= (\on{R}_-^*)^{\boxtimes 2} \circ \overline{\Delta}_\Pi \circ 
(\on{R}_-^{(-1)})^*.  
\end{equation}
Let $m_a\in \hbox{\boldmath$\Pi$\unboldmath}
(S\underline\boxtimes {\mathfrak 1}, 
(S\boxtimes S)\underline\boxtimes {\mathfrak 1})$
be the value of both sides of the first identity of (\ref{pre:ma}), 
and let $\Delta_a\in \hbox{\boldmath$\Pi$\unboldmath}
((S\boxtimes S)\underline\boxtimes {\mathfrak 1}, 
S\underline\boxtimes {\mathfrak 1})$ be the common value of both sides of 
the second identity. Then $m_a,\Delta_a$ satisfy (\ref{m:assoc}) 
and (\ref{other:ids}). 

Then there is a unique morphism $Q: \on{Bialg} \to S({\bf LBA})$, 
such that $Q(m) =$ the element of ${\bf LBA}(S^{\otimes 2},S)$
corresponding to $m_a$, $Q(\Delta) =$ the element of 
${\bf LBA}(S,S^{\otimes 2})$ corresponding to $\Delta_a$, 
$Q(\eps) =$ the element of ${\bf LBA}(S,{\bf 1})$ corresponding to 
$1\in \kk$, $Q(\eta) =$ the element of ${\bf LBA}({\bf 1},S)$ corresponding 
to $1\in \kk$.  
\end{proposition}

{\em Proof.} The proof follows that of the following
statement: let $\cS$ be a symmetric tensor category, let 
$A,X,B\in\on{Ob}(\cS)$. Assume that $m_A\in \cS(A^{\otimes 2},A)$,
$\Delta_A\in \cS(A,A^{\otimes 2})$,...  is a bialgebra structure
on $A$ in the category $\cC$. Let similarly $(m_X,\Delta_X,...)$ 
and $(m_B,\Delta_B,...)$ be $\cS$-bialgebra structures on $X$ and $B$. 
Let $\ell_{AX}\in \cS(A,X)$ and $\ell_{XB}\in \cS(X,B)$ be morphisms 
of $\cS$-bialgebras, such that $\ell_{AX}$ is right invertible and 
$\ell_{XB}$ is left invertible: let $\ell_{XA}\in \cS(X,A)$ and 
$\ell_{BX}\in\cS(X,B)$ be such that $\ell_{AX} \circ \ell_{XA} = 
\on{id}_X$, and $\ell_{BX} \circ \ell_{XB} = \on{id}_X$. Then 
$\ell_{AX} \circ m_A \circ \ell_{XA}^{\otimes 2}
= \ell_{BX} \circ m_B \circ \ell_{XB}^{\otimes 2}$, 
$\ell^{\otimes 2}_{AX} \circ \Delta_A \circ \ell_{XA}
= \ell_{BX}^{\otimes 2} \circ \Delta_B \circ \ell_{XB}$, etc. If we call 
$m_X\in \cS(X^{\otimes 2},X)$ (resp., $\Delta_X\in \cS(X,X^{\otimes 2})$, 
etc.)
the common value of both sides of the first (resp., second) identity, then 
$(m_X,\Delta_X,...)$ is a $\cS$-bialgebra structure on $X$ . 
\hfill \qed \medskip 
 
According to \cite{Enr:coh}, $\on{J}$ is uniquely determined by 
(\ref{eq:J}) only up to a gauge transformation $\on{J} \mapsto 
{}^u\on{J} = 
u^1 u^2 \on{J} (u^{12})^{-1}$, where $u\in 
(\wh{\bf U}_1)^\times_1$. 

\begin{lemma}
Quantization functors corresponding to $\on{J}$ and to ${}^u\on{J}$
are equivalent. 
\end{lemma} 

{\em Proof.} We have $u,u^{-1} \in \wh{{\bf U}}_1 \simeq 
$\boldmath$\Pi$\unboldmath$({\mathfrak 1}\underline\boxtimes {\mathfrak 1},
S\underline\boxtimes S)$. 
Let us set 
$$
\on{Ad}(u) := m_\Pi^{(2)}(u \boxtimes
\on{id}_{S\underline\boxtimes S} 
\boxtimes u^{-1}) \in \hbox{\boldmath$\Pi$\unboldmath}(S\underline\boxtimes
S,S\underline\boxtimes S)^\times
$$
(one checks that the r.h.s. makes sense). 

Let us view $\on{Ad}(u)$ as an element of \boldmath$\Pi$\unboldmath$
(F\underline\boxtimes G,F'\underline\boxtimes G')$, then a graph for 
$\on{Ad}(u)$ and $\on{Ad}(u)^{-1}(= \on{Ad}(u^{-1})$) is 
$F\to F'$, $G'\to G$, $G'\to F'$. 

In the same way, $\on{Ad}(u)^*, (\on{Ad}(u)^{-1})^* \in 
$\boldmath$\Pi$\unboldmath$(S\underline\boxtimes
S,S\underline\boxtimes S)^\times$, and a graph for these elements 
is $F\to F'$, $G'\to G$, $F\to G$. 

Let us denote by ${}^u\on{R},{}^u\ell,...,{}^uQ$ the analogues of 
$\on{R},\ell,...,Q$, with $\on{J}$ replaced by ${}^u\on{J}$. 
These analogues can be expressed as follows: ${}^u m_\Pi = m_\Pi$, 
${}^u\Delta_\Pi = \on{Ad}(u)^{\boxtimes 2} \circ \Delta_\Pi \circ 
\on{Ad}(u)^{-1}$ (one checks that the r.h.s. is well-defined), 
${}^u\on{R} = u^1 u^2 \on{R} (u^1u^2)^{-1}$, therefore 
${}^u\ell = \on{Ad}(u) \circ \ell \circ \on{Ad}(u)^*$ (one checks that 
the r.h.s. is well-defined). We then make the following choices for 
${}^u\on{R}_\pm$: ${}^u\on{R}_\pm = \on{Ad}(u) \circ \on{R}_\pm$ 
(one checks that
both r.h.s. are well-defined). 

We then have ${}^u\on{R}_\pm^{(-1)} = \on{R}_\pm^{(-1)} \circ \on{Ad}(u)^{-1}$. 
Then: 
$$
{}^um_a = ({}^u\on{R}_+^{(-1)})^{\boxtimes 2} \circ {}^um_\Pi \circ 
{}^u\on{R}_+ 
= (\on{R}_+^{(-1)})^{\boxtimes 2} \circ 
(\on{Ad}(u)^{-1})^{\boxtimes 2} \circ m_\Pi \circ \on{Ad}(u) \circ 
\on{R}_+ 
= (\on{R}_+^{(-1)})^{\boxtimes 2} \circ m_\Pi \circ \on{R}_+ = m_a, 
$$
and 
$$
{}^u\Delta_a = ({}^u\on{R}_+^{(-1)})^{\boxtimes 2} 
\circ {}^u\Delta_\Pi \circ {}^u\on{R}_+
= (\on{R}_+^{(-1)})^{\boxtimes 2} \circ (\on{Ad}(u)^{-1})^{\boxtimes 2}
\circ {}^u\Delta_\Pi \circ \on{Ad}(u) \circ \on{R}_+
= (\on{R}_+^{(-1)})^{\boxtimes 2} \circ \Delta_\Pi \circ \on{R}_+
= \Delta_a, 
$$
so ${}^uQ(m) = Q(m)$ and ${}^uQ(\Delta) = Q(\Delta)$, so ${}^uQ = Q$. 
\hfill \qed \medskip 

Here are pictures of the main graphs of the above construction. 
The object $S$ is represented by black vertices, and the  
object ${\mathfrak 1}\in{\bf Sch}_{(1)}$ is represented by 
white vertices.

{\large\bf\boldmath

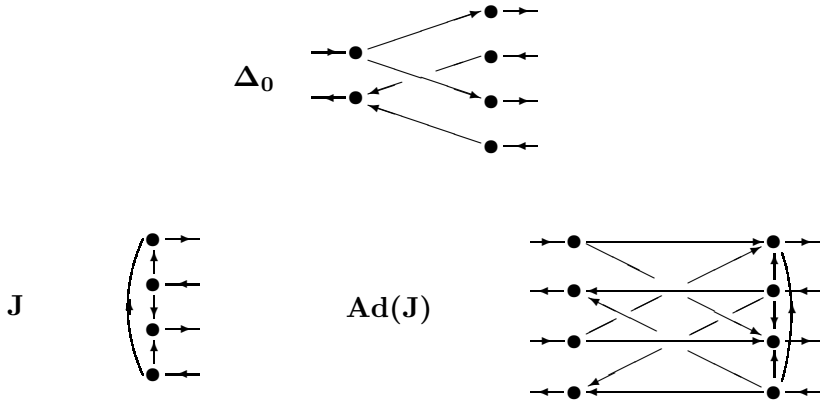
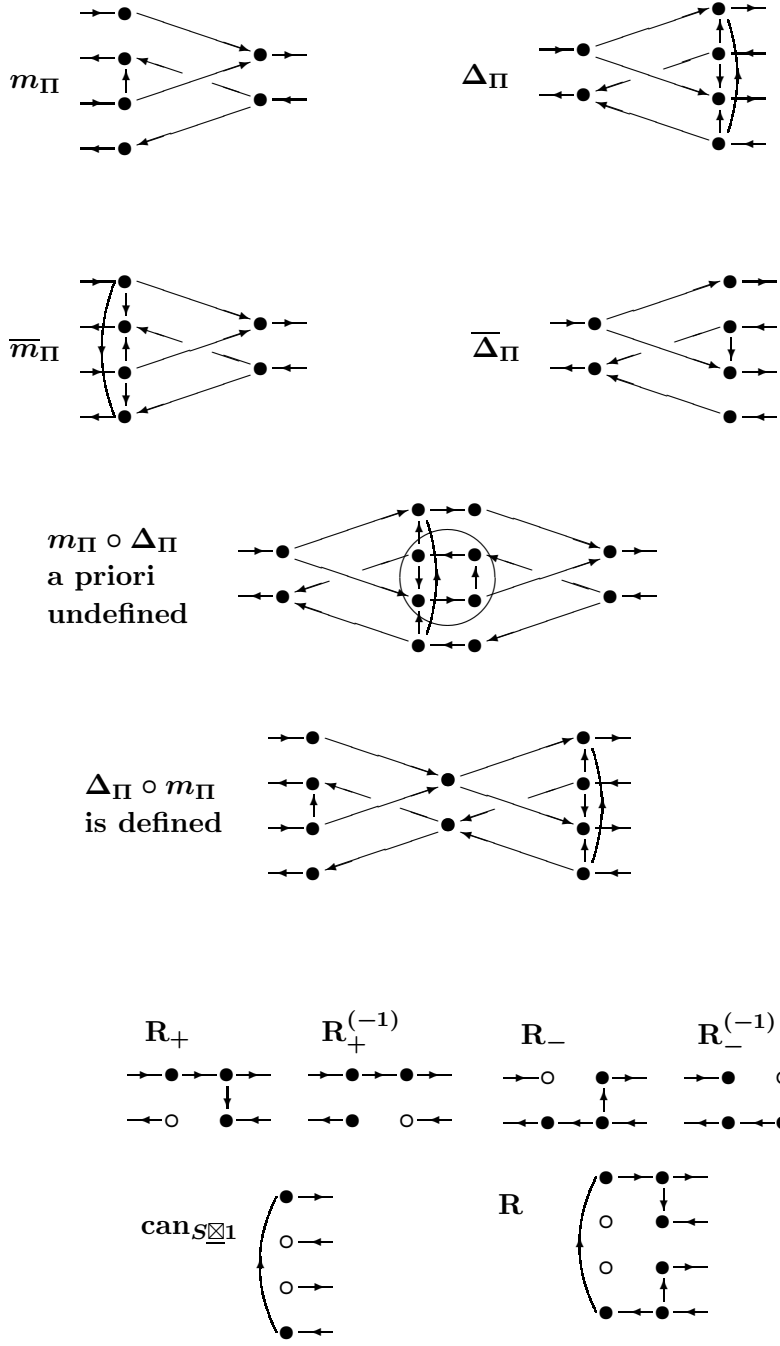
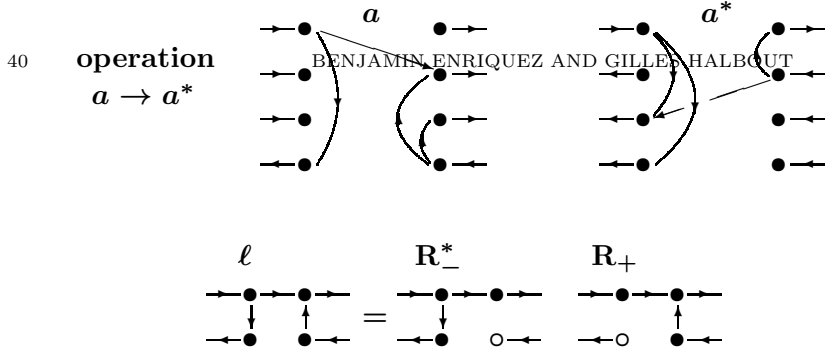
\setlength{\unitlength}{0.5cm}

$\hbox{}$

\begin{picture}(10,5)(-6,0)
\put(-1,2){$\Delta_0$}
\put(7.1,2.8){\vector(-1,0){.6}}
\put(6.5,2.8){\line(-1,0){0.3}}
\put(5.65,3.8){$\bullet$}

\put(6.2,4){\vector(1,0){0.7}}
\put(6.9,4){\line(1,0){0.2}}
\put(5.65,2.6){$\bullet$}

\put(7.1,.4){\vector(-1,0){.6}}
\put(6.5,.4){\line(-1,0){0.3}}
\put(5.65,1.4){$\bullet$}

\put(6.2,1.6){\vector(1,0){0.7}}
\put(6.9,1.6){\line(1,0){0.2}}
\put(5.65,.2){$\bullet$}


\put(5.6,2.8){\line(-3,-1){1.2}}
\put(3.5,2.1){\vector(-3,-1){0.9}}
\put(2.05,2.7){$\bullet$}
\put(2,1.7){\vector(-1,0){.6}}
\put(1.4,1.7){\line(-1,0){0.3}}

\put(5.6,.5){\vector(-3,1){3}}

\put(2.6,3){\vector(3,1){3}}
\put(2.05,1.5){$\bullet$}
\put(1.1,2.9){\vector(1,0){0.7}}
\put(1.8,2.9){\line(1,0){0.2}}

\put(2.6,2.7){\vector(3,-1){3}}

\end{picture}

\vskip0cm

\begin{picture}(10,5)(0,1)
\put(-1,2){$\operatorname{J}$}
\put(4.1,2.8){\vector(-1,0){.6}}
\put(3.5,2.8){\line(-1,0){0.3}}
\put(2.65,3.8){$\bullet$}

\put(3.2,4){\vector(1,0){0.7}}
\put(3.9,4){\line(1,0){0.2}}
\put(2.65,2.6){$\bullet$}

\put(4.1,.4){\vector(-1,0){.6}}
\put(3.5,.4){\line(-1,0){0.3}}
\put(2.65,1.4){$\bullet$}

\put(3.2,1.6){\vector(1,0){0.7}}
\put(3.9,1.6){\line(1,0){0.2}}
\put(2.65,.2){$\bullet$}

\put(2.9,.7){\vector(0,1){0.6}}
\put(2.9,3.1){\vector(0,1){0.6}}
\put(2.9,2.5){\vector(0,-1){0.6}}

\qbezier(2.6,4)(1.8,2.2)(2.6,.4)
\put(2.25,2){\vector(0,1){0.4}}
\end{picture}

\begin{picture}(10,5)(-12,-4)
\put(-4,2){Ad$(\operatorname{J})$}
\put(0.9,4){\vector(1,0){.6}}
\put(1.5,4){\line(1,0){0.3}}
\put(1.85,3.8){$\bullet$}

\put(1.8,2.7){\vector(-1,0){0.7}}
\put(1.1,2.7){\line(-1,0){0.2}}
\put(1.85,2.5){$\bullet$}

\put(0.9,1.35){\vector(1,0){.6}}
\put(1.5,1.35){\line(1,0){0.3}}
\put(1.85,1.15){$\bullet$}

\put(1.8,0){\vector(-1,0){0.7}}
\put(1.1,0){\line(-1,0){0.2}}
\put(1.85,-0.2){$\bullet$}

\put(2.4,4){\vector(1,0){4.7}}
\put(7.1,2.7){\vector(-1,0){4.7}}
\put(2.4,1.35){\vector(1,0){4.7}}
\put(7.1,0){\vector(-1,0){4.7}}

\put(2.4,3.95){\line(2,-1){1.8}}
\put(5.1,2.5){\vector(2,-1){1.95}}

\put(4.26,1.62){\vector(-2,1){1.8}}
\put(5.1,1.1){\line(2,-1){1.95}}

\put(2.4,1.55){\line(2,1){0.8}}
\put(3.6,2.15){\line(2,1){0.8}}
\put(5.1,2.9){\vector(2,1){1.95}}

\put(7.1,2.5){\line(-2,-1){0.8}}
\put(5.9,1.9){\line(-2,-1){0.75}}
\put(4.4,1.15){\vector(-2,-1){1.95}}

\put(8.6,2.7){\vector(-1,0){.6}}
\put(8,2.7){\line(-1,0){0.3}}
\put(7.15,3.8){$\bullet$}

\put(7.7,4){\vector(1,0){0.7}}
\put(8.4,4){\line(1,0){0.2}}
\put(7.15,2.5){$\bullet$}

\put(8.6,0){\vector(-1,0){.6}}
\put(8,0){\line(-1,0){0.3}}
\put(7.15,1.15){$\bullet$}

\put(7.7,1.35){\vector(1,0){0.7}}
\put(8.4,1.35){\line(1,0){0.2}}
\put(7.15,-0.2){$\bullet$}

\put(7.4,.4){\vector(0,1){0.7}}
\put(7.4,3){\vector(0,1){0.7}}
\put(7.4,2.4){\vector(0,-1){0.7}}

\qbezier(7.6,3.7)(8.1,2.2)(7.6,.3)
\put(7.87,2){\vector(0,1){0.4}}

\end{picture}


\begin{picture}(10,5)(0,1)
\put(-1,2){$m_\Pi$}
\put(0.9,4){\vector(1,0){.6}}
\put(1.5,4){\line(1,0){0.3}}
\put(1.85,3.8){$\bullet$}

\put(1.8,2.8){\vector(-1,0){0.7}}
\put(1.1,2.8){\line(-1,0){0.2}}
\put(1.85,2.6){$\bullet$}

\put(0.9,1.6){\vector(1,0){.6}}
\put(1.5,1.6){\line(1,0){0.3}}
\put(1.85,1.4){$\bullet$}

\put(1.8,.4){\vector(-1,0){0.7}}
\put(1.1,.4){\line(-1,0){0.2}}
\put(1.85,.2){$\bullet$}

\put(2.1,1.9){\vector(0,1){0.6}}

\put(2.4,4){\vector(3,-1){3}}
\put(5.45,2.7){$\bullet$}
\put(6,2.9){\vector(1,0){.6}}
\put(6.6,2.9){\line(1,0){0.3}}

\put(2.4,1.7){\vector(3,1){3}}

\put(3.6,2.4){\vector(-3,1){1.2}}
\put(5.4,1.8){\line(-3,1){1}}
\put(5.45,1.5){$\bullet$}
\put(6.9,1.7){\vector(-1,0){0.7}}
\put(6.2,1.7){\line(-1,0){0.2}}

\put(5.4,1.5){\vector(-3,-1){3}}

\end{picture}

\begin{picture}(10,5)(-12,-4.2)
\put(-1,2){$\Delta_\Pi$}


\put(5.6,2.8){\line(-3,-1){1.2}}
\put(3.5,2.1){\vector(-3,-1){0.9}}
\put(2.05,2.7){$\bullet$}
\put(2,1.7){\vector(-1,0){.6}}
\put(1.4,1.7){\line(-1,0){0.3}}

\put(5.6,.5){\vector(-3,1){3}}

\put(2.6,3){\vector(3,1){3}}
\put(2.05,1.5){$\bullet$}
\put(1.1,2.9){\vector(1,0){0.7}}
\put(1.8,2.9){\line(1,0){0.2}}

\put(2.6,2.7){\vector(3,-1){3}}

\put(7.1,2.8){\vector(-1,0){.6}}
\put(6.5,2.8){\line(-1,0){0.3}}
\put(5.65,2.6){$\bullet$}

\put(6.2,4){\vector(1,0){0.7}}
\put(6.9,4){\line(1,0){0.2}}
\put(5.65,3.8){$\bullet$}

\put(7.1,.4){\vector(-1,0){.6}}
\put(6.5,.4){\line(-1,0){0.3}}
\put(5.65,.2){$\bullet$}

\put(6.2,1.6){\vector(1,0){0.7}}
\put(6.9,1.6){\line(1,0){0.2}}
\put(5.65,1.4){$\bullet$}

\put(5.9,.7){\vector(0,1){0.6}}
\put(5.9,3.1){\vector(0,1){0.6}}
\put(5.9,2.5){\vector(0,-1){0.6}}

\qbezier(6.1,3.7)(6.6,2.2)(6.1,.7)
\put(6.37,2){\vector(0,1){0.4}}

\end{picture}

\vskip-1cm

\begin{picture}(10,4)(0,-0)
\put(-1,2){$\overline{m}_\Pi$}

\put(0.9,4){\vector(1,0){.6}}
\put(1.5,4){\line(1,0){0.3}}
\put(1.85,3.8){$\bullet$}

\put(1.8,2.8){\vector(-1,0){0.7}}
\put(1.1,2.8){\line(-1,0){0.2}}
\put(1.85,2.6){$\bullet$}

\put(0.9,1.6){\vector(1,0){.6}}
\put(1.5,1.6){\line(1,0){0.3}}
\put(1.85,1.4){$\bullet$}

\put(1.8,.4){\vector(-1,0){0.7}}
\put(1.1,.4){\line(-1,0){0.2}}
\put(1.85,.2){$\bullet$}

\put(2.1,1.9){\vector(0,1){0.6}}
\put(2.1,1.3){\vector(0,-1){0.6}}
\put(2.1,3.7){\vector(0,-1){0.6}}

\qbezier(1.8,4)(1.1,2.2)(1.8,.4)
\put(1.46,2.4){\vector(0,-1){0.4}}

\put(2.4,4){\vector(3,-1){3}}
\put(5.45,2.7){$\bullet$}
\put(6,2.9){\vector(1,0){.6}}
\put(6.6,2.9){\line(1,0){0.3}}

\put(2.4,1.7){\vector(3,1){3}}

\put(3.6,2.4){\vector(-3,1){1.2}}
\put(5.4,1.8){\line(-3,1){1}}
\put(5.45,1.5){$\bullet$}
\put(6.9,1.7){\vector(-1,0){0.7}}
\put(6.2,1.7){\line(-1,0){0.2}}

\put(5.4,1.5){\vector(-3,-1){3}}

\end{picture}
\nopagebreak
\begin{picture}(10,5)(-2,0)

\put(-1,2){$\overline{\Delta}_\Pi$}
\put(7.1,2.8){\vector(-1,0){.6}}
\put(6.5,2.8){\line(-1,0){0.3}}
\put(5.65,3.8){$\bullet$}

\put(6.2,4){\vector(1,0){0.7}}
\put(6.9,4){\line(1,0){0.2}}
\put(5.65,2.6){$\bullet$}
\put(5.9,2.5){\vector(0,-1){0.6}}

\put(7.1,.4){\vector(-1,0){.6}}
\put(6.5,.4){\line(-1,0){0.3}}
\put(5.65,1.4){$\bullet$}

\put(6.2,1.6){\vector(1,0){0.7}}
\put(6.9,1.6){\line(1,0){0.2}}
\put(5.65,.2){$\bullet$}


\put(5.6,2.8){\line(-3,-1){1.2}}
\put(3.5,2.1){\vector(-3,-1){0.9}}
\put(2.05,2.7){$\bullet$}
\put(2,1.7){\vector(-1,0){.6}}
\put(1.4,1.7){\line(-1,0){0.3}}

\put(5.6,.5){\vector(-3,1){3}}

\put(2.6,3){\vector(3,1){3}}
\put(2.05,1.5){$\bullet$}
\put(1.1,2.9){\vector(1,0){0.7}}
\put(1.8,2.9){\line(1,0){0.2}}

\put(2.6,2.7){\vector(3,-1){3}}

\end{picture}

\begin{picture}(10,5)(-4,1)
\put(-4,3){$m_\Pi \circ \Delta_\Pi$}
\put(-4,2){a priori}
\put(-4,1){undefined}

\put(5.6,2.8){\line(-3,-1){1.2}}
\put(3.5,2.1){\vector(-3,-1){0.9}}
\put(2.05,2.7){$\bullet$}
\put(2,1.7){\vector(-1,0){.6}}
\put(1.4,1.7){\line(-1,0){0.3}}

\put(5.6,.5){\vector(-3,1){3}}

\put(2.6,3){\vector(3,1){3}}
\put(2.05,1.5){$\bullet$}
\put(1.1,2.9){\vector(1,0){0.7}}
\put(1.8,2.9){\line(1,0){0.2}}

\put(2.6,2.7){\vector(3,-1){3}}

\put(7.1,2.8){\vector(-1,0){.6}}
\put(6.5,2.8){\line(-1,0){0.3}}
\put(5.65,2.6){$\bullet$}

\put(6.2,4){\vector(1,0){0.7}}
\put(6.9,4){\line(1,0){0.2}}
\put(5.65,3.8){$\bullet$}

\put(7.1,.4){\vector(-1,0){.6}}
\put(6.5,.4){\line(-1,0){0.3}}
\put(5.65,.2){$\bullet$}

\put(6.2,1.6){\vector(1,0){0.7}}
\put(6.9,1.6){\line(1,0){0.2}}
\put(5.65,1.4){$\bullet$}

\put(5.9,.7){\vector(0,1){0.6}}
\put(5.9,3.1){\vector(0,1){0.6}}
\put(5.9,2.5){\vector(0,-1){0.6}}

\qbezier(6.1,3.7)(6.6,2.2)(6.1,.7)
\put(6.37,2){\vector(0,1){0.4}}

\put(6.65,2.2){\circle{2.5}}

\put(7.15,3.8){$\bullet$}

\put(7.15,2.6){$\bullet$}

\put(7.15,1.4){$\bullet$}

\put(7.15,.2){$\bullet$}

\put(7.4,1.9){\vector(0,1){0.6}}

\put(7.7,4){\vector(3,-1){3}}
\put(10.75,2.7){$\bullet$}
\put(11.3,2.9){\vector(1,0){.6}}
\put(11.9,2.9){\line(1,0){0.3}}

\put(7.7,1.7){\vector(3,1){3}}

\put(8.9,2.4){\vector(-3,1){1.2}}
\put(10.7,1.8){\line(-3,1){1}}
\put(10.75,1.5){$\bullet$}
\put(12.2,1.7){\vector(-1,0){0.7}}
\put(11.5,1.7){\line(-1,0){0.2}}

\put(10.7,1.5){\vector(-3,-1){3}}

\end{picture}

\vskip1cm

\begin{picture}(10,5)(-5,0)
\put(-4,2.5){$\Delta_\Pi\circ m_\Pi $}
\put(-4,1.5){is defined}
\put(0.9,4){\vector(1,0){.6}}
\put(1.5,4){\line(1,0){0.3}}
\put(1.85,3.8){$\bullet$}

\put(1.8,2.8){\vector(-1,0){0.7}}
\put(1.1,2.8){\line(-1,0){0.2}}
\put(1.85,2.6){$\bullet$}

\put(0.9,1.6){\vector(1,0){.6}}
\put(1.5,1.6){\line(1,0){0.3}}
\put(1.85,1.4){$\bullet$}

\put(1.8,.4){\vector(-1,0){0.7}}
\put(1.1,.4){\line(-1,0){0.2}}
\put(1.85,.2){$\bullet$}

\put(2.1,1.9){\vector(0,1){0.6}}

\put(2.4,4){\vector(3,-1){3}}
\put(5.45,2.7){$\bullet$}

\put(2.4,1.7){\vector(3,1){3}}

\put(3.6,2.4){\vector(-3,1){1.2}}
\put(5.4,1.8){\line(-3,1){1}}
\put(5.45,1.5){$\bullet$}

\put(5.4,1.5){\vector(-3,-1){3}}

\put(9,2.8){\line(-3,-1){1.2}}
\put(6.9,2.1){\vector(-3,-1){0.9}}
\put(5.45,2.7){$\bullet$}

\put(9,.5){\vector(-3,1){3}}

\put(6,3){\vector(3,1){3}}
\put(5.45,1.5){$\bullet$}

\put(6,2.7){\vector(3,-1){3}}

\put(10.5,2.8){\vector(-1,0){.6}}
\put(9.9,2.8){\line(-1,0){0.3}}
\put(9.05,2.6){$\bullet$}

\put(9.6,4){\vector(1,0){0.7}}
\put(10.3,4){\line(1,0){0.2}}
\put(9.05,3.8){$\bullet$}

\put(10.5,.4){\vector(-1,0){.6}}
\put(9.9,.4){\line(-1,0){0.3}}
\put(9.05,.2){$\bullet$}

\put(9.6,1.6){\vector(1,0){0.7}}
\put(10.3,1.6){\line(1,0){0.2}}
\put(9.05,1.4){$\bullet$}

\put(9.3,.7){\vector(0,1){0.6}}
\put(9.3,3.1){\vector(0,1){0.6}}
\put(9.3,2.5){\vector(0,-1){0.6}}

\qbezier(9.5,3.7)(10,2.2)(9.5,.7)
\put(9.77,2){\vector(0,1){0.4}}

\end{picture}

\vskip1.5cm 

\begin{picture}(10,3)(0,2)
\put(2.6,4){$\operatorname{R}_+$}
\put(3.1,2.9){$\bullet$}
\put(3.1,1.7){$\circ$}
\put(3.6,3.1){\vector(1,0){.6}}
\put(4.2,3.1){\line(1,0){0.3}}
\put(2.15,3.1){\vector(1,0){.6}}
\put(2.75,3.1){\line(1,0){0.3}}
\put(3.05,1.9){\vector(-1,0){.6}}
\put(2.45,1.9){\line(-1,0){0.3}}
\put(5.05,3.1){\vector(1,0){.6}}
\put(5.65,3.1){\line(1,0){0.3}}
\put(5.95,1.9){\vector(-1,0){.6}}
\put(5.35,1.9){\line(-1,0){0.3}}
\put(4.55,2.9){$\bullet$}
\put(4.55,1.7){$\bullet$}
\put(4.8,2.8){\vector(0,-1){0.6}}

\put(7.3,4){$\operatorname{R}_+^{(-1)}$}
\put(7.9,2.9){$\bullet$}
\put(7.9,1.7){$\bullet$}
\put(6.95,3.1){\vector(1,0){.6}}
\put(7.55,3.1){\line(1,0){0.3}}
\put(7.85,1.9){\vector(-1,0){.6}}
\put(7.25,1.9){\line(-1,0){0.3}}
\put(9.85,3.1){\vector(1,0){.6}}
\put(10.45,3.1){\line(1,0){0.3}}
\put(10.75,1.9){\vector(-1,0){.6}}
\put(10.15,1.9){\line(-1,0){0.3}}

\put(8.4,3.1){\vector(1,0){.6}}
\put(9,3.1){\line(1,0){0.3}}
\put(9.35,2.9){$\bullet$}
\put(9.35,1.7){$\circ$}
\end{picture}

\begin{picture}(10,3)(-10,-1)
\put(2.6,4){$\operatorname{R}_-$}
\put(3.1,2.9){$\circ$}
\put(3.1,1.7){$\bullet$}
\put(2.15,3.1){\vector(1,0){.6}}
\put(2.75,3.1){\line(1,0){0.3}}
\put(3.05,1.9){\vector(-1,0){.6}}
\put(2.45,1.9){\line(-1,0){0.3}}
\put(5.05,3.1){\vector(1,0){.6}}
\put(5.65,3.1){\line(1,0){0.3}}
\put(5.95,1.9){\vector(-1,0){.6}}
\put(5.35,1.9){\line(-1,0){0.3}}
\put(4.5,1.9){\vector(-1,0){.6}}
\put(3.9,1.9){\line(-1,0){0.3}}
\put(4.55,2.9){$\bullet$}
\put(4.55,1.7){$\bullet$}
\put(4.8,2.2){\vector(0,1){0.6}}

\put(7.3,4){$\operatorname{R}_-^{(-1)}$}
\put(7.9,2.9){$\bullet$}
\put(7.9,1.7){$\bullet$}
\put(6.95,3.1){\vector(1,0){.6}}
\put(7.55,3.1){\line(1,0){0.3}}
\put(7.85,1.9){\vector(-1,0){.6}}
\put(7.25,1.9){\line(-1,0){0.3}}
\put(9.85,3.1){\vector(1,0){.6}}
\put(10.45,3.1){\line(1,0){0.3}}
\put(10.75,1.9){\vector(-1,0){.6}}
\put(10.15,1.9){\line(-1,0){0.3}}
\put(9.35,1.9){\vector(-1,0){.6}}
\put(8.75,1.9){\line(-1,0){0.3}}

\put(9.35,2.9){$\circ$}
\put(9.35,1.7){$\bullet$}

\end{picture}

 \vskip -1cm

\begin{picture}(10,5)(-2.5,0)

\put(0,3){$\operatorname{can}_{S \underline{\boxtimes} 1}$}
\qbezier(3.6,0.4)(2.7,2.2)(3.6,4)
\put(3.17,2.1){\vector(0,1){.3}}
\put(5.1,2.8){\vector(-1,0){.6}}
\put(4.5,2.8){\line(-1,0){0.3}}
\put(3.65,3.8){$\bullet$}

\put(4.2,4){\vector(1,0){0.7}}
\put(4.9,4){\line(1,0){0.2}}
\put(3.65,2.6){$\circ$}

\put(5.1,.4){\vector(-1,0){.6}}
\put(4.5,.4){\line(-1,0){0.3}}
\put(3.65,1.4){$\circ$}

\put(4.2,1.6){\vector(1,0){0.7}}
\put(4.9,1.6){\line(1,0){0.2}}
\put(3.65,.2){$\bullet$}
\end{picture}

 \vskip -0.3cm
\begin{picture}(10,5)(-11,-5)
\put(1,3){$\operatorname{R}$}
\qbezier(3.6,0.4)(2.7,2.2)(3.6,4)
\put(3.17,2.1){\vector(0,1){.3}}
\put(3.65,3.8){$\bullet$}

\put(4.2,4){\vector(1,0){0.7}}
\put(4.9,4){\line(1,0){0.2}}
\put(3.65,2.6){$\circ$}

\put(5.1,.4){\vector(-1,0){.6}}
\put(4.5,.4){\line(-1,0){0.3}}
\put(3.65,1.4){$\circ$}

\put(3.65,.2){$\bullet$}

\put(5.65,4){\vector(1,0){.6}}
\put(6.25,4){\line(1,0){0.3}}
\put(6.55,2.8){\vector(-1,0){.6}}
\put(5.95,2.8){\line(-1,0){0.3}}
\put(5.15,3.8){$\bullet$}
\put(5.15,2.6){$\bullet$}
\put(5.4,3.7){\vector(0,-1){0.6}}

\put(5.65,1.6){\vector(1,0){.6}}
\put(6.25,1.6){\line(1,0){0.3}}
\put(6.55,.4){\vector(-1,0){.6}}
\put(5.95,.4){\line(-1,0){0.3}}
\put(5.15,1.4){$\bullet$}
\put(5.15,.2){$\bullet$}
\put(5.4,.7){\vector(0,1){0.6}}

\end{picture}

\vskip-1cm
\begin{picture}(10,5)(-5,-3)
\put(-4,3){operation}
\put(-3.6,2){$a \to a^*$}

\put(0.9,4){\vector(1,0){.6}}
\put(1.5,4){\line(1,0){0.3}}
\put(1.85,3.8){$\bullet$}

\put(0.9,2.8){\vector(1,0){.6}}
\put(1.5,2.8){\line(1,0){0.3}}
\put(1.85,2.6){$\bullet$}

\put(0.9,1.6){\vector(1,0){.6}}
\put(1.5,1.6){\line(1,0){0.3}}
\put(1.85,1.4){$\bullet$}

\put(1.8,.4){\vector(-1,0){0.7}}
\put(1.1,.4){\line(-1,0){0.2}}
\put(1.85,.2){$\bullet$}

\put(6,4){\vector(1,0){0.7}}
\put(6.7,4){\line(1,0){0.2}}
\put(5.45,3.8){$\bullet$}

\put(6,2.8){\vector(1,0){0.7}}
\put(6.7,2.8){\line(1,0){0.2}}
\put(5.45,2.6){$\bullet$}

\put(6.9,.4){\vector(-1,0){.6}}
\put(6.3,.4){\line(-1,0){0.3}}
\put(5.45,.2){$\bullet$}

\put(6,1.6){\vector(1,0){0.7}}
\put(6.7,1.6){\line(1,0){0.2}}
\put(5.45,1.4){$\bullet$}

\put(9.9,4){\vector(1,0){.6}}
\put(10.5,4){\line(1,0){0.3}}
\put(10.85,3.8){$\bullet$}

\put(10.8,2.8){\vector(-1,0){.7}}
\put(10.1,2.8){\line(-1,0){0.2}}
\put(10.85,2.6){$\bullet$}

\put(10.8,1.6){\vector(-1,0){.7}}
\put(10.1,1.6){\line(-1,0){0.2}}
\put(10.85,1.4){$\bullet$}

\put(10.8,.4){\vector(-1,0){0.7}}
\put(10.1,.4){\line(-1,0){0.2}}
\put(10.85,.2){$\bullet$}

\put(15,4){\vector(1,0){0.7}}
\put(15.7,4){\line(1,0){0.2}}
\put(14.45,3.8){$\bullet$}

\put(15.9,2.8){\vector(-1,0){0.6}}
\put(15.3,2.8){\line(-1,0){0.3}}
\put(14.45,2.6){$\bullet$}

\put(15.9,.4){\vector(-1,0){.6}}
\put(15.3,.4){\line(-1,0){0.3}}
\put(14.45,.2){$\bullet$}

\put(15.9,1.6){\vector(-1,0){0.6}}
\put(15.3,1.6){\line(-1,0){0.3}}
\put(14.45,1.4){$\bullet$}

\qbezier(2.4,.4)(3.5,2.1)(2.4,3.9)
\put(2.95,2.28){\vector(0,-1){0.4}}

\qbezier(11.35,1.7)(12.4,2.9)(11.35,3.85)
\qbezier(11.4,.4)(13.5,2.1)(11.4,3.9)

\put(11.87,3){\vector(0,-1){0.4}}
\put(12.45,2.15){\vector(0,-1){0.4}}

\put(14.45,2.65){\line(-3,-1){1.6}}
\put(12.15,1.9){\vector(-3,-1){.8}}

\qbezier(14.4,2.75)(13.7,3.28)(14.4,3.9)
\put(14.1,3.18){\vector(0,1){0.3}}

\put(2.5,3.95){\vector(3,-1){3}}

\qbezier(5.35,.4)(3.7,1.6)(5.35,2.75)

\qbezier(5.4,.45)(4.9,1)(5.4,1.55)

\put(4.57,1.45){\vector(0,1){0.3}}
\put(5.20,.88){\vector(0,1){0.3}}

\put(3.6,4.2){$a$}
\put(12.6,4.2){$a^*$}

\end{picture}

\vskip-1cm

\begin{picture}(10,5.5)(-5,0)
\put(0.3,4.3){$\ell$}
\put(.4,3.3){$\bullet$}
\put(.4,2.1){$\bullet$}
\put(0.9,3.5){\vector(1,0){.6}}
\put(1.5,3.5){\line(1,0){0.3}}
\put(-.55,3.5){\vector(1,0){.6}}
\put(0.05,3.5){\line(1,0){0.3}}
\put(.35,2.3){\vector(-1,0){.6}}
\put(-0.25,2.3){\line(-1,0){0.3}}
\put(2.35,3.5){\vector(1,0){.6}}
\put(2.95,3.5){\line(1,0){0.3}}
\put(3.25,2.3){\vector(-1,0){.6}}
\put(2.65,2.3){\line(-1,0){0.3}}
\put(1.85,3.3){$\bullet$}
\put(1.85,2.1){$\bullet$}
\put(0.65,3.2){\vector(0,-1){0.6}}
\put(2.1,2.6){\vector(0,1){0.6}}

\put(3.6,2.7){$=$}

\put(5,4.3){$\operatorname{R}_{-}^*$}
\put(5.5,3.3){$\bullet$}
\put(5.5,2.1){$\bullet$}
\put(6,3.5){\vector(1,0){.6}}
\put(6.6,3.5){\line(1,0){0.3}}
\put(4.55,3.5){\vector(1,0){.6}}
\put(5.15,3.5){\line(1,0){0.3}}
\put(5.45,2.3){\vector(-1,0){.6}}
\put(4.85,2.3){\line(-1,0){0.3}}
\put(7.45,3.5){\vector(1,0){.6}}
\put(8.05,3.5){\line(1,0){0.3}}
\put(8.35,2.3){\vector(-1,0){.6}}
\put(7.75,2.3){\line(-1,0){0.3}}
\put(6.95,3.3){$\bullet$}
\put(6.95,2.1){$\circ$}
\put(5.75,3.2){\vector(0,-1){0.6}}

\put(9.7,4.3){$\operatorname{R}_+$}
\put(10.3,3.3){$\bullet$}
\put(10.3,2.1){$\circ$}
\put(9.35,3.5){\vector(1,0){.6}}
\put(9.95,3.5){\line(1,0){0.3}}
\put(10.25,2.3){\vector(-1,0){.6}}
\put(9.65,2.3){\line(-1,0){0.3}}
\put(12.25,3.5){\vector(1,0){.6}}
\put(12.85,3.5){\line(1,0){0.3}}
\put(13.15,2.3){\vector(-1,0){.6}}
\put(12.55,2.3){\line(-1,0){0.3}}

\put(10.8,3.5){\vector(1,0){.6}}
\put(11.4,3.5){\line(1,0){0.3}}
\put(11.75,3.3){$\bullet$}
\put(11.75,2.1){$\bullet$}
\put(12,2.6){\vector(0,1){0.6}}

\end{picture}
}


\section{Compatibility of quantization functors with twists}
\label{sect:twists}
    
\subsection{The category $\cY$}

We define $\cY$  as the category where objects are integer numbers
$\geq 0$, and $\cY(n,m)$ is the set of pairs $(\phi,o)$, where 
$\phi : [m] \to [n]$ is a partially defined function and 
$o = (o_1,...,o_n)$, where $o_i$ is a total order on $\phi^{-1}(i)$. 
If $(\phi,o)\in\cY(n,m)$ and $(\phi',o')\in\cY(m,p)$, then their composition 
is $(\phi'',o'')\in \cY(n,p)$, where $\phi'' = \phi \circ \phi'$ and 
$o'' = (o''_1,...,o''_n)$, where $o''_i$ is te lexicographic order on
$(\phi'')^{-1}(i) = \sqcup_{j\in \phi^{-1}(i)} (\phi')^{-1}(j)$. 

A $\cY$-vector space, (resp., a $\cY$-algebra) is a functor $\cY\to \on{Vect}$
(resp., $\cY\to \on{Alg}$). The forgetful morphism $\cY \to \cX$ gives rise 
to functors $\{\cX$-vector spaces$\} \to \{\cY$-vector spaces$\}$ 
and $\{\cX$-algebras$\} \to \{\cY$-algebras$\}$. A $\cY$-vector 
space is therefore a collection of vector spaces $(V_n)_{n\geq 0}$
and of maps $V_n \to V_m$, $x\mapsto x^{\phi,o}$ for $(\phi,o)\in 
\cY(n,m)$. 

If $H$ is a (non-necessarily cocommutative) coalgebra (resp., bialgebra), 
then $(H^{\otimes n})_{n\geq 0}$ is a $\cY$-vector space 
(resp., $\cY$-algebra).

\subsection{$\cY$-algebra structures on $\wh{\bf U}_n$, 
$\wh{\bf U}_{n,f}$ associated with 
$\on{J}$} \label{new:C}

A solution $\on{J}$ of (\ref{eq:J}) gives rise to a $\cY$-algebra 
structures on $(\wh{\bf U}_n)_{n\geq 0}$, $(\wh{\bf U}_{n,f})_{n\geq 0}$, 
which we now define (we will call them the $\on{J}$-twisted structures). 

For $(\phi,o)\in \cY(n,m)$, define $\on{J}^{\phi,o}\in 
\wh{\bf U}_n^\times$
as follows. For $\psi : [k] \to [m]$ an injective map, we set 
$$
\on{J}_\psi = 
\on{J}^{\psi(1),\psi(2)}
...
\on{J}^{\psi(1)...\psi(k-2),\psi(k-1)}
\on{J}^{\psi(1)...\psi(k-1),\psi(k)} 
$$ 
and 
$$
\on{J}^{\phi,o} = \on{J}_{\psi_1} ... \on{J}_{\psi_n}, 
$$
where $\psi_i : [|\phi^{-1}(i)|] \to \phi^{-1}$ is the unique 
order-preserving bijection. 

The $\cY$-vector space structure on $({\bf U}_n)_{n\geq 0}$ is then 
defined by $x\mapsto (x)_{\on{J}}^{\phi,o} := (x)_{\on{J}}^{\phi,o} 
:= \on{J}^{\phi,o} x^\phi (\on{J}^{\phi,o})^{-1}$; 
the algebra structure is unchanged. 

In the case of $\wh{\bf U}_{n,f}$, the $\cY$-algebra structure is defined
by $x\mapsto (x)_{\on{J}}^{\phi,o} 
:= \kappa_1^\Pi(\on{J}^{\phi,o}) x^\phi \kappa_1^\Pi(\on{J}^{\phi,o})^{-1}$. 

Both $(\wh{\bf U}_n)_{n\geq 0}$ and $(\wh{\bf U}_{n,f})_{n\geq 0}$
are $\cY$-algebras, equipped with decreasing $\cY$-algebra filtrations 
(where the $N$th step consists of the elements of degree $\geq N$). 

\subsection{$\cY$-algebra structure on ${\bf P}({\bf 1},S^{\otimes n})$}
\label{C:algebras} \label{7:3}

Let ${\bf P}$ be a topological prop and let 
$\overline Q : \on{Bialg} \to S({\bf P})$ be a prop morphism. 
Recall that if $H$ is a coalgebra, then $(H^{\otimes n})_{n\geq 0}$
is a $\cY$-vector space. Let us denote by $\Delta^{\phi,o}_H : H^{\otimes n}
\to H^{\otimes m}$ the map corresponding to $(\phi,o)\in \cY(n,m)$. 
The propic versions of the maps $\Delta^{\phi,o}_H$ are elements 
$\Delta^{\phi,o} \in \on{Coalg}(T_n,T_m)$, where $\on{Coalg}$ is the prop 
of algebras (with generators $\Delta,\eta$ with the same relations as in
$\on{Bialg}$). We also denote by $\Delta^{\phi,o}\in \on{Bialg}(T_n,T_m)$ 
the images of these elements
under the prop morphism $\on{Coalg} \to \on{Bialg}$. 

Then $({\bf P}({\bf 1},S^{\otimes n}))_{n\geq 0}$ 
is a $\cY$-vector 
space: the map ${\bf P}({\bf 1},S^{\otimes n}) \to 
{\bf P}({\bf 1},S^{\otimes m})$
corresponding to $(\phi,o)\in \cD(n,m)$ is $x\mapsto 
(x)_{\overline Q}^{\phi,o} := \overline Q(\Delta^{\phi,o}) \circ x$. 

Each ${\bf P}({\bf 1},S^{\otimes n})$ 
is equipped with the algebra structure
$$
x\otimes y \mapsto x \ast_{\overline Q} y :=  
\overline Q(m)^{\otimes n} \circ 
(1,n+1,2,n+2,...) \circ (x\otimes y). 
$$   
The unit for this algebra is $\overline Q(\eta^{\otimes n})$. 

Then this family of algebra structures is compatible with the 
$\cY$-structure, so $({\bf P}({\bf 1},S^{\otimes n}))_{n\geq 0}$ is 
a $\cY$-algebra. 

In particular, the morphism $\kappa_1^\Pi \circ Q : \on{Bialg}
\to S({\bf LBA}_f)$ gives rise to a $\cY$-algebra structure on 
${\bf LBA}_f({\bf 1},S^{\otimes n})$. Using the identification 
${\bf LBA}_f({\bf 1},S^{\otimes n}) \simeq \hbox{\boldmath$\Pi$\unboldmath}_f
({\mathfrak 1}\underline\boxtimes {\mathfrak 1},
(S \underline\boxtimes {\mathfrak 1})^{\otimes n})$, the algebra structure
is given by 
$$
x\ast_{Q} y := \kappa_1^\Pi(m_a)^{\boxtimes n} \circ (1,n+1,2,n+2,...)
\circ (x\boxtimes y), 
$$
and the $\cY$-vector space structure by 
$$
(x)^{\phi,o}_Q := \kappa_1^\Pi(\Delta_a^{\phi,o}) \circ x. 
$$

\subsection{A $\cY$-algebra morphism $I_n : {\bf LBA}_f({\bf 1},S^{\otimes n})
\to \wh{\bf U}_{n,f}$} \label{sect:ass:gr}

Define a linear map 
$$
I_n : {\bf LBA}_f({\bf 1},S^{\otimes n}) \simeq (\wh{S^{\otimes n} \boxtimes 
{\bf 1}})_f \simeq \hbox{\boldmath$\Pi$\unboldmath}_f
({\mathfrak 1}\underline\boxtimes {\mathfrak 1},
(S \underline\boxtimes {\mathfrak 1})^{\otimes n})\to \wh{\bf U}_{n,f}, 
$$
$$
\hbox{\boldmath$\Pi$\unboldmath}_f
({\mathfrak 1}\underline\boxtimes {\mathfrak 1},
(S \underline\boxtimes {\mathfrak 1})^{\otimes n})
\ni x\mapsto \kappa_1^\Pi(\on{R}_+)^{\boxtimes n} \circ x. 
$$
This is a morphism of $\cY$-algebras, where 
${\bf LBA}_f({\bf 1},S^{\otimes n})$ 
is equipped with the structure correponding to $S(\kappa_1)\circ Q$ and 
$\wh{\bf U}_{n,f}$ is equipped with its $\on{J}$-twisted structure. 

Then $I_n$ is a filtered map, and the associated graded is the inclusion  
$\on{LBA}_f({\bf 1},S^{\otimes n}) \simeq (S^{\otimes n} \boxtimes {\bf 1})_f
\hookrightarrow (\Delta(S^{\otimes n}))_f \simeq \Pi_f({\mathfrak 1}
\underline\boxtimes {\mathfrak 1},(S\underline\boxtimes S)^{\otimes n}) 
\simeq {\bf U}_{n,f}$.

\subsection{Construction of $(v,\on{F})$}

\begin{theorem} \label{thm:twists}
There exists a pair $(v,\on{F})$, 
where $v\in (\wh{\bf U}_{1,f})^\times_1$ and $\on{F}\in 
(\wh{(S^{\otimes 2} \boxtimes {\bf 1})_f})^\times_1$, such that 
\begin{equation} \label{gauge:J}
{\on J}(r+f) = v^1 v^2 I_2(\on{F}) 
{\on J}(r) (v^{12})^{-1} 
\end{equation}
(equality in $\wh{\bf U}_{2,f}$, where $v^{12}$ is defined using the 
$\cC$-algebra structure on $\wh{\bf U}_{n,f}$). 

Then 
\begin{equation} \label{F=twist}
({\on F})_Q^{1,2} \ast_Q ({\on F})_Q^{12,3}
= ({\on F})_Q^{2,3} \ast_Q ({\on F})_Q^{1,23} . 
\end{equation}
\end{theorem}

{\em Proof.} Write $v = 1 + v_1 + ...$, where $v_i\in {\bf U}_{1,f}$ 
has degree $i$ and $\on{F} = 1 + \on{F}_1 + \on{F}_2+...$, 
where $\on{F}_i \in (S^{\otimes n} \boxtimes {\bf 1})_f$ has degree $i$. 

If we set $\on{F}_1 = -f/2$, $v_1 = 0$, then (\ref{gauge:J}) holds 
modulo terms of degree $\geq 2$. 
 
Assume that we have found $v_1,...,v_{n-1}$ and $\on{F}_1,...,\on{F}_{n-1}$
such that (\ref{gauge:J}) holds modulo terms of degree $\geq n$.  

Let us set $v_{<n} = 1 + v_1 + ... + v_{n-1}$, 
$\on{F}_{<n} = 1 + \on{F}_1 + ... + \on{F}_{n-1}$. We then have 
\begin{equation} \label{gauge:1}
(v_{<n}^1 v_{<n}^2)^{-1} \on{J}(r+f) v_{<n}^{12} \on{J}(r)^{-1}
= I_2(\on{F}_{<n}) + \psi, 
\end{equation}
where $\psi = \psi_n + \psi_{n+1} + ...$ is an element of 
$\wh{\bf U}_{2,f}$ of degree $\geq n$. Let us denote by $\on{K} \in 
(\wh{\bf U}_{2,f})^\times_1$ the l.h.s. of (\ref{gauge:1}), then $\on{K}$
satisfies 
$$
\on{K}^{1,2} \on{J}(r)^{1,2}\on{K}^{12,3} (\on{J}(r)^{1,2})^{-1} 
= \on{K}^{2,3} \on{J}(r)^{2,3}\on{K}^{1,23} (\on{J}(r)^{2,3})^{-1}. 
$$
This implies that  
\begin{equation} \label{pre:twist}
I_3 \big( (\on{F}_{<n})_Q^{1,2} \ast_Q (\on{F}_{<n})_Q^{12,3}
\ast_Q ((\on{F}_{<n})_Q^{2,3} \ast_Q (\on{F}_{<n})_Q^{1,23})^{-1} \big) 
= 1 + \psi^{2,3} +  \psi^{1,23} -  (\psi^{1,2} + \psi^{12,3})  
\end{equation} 
modulo degree $>n$. The associated graded of $I_3$ is the 
composed map $(S^{\otimes 3} \boxtimes {\bf 1})_f \to 
(\Delta(S^{\otimes 3}))_f \simeq {\bf U}_{n,f}$, which is 
injective; hence so is $I_3$. 
Therefore $(\on{F}_{<n})_Q^{1,2} \ast_Q (\on{F}_{<n})_Q^{12,3}
\ast_Q ((\on{F}_{<n})_Q^{2,3} \ast_Q (\on{F}_{<n})_Q^{1,23})^{-1} = 1$ 
modulo degree $\geq n$. 

Moreover, $(S^{\otimes 3} \boxtimes {\bf 1})_f \to 
(\Delta(S^{\otimes 3}))_f \simeq {\bf U}_{3,f}$ 
is the linear isomorphism $(S^{\otimes 3} \boxtimes {\bf 1})_f 
\stackrel{\sim}{\to}  {\bf U}_{3,f}^{aaa}$, so 
$d(\psi_n) := \psi_n^{2,3} +  \psi_n^{1,23} -  (\psi_n^{1,2} + \psi_n^{12,3})
\in {\bf U}^{aaa}_{3,f}$.

Now $d(d(\psi_n)) = 0$ and $\on{Alt}(d(\psi_n)) = 0$, so the 
computation of the co-Hochschild cohomology of 
${\bf U}^{a...a}_{*,f}$ in Subsection \ref{cohoch} implies that 
$d(\psi_n) = d(\bar{\on{F}}'_n)$, where $\bar{\on{F}}'_n\in {\bf U}_{2,f}^{aa}$. 
The computation of the co-Hochschild cohomology for 
${\bf U}_{*,f}$  then implies that $\psi_n = \bar{\on{F}}'_n 
+ (v_n^{12} - v_n^1 - v_n^2) + \lambda'$, where 
$v_n\in {\bf U}_{1,f}$ and $\lambda'\in (\Delta(\wedge^2))_f$
all have degree $n$. 

Now $(\Delta(\wedge^2))_f = (\wedge^2\boxtimes {\bf 1})_f \oplus 
({\bf id} \boxtimes {\bf id})_f \oplus ({\bf 1} \boxtimes \wedge^2)_f$. 
Since $({\bf 1} \boxtimes \wedge^2)_f = 0$, we decompose $\lambda'$
as $\lambda'' + \lambda - \lambda^{2,1}$, where $\lambda''\in 
(\wedge^2\boxtimes {\bf 1})_f$ and $\lambda \in ({\bf id} \boxtimes 
{\bf id})_f$. 
Set $\bar{\on{F}}_n:= \bar{\on{F}}'_n + \lambda'' \in {\bf U}^{aa}_{2,f}$. 

Then 
$$
\psi_n = (v_n^{12} - v_n^1 - v_n^2) + \bar{\on{F}}_n + \lambda - \lambda^{2,1}.
$$ 

Let $\on{F}_n \in (S^{\otimes 2} \boxtimes {\bf 1})_f$ be the preimage of 
$\bar{\on{F}}_n$ under the symmetrization map 
$(S^{\otimes 2} \boxtimes {\bf 1})_f \to {\bf U}_{2,f}^{aa}$.  

Let us set $v_{\leq n} = (1 + v_n)(1+v_{<n})$, $\on{F}_{\leq n} 
= \on{F}_{<n} + \on{F}_n$. 
Then (\ref{gauge:1}) is rewritten as 
\begin{equation} \label{gauge:2}
(v^1_{\leq n} v^2_{\leq n})^{-1} \on{J}(r+f) v^{12}_{\leq n} \on{J}(r)^{-1}
= I_2(\on{F}_{\leq n}) + \lambda - \lambda^{2,1} + \psi', 
\end{equation}
where $\psi' = \psi'_{n+1} + ... \in \wh{\bf U}_{2,f}$ has degree $\geq n+1$. 

As above, we denote by $\on{K}'$ the l.h.s. of (\ref{gauge:2}).
We have again 
$$
(\on{K}')^{1,2} \on{J}(r)^{1,2} (\on{K}')^{12,3} (\on{J}(r)^{1,2})^{-1}
= (\on{K}')^{1,2} \on{J}(r)^{2,3} (\on{K}')^{1,23} (\on{J}(r)^{2,3})^{-1}, 
$$
which according to (\ref{gauge:2}) can be rewritten as 
\begin{align} \label{gauge:3}  
& I_3((\on{F}_{\leq n})_Q^{1,2} \ast_Q 
(\on{F}_{\leq n})_Q^{12,3} \ast_Q ((\on{F}_{\leq n})_Q^{2,3}
\ast_Q (\on{F}_{\leq n})_Q^{1,23})^{-1}) 
\\ & \nonumber 
= 1 + (\psi')^{2,3} + (\psi')^{1,23} 
- ((\psi')^{1,2} + (\psi')^{12,3} ) 
+ [r_1^{1,2},\lambda^{12,3} - \lambda^{3,12}]/2 
- [r_1^{2,3},\lambda^{1,23} - \lambda^{23,1}]/2 
\end{align}
modulo terms of degree $>n+1$. Here $r_1 = \kappa_1^\Pi(r)\in {\bf
U}_{2,f}^{ab}$. 

As above, this equation implies that 
$(\on{F}_{\leq n})_Q^{1,2} \ast_Q (\on{F}_{\leq n})_Q^{12,3} 
\ast_Q ((\on{F}_{\leq n})_Q^{2,3} \ast_Q
(\on{F}_{\leq n})_Q^{1,23})^{-1}$ has the form $1 + \bar g_{n+1} + ...$, 
where $\bar g_{n+1},...$ have degree $\geq n+1$. The degree $n+1$ part of 
(\ref{gauge:3}) yields 
$$
g_{n+1} = (\psi'_{n+1})^{2,3} + (\psi'_{n+1})^{1,23} 
- ((\psi'_{n+1})^{1,2} + (\psi'_{n+1})^{12,3} ) 
+ [r_1^{1,2},\lambda^{12,3} - \lambda^{3,12}]/2 
- [r_1^{2,3},\lambda^{1,23} - \lambda^{23,1}]/2 
$$   
where $g_{n+1} \in {\bf U}^{aaa}_{3,f}$ is the image of $g_{n+1}$
by the isomorphism $(S^{\otimes 3} \boxtimes {\bf 1})_f \simeq 
{\bf U}_{2,f}^{aaa}$. Applying $\on{Alt}$ to this equation, we get 
$$
\on{Alt}([r_1^{1,2},\lambda^{12,3} - \lambda^{3,12}] 
- [r_1^{2,3},\lambda^{1,23} - \lambda^{23,1}]) \in 
{\bf U}_{3,f}^{aaa}. 
$$ 
Now the terms under Alt belong to ${\bf U}_{3,f}^{c_1 c_2 c_3}$
where $c_1 c_2 c_3$ is respectively $aab$, $bba$, $abb$, $baa$. 
These terms are antisymmetric w.r.t. the pairs of repeated indices, 
and $- [r_1^{2,3},\lambda^{1,23}] = ([r_1^{1,2},-\lambda^{3,12}])^{2,3,1}$, 
$[r_1^{2,3}, \lambda^{23,1}] = ([r_1^{1,2},\lambda^{12,3}])^{2,3,1}$. 
Hence 
$$
[r_1^{1,2},\lambda^{12,3} - \lambda^{3,12}] + \on{cyclic\ permutations} = 0. 
$$
Since the spaces
${\bf U}^{c_1 c_2 c_3}_{3,f}$ are in direct sum for distinct 
$c_1 c_2 c_3$, we get  
$$
[r_1^{1,2},\lambda^{12,3}] = [r_1^{2,3},\lambda^{1,23}] = 0. 
$$  
We will show that the second equality implies that 
$\lambda =0$. This will prove the induction step, because (\ref{gauge:2}) 
then means that (\ref{gauge:J}) holds at step $n+1$.  

So it remains to prove: 
\begin{lemma}
The composition 
\begin{equation} \label{map:star}
({\bf id} \boxtimes {\bf id})_f \hookrightarrow 
{\bf U}^{ab}_{2,f} \to {\bf U}_{3,f} 
\end{equation} 
is injective, where the first map is 
$({\bf id} \boxtimes {\bf id})_f \simeq \Pi_f({\mathfrak 1}\underline\boxtimes
{\mathfrak 1},{\bf id}\underline\boxtimes {\bf id}) 
\hookrightarrow \Pi_f({\mathfrak 1}\underline\boxtimes {\mathfrak 1},
S\underline\boxtimes S) \simeq {\bf U}_{2,f}^{ab}$
and the second map is ${\bf U}_{2,f}^{ab} \ni \lambda \mapsto 
[r_1^{2,3},\lambda^{1,23}] \in {\bf U}_{3,f}$. 
\end{lemma}

{\em Proof of Lemma.} It follows from Subsection \ref{def:r} 
that (\ref{map:star}) coincides with the composition 
$({\bf id} \boxtimes {\bf id})_f \to 
({\bf id} \boxtimes \wedge^2)_f \hookrightarrow 
({\bf id} \boxtimes {\bf id}^{\otimes 2})_f \hookrightarrow 
{\bf U}^{abb}_{3,f} \hookrightarrow {\bf U}_{3,f}$, where the 
first map is 
$$
({\bf id}\boxtimes {\bf id})_f \simeq \on{LBA}_f({\bf id},{\bf id}) 
\stackrel{- \circ \mu}{\to} \on{LBA}_f(\wedge^2,{\bf id})\simeq 
({\bf id}\boxtimes \wedge^2)_f, 
$$
the second and the fourth maps are the natural injections, and the 
third map is the injection 
$$
({\bf id}\boxtimes {\bf id}^{\otimes 2})_f \simeq 
\Pi_f({\mathfrak 1}\underline\boxtimes {\mathfrak 1},
({\bf id}\underline\boxtimes {\mathfrak 1}) \otimes 
({\mathfrak 1} \underline\boxtimes {\bf id})^{\otimes 2}) 
\hookrightarrow \Pi_f({\mathfrak 1}\underline\boxtimes {\mathfrak 1},
(S\underline\boxtimes {\mathfrak 1}) 
\otimes ({\mathfrak 1} \underline\boxtimes S)^{\otimes 2}) 
\simeq {\bf U}_{3,f}^{abb}. 
$$
It follows from Proposition \ref{prop:inject} that the first map 
is also injective. Therefore the map 
$({\bf id}\boxtimes {\bf id})_f \to {\bf U}_{3,f}$
given by (\ref{map:star}) is injective. This proves the Lemma. 
\hfill \qed \medskip 

This ends the proof of the first part of Theorem \ref{thm:twists}. 
Equation (\ref{F=twist}) is then obtained by taking the limit 
$n\to\infty$ in (\ref{pre:twist}). This proves Theorem \ref{thm:twists}. 
\hfill \qed \medskip 

We prove that pairs $(v,\on{F})$ are unique up to gauge
(this fact will not be used in the sequel). 

\begin{lemma}
The set of pairs $(v,\on{F})$ as in Theorem \ref{thm:twists} 
is a torsor under the action of 
$(\wh{(S\boxtimes {\bf 1})_f})^\times_1$: an element 
$g\in (\wh{(S\boxtimes {\bf 1})_f})^\times_1$
transforms $(v,\on{F})$ into $(v I_1(g),((g)_Q^1 \ast_Q (g)_Q^2)^{-1} 
\ast_Q \on{F} \ast_Q (g)^{12}_Q)$. 
\end{lemma}

{\em Proof.} Since $I_2(g^{12}_Q) = \on{J}(r) I_1(g)^{12}
\on{J}(r)^{-1}$, the pair $(v I_1(g),((g)_Q^1 \ast_Q (g)_Q^2)^{-1} 
\ast_Q \on{F} \ast_Q (g)^{12}_Q)$ is also a solution 
of the equation of Theorem \ref{thm:twists}. 
Conversely, let 
$(v_1,\on{F_1})$ and $(v_2,\on{F_2})$ be solutions of this equation. 
Then $v_1^1 v_1^2 I_2(\on{F}_1) \on{J}(r) (v_1^{12})^{-1}
= v_2^1 v_2^2 I_2(\on{F}_2) \on{J}(r) (v_2^{12})^{-1}$. 
Let $n$ be the smallest index such that the degree $n$ components 
of $(v_1,\on{F}_1)$ and $(v_2,\on{F}_2)$ are different. We denote with an
additional index $n$ these components. Then we have 
$$
(v_{2,n} - v_{1,n})^{12} - (v_{2,n} - v_{1,n})^1 - (v_{2,n} - v_{1,n})^2
= \on{sym}_2(\on{F}_{2,n} - \on{F}_{1,n}), 
$$
where $\on{sym}_2 : (S^{\otimes 2}\boxtimes {\bf 1})_f \to 
(\Delta(S^{\otimes 2}))_f \simeq {\bf U}_{2,f}$ is the 
canonical injection. So $d(v_{2,n} - v_{1,n}) \in {\bf U}^{aa}_{2,f}$. 
As above, we
obtain the existence of $w\in {\bf U}^a_{1,f}$ of degree $n$, 
such that $d(v_{2,n} - v_{1,n}) = d(w)$. Therefore 
$v_{2,n} - v_{1,n} - w\in (\Delta({\bf id}))_{f} = 
({\bf id} \boxtimes {\bf 1})_f \oplus ({\bf 1} \boxtimes {\bf id})_f$. 
Now $({\bf 1} \boxtimes {\bf id})_f = 0$, so 
$v_{2,n} - v_{1,n} - w\in ({\bf id} \boxtimes {\bf 1})_f \subset 
{\bf U}^a_{1,f}$. Therefore $w' := v_{2,n} - v_{1,n} \in {\bf U}^a_{1,f}$. 
Replacing $(v_2,\on{F}_2)$ by $(v_2 I_1(1-w'),
((1-w')_Q^1 \ast_Q (1-w')_Q^2)^{-1} \ast_Q \on{F}_2 
\ast_Q (1-w')^{12}_Q)$, we obtain
a solution equal to $(v_1,\on{F}_1)$ up to degree $n$. Proceeding 
inductively, we see that 
$(v_1,\on{F}_1)$ and $(v_2,\on{F}_2)$ are related by the action
of an element of $(\wh{(S\boxtimes {\bf 1})_f})^\times_1$.   
\hfill \qed \medskip

\subsection{Compatibility of quantization functors with twists}

Let ${\bf P}$ be a prop and 
$\overline Q : \on{Bialg} \to S({\bf P})$ be a 
prop morphism. 
Using $\overline Q$, we equip the collection of all 
${\bf P}({\bf 1},S^{\otimes n})$ with 
the structure of a $\cY$-algebra (see Subsection \ref{C:algebras})
with unit $\overline Q(\eta^{\otimes n})$. 

We define a twist of $\overline Q$ to be an element $\on{F}$ of 
${\bf P}({\bf 1},S^{\otimes 2})^\times$, 
such that the relations
$$
(\on{F})_{\overline Q}^{1,2} \ast_{\overline{Q}}
(\on{F})_{\overline Q}^{12,3} = 
(\on{F})_{\overline Q}^{2,3} \ast_{\overline{Q}}
(\on{F})_{\overline Q}^{1,23}, \quad 
(\on{F})_{\overline Q}^{\emptyset,1} 
= (\on{F})_{\overline Q}^{1,\emptyset} = \overline Q(\eta)
$$
hold in ${\bf P}({\bf 1},S^{\otimes 3})$ 
(here we use the $\cY$-algebra 
structure on ${\bf P}({\bf 1},S^{\otimes n})$ 
given by $\overline Q$). 
Then get a new prop morphism ${}^{\on{F}}\overline Q : \on{Bialg} \to 
S({\bf P})$, 
defined by ${}^{\on{F}}\overline Q(m) = \overline Q(m)$, 
${}^{\on{F}}\overline Q(\Delta) = \underline{\on{Ad}}(\on{F}) 
\circ \overline Q(\Delta)$, ${}^{\on{F}}\overline Q(\eps) 
= \overline Q(\eps)$, ${}^{\on{F}}\overline Q(\eta) = \overline Q(\eta)$. 
Here $\underline{\on{Ad}}({\on F}) \in S({\bf P})({\bf id}^{\otimes 2},
{\bf id}^{\otimes 2})$ is given by 
$$
\underline{\on{Ad}}(\on{F}) = \overline Q( m^{(2)} \otimes 
m^{(2)}) \circ (142536) \circ 
(\on{F} \otimes \on{id}_{S^{\otimes 2}} 
\otimes \on{F}^{-1})
\in {\bf P}(S^{\otimes 2},S^{\otimes 2});   
$$ 
here $m^{(2)} = m \circ (m\otimes \on{id}_{{\bf id}})
\in \on{Bialg}(T_3,{\bf id})$. 

We say that the prop morphisms $\overline Q,\overline Q' : \on{Bialg} \to 
S({\bf P})$ are equivalent if 
$\overline Q' = \theta(\xi) \circ \overline Q$, 
where $\xi\in S({\bf P})({\bf id},{\bf id})^\times$
and $\theta(\xi)$ is the corresponding inner automorphism of 
$S({\bf P})$. 

\begin{theorem} Let $Q : \on{Bialg} \to S({\bf LBA})$ 
be an Etingof-Kazhdan quantization functor. Then  
$S(\kappa_i) \circ Q : \on{Bialg} \to S({\bf LBA}_f)$ 
($i=1,2$) are prop morphisms. 
There exists $\on{i}\in S({\bf LBA}_f)
({\bf id},{\bf id})^\times$, 
such that $\kappa_0(\on{i}) = S({\bf LBA})(\on{id}_{{\bf id}})$,   
and a twist ${\on F}$ of $S(\kappa_1) \circ Q$, such that 
$$
S(\kappa_2) \circ Q = \theta(\on{i}) \circ {}^{\on{F}}(S(\kappa_1) \circ Q). 
$$
\end{theorem}

{\em Proof.} We will construct $\on{i}$, such that 
$$
\kappa_2^\Pi(m_a) = \on{i} \circ \kappa_1^\Pi(m_a) \circ 
(\on{i}^{\boxtimes 2})^{-1}, \quad 
\kappa_2^\Pi(\Delta_a) = \on{i}^{\boxtimes 2} \circ 
\underline{\on{Ad}}(\on{F})
\circ \kappa_1^\Pi(\Delta_a) \circ \on{i}^{-1},  
$$
where as before 
$$
\underline{\on{Ad}}(\on{F}) = (m_a^{(2)} \boxtimes m_a^{(2)}) \circ 
(142536) \circ \big( \on{F} \boxtimes 
\on{id}_{(S\underline\boxtimes
{\mathfrak 1})^{\boxtimes 2}} \boxtimes \on{F}^{-1}\big) \in
\hbox{\boldmath$\Pi$\unboldmath}_f
((S\underline\boxtimes {\mathfrak 1})^{\otimes 2}, 
(S\underline\boxtimes {\mathfrak 1})^{\otimes 2})^\times. 
$$

Let us relate the $\kappa_i^\Pi(m_\Pi)$, $i=1,2$. (\ref{24}) implies that 
$$
\kappa^\Pi_2\big( \on{Ad}(\on{J}) \big) = \Xi_f^{\boxtimes 2} \circ
\on{Ad}(v)^{\otimes 2}\circ \kappa^\Pi_1(
\on{Ad}(\on{R}_+^{\boxtimes 2} \circ \on{F}) \circ  
\on{Ad}(\on{J})) \circ \on{Ad}(v^{12})^{-1}
\circ (\Xi_f^{-1})^{\boxtimes 2},  
$$
and therefore (\ref{25}) implies that  
$$
\kappa^\Pi_2(\Delta_\Pi) = \Xi_f^{\boxtimes 2} \circ 
\on{Ad}(v)^{\boxtimes 2}\circ \kappa^\Pi_1(
\on{Ad}(\on{R}_+^{\boxtimes 2} \circ \on{F}) \circ  \Delta_\Pi) 
\circ \on{Ad}(v)^{-1}\circ \Xi_f^{-1}. 
$$

Now $\on{R}(r+f) = (m_\Pi^{(2)} \boxtimes m_\Pi^{(2)}) \circ (142536) 
\circ \Big(  
( \on{R}_+^{\boxtimes 2} \circ \on{F}^{2,1}) 
\boxtimes \on{R}(r) \boxtimes 
(\on{R}_+^{\boxtimes 2} \circ \on{F}^{-1}) \Big)$, where $\on{F}
\in \hbox{\boldmath$\Pi$\unboldmath}_f({\mathfrak 1}\underline\boxtimes 
{\mathfrak 1}, (S\underline\boxtimes {\mathfrak 1})^{\otimes 2})$. 
For $X\in \hbox{\boldmath$\Pi$\unboldmath}_f
({\mathfrak 1}\underline\boxtimes {\mathfrak 1}, (S\underline\boxtimes 
{\mathfrak 1})^{\otimes 2})$, set $\underline X :=
(\on{can}_{{\mathfrak 1}\underline\boxtimes S}^* 
\boxtimes \on{id}_{S\underline\boxtimes {\mathfrak 1}}) 
\circ (\on{id}_{{\mathfrak 1}\underline\boxtimes S} \boxtimes X)$. 

Then 
\begin{align*}
& (\kappa_2^\Pi(\on{R}_+)\boxtimes \kappa_2^\Pi(\on{R}_-)) \circ 
\on{can}_{S\underline\boxtimes{\mathfrak 1}} = \kappa_2^\Pi(\on{R}) 
= \Xi_f^{\boxtimes 2} \circ \on{R}(r+f) 
= 
\Big( \big( \Xi_f \circ \on{Ad}(v) \circ \kappa_1^\Pi(\on{R}_+) \big) 
\\ & \boxtimes 
\big( \Xi_f \circ \on{Ad}(v) \circ 
\kappa_1^\Pi\big( m_\Pi^{(2)*} \circ (\on{R}_+ \boxtimes \on{R}_- 
\boxtimes \on{R}_+)\big) \circ \big( \underline{\on{F}}^{2,1} \boxtimes 
\on{id}_{{\mathfrak 1}\underline\boxtimes S} \boxtimes 
\underline{\on{F}}^{-1} \big) \big) \Big) \circ 
\on{can}_{S\underline\boxtimes {\mathfrak 1}}.  
\end{align*}
It follows that for some $\on{i}\in 
\hbox{\boldmath$\Pi$\unboldmath}_f(S\underline\boxtimes 
{\mathfrak 1},S \underline\boxtimes{\mathfrak 1})^\times$, we have 
$$
\kappa_2^\Pi(\on{R}_+) = \Xi_f \circ \on{Ad}(v) 
\circ \kappa_1^\Pi(\on{R}_+) \circ \on{i}^{-1},  
$$
therefore $\kappa_2^\Pi(\on{R}_+^{(-1)}) = \on{i} \circ 
\kappa_1^\Pi(\on{R}_+^{(-1)}) \circ \on{Ad}(v^{-1}) \circ \Xi_f^{-1}$. 

Now 
\begin{align*}
& \kappa_2^\Pi(m_a) = \kappa_2^\Pi(\on{R}_+^{(-1)}) \circ \kappa_2^\Pi(m_\Pi) 
\circ \kappa_2^\Pi(\on{R}_+)^{\boxtimes 2} = \on{i} \circ 
\kappa_1^\Pi(\on{R}_+^{(-1)}) \circ \kappa_1^\Pi(m_\Pi) \circ 
\kappa_1^\Pi(\on{R}_+^{\boxtimes 2}) \circ (\on{i}^{-1})^{\boxtimes 2} 
\\ & = \on{i} \circ \kappa_1^\Pi(m_a) \circ (\on{i}^{-1})^{\boxtimes 2},   
\end{align*}
and 
$$
\kappa_2^\Pi(\Delta_a) = \kappa_2^\Pi(\on{R}_+^{(-1)})^{\boxtimes 2} 
\circ \kappa_2^\Pi(\Delta_\Pi) \circ \kappa_2^\Pi(\on{R}_+) = 
\on{i}^{\boxtimes 2} \circ \kappa_1^\Pi(\on{R}_+^{(-1)})^{\boxtimes 2} \circ 
\kappa_1^\Pi(\on{Ad}(\on{R}_+^{\boxtimes 2} \circ \on{F}) \circ \Delta_\Pi)
\circ \kappa_1^\Pi(\on{R}_+) \circ \on{i}^{-1}. 
$$

We first prove that 
\begin{equation} \label{Delta:a}
(\kappa_1^\Pi(\on{R}_+^{(-1)}))^{\boxtimes 2} \circ 
\on{Ad}(\kappa_1^\Pi(\on{R}_+)^{\boxtimes 2}
\circ \on{F}) \circ \kappa_1^\Pi(\Delta_\Pi \circ \on{R}_+) 
= \underline{\on{Ad}}(\on{F}) \circ \kappa_1^\Pi(\Delta_a). 
\end{equation}

One checks that $\on{Ad}(\kappa_1^\Pi(\on{R}_+)^{\boxtimes 2} \circ 
\on{F}) \circ \kappa_1^\Pi(\Delta_\Pi \circ \on{R}_+) = 
\kappa_1^\Pi(\on{R}_+)^{\boxtimes 2} \circ \underline{\on{Ad}}(\on{F})$. 
We have 
\begin{align*}
& \Delta_\Pi \circ \on{R}_+ = \Delta_\Pi \circ \on{R}_+ \circ 
\on{R}_-^* \circ (\on{R}_-^{(-1)})^* = \Delta_\Pi \circ \ell \circ 
(\on{R}_-^{(-1)})^*  
\\ & 
= \ell^{\boxtimes 2} \circ \overline{\Delta}_\Pi \circ (\on{R}_-^{(-1)})^*
= \on{R}_+^{\boxtimes 2} \circ (\on{R}_-^*)^{\boxtimes 2} \circ 
\overline{\Delta}_\Pi 
\circ (\on{R}_-^{(-1)})^* = \on{R}_+^{\boxtimes 2} \circ \Delta_a. 
\end{align*}
Applying $\kappa_1^\Pi$ and composing from the left with 
$\on{Ad}(\kappa_1^\Pi(\on{R}_+)^{\boxtimes 2} 
\circ \on{F})$, we get $\on{Ad}(\kappa_1^\Pi(\on{R}_+)^{\boxtimes 2} 
\circ \on{F}) \circ \kappa_1^\Pi(\Delta_\Pi \circ \on{R}_+) =
\kappa_1^\Pi(\on{R}_+)^{\boxtimes 2} \circ 
\underline{\on{Ad}}(\on{F})\circ \kappa_1^\Pi(\Delta_a)$. 
Composing from the left with $\kappa_1^\Pi(\on{R}_+^{(-1)})^{\boxtimes 2}$, 
we get (\ref{Delta:a}). 

It follows that 
$$
\kappa_2^\Pi(\Delta_a) 
= \on{i}^{\boxtimes 2} \circ \underline{\on{Ad}}(\on{F}) 
\circ \kappa_1^\Pi(\Delta_a) \circ \on{i}^{-1}, 
$$
as wanted. \hfill \qed \medskip

\section{Quantization of coboundary Lie bialgebras} \label{sect:cob}

\subsection{Compatibility with coopposite}

Let $\Phi$ be an associator. Then $\Phi' := \Phi(-A,-B)$ is also an 
associator. Let $Q,Q'$ be the Etingof-Kazhdan quantization functors
corresponding to $\Phi,\Phi'$. 

Recall that $\tau_{\on{LBA}}\in \on{Aut}(\on{LBA})$ 
is defined by $\mu\mapsto \mu$, 
$\delta \mapsto -\delta$ and let $\tau_{\on{Bialg}}\in 
\on{Aut}(\on{Bialg})$ be defined by $m \mapsto m$, $\Delta \mapsto 
(21) \circ \Delta$. 

\begin{proposition} \label{prop:coop}
There exists $\xi_\tau\in S({\bf LBA})({\bf id},{\bf id})^\times$, 
with $\xi_\tau = \on{id}_{{\bf id}}$ + terms of 
positive degree in both $\mu$ and $\delta$, such that 
$$
Q'\circ \tau_{\on{Bialg}} = \theta(\xi_\tau) \circ S(\tau_{\on{LBA}}) \circ Q.  
$$
\end{proposition}

{\em Proof.} This means that 
$$
Q'(m) = \xi_\tau^{\boxtimes 2} \circ S(\tau_{\on{LBA}})(Q(m)) \circ
\xi_\tau^{-1}, \quad 
Q'(\Delta) \circ (21) = \xi_\tau \circ S(\tau_{\on{LBA}})(Q(\Delta))
\circ (\xi_\tau^{-1})^{\boxtimes 2}. 
$$
We will therefore construct 
$\tilde \xi_\tau\in \hbox{\boldmath$\Pi$\unboldmath}
(S\underline\boxtimes {\mathfrak 1}, 
S\underline\boxtimes {\mathfrak 1})^\times$, such that 
$$
m'_a = \tilde \xi_\tau^{\boxtimes 2} \circ \tau_\Pi(m_a) \circ 
\tilde\xi_\tau^{-1}, 
\quad
\Delta'_a \circ (21) = \tilde \xi_\tau \circ \tau_\Pi(\Delta_a) 
\circ (\tilde\xi_\tau^{\boxtimes 2})^{-1}, 
$$
where $m'_a,\Delta'_a$ are the analogues of $m_a,\Delta_a$ for $\Phi'$. 

\begin{lemma} \label{tau:adJ}
$\tau_\Pi(\on{Ad}(\on{J})) = 
(\on{id}_S \underline\boxtimes \omega_S)^{\boxtimes 2} \circ 
\on{Ad}(\on{J}(-r)) \circ 
((\on{id}_S \underline\boxtimes \omega_S)^{\boxtimes 2})^{-1}$. 
\end{lemma}

This follows from Lemma \ref{tau:pi:m:pi}. 

\begin{lemma}
Let $\on{J}'$ be the analogue of $\on{J}$ for $\Phi'$. 
There exists $u\in \wh{{\bf U}}_1$ of the form $u = 1$ + terms of 
degree $\geq 1$, such that  
\begin{equation} \label{twist:J:J'}
(\on{J}')^{2,1} = u^1 u^2 \on{J}(-r) (u^{12})^{-1}.  
\end{equation}
\end{lemma}

{\em Proof of Lemma.}
We have 
$$
\on{J}(-r)^{1,2} \on{J}(-r)^{12,3} = \on{J}(-r)^{2,3} \on{J}(-r)^{1,23} 
\Phi(-t_{12},-t_{23}) = \on{J}(-r)^{2,3} \on{J}(-r)^{1,23}  \Phi' 
$$
(equality in $\wh{\bf U}_2$, where we use the 
$\cX$-algebra structure on $\wh{\bf U}_n$). 
Let us set $u_0 = 1$ + class of $a_1b_1/2$ and $\overline{\on{J}}
:= u_0^1u_0^2 \on{J}(-r)^{2,1} (u_0^{12})^{-1}$.  Then $\overline{\on{J}}$ 
satisfies 
$\overline{\on{J}}^{1,2}\overline{\on{J}}^{12,3} = 
\overline{\on{J}}^{2,3}\overline{\on{J}}^{1,23}\Phi'$
(since $(\Phi')^{3,2,1} = (\Phi')^{-1}$)
and $\overline{\on{J}} = 1 - r/2$ + terms of degree $>1$, and 
$\on{J}'$
satisfies the same conditions. According to \cite{Enr:coh}, this implies the
existence of $u_1\in \wh{\bf U}_1$ of the form $u_1 = 1$ + terms of
degree $>1$, such that $\on{J}' = u_1^1 u_1^2\overline{\on{J}}
(u_1^{12})^{-1}$, 
so if we set $u = u_1u_0\in \wh{\bf U}_1$, then $u$ has the form 
$u = 1$ + class of $a_1 b_1/2$ + terms of degree $>1$, and 
satisfies (\ref{twist:J:J'}). \hfill \qed \medskip 

\begin{lemma} We have $\tau_\Pi(\Delta_\Pi) = 
(\on{id}_S \underline\boxtimes \omega_S)^{\boxtimes 2}
\circ \on{Ad}(u^{-1})^{\boxtimes 2} \circ \Big( (21) \circ \Delta'_\Pi\Big) 
\circ \on{Ad}(u) \circ (\on{id}_S \underline\boxtimes \omega_S)^{-1}$. 
\end{lemma}

This follows from Lemmas \ref{tau:pi:m:pi} and \ref{tau:adJ}.

\begin{lemma}
$\tau_\Pi(\on{R}) =  
(\on{id}_S \underline\boxtimes \omega_S)^{\boxtimes 2}
\circ \on{Ad}(u^{-1})^{\boxtimes 2} \circ (\on{R}')^{-1}$.  
\end{lemma}

{\em Proof of Lemma.} We have $\on{R} = \sum_{n\geq 0}
(n!)^{-1} (m_\Pi^{(n+1)} \boxtimes m_\Pi^{(n+1)}) \circ (1,n+3,2,n+4,...)\circ
\big(\on{J}^{2,1} \boxtimes t^{\boxtimes n} \circ \on{J}^{-1}\big)$, 
so 
\begin{align*}
& \tau_\Pi(\on{R}) = \sum_{n\geq 0} (n!)^{-1}
(\on{id}_S \underline\boxtimes
\omega_S)^{\boxtimes 2} \circ (m_\Pi^{(n+1)} \boxtimes m_\Pi^{(n+1)}) \circ 
(\on{id}_S \underline\boxtimes
\omega_S)^{\boxtimes 2(n+2)})^{-1} 
\\ & \circ \big( \tau_\Pi(\on{J}^{2,1})
\boxtimes \tau_\Pi(t)^{\boxtimes n} \boxtimes \tau_\Pi(\on{J}^{-1}) \big) 
\\ & 
= \sum_{n\geq 0} (n!)^{-1} 
(\on{id}_S \underline\boxtimes
\omega_S)^{\boxtimes 2} \circ (m_\Pi^{(n+1)} \boxtimes m_\Pi^{(n+1)}) \circ 
\big( \on{J}(-r)^{2,1}
\boxtimes (-t)^{\boxtimes n} \boxtimes \on{J}(-r)^{-1} \big) 
\\ & 
= (\on{id}_S \underline\boxtimes
\omega_S)^{\boxtimes 2} \circ \big( \on{J}(-r)^{2,1} e^{-t/2} \on{J}(-r)\big) 
= (\on{id}_S \underline\boxtimes
\omega_S)^{\boxtimes 2} \circ \on{R}(-r).  
\end{align*}

Now $\on{R} = (u^1u^2) (\on{R}')^{-1}(u^1u^2)^{-1}$, whence the result. 
\hfill \qed \medskip

\begin{lemma}
There exists $\sigma\in \hbox{\boldmath$\Pi$\unboldmath}({\mathfrak 1}
\underline\boxtimes S, {\mathfrak 1}\underline\boxtimes S)^\times$, such that 
$$
\on{R}^{-1} = \Big( \on{R}_+ \boxtimes (\on{R}_- \circ \sigma)\Big) \circ 
\on{can}_{S\underline\boxtimes {\mathfrak 1}}.  
$$ 
\end{lemma}

{\em Proof.} Set $\on{can} := \on{can}_{S\underline\boxtimes {\mathfrak 1}}$
and $\on{can}_+ := \sum_{i\geq 1} \on{can}_{S^i \underline\boxtimes 
{\mathfrak 1}}$. Set $m_b := \Delta_a^*$. Then $m_\Pi \circ 
\on{R}_-^{\boxtimes 2} = \on{R}_- \circ m_b$. The series 
$\on{can}' := \on{can}_{{\bf 1}\underline\boxtimes{\mathfrak 1}}
+ \sum_{i>0} (m_a^{(i)} \boxtimes m_a^{(i)}) \circ (1,i+1,2,i+2,...) \circ 
(-\on{can}_+)^{\boxtimes i}$
is convergent and has the form $(\on{id}_{S\underline\boxtimes{\mathfrak 1}}
\boxtimes \sigma)\circ \on{can}_{S\underline\boxtimes
{\mathfrak 1}}$ for a suitable invertible $\sigma$. 
We then have $(m_a \boxtimes m_b) \circ (32) \circ
(\on{can}_{S\underline\boxtimes {\mathfrak 1}} \boxtimes \on{can}') = 
\on{can}_{{\mathfrak 1}\underline\boxtimes {\mathfrak 1}}$. 

It follows that $\on{R}^{-1} = (\on{R}_+ \boxtimes \on{R}_-) \circ 
\on{can}' = \Big( \on{R}_+ \boxtimes (\on{R}_- \circ \sigma)\Big) \circ 
\on{can}_{S\underline\boxtimes {\mathfrak 1}}$. \hfill \qed \medskip 

{\em End of proof of Proposition \ref{prop:coop}.} 
The above lemmas imply  
$$
\tau_\Pi(\on{R}) = \Big( (\on{id}_S \boxtimes \omega_S) 
\circ \on{Ad}(u^{-1}) \Big)^{\boxtimes 2}
\circ \Big( \on{R}'_+ \boxtimes (\on{R}'_- \circ \sigma')\Big) \circ 
\on{can}_{S\underline\boxtimes {\bf 1}}, 
$$
where $\sigma'$ is the analogue of $\sigma$ for $\Phi'$. 
Since $\tau_\Pi(\on{R}) = \Big( \tau_\Pi(\on{R}_+) \boxtimes 
\tau_\Pi(\on{R}_-)\Big) \circ \on{can}_{S\underline\boxtimes 
{\mathfrak 1}}$, there exists $\xi_\tau\in \hbox{\boldmath$\Pi$\unboldmath}
(S\underline\boxtimes{\mathfrak 1}, S\underline\boxtimes{\mathfrak 1})^\times$, 
such that 
$$
\tau_\Pi(\on{R}_+) = (\on{id}_S \underline\boxtimes \omega_S)
\circ \on{Ad}(u)^{-1} \circ \on{R}'_+ \circ \xi_\tau. 
$$
It follows that 
$$
\tau_\Pi(\on{R}_+^{(-1)}) = \xi_\tau^{-1} \circ \on{R}_+^{\prime(-1)}
\circ \on{Ad}(u) \circ (\on{id}_S \underline\boxtimes \omega_S)^{-1}. 
$$

Then 
\begin{align*}
& \tau_\Pi(m_a) = \tau_\Pi(\on{R}_+^{(-1)} \circ m_\Pi \circ 
\on{R}_+^{\boxtimes 2}) = \xi_\tau^{-1} \circ \on{R}_+^{\prime (-1)}
\circ \on{Ad}(u) \circ m_\Pi \circ (\on{Ad}(u)^{-1})^{\boxtimes 2}
\circ (\on{R}'_+)^{\boxtimes 2} \circ \xi_\tau^{\boxtimes 2}
\\ & = \xi_\tau^{-1} \circ m'_a \circ \xi_\tau^{\boxtimes 2} 
\end{align*}
and 
\begin{align*}
& \tau_\Pi(\Delta_a) = \tau_\Pi\Big( (\on{R}_+^{(-1)})^{\boxtimes 2} \circ 
\Delta_\Pi \circ \on{R}_+ \Big) = \Big( \xi_\tau^{-1} \circ 
\on{R}_+^{\prime(-1)}\Big)^{\boxtimes 2} \circ (21)
\circ \Delta'_\Pi \circ (\on{R}'_+)^{\boxtimes 2} \circ \xi_\tau 
\\ & = (\xi_\tau^{-1})^{\boxtimes 2} \circ \Delta'_a \circ \xi_\tau.  
\end{align*}

Moreover, the image $(\xi_\tau)_{|\mu = \delta = 0} 
\in S({\bf Sch})({\bf id},{\bf id})$ of $\xi_\tau$ by the morphism 
${\bf LBA} \to {\bf Sch}$, $\mu,\delta\mapsto 0$ is equal to 
$\on{id}_{{\bf id}}$. 
Set $\xi' := (\xi_\tau)_{|\delta=0}$, then $\xi'\in {\bf LA}(S,S)$. 
We have $\on{LA}(S^p,S^q) = 0$ unless $p=q$, and $\on{LA}(S^p,S^p)
= \kk \on{id}_{S^p}$. So $\xi' = \on{id}_S$. Now  $\xi_\tau
= \xi'$ + terms of positive degree in $\delta$, so
$\xi_\tau = \on{id}_S$ + terms of positive degree in $\delta$. 
In the same way, $\xi_\tau = \on{id}_S$ + terms of positive 
degree in $\mu$. So $\xi_\tau = \on{id}_S$ + terms of positive
degree in both $\delta$ and $\mu$.    
 
This ends the proof of Proposition \ref{prop:coop}. \hfill \qed \medskip 

\begin{remark}
We take this opportunity to correct a mistake in Theorem 2.1 in  
\cite{Enr:coh}. Let $\on{J} = 1 - r/2 + ...$ be a solution 
of (\ref{eq:J}). Then the set of solutions of (\ref{eq:J}) of the form 
$1 +$ terms of degree $\geq 1$ consists in the disjoint union of 
{\it two} gauge orbits (and not one), that of $\on{J}$ and 
that of $\on{J}^{2,1}$. The degree one term of the solution has
the form $\alpha r + \beta r^{2,1}$, where $\alpha - \beta = \pm 1/2$; 
the solution is in the gauge class of $\on{J}$ (resp., $\on{J}^{2,1}$)
iff $\alpha - \beta = -1/2$ (resp., $1/2$). This follows from a more careful
analysis in degree one in the proof of Theorem 2.1 in \cite{Enr:coh}.  
\end{remark}

\subsection{Quantization functors for coboundary Lie bialgebras}

A quantization functor of coboundary Lie bialgebras is 
a prop morphism $\overline Q : \on{COB} \to S({\bf Cob})$, such that: 

(a) the composed morphism $\on{Bialg}\to\on{COB} \stackrel{\overline Q}{\to} 
S({\bf Cob}) \stackrel{\mu,r\mapsto 0}{\to}
S({\bf Sch})$ is the propic version of the bialgebra structure on 
the symmetric algebras, 

(b) $\overline Q(R) = \on{inj}_0^{\otimes 2}$ + terms of degree 
$\geq 1$ in $\rho$, and 
$\overline Q(R) - (21) \circ \overline Q(R) = \on{inj}_1^{\otimes 2} 
\circ \rho$ + terms of degree $\geq 2$ in $\rho$, where 
$\on{inj}_0\in {\bf Sch}({\bf 1},S)$ and $\on{inj}_1 \in 
{\bf Sch}({\bf id},S)$ are the canonical injection maps
(recall that $\on{Cob}$ has a grading where $\mu$ has degree $0$ and $\rho$
has degree $1$). 

As in the case of quantization functors of Lie bialgebras, 
$\overline Q$ necessarily satisfies $\overline Q(\eta)
= \on{inj}_0$, $\overline Q(\eps) = \on{pr}_0$. 
As we explained, each such morphism $\overline Q$ yields a solution of the 
quantization problem of coboundary Lie bialgebras. 

\subsection{Construction of quantization functors of 
coboundary Lie bialgebras}

\begin{theorem} \label{thm:q:cob}
Any even associator defined over $\kk$ 
gives rise to a quantization functor of coboundary Lie bialgebras. 
\end{theorem}

\begin{remark} In \cite{BN}, the existence of rational even associators 
is proved. This implies the existence of quantization functors
of coboundary Lie bialgebras over any field $\kk$ of characteristic $0$. 
\end{remark}

{\em Proof.} 
There is a unique automorphism $\tau_{\on{Cob}}$ of $\on{Cob}$, 
defined by $\mu\mapsto \mu$ and $\rho\mapsto -\rho$. Then the
following diagrams of prop morphisms commute 
\begin{equation} \label{diagrams}
\begin{matrix}
\on{LBA} & \stackrel{\kappa\circ\kappa_1}{\to} & \on{Cob}\\ 
\scriptstyle{\tau_{\on{LBA}}}\downarrow & & 
\downarrow\scriptstyle{\tau_{\on{Cob}}}\\
\on{LBA} & \stackrel{\kappa\circ\kappa_1}{\to} & \on{Cob}
\end{matrix}
\quad \on{and} \quad 
\begin{matrix}
 & \scriptstyle{\kappa_1}\nearrow & \on{LBA}_f & \stackrel{\kappa}{\to}& 
 \on{Cob}\\
\on{LBA} & & & & \downarrow\scriptstyle{\tau_{\on{Cob}}} \\
 & \scriptstyle{\kappa_2}\searrow & \on{LBA}_f & \stackrel{\kappa}{\to} 
 & \on{Cob}
\end{matrix}
\end{equation}

Let $Q: \on{Bialg} \to S({\bf LBA})$ be a quantization functor 
corresponding to an even associator. 
Then $\overline Q := S(\kappa\circ\kappa_1) \circ Q : \on{Bialg} \to 
S({\bf Cob})$ is a prop morphism.  
We have 
\begin{align*}
& S(\tau_{\on{Cob}}) \circ S(\kappa\circ\kappa_1) \circ Q = 
S(\kappa\circ\kappa_1) \circ S(\tau_{\on{LBA}}) \circ Q = 
S(\kappa\circ\kappa_1) \circ \theta(\xi_\tau^{-1}) \circ Q \circ 
\tau_{\on{Bialg}}
\\ & = 
\theta(S(\kappa\circ\kappa_1)(\xi_\tau^{-1})) \circ S(\kappa\circ\kappa_1) 
\circ Q \circ \tau_{\on{Bialg}}
\end{align*} 
(the first equality uses the first diagram of (\ref{diagrams}), and the 
second equality uses Proposition \ref{prop:coop}), 
so 
$$
S(\tau_{\on{Cob}}) \circ \overline Q = \theta(\xi'_\tau) \circ \overline Q 
\circ \tau_{\on{Bialg}}, 
$$ 
where $\xi_\tau = S(\kappa\circ \kappa_1)(\xi_\tau^{-1})$. 

On the other hand, there exists $\on{F}\in 
{\bf LBA}_f({\bf 1},S^{\otimes 2})^\times$
and $\on{i}\in S({\bf LBA}_f)({\bf id},{\bf id})^\times$, such that  
$S(\kappa_2) \circ Q = \theta(\on{i}) \circ {}^{\on{F}}(S(\kappa_1) \circ Q)$. 
Composing this equality with $S(\kappa)$, we get: 
$S(\kappa\circ\kappa_2) \circ Q = \theta(S(\kappa)(\on{i})) \circ 
{}^{S(\kappa)(\on{F})}(\overline Q)$. 
Now $S(\kappa\circ\kappa_2) \circ Q = S(\tau_{\on{Cob}}) \circ \overline Q$
(using the second diagram in (\ref{diagrams})), so 
$$
S(\tau_{\on{Cob}}) \circ \overline Q = \theta(S(\kappa)(\on{i})) \circ 
{}^{S(\kappa)(\on{F})}(\overline Q). 
$$ 

We therefore get $\theta(\xi'_\tau) \circ \overline Q \circ \tau_{\on{Bialg}} 
= \theta(S(\kappa)(\on{i})) \circ 
{}^{S^{\otimes 2}(\kappa)(\on{F})}\overline Q$, 
so 
\begin{equation} \label{starting:point}
{}^{\on{F}'}\overline Q = \theta(\xi'') \circ \overline Q 
\circ \tau_{\on{Bialg}}, 
\end{equation} 
where $\xi'' = S(\kappa)(\on{i})^{-1} \circ \xi'_\tau$, and 
$\on{F}' = S^{\otimes 2}(\kappa)(\on{F})$. Here $\xi''\in 
S({\bf Cob})({\bf id},{\bf id})^\times$ has the form $\on{id}_S$ + 
terms of positive degree in $\rho$, and $\on{F}'\in {\bf Cob}({\bf
1},S^{\otimes 2})$ satisfies 
\begin{equation} \label{twist:rel}
(\on{F}')_{\overline Q}^{1,2} \ast_{\overline Q} 
(\on{F}')_{\overline Q}^{12,3} 
= (\on{F}')_{\overline Q}^{2,3} \ast_{\overline Q}
(\on{F}')_{\overline Q}^{1,23}, \quad 
(\on{F}')_{\overline Q}^{\emptyset,1} = (\on{F}')_{\overline Q}^{1,\emptyset} 
= \on{inj}_0, 
\end{equation}
and 
\begin{equation} \label{limit:F}
\on{F}' = \on{inj}_0^{\otimes 2} + \rho  + \on{terms\ of\ degree\ }\geq 2
\on{\ in\ }\rho.
\end{equation}
We will prove: 

\begin{proposition} \label{constr:G}
There exists $\on{G}\in {\bf Cob}({\bf 1},S^{\otimes 2})$, satisfying 
(\ref{twist:rel}) (where $\overline Q$ is also used), (\ref{limit:F}), and
$$
\on{G} \ast_{\overline Q} \on{G}^{2,1} 
= \on{G}^{2,1} 
\ast_{\overline Q} \on{G} = \on{inj}_0^{\boxtimes 2}, 
\quad 
(21) \circ \overline Q(\Delta) = \underline{\on{Ad}}(\on{G}) 
\circ \overline Q(\Delta)
$$ 
(recall that the definition of $\underline{\on{Ad}}(\on{G})$
involves $\overline Q(m)$).
\end{proposition}

This proposition implies the theorem, since we now have a prop morphism 
$\on{COB} \to S({\bf Cob})$, obtained by extending $\overline Q : 
\on{Bialg} \to S({\bf Cob})$ by $R\mapsto \on{G}$. 

Let us now prove Proposition \ref{constr:G}. We start by making 
(\ref{starting:point}) explicit: this means that 
$$
\overline Q(m) = (\xi'')^{\boxtimes 2} \circ \overline Q(m) \circ (\xi'')^{-1}, 
\quad (\xi'')^{\boxtimes 2} \circ (21)
\circ \overline Q(\Delta) \circ (\xi'')^{-1} 
= \underline{\on{Ad}}(\on{F}') \circ \overline Q(\Delta),  
$$
$\overline Q(\eta) = \xi'' \circ \overline Q(\eta)$, 
$\overline Q(\eps) = \overline Q(\eps) \circ \xi''$. 

We first prove: 
 
\begin{lemma} \label{lemma:H}
There exists a unique $\on{H}\in {\bf Cob}({\bf 1},S^{\otimes 2})^\times$
such that $\on{H} = \on{inj}_0^{\otimes 2}$ + terms of degree 
$\geq 1$ in $\rho$, and $((\xi'')^{\boxtimes 2} \circ \on{H}) 
\ast_{\overline Q} \on{H} =  
((\xi'')^{\boxtimes 2} \circ (\on{F}')^{2,1}) 
\ast_{\overline Q} \on{F}'$. 
Then $(\xi'')^{\boxtimes 2}\circ \overline Q(\Delta) \circ (\xi'')^{-1}
= \underline{\on{Ad}}(\on{H}) \circ \overline Q(\Delta)$, 
$\on{H}$ satisfies the
identities (\ref{twist:rel}), and
$\on{H} = \on{inj}_0^{\otimes 2}$ + terms of degree $\geq 2$ in $\rho$.  
\end{lemma}  

{\em Proof of Lemma.} The existence of $\on{H}$ is a consequence of 
the following statement. Let $A = A^0 \supset A^1 \supset ...$ 
be a filtered algebra, complete and separated for this filtration. 
Let $\theta$ be a topological automorphism of $A$, such that 
$(\theta - \on{id}_A)(A^n) \subset A^{n+1}$ for any $n$.   
Let $u\in A$ be such that $u\equiv 1$ modulo $A^1$. Then there exists a 
unique $v\in A$, with $v\equiv 1$ modulo $A^1$ and $v\theta(v) = u$. 
We will apply this statement to $A = {\bf Cob}({\bf 1},S^{\otimes 2})$
equipped with the product given by $\overline Q(m)$. The filtration is 
given by the degree in $\rho$, and $\theta(\on{F}) = (\xi'')^{\boxtimes 2}
\circ \on{F}$. 

To prove the existence of $v$, we construct inductively the class $[v]_n$ 
of $v$ in $A/A^n$: assume that $[v]_n$ has been found such that 
$[v]_n \theta([v]_n) = [u]_n$ in $A/A^n$, and let $v'$ be a lift
of $[v]_n$ to $A/A^{n+1}$, then $v'\theta(v') \equiv [u]_{n+1}$
modulo $A^n/A^{n+1}$. Then we set $[v]_{n+1} = v' - 
(1/2) (v'\theta(v') - [u_{n+1}])$ (in $A/A^{n+1}$). 

Let us prove the uniqueness of $v$. Let $v$ and $v'$ be solutions; 
let us prove by induction on $n$ that $[v]_n = [v']_n$.
Assume that this has been proved up to order $n-1$ and let us prove it at 
order $n$. We have $v\theta(v) - v'\theta(v') = (v-v')\theta(v) 
+ v'(\theta(v) - \theta(v'))$. Then we have $v-v'\in A^{n-1}$, 
$\theta(v) - \theta(v') \in A^{n-1}$, and the classes of these
elements are equal in $A^{n-1}/A^n$. So the class of
$v\theta(v) - v'\theta(v')$ in $A^{n-1}/A^n$ is equal to twice the class of 
$v-v'$ in $A^{n-1}/A^n$. Since $v\theta(v) = v'\theta(v')$, the latter
class is $0$, so $v-v'\in A^n$. 
 
Before we prove the properties of $\on{H}$, we construct the following 
propic version of the theory of twists.  Let us denote by 
\boldmath$\Delta$\unboldmath \ 
the set of all $\overline\Delta\in 
{\bf Cob}(S,S^{\otimes 2})$, such that there exists a prop morphism 
$\overline Q_{\overline\Delta} : \on{Bialg} \to S({\bf Cob})$, such that 
$\overline Q_{\overline\Delta}(\Delta) = \overline\Delta$ and 
$\overline Q_{\overline\Delta}(m)
= \overline Q(m)$, $\overline Q_{\overline\Delta}(\eps) = \overline Q(\eps)$,
$\overline Q_{\overline\Delta}(\eta) = \overline Q(\eta)$.
For $\overline\Delta_1\in$\boldmath$\Delta$\unboldmath, 
we denote by $\on{Tw}(\overline\Delta_1,-)$ the set of all 
$\on{F}_1\in {\bf Cob}({\bf 1},S^{\otimes 2})^\times$ satisfying 
(\ref{twist:rel}) where the underlying structure is that given by 
$\overline Q_{\overline\Delta_1}$. Then if $\overline\Delta_2 := 
\on{Ad}(\on{F}_{1}) \circ \overline\Delta_1$, we have 
$\overline\Delta_2\in$\boldmath$\Delta$\unboldmath. 
If $\overline\Delta_1,\overline\Delta_2\in$\boldmath$\Delta$\unboldmath, 
let us denote by $\on{Tw}(\overline\Delta_1,\overline\Delta_2) \subset 
\on{Tw}(\overline\Delta_1,-)$ the set of all $\on{F}_1$ such that 
$\overline\Delta_2 = \underline{\on{Ad}}(\on{F}_{1}) 
\circ \overline\Delta_1$. 

Then if $\overline\Delta_i\in$\boldmath$\Delta$\unboldmath \ ($i=1,2,3$), 
the map $(\on{F}_1,\on{F}_2) \mapsto \on{F}_2\on{F}_1$ (product 
in ${\bf Cob}({\bf 1},S^{\otimes 2})$ using $\overline Q(m)$)   
defines a map $\on{Tw}(\overline\Delta_1,\overline\Delta_2)
\times \on{Tw}(\overline\Delta_2,\overline\Delta_3)
\to \on{Tw}(\overline\Delta_1,\overline\Delta_3)$. 

Let us now prove the properties of $\on{H}$. We have 
$\on{F}'\in \on{Tw}(\overline Q(\Delta),(\xi'')^{\boxtimes 2}
\circ (21) \circ \overline Q(\Delta) \circ (\xi'')^{-1})$
and $(\xi'')^{\boxtimes 2} \circ (\on{F}')^{2,1} \in 
\on{Tw} \big(   (\xi'')^{\boxtimes 2}
\circ (21) \circ \overline Q(\Delta) \circ (\xi'')^{-1} , 
((\xi'')^2)^{\boxtimes 2} \circ \overline Q(\Delta) \circ (\xi'')^{-2} \big)$, 
therefore 
$$
\cF := \big( (\xi'')^{\boxtimes 2} \circ (\on{F}')^{2,1} \big) 
\ast_{\overline Q} \on{F}'\in \on{Tw} \big( \overline Q(\Delta),
((\xi'')^2)^{\boxtimes 2} \circ \overline Q(\Delta) \circ (\xi'')^{-2} \big). 
$$
In particular, we have 
$$
((\xi'')^2)^{\boxtimes 2} \circ \overline Q(\Delta) \circ (\xi'')^{-2}
= \underline{\on{Ad}}(\cF) \circ \overline Q(\Delta). 
$$
Then if we set 
$$
\cF(n) = ((\xi'')^{\boxtimes 2(n-1)} \circ \cF) \ast_{\overline Q} ... 
\ast_{\overline Q} ((\xi'')^{\boxtimes 2} \circ \cF) \ast_{\overline Q} 
\cF, 
$$
we have for $n$ integer $\geq 0$, 
\begin{equation} \label{equality:cF}
((\xi'')^{2n})^{\boxtimes 2} \circ \overline Q(\Delta) \circ (\xi'')^{-2n}
= \underline{\on{Ad}}(\cF(n)) \circ \overline Q(\Delta). 
\end{equation}

${\bf Cob}(S,S)$ is the completion of a $\NN$-graded algebra, where 
the degree is given by $\on{deg}(\rho) = 1$, $\on{deg}(\mu)=0$.  
$\xi''\in {\bf Cob}(S,S)$ is equal to the identity modulo terms of 
positive degree. We have therefore a unique formal map 
$t\mapsto (\xi'')^t$, inducing a polynomial map $\kk \to {\bf Cob}(S,S)/\{$its
part of degree $>k\}$ for each $k\geq 0$, which coincides with 
the map $n\mapsto$ (class of $(\xi'')^n$) for $t\in\NN$. 
 
On the other hand, one checks that there is a unique formal map 
$t\mapsto \cF(t)$ with values in ${\bf Cob}({\bf 1},S^{\otimes 2})$, 
such that the induced map $\kk \to {\bf Cob}({\bf 1},S^{\otimes 2}) / \{$its
part of degree $>k\}$ is polynomial for any $k\geq 0$ and coincides with the 
maps $n\mapsto$ (class of $\cF(n)$) for $t\in \NN$.  

It follows that (\ref{equality:cF}) also holds when $n$ is replaced
by the formal variable $t$. The resulting identity can be specialized 
for $n=1/2$. The specialization of $(\xi'')^{2t}$ for $t=1/2$ is $\xi''$. 

We now prove that $\cF(1/2) = \on{H}$. 

Let us set $\on{H}(n) = ((\xi'')^{\boxtimes 2(n-1)} \circ \on{H}) 
\ast_{\overline Q} ... \ast_{\overline Q} 
((\xi'')^{\boxtimes 2} \circ \on{H}) \ast_{\overline Q} \on{H}$. 
Then we have a unique formal map 
$t\mapsto \on{H}(t)$ with values in ${\bf Cob}({\bf 1},S^{\otimes 2})$, 
such that the induced map $\kk \to {\bf Cob}({\bf 1},S^{\otimes 2}) / \{$its
part of degree $>k\}$ is polynomial for any $k\geq 0$ and coincides with the 
maps $n\mapsto$ (class of $\on{H}(n)$) for $t\in \NN$.  

We have $\on{H}(2n) = \cF(n)$ for any integer $n$, so this identity also 
holds when $n$ is replaced by the formal variable $t$. Specializing the
resulting identity for $t=1/2$, we get $\cF(1/2) = \on{H}(1) = \on{H}$. 

The specialization of the formal version of (\ref{equality:cF}) for 
$t=1/2$ then gives
\begin{equation} \label{first:twist}
(\xi'')^{\boxtimes 2} \circ \overline Q(\Delta) \circ (\xi'')^{-1}
= \underline{\on{Ad}}(\on{H}) \circ \overline Q(\Delta). 
\end{equation} 

Let us now prove the identities (\ref{twist:rel}) in $\on{H}$.
We have $(\xi''\circ (\on{H}^{\emptyset,1})) \ast_{\overline Q}
\on{H}^{\emptyset,1} = \on{inj}_0$ and 
$\on{H}^{\emptyset,1} = \on{inj}_0$ + 
terms of positive degree in $\rho$, hence by the uniqueness result 
proved above $\on{H}^{\emptyset,1} = \overline Q(\eta)$. In the same way, 
$\on{H}^{1,\emptyset} = \on{inj}_0$.  

We now prove that 
\begin{equation} \label{second:twist}
(\on{H})_{\overline Q}^{1,2} \ast_{\overline Q} 
(\on{H})_{\overline Q}^{12,3} = 
(\on{H})_{\overline Q}^{2,3} \ast_{\overline Q} 
(\on{H})_{\overline Q}^{1,23}. 
\end{equation}
Let $\Psi\in {\bf Cob}({\bf 1},S^{\otimes 3})^\times$ be such that 
$(\on{H})_{\overline Q}^{1,2} \ast_{\overline Q} 
(\on{H})_{\overline Q}^{12,3} 
= (\on{H})_{\overline Q}^{2,3} \ast_{\overline Q} 
(\on{H})_{\overline Q}^{1,23} \ast_{\overline Q} \Psi$. 

We have $(\cF)_{\overline Q}^{1,2} \ast_{\overline Q} 
(\cF)_{\overline Q}^{12,3} = (\cF)_{\overline Q}^{2,3} 
\ast_{\overline Q} (\cF)_{\overline Q}^{1,23}$, that is 
$$
(\cF \boxtimes \overline Q(\eta)) \ast_{\overline Q}
(\overline Q(\Delta \boxtimes \eta) \circ \cF) 
= (\overline Q(\eta) \boxtimes \cF) \ast_{\overline Q}
(\overline Q(\eta\boxtimes \Delta) \circ \cF)  
$$
This is rewritten 
\begin{align*}
& \big( (\xi'')^{\boxtimes 3} \circ (\on{H} \boxtimes \overline Q(\eta)) \big) 
\ast_{\overline Q}
(\on{H} \boxtimes \overline Q(\eta)) 
\ast_{\overline Q}
(\overline Q(\Delta \boxtimes \eta) \circ 
\xi^{\prime\prime\boxtimes 2} \circ \on{H} )
\ast_{\overline Q} 
(\overline Q(\Delta \boxtimes \eta) \circ \on{H} ) 
\\ & = 
\big( (\xi'')^{\boxtimes 3} \circ (\overline Q(\eta) \boxtimes \on{H} ) \big) 
\ast_{\overline Q}
(\overline Q(\eta) \boxtimes \on{H} )
\ast_{\overline Q} 
(\overline Q(\eta\boxtimes \Delta) \circ 
\xi^{\prime\prime\boxtimes 2} \circ \on{H} )  
\ast_{\overline Q}
(\overline Q(\eta\boxtimes \Delta) \circ \on{H})  
\end{align*}
Using $\xi''\circ \on{inj}_0 = \on{inj}_0$ and (\ref{first:twist}), we get  
\begin{align*}
& \big( (\xi'')^{\boxtimes 3} \circ (\on{H} \boxtimes \overline Q(\eta)) \big) 
\ast_{\overline Q}
( \xi^{\prime\prime\boxtimes 3} 
\circ \overline Q(\Delta \boxtimes \eta) \circ \on{H} ) 
\ast_{\overline Q}
(\on{H} \boxtimes \overline Q(\eta)) 
\ast_{\overline Q}
(\overline Q(\Delta \boxtimes \eta) \circ \on{H} ) 
\\ & = 
\big( (\xi'')^{\boxtimes 3} \circ (\overline Q(\eta) \boxtimes \on{H} ) \big) 
\ast_{\overline Q}
( \xi^{\prime\prime\boxtimes 3} \circ 
\overline Q(\eta\boxtimes \Delta) 
\circ \on{H} ) 
\ast_{\overline Q} 
(\overline Q(\eta) \boxtimes \on{H} ) 
\ast_{\overline Q}
(\overline Q(\eta\boxtimes \Delta) \circ \on{H})  
\end{align*}
and since $X\mapsto (\xi'')^{\boxtimes 3} \circ X$ is an automorphism of 
${\bf Cob}({\bf 1},S^{\otimes 3})$, we get 
\begin{align*}
& \Big( (\xi'')^{\boxtimes 3} \circ 
\Big( (\on{H} \boxtimes \overline Q(\eta))  
(\overline Q(\Delta \boxtimes \eta) \circ \on{H} ) \Big) \Big) 
\ast_{\overline Q}
(\on{H} \boxtimes \overline Q(\eta)) \ast_{\overline Q}
(\overline Q(\Delta \boxtimes \eta) \circ \on{H} ) 
\\ & = 
\Big( (\xi'')^{\boxtimes 3} \circ 
\Big( (\overline Q(\eta) \boxtimes \on{H} )  
( \overline Q(\eta\boxtimes \Delta) 
\circ \on{H} )  \Big) \Big) 
\ast_{\overline Q}
(\overline Q(\eta) \boxtimes \on{H} ) 
\ast_{\overline Q}
(\overline Q(\eta\boxtimes \Delta) \circ \on{H}) ,  
\end{align*}
i.e., 
$$
((\xi'')^{\boxtimes 3} \circ \Psi)
\ast_{\overline Q}
\big((\on{H})_{\overline Q}^{2,3} \ast_{\overline Q} 
(\on{H})_{\overline Q}^{1,23}\big)
\ast_{\overline Q} \Psi = 
((\on{H})_{\overline Q}^{2,3} \ast_{\overline Q}
(\on{H})_{\overline Q}^{1,23}). 
$$
We now prove that this implies that $\Psi = \on{inj}_0^{\boxtimes 3}$. 

For this, we apply the following general statement. Let $A = A^0 \supset A^1
\supset ...$ be an algebra equipped with a decreasing filtration, complete and
separated for this filtration. Let $\theta$ be a topological automorphism of
$A$, such that $(\theta - \on{id}_A)(A^n) \subset A^{n+1}$. 
Let $X\in A$ be such that $X \equiv 1$
modulo $A^1$. Let $x\in A$ be such that $x\equiv 1$ modulo $A^1$, and 
$$
\theta(x) X x = X. 
$$
Then $x = 1$. This is proved by induction. Assume that we have proved
that $x \equiv 1$ modulo $A^{n-1}$.  Then $\theta(x) X x X^{-1} \equiv 
1 + 2(x-1)$ modulo $A^n$. Therefore $x\equiv 1$ modulo $A^n$. Finally $x=1$. 

We then apply the general statement to 
$A = {\bf Cob}({\bf 1},S^{\otimes 3})$ and $\theta : X\mapsto 
(\xi'')^{\boxtimes 3} \circ X$ and get 
$\Psi = \on{inj}_0^{\boxtimes 3}$. This implies (\ref{second:twist}). 
This ends the proof of Lemma \ref{lemma:H}. \hfill \qed \medskip

We now end the proof of Proposition \ref{constr:G}. 
Lemma \ref{lemma:H} says that 
$\on{H} \in \on{Tw}(\overline Q(\Delta),(\xi'')^{\boxtimes 2} \circ 
\overline Q(\Delta)\circ (\xi'')^{-1})$, and since 
$\on{F}' \in \on{Tw}(\overline Q(\Delta),(\xi'')^{\boxtimes 2} \circ 
(21) \circ \overline Q(\Delta)\circ (\xi'')^{-1})$, we have 
$$
\on{G}' := \big( (21) \circ \on{H}^{-1} \big) 
\ast_{\overline Q} \on{F}' \in 
\on{Tw}(\overline Q(\Delta),(21) \circ \overline Q(\Delta)). 
$$

Let us set $\cG' := (\on{G}')^{2,1} \ast_{\overline Q}\on{G}'$, then  
$\cG' \in \on{Tw}(\overline Q(\Delta),\overline Q(\Delta))$. 
For any integer $n\geq 0$, we then have $(\cG')^n\in 
\on{Tw}(\overline Q(\Delta),\overline Q(\Delta))$. 
As before, there exists a unique formal map $t\mapsto (\cG')^t$, 
such that the map $t\mapsto$ (class of $(\cG')^t$) in 
${\bf Cob}({\bf
1},S^{\otimes n}) / \{$its part of degree $\geq k\}$ is polynomial and
extends $n\mapsto (\cG')^n$. Specializing for $t=-1/2$, we get 
$(\cG')^{-1/2}\in 
\on{Tw}(\overline Q(\Delta),\overline Q(\Delta))$. Set 
$$
\on{G} := \on{G}' \ast_{\overline Q} (\cG')^{-1/2} 
= \on{G}' \ast_{\overline Q} 
(\on{G}^{\prime 2,1} \ast_{\overline Q}\on{G}')^{-1/2}, 
$$
then 
$$
\on{G}\in \on{Tw}(\overline Q(\Delta),(21)\circ \overline Q(\Delta)). 
$$

Then we have: $\on{G}' \ast_{\overline Q}
(\on{G}^{\prime 2,1} \ast_{\overline Q}\on{G}')^n
= (\on{G}' \ast_{\overline Q}\on{G}^{\prime 2,1})^n \ast_{\overline Q}
\on{G}'$ for any integer $n\geq 0$, 
so this identity also holds when $n$ is replaced by a formal variable 
$t$. Specializing the latter identity to $t=1/2$, we get 
$\on{G} = (\on{G}' \ast_{\overline Q}\on{G}^{\prime 2,1})^{-1/2}
\ast_{\overline Q} \on{G}'$. 
Then 
$\on{G}\ast_{\overline Q} \on{G}^{2,1} = \on{G}' 
\ast_{\overline Q} (\on{G}^{\prime 2,1} \ast_{\overline Q}
\on{G}')^{-1} \ast_{\overline Q}
\on{G}^{\prime 2,1} = \on{inj}_0^{\boxtimes 2}$, so we also have 
$\on{G}^{2,1} \ast_{\overline Q} \on{G} = \on{inj}_0^{\boxtimes 2}$.

This ends the proof of Proposition \ref{constr:G}, and therefore also 
of Theorem \ref{thm:q:cob}. \hfill \qed \medskip 

\begin{remark} The proof of Proposition \ref{constr:G} is a propic version of
the proof of the following statement. Let $(U,m_U,\eta_U)$ be a formal 
deformation over $\kk[[\hbar]]$ of an enveloping algebra $U(\a)$
(as an algebra). Let \boldmath$\Delta$\unboldmath$_U$ be the set of all
morphisms $\Delta : U \to
U^{\otimes 2}$ such that $(U,m_U,\Delta,\eta_U,\eps_U)$ is a QUE 
algebra formally deforming the bialgebra $U(\a)$. If $\Delta_1,\Delta_2
\in $\boldmath$\Delta$\unboldmath$_U$, we say that $F_U \in
\on{Tw}(\Delta_1,\Delta_2)$ iff $F_U\in (U^{\otimes 2})^\times$, 
$(\eps_U\otimes\on{id}_U)(F_U) = (\on{id}_U\otimes \eps_U)(F_U) = 1_U$, 
$(F_U\otimes 1_U)(\Delta_1 \otimes \on{id}_U)(F_U) = 
(1_U\otimes F_U)(\on{id}_U \otimes \Delta_1)(F_U)$ and $\Delta_2 = \on{Ad}(F_U)
\circ \Delta_1$, where $\on{Ad}(F_U) : U^{\otimes 2}\to U^{\otimes 2}$ 
is given by $x\mapsto F_U x F_U^{-1}$ (and $1_U = \eta_U(1)$). 
For $\Delta \in $\boldmath$\Delta$\unboldmath$_U$
and $\theta_U\in \on{Aut}(U,m_U,\eta_U)$ such that
$\theta_U = \on{id}_U + O(\hbar)$, we also have 
$\Delta^{21}\in $\boldmath$\Delta$\unboldmath$_U, 
\theta_U^{\otimes 2} \circ \Delta^{21} \circ \theta_U^{-1}\in 
$\boldmath$\Delta$\unboldmath$_U$. The statement is that if for such 
$\Delta,\theta_{U}$, there exists 
$F_U\in \on{Tw}(\Delta,\theta_U^{\otimes 2} \circ \Delta^{21} 
\circ \theta_U^{-1})$, then there exists $G_U\in \on{Tw}(\Delta,\Delta^{21})$ 
such that $G_UG_U^{2,1}=1_U^{\otimes 2}$. 
\end{remark}

\subsection{Relation with quasi-Poisson manifolds} \label{sec:qpoisson}

Define a coboundary quasi-Lie bialgebra (QLBA) as a set $(\g,\mu_\g,\delta_\g,
Z_\g,r_\g)$, where $(\g,\mu_\g,\delta_\g,Z_\g)$ is a quasi-Lie bialgebra, 
and $r_\g\in\wedge^2(\g)$ is such that $\delta_\g(x) = [r_\g,x\otimes 1 + 1
\otimes x]$. In \cite{Dr:QH}, coboundary QUE quasi-Hopf algebras were 
introduced; the classical limit of this structure is a coboundary QLBA. 

According to \cite{Dr:QH}, Proposition 3.13, a coboundary QUE quasi-Hopf
algebra with classical limit the coboundary QLBA $(\g,\mu_\g,\delta_\g,
Z_\g,r_\g)$ is twist-equivalent to a coboundary QUE Hopf algebra of the form 
$(U(\g_\hbar),m_0,\Delta_0,R = 1, \Phi = {\cal E}(\hbar^2 Z_\hbar))$, 
where $\g_\hbar$ is a deformation of $\g$ (as a Lie algebra) 
in the category of topologically free $\kk[[\hbar]]$-modules, 
$Z_\hbar \in \wedge^3(\g_\hbar)^{\g_\hbar}$ is a deformation of 
$Z_\g + (\delta_\g\otimes\on{id})(r_\g) + \on{c.p.} - \on{CYB}(r_\g)$, 
and ${\cal E}(Z) = 1 + Z/6 + ...$ is a series introduced in 
\cite{Dr:QH} ($m_0,\Delta_0$ are the undeformed operations). 

Let now $(\a,r_\a)$ be a coboundary Lie bialgebra. Let $(U_\hbar(\a),R_\a)$
be a quantization of it: this is a coboundary QUE Hopf algebra. Applying
to it the above result, we obtain: 

a) there exists a deformation $\a_\hbar$ of $\a$ in the category of 
topologically free $\kk[[\hbar]]$-Lie algebras, such that 
$U_\hbar(\a)$ is isomorphic to $U(\a_\hbar)$ as an algebra; 

b) there exists $J\in U(\a_\hbar)^{\otimes 2}$ of the form 
$J = 1 + \hbar r_\a/2 + O(\hbar^2)$ such that 
$J^{2,3}J^{1,23}{\cal E}(\hbar^2 Z_\hbar) = J^{1,2}J^{12,3}$, 
where $Z_\hbar\in \wedge^3(\a_\hbar)^{\a_\hbar}$
is a deformation of $Z_\a$. 

If $A$ is a Lie group with Lie algebra $\a$, then $r_\a$ induces 
a quasi-Poisson homogeneous structure on $A$ under the action of the 
quasi-Lie bialgebra $(\a,\delta_\a = 0,Z_\a)$: the action of $\a$
is the regular left action, and the quasi-Poisson structure is 
$\{f,g\} = m \circ {\bf L}^{\otimes 2}(r_\a)(f\otimes g)$, where 
$f,g$ are functions on $A$ (and $m$ is the product of functions). 
As explained in \cite{EE}, $J$ constructed above gives rise to a 
quantization of this quasi-Poisson homogeneous space, compatible 
with the quasi-Hopf algebra $(U(\a_\hbar),m_0,\Delta_0,R = 1,\Phi = 
{\cal E}(\hbar^2 Z_\hbar))$. 

Notice that the deformation class of $(\a_\hbar,Z_\hbar)$ is a priori 
dictated by $r_\a$. 


\end{document}